\documentclass{gtart_a}
\pdfoutput=1

\usepackage{xy,graphicx}
\xyoption{all}
\usepackage[dvipsnames]{color}

\title{A cylindrical reformulation of Heegaard Floer homology}
\author{Robert Lipshitz}
\givenname{Robert}
\surname{Lipshitz}
\address{Department of Mathematics\\Stanford University\\
Stanford, CA 94305-2125\\USA}
\email{lipshitz@math.stanford.edu}
\urladdr{}

\volumenumber{10}
\issuenumber{}
\publicationyear{2006}
\papernumber{25}
\lognumber{0609}
\startpage{955}
\endpage{1096}

\doi{}
\MR{}
\Zbl{}


\received{14 May 2005}
\revised{9 October 2005}
\accepted{3 January 2006}
\proposed{Peter Ozsv\'ath}
\seconded{John Morgan, Ronald Fintushel}
\published{9 August 2006}
\publishedonline{9 August 2006}
\corresponding{}
\editor{}
\version{}

\arxivreference{math.SG/0502404}

\let\xysavmatrix\xymatrix
\def\xymatrix{\disablesubscriptcorrection\xysavmatrix}
\AtBeginDocument{\let\bar\wbar\let\tilde\wtilde\let\hat\what}
\def\SetFigFont#1#2#3#4#5{\small}
\makeop{Sym}
\makeop{Spin}
\makeop{Map}

\makeautorefname{subsection}{Subsection}
\makeautorefname{subsubsection}{Subsubsection}

\newcommand{\bZ}{\mathbb{Z}}
\newcommand{\vx}{{\vec{x}}}
\newcommand{\vy}{{\vec{y}}}
\newcommand{\vz}{{\vec{z}}}
\newcommand{\vw}{{\vec{w}}}
\newcommand{\frako}{\mathfrak{o}}
\newcommand{\bN}{\mathbb{N}}
\newcommand{\bC}{\mathbb{C}}
\newcommand{\bD}{\mathbb{D}}
\newcommand{\bR}{\mathbb{R}}
\newcommand{\bH}{\mathbb{H}}
\newcommand{\Ss}{\mathfrak{s}}
\renewcommand{\SS}{\mathfrak{S}}
\newcommand{\St}{\mathfrak{t}}
\newcommand{\sm}{\setminus}

\newcommand{\ba}{\boldsymbol{\alpha}}
\newcommand{\bb}{\boldsymbol{\beta}}
\newcommand{\bc}{\boldsymbol{\gamma}}
\newcommand{\bd}{\boldsymbol{\delta}}
\newcommand{\hcM}{\widehat{\mathcal{M}}}
\newcommand{\cM}{\mathcal{M}}
\newcommand{\cX}{\mathcal{X}}
\newcommand{\cE}{\mathcal{E}}

\newcommand{\cB}{\mathcal{B}}
\newcommand{\cJ}{\mathcal{J}}
\newcommand{\tD}{\tilde{D}}
\newcommand{\tE}{\tilde{E}}
\newcommand{\coker}{\mathrm{coker}}
\newcommand{\ind}{\mathrm{ind}}
\newcommand{\sign}{\mathrm{sign}}
\newcommand{\cL}{\mathcal{L}}
\newcommand{\cT}{\mathcal{T}}
\newcommand{\cH}{\mathcal{H}}
\newcommand{\va}{\vec{\alpha}}
\newcommand{\vb}{\vec{\beta}}
\newcommand{\vc}{\vec{\gamma}}
\newcommand{\vd}{\vec{\delta}}
\newcommand{\Fz}{\mathfrak{z}}
\newcommand{\Yab}{Y_{\alpha,\beta}}
\newcommand{\Ybc}{Y_{\beta,\gamma}}
\newcommand{\Yac}{Y_{\alpha,\gamma}}
\newcommand{\End}{\mathrm{End}}
\newcommand{\Hom}{\mathrm{Hom}}
\newcommand{\Span}{\mathrm{Span}}
\newcommand{\rembds}{\hookrightarrow}
\newcommand{\bdy}{\partial}
\newcommand{\dbar}{\overline\partial}
            
\newcommand{\ol}{\overline}

\newcommand{\mapsinto}{\hookrightarrow}

\makeatletter
\def\cnewtheorem#1[#2]#3{\newtheorem{#1}{#3}[section]
\expandafter\let\csname c@#1\endcsname\c@Def}

\theoremstyle{definition}
\newtheorem{Def}{Definition}[section]
\theoremstyle{plain}
\cnewtheorem{Lem}[Def]{Lemma}
\cnewtheorem{Cor}[Def]{Corollary}
\cnewtheorem{Con}[Def]{Conjecture}
\newtheorem{Thm}{Theorem}
\cnewtheorem{Prop}[Def]{Proposition}
\cnewtheorem{Construct}[Def]{Construct}
\cnewtheorem{Sublem}[Def]{Sublemma}

\makeatother  

\makeautorefname{Construct}{Construct}
\makeautorefname{Cor}{Corollary}
\makeautorefname{Lem}{Lemma}
\makeautorefname{Def}{Definition}
\makeautorefname{Thm}{Theorem}
\makeautorefname{Sublem}{Sublemma}
\makeautorefname{Prop}{Proposition}

\newcommand{\Ah}{\Xi}
\newcommand{\Ag}{\Upsilon}

\begin{document}

\begin{asciiabstract}
We reformulate Heegaard Floer homology in terms of holomorphic curves
in the cylindrical manifold Sigma x [0,1] x R, where Sigma is the
Heegaard surface, instead of Sym^g(Sigma). We then show that the
entire invariance proof can be carried out in our setting. In the
process, we derive a new formula for the index of the dbar-operator in
Heegaard Floer homology, and shorten several proofs. After proving
invariance, we show that our construction is equivalent to the
original construction of Ozsvath-Szabo. We conclude with a discussion
of elaborations of Heegaard Floer homology suggested by our
construction, as well as a brief discussion of the relation with a
program of C Taubes.
\end{asciiabstract}

\begin{htmlabstract}
We reformulate Heegaard Floer homology in terms of holomorphic curves
in the cylindrical manifold &Sigma;&times;[0,1]&times;R, where
&Sigma; is the Heegaard surface, instead of Sym<sup>g</sup>(&Sigma;). We
then show that the entire invariance proof can be carried out in our
setting. In the process, we derive a new formula for the index of the
<span style="text-decoration:overline">&part;</span>&ndash;operator
in Heegaard Floer homology, and shorten several proofs. After proving
invariance, we show that our construction is equivalent to the original
construction of Ozsv&aacute;th&ndash;Szab&oacute;. We conclude with
a discussion of elaborations of Heegaard Floer homology suggested by
our construction, as well as a brief discussion of the relation with a
program of C Taubes.
\end{htmlabstract}

\begin{abstract}
We reformulate Heegaard Floer homology in terms of holomorphic curves
in the cylindrical manifold $\Sigma \times [0,1] \times R$, where
$\Sigma$ is the Heegaard surface, instead of $\mathrm{Sym}^g(\Sigma)$. We then
show that the entire invariance proof can be carried out in our
setting. In the process, we derive a new formula for the index of the
$\overline\partial$--operator in Heegaard Floer homology, and shorten several
proofs. After proving invariance, we show that our construction is
equivalent to the original construction of Ozsv\'ath--Szab\'o. We
conclude with a discussion of elaborations of Heegaard Floer homology
suggested by our construction, as well as a brief discussion of the
relation with a program of C Taubes.
\end{abstract}

\maketitle

In~\cite{OS1}, P Ozsv\'ath and Z Szab\'o associated to a
three--manifold $Y$ and a $\Spin^\bC$--structure $\Ss$ on $Y$ a
collection of abelian groups, known together as Heegaard Floer
homology.  These groups, which are believed to be isomorphic to
certain Seiberg--Witten Floer homology groups (Ozsv\'ath--Szab\'o
\cite{OS2} and Kronheimer--Manolescu \cite{Ciprian}), fit into the
framework of a $(3+1)$--dimensional topological quantum field theory.
Since its discovery around the turn of the millennium, Heegaard Floer
homology has been applied by Ozsv\'ath, Rasmussen and Szab\'o to the
study of knots and surgery \cite{OSKnots,Rasmussen,OSKnots2}, contact
structures \cite{OSContact} and symplectic structures
\cite{OSSymplectic}, and is strong enough to reprove most results
about smooth four--manifolds originally proved by gauge theory
\cite{OSTri}.  In this paper we give an alternate definition of the
Heegaard Floer homology groups.

Rather than being associated directly to a three--manifold $Y$, the
Heegaard Floer homology groups defined in~\cite{OS1} and in this paper
are associated to a Heegaard diagram for $Y$, as well as a
$\Spin^\bC$--structure $\Ss$ and
some additional structure.  A \emph{Heegaard diagram} is a closed,
orientable surface $\Sigma$ of genus $g$, together with two
$g$--tuples of pairwise disjoint, homologically linearly
independent, simple closed curves $\va=\{\alpha_1,\cdots,\alpha_g\}$
and $\vb=\{\beta_1,\cdots,\beta_g\}$ in $\Sigma$.  A Heegaard
diagram specifies a three--manifold as follows.  Thicken $\Sigma$ to
$\Sigma\times[0,1]$.  Glue thickened disks along the
$\alpha_i\times\{0\}$ and along the
$\beta_j\times\{1\}$.  The resulting space has two boundary
components, each homeomorphic to $S^2$.  Cap each with a
three--ball.  The result is the three--manifold specified by
$(\Sigma,\va,\vb)$.

Different Heegaard diagrams can specify the same three--manifold.  Two different Heegaard diagrams specify the same three--manifold if and only if they agree after a sequence of moves of the following three kinds:
\begin{itemize}
\item Isotopies of the $\alpha$-- or $\beta$--circles.
\item Handleslides among the $\alpha$-- or $\beta$--circles.  These correspond to pulling one $\alpha$-- (or $\beta$--) circle over another.
\item Stabilization, which corresponds to taking the connect sum of the Heegaard diagram with the standard genus--one Heegaard diagram for $S^3$.
\end{itemize}
See Gompf and Stipsicz~\cite[Sections 4.3 and 5.1]{GS} or Ozsv\'ath
and Szab\'o~\cite[Section 2]{OS1} for more details.

So, after associating the Heegaard Floer homology groups to a Heegaard
diagram, one must prove they are unchanged by these three kinds of
Heegaard moves (as well as deforming the additional structure involved in their definition).  Doing so comprises most of~\cite{OS1} for the original definition.  Similarly, for our definition, most of this paper is involved in proving:
\begin{Thm}\label{Theorem:Main}The Heegaard Floer homology groups
  $$HF^\infty(\Sigma,\va,\vb,\Ss),\qua HF^+(\Sigma,\va,\vb,\Ss),\qua
  HF^-(\Sigma,\va,\vb,\Ss)\qua\text{and}\qua
  \widehat{HF}(\Sigma,\va,\vb,\Ss)$$ associated to a Heegaard
  diagram $(\Sigma,\va,\vb)$ and $\Spin^\bC$--structure $\Ss$ are
  in fact invariants of the pair $(Y,\Ss)$.
\end{Thm}

We are also able to prove:
\begin{Thm}\label{Theorem:Equivalent} The Heegaard Floer homology groups defined in this paper are isomorphic to the corresponding groups defined in~\cite{OS1}.
\end{Thm}
\fullref{Theorem:Equivalent} is proved in
\fullref{Section:Comparison}.  The proof does not rely on the
invariance results proved in this paper; it could be carried out immediately after \fullref{Section:ChainComplexes}.  (We defer the proof to the end to avoid interrupting the narrative flow.)  Clearly, \fullref{Theorem:Equivalent} implies \fullref{Theorem:Main}.  However, one key goal of this paper is to demonstrate that the entire invariance proof can be carried out in our setting, and to develop the tools necessary to do so.

The only esentially new results in this paper are in
\fullref{Section:Index}, where we give a nice formula for the
index of the $\dbar$ operator in our setup, and hence also the
Maslov index in the traditional setting, and in the discussion of
elaborations of Heegaard Floer homology in the last section (\fullref{Section:OtherRemarks}).  The casual reader might also be interested in
looking at the elaboration and speculation in \fullref{Section:OtherRemarks}.

Although this paper is essentially self contained, it is probably
most useful to read it in parallel with~\cite{OS1}.  To facilitate this, the paper is
organized similarly to~\cite{OS1}, and throughout there are precise
references to corresponding results in their original forms.  In addition, the last
appendix is a table cross referencing most of the results in this
paper with those of~\cite{OS1}.

A more technical discussion of the difference between our setup and that of~\cite{OS1} follows.

The original definition of Heegaard Floer homology involves holomorphic disks in $\Sym^g(\Sigma)$.  In this paper, we consider holomorphic curves in $\Sigma\times[0,1]\times\bR$.  For instance, for us the chain complex $\widehat{CF}$ is generated by $g$--tuples of Reeb chords $\{x_i\times [0,1]\textrm{ }|\textrm{
}x_i\in\alpha_i\cap\beta_{\sigma(i)}\}$.  For an appropriate almost
complex structure $J$ on $\Sigma\times[0,1]\times\bR$, the coefficient of
$\{y_i\times[0,1]\}$ in $\bdy \left(\{x_i\times[0,1]\}\right)$ is
given by counting holomorphic curves in $\Sigma\times[0,1]\times\bR$ asymptotic to
$\{x_i\times[0,1]\}$ at $-\infty$ and to $\{y_i\times[0,1]\}$ at
$\infty$, with boundary mapped to the Lagrangian cylinders
$\{\alpha_j\times\{1\}\times\bR\}$ and $\{\beta_j\times\{0\}\times\bR\}$.
(We impose a few further technical conditions on the curves that we count;
see \fullref{Section:Notation}.)

If $J$ is the split complex structure $j_\Sigma\times j_\bD$ then
a holomorphic curve in $\Sigma\times[0,1]\times\bR$ is just a surface
$S$ and a pair of holomorphic maps $u_\Sigma:(S,\bdy S)\to
(\Sigma,\alpha_1\cup\cdots\cup\alpha_g\cup\beta_1\cup\cdots\cup\beta_g)$
and $u_\bD:(S,\bdy S)\to(\bD,\bdy\bD)$.
If the map $u_\bD$ is a $g$--fold branched covering then this data
specifies a map $\bD\to \Sym^g(\Sigma)$ as follows.  For $p\in\bD$,
let $p_1,\cdots,p_g$ be the preiamges of $p$ under $\pi_\bD\circ
u$, listed with multiplicity.  Then the map $\bD\to \Sym^g(\Sigma)$
sends $p$ to $\{u_\Sigma(p_1),\cdots,u_\Sigma(p_g)\}$.

 Note that the idea of viewing a map to $\Sym^g(\Sigma)$ as a pair
$$(\textrm{a $g$--fold covering $S\to\bD$},\textrm{ a map
  $S\to\Sigma$})$$
is already implicit in~\cite{OS1}, although they use this idea mainly for calculations in special cases.

Working in $\Sigma\times[0,1]\times\bR$ has several advantages.  A
main advantage is that, unlike a $g$--fold symmetric product, one can
actually visualize $\Sigma\times[0,1]\times\bR$.  A second advantage
is that a number of the technical details become somewhat simpler.
The main disadvantage is that we must now consider higher genus
holomorphic curves, not just disks.  Another difficulty is that our
setup requires compactness for holomorphic curves in manifolds with
cylindrical ends, proved by Bourgeois et al in~\cite{Ya2}.  I also
borrow from the language of symplectic field theory.  Fortunately,
much of the subtle machinery of symplectic field theory, like virtual
cycles or the operator formalism, is unnecessary for this paper.

The paper is organized as follows.  The first two sections are devoted to basic definitions and notation, and certain algebro--topological considerations.  The third section proves transversality results necessary for the rest of the paper.  These results should be standard, but I am unaware of a reference that applies to our setting.

The fourth section discusses the index of the $\dbar$--operator in our context.  We prove this index is the same as the Maslov index in the traditional setting, and obtain a combinatorial formula for it.  The fifth section discusses so--called admissibility criteria necessary for the case $b_1(Y)>0$.  The definitions and results are completely analogous with~\cite{OS1}.  The sixth section discusses coherent orientations of the moduli spaces.  Again, our treatment is close to~\cite{OS1}.

The seventh section rules out undesirable codimension--one degenerations of our holomorphic curves.  After doing so, we are finally ready to define the Heegaard Floer chain complexes in \fullref{Section:ChainComplexes}, and turn to the invariance proof.  The ninth section proves isotopy independence.  Before proving handleslide independence, we introduce triangle maps in \fullref{Section:Triangles}.  (As in~\cite{OS1}, to a \emph{Heegaard triple--diagram} $(\Sigma,\alpha_1,\cdots,\alpha_g,\beta_1,\cdots,\beta_g,\gamma_1,\cdots,\gamma_g)$ is associated maps $$HF(\Sigma,\va,\vb)\otimes HF(\Sigma,\vb,\vc)\to HF(\Sigma,\va,\vc),$$ for various decorations of $HF$.)  Using these triangle maps and a model computation, we prove handleslide invariance in \fullref{Section:Handleslides}.

Finally, in section twelve we prove stabilization invariance, completing the invariance proof.  After this, we devote a section to proving equivalence with traditional Heegaard Floer homology and a section to elaborations and speculation.

There are also two appendices.  The first is devoted to the gluing results used throughout the paper.  The second cross references our results with those in~\cite{OS1}.

For technical results about holomorphic curves, this paper sometimes
cites recent sources when older ones would suffice.  This generally
reflects either that the newer results are more broadly applicable or
that I found the newer exposition significantly clearer.

I thank Ya Eliashberg, who is responsible for
communicating to me most of the ideas in this paper.  I also
thank Z Szab\'o for a helpful conversation about the index;
P Ozsv\'ath for a helpful conversation about annoying curves (see \fullref{Section:ChainComplexes} below) and pointing out a serious omission in \fullref{Section:Comparison}; P Melvin for a stimulating conversation about the index;
M Hutchings for a discussion clarifying the relation between the
$H_1(Y)$--action and twisted coefficients; M Hedden and C Wendl for pointing out errors, both typographical and otherwise, in a previous version; and
W Hsiang, C Manolescu, L Ng and B Parker for comments that have improved the exposition.  Finally, I thank the referees for finding several errors and making many helpful suggestions.

This work was partially supported by the NSF Graduate Research Fellowship Program, and partly by the NSF Focused Research Group grant DMS--0244663.

\section{Basic definitions and notation}
\label{Section:Notation}

By a \emph{pointed Heegaard diagram} we mean a Heegaard diagram (as discussed in the introduction) together with a chosen point $\Fz$ of the Heegaard surface in the complement of the $\alpha$-- and $\beta$--circles.  Fix a pointed Heegaard diagram $\cH=(\Sigma_g,\vec{\alpha}=\{\alpha_1,\cdots,\alpha_g\},
\vec{\beta}=\{\beta_1,\cdots,\beta_g\},\Fz)$.  Let $\ba=\alpha_1\cup\cdots\cup\alpha_g\subset\Sigma$ and $\bb=\beta_1\cup\cdots\cup\beta_g\subset\Sigma$.  Consider the manifold
$W=\Sigma\times[0,1]\times\bR$.  We let $(p,s,t)$
denote a point in $W$ (so $p\in\Sigma$, $s\in[0,1]$ and
 $t\in\bR$).  Let $\pi_\bD\co W\to[0,1]\times\bR$, $\pi_\bR\co W\to\bR$
and $\pi_\Sigma\co W\to\Sigma_g$ denote the obvious projections.  Consider the
cylinders $C_\alpha = \ba\times\{1\}\times\bR$ and
$C_\beta = \bb\times\{0\}\times\bR.$  We will obtain Heegaard Floer
 homology by constructing a boundary map counting holomorphic curves
 with boundary on $C_\alpha\cup C_\beta$ and appropriate asymptotics
 at $\pm\infty$.

We shall always assume $g>1$, as the $g=1$ case is slightly
different technically.  Since we can stabilize any Heegaard diagram,
this does not restrict the class of manifolds under consideration.

Fix a point $\Fz_i$ in each component $D_i$ of
$\Sigma_g\sm(\ba\cup\bb)$.  Let $dA$ be an area form
on $\Sigma$, and $j_\Sigma$ a complex structure on $\Sigma$
tamed by $dA$.  Let $\omega=ds\wedge dt + dA$, a split
symplectic form on $W$.  Let $J$ be an almost complex
structure on $W$ such that
\begin{enumerate}
\item[(\textbf{J1})] $J$ is tamed by $\omega$.
\item[(\textbf{J2})] In a cylindrical neighborhood $U_{\{\Fz_i\}}$ of
 $\{\Fz_i\}\times[0,1]\times\bR$, 
 $J = j_\Sigma\times j_\bD$ is split.  (Here, $U_{\{\Fz_i\}}$ is small
 enough that its closure does not intersect
 $(\ba\cup\bb)\times[0,1]\times\bR$).
\item[(\textbf{J3})] $J$ is translation invariant in the $\bR$--factor.
\item[(\textbf{J4})] $J(\partial/\partial t) = \partial/\partial s$
\item[(\textbf{J5})] $J$ preserves $T(\Sigma\times\{(s,t)\})\subset
  TW$ for all 
  $(s,t)\in [0,1]\times\bR$.
\end{enumerate}
The first requirement is in order to obtain compactness of the moduli
spaces.  The second is for ``positivity of domains'' (see twelve paragraphs below).  The
third and fourth make $W$ \emph{cylindrical} as defined in
\cite[Section 2.1]{Ya2}.  The fifth ensures that our complex structure is
\emph{symmetric} and \emph{adjusted to $\omega$} in the sense 
of \cite[Section 2.1 and Section 2.2]{Ya2}.  (Note that $W$ is Levi--flat as defined there.
The vector field $\mathbf{R}$ introduced there is
$\partial/\partial s$.  The form $\lambda$ is just $ds $.) 

Note that we can view $J$ as a path $J_s$ of
complex structures on $\Sigma$.  Also notice that $C_\alpha$ and
$C_\beta$ are Lagrangian with respect to $\omega$.

At one point later -- the proof of~\ref{Lemma:InfinityDefined} -- we
need to consider almost complex structures which, instead of satisfying (\textbf{J5}),  satisfy the slightly less restrictive condition
\begin{enumerate}\leftskip 3pt
\item[(\textbf{J5$'$})] there is a 2--plane distribution $\xi$ on
  $\Sigma\times[0,1]$ such that the restriction of $\omega$ to
  $\xi$ is non--degenerate, $J$ preserves $\xi$, and the
  restriction of $J$ to $\xi$ is compatible with $\omega$.  We
  further assume that $\xi$ is tangent to $\Sigma$ near
  $(\ba\cup\bb)\times[0,1]$ and near $\Sigma\times\left(\bdy[0,1]\right)$.
\end{enumerate}
This still guarantees that $J$ is symmetric and adjusted to
$\omega$.

By an \emph{intersection point} we mean a set of $g$ distinct points 
$\vec{x}=\{x_1,\ldots,x_g\}$
in $\ba\cap \bb$ such that exactly one $x_i$ lies on each
$\alpha_j$ and exactly one $x_i$ lies on each
$\beta_k$.  (This corresponds to an
  intersection point of the $\alpha$-- and $\beta$--tori
  in~\cite{OS1}.)

Observe that the characteristic foliation on $\Sigma\times[0,1]$
induced by $\omega$ has leaves $\{p\}\times[0,1]\times\{t\}$.  So,an intersection point
$\vec{x}$ specifies a $g$--tuple of distinct ``Reeb chords'' (with respect
to the characteristic foliation on $\Sigma\times[0,1]$ induced by
$\omega$) in $\Sigma\times[0,1]$ with boundaries on $\ba\times\{1\}\cup
\bb\times\{0\}$.  (The collection of Reeb chords is just $\{x_i\}\times[0,1]$.)  We will call a $g$--tuple of Reeb chords at $\pm\infty$ specified by an intersection point an \emph{I--chord collection}.  (I stands for ``intersection.'')  We will abuse notation and also use $\vx$ to denote the I--chord collection specified by $\vx$.

Let $\cM$ denote the moduli space of Riemann surfaces $S$ with boundary,
$g$ ``negative''
punctures $\vec{p}=\{p_1,\cdots,p_g\}$ and $g$ ``positive''
punctures $\vec{q}=\{q_1,\cdots,q_g\}$, all on the boundary of
$S$, and such that $S$ is compact away from the
punctures.

For $J$ satisfying (\textbf{J1})--(\textbf{J5}), we will consider $J$--holomorphic maps
$u\co S\to W$ such that
\begin{itemize}\leftskip 3pt
\item[(\textbf{M0})] The source $S$ is smooth.
\item[(\textbf{M1})] $u(\bdy(S))\subset C_\alpha\cup C_\beta$.  
\item[(\textbf{M2})]  There are no components of $S$ on which
  $\pi_\bD\circ u$ is constant.
\item[(\textbf{M3})]  For each $i$,
  $u^{-1}(\alpha_i\times\{1\}\times\bR)$ and
  $u^{-1}(\beta_i\times\{0\}\times\bR)$ each consist
  of exactly one component of $\bdy
  S\setminus\{p_1,\cdots,p_g,q_1,\cdots,q_g\}$.
\item[(\textbf{M4})] $\lim_{w\to p_i}\pi_\bR\circ u(w)=-\infty$ and
    $\lim_{w\to q_i}\pi_\bR\circ u(w)=\infty$.
\item[(\textbf{M5})]  The \emph{energy} of $u$, as defined in~\cite[Section
  5.3]{Ya2}, is finite.  (For the moduli spaces defined later
    in the paper we shall always assume this technical condition is
    satisfied, but shall not usually state it.)
\item[(\textbf{M6})] $u$ is an embedding.
\end{itemize}

Note that condition (\textbf{M3}) implies that $\bdy
S\setminus\{p_1,\cdots,p_g,q_1,\cdots,q_g\}$ consist of exactly
$2g$ components, none of them compact.  Also note that we allow
holomorphic curves to be disconnected.

\begin{figure}
\centering
\begin{picture}(0,0)%
\includegraphics[scale=0.55]{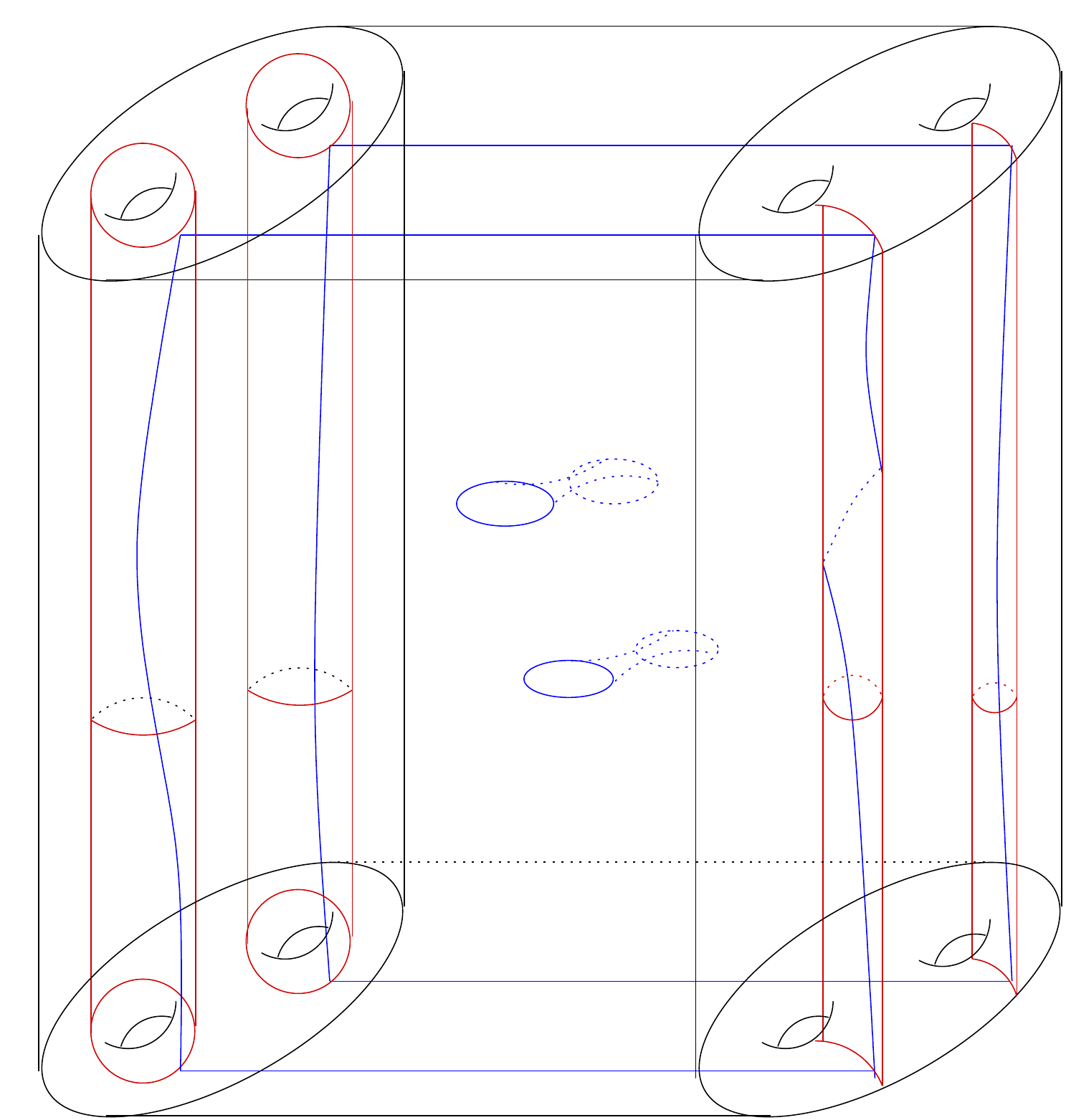}%
\end{picture}%
\setlength{\unitlength}{2170sp}%
\begingroup\makeatletter\ifx\SetFigFont\undefined%
\gdef\SetFigFont#1#2#3#4#5{%
  \reset@font\fontsize{#1}{#2pt}%
  \fontfamily{#3}\fontseries{#4}\fontshape{#5}%
  \selectfont}%
\fi\endgroup%
\begin{picture}(7112,7482)(226,-7449)
\put(3001,-7061){\makebox(0,0)[lb]{\smash{\SetFigFont{12}{14.4}{\rmdefault}{\mddefault}{\updefault}{\color[rgb]{0,0,0}{\(\color[named]{Blue}u(S)\)}}%
}}}
\put(977,-536){\makebox(0,0)[lb]{\smash{\SetFigFont{9}{10.8}{\rmdefault}{\mddefault}{\updefault}{\color[rgb]{0,0,0}\(\Sigma\)}%
}}}
\put(3926,-27){\makebox(0,0)[lb]{\smash{\SetFigFont{9}{10.8}{\rmdefault}{\mddefault}{\updefault}{\color[rgb]{0,0,0}\([0,1]\)}%
}}}
\put(176,-3386){\makebox(0,0)[lb]{\smash{\SetFigFont{9}{10.8}{\rmdefault}{\mddefault}{\updefault}{\color[rgb]{0,0,0}\(\mathbb{R}\)}%
}}}
\put(1526,-786){\makebox(0,0)[lb]{\smash{\SetFigFont{9}{10.8}{\rmdefault}{\mddefault}{\updefault}{\color[rgb]{0,0,0}{\(\color[named]{Red}\beta_2\)}}%
}}}
\put(525,-1385){\makebox(0,0)[lb]{\smash{\SetFigFont{9}{10.8}{\rmdefault}{\mddefault}{\updefault}{\color[rgb]{0,0,0}{\(\color[named]{Red}\beta_1\)}}%
}}}
\put(5977,-1386){\makebox(0,0)[lb]{\smash{\SetFigFont{9}{10.8}{\rmdefault}{\mddefault}{\updefault}{\color[rgb]{0,0,0}{\(\color[named]{Red}\alpha_2\)}}%
}}}
\put(6876,-785){\makebox(0,0)[lb]{\smash{\SetFigFont{9}{10.8}{\rmdefault}{\mddefault}{\updefault}{\color[rgb]{0,0,0}{\(\color[named]{Red}\alpha_1\)}}%
}}}
\end{picture}
\caption{A curve in $W$ we might consider.  Note that our curves can also be disconnected.}
\label{Figure:PrettyPicture}
\end{figure}
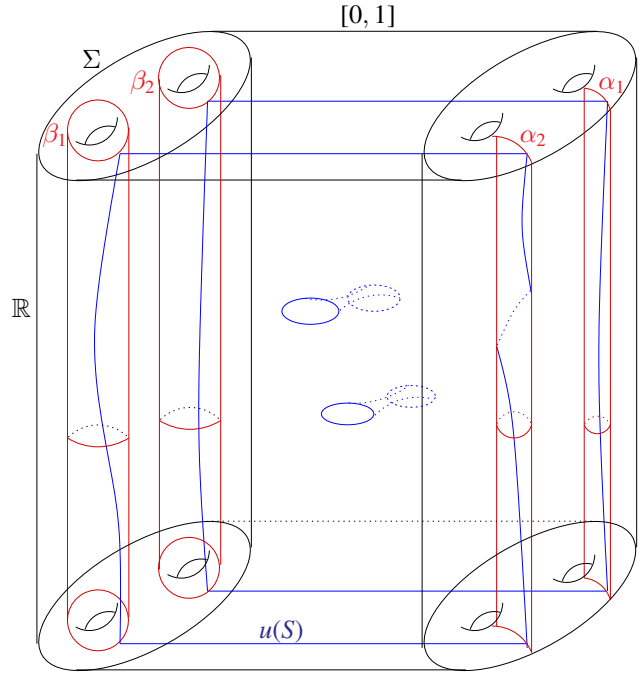

It follows from~\cite[Proposition 5.8]{Ya2} that near each negative puncture
(respectively positive puncture), a
holomorphic curve satisfying (\textbf{M0})--(\textbf{M6})
converges exponentially (in $t$) to an I--chord collection $\vx$ (respectively
$\vy$) at
$-\infty$ (respectively $\infty$).  We say the holomorphic curve
connects $\vx$ to $\vy$.  It follows from this asymptotic
convergence to Reeb chords that $\pi_\bD\circ u$ is a $g$--fold
branched covering map. 

Consider the space $\overline{W}=\Sigma\times[0,1]\times[-1,1]$
as a compactification of $W$.  Let $\overline{C_\alpha}$,
$\overline{C_\beta}$ denote the closures of the images of
$C_\alpha$ and $C_\beta$ in $\overline{W}$.  Let
$\overline{S}$ denote the surface obtained by blowing up $S$ at
the punctures.  Then, the asymptotic
convergence to Reeb orbits mentioned earlier implies that $u$ can be
extended to a continuous map
$\overline{u}\co \overline{S}\to\overline{W}$.  (Compare, for
example,~\cite[Proposition 6.2]{Ya2}.)

Let $\pi_2(\vec{x},\vec{y})$ denote the set of homology classes of
continuous maps $(S,\bdy S)\to
(W,C_\alpha\cup C_\beta)$
which converge to $\vec{x}$ (respectively 
$\vec{y}$) near the negative (respectively positive) punctures of
$S$.
  That is, two such maps are equivalent if they induce the same
  element in
$H_2(\overline{W},\overline{C_\alpha}\cup\overline{C_\beta}\cup
\left(\{x_i\}\times[0,1]\times\{-1\}\right)\cup\left(\{y_i\}\times[0,1]\times\{1\}\right))$.  (The notation is chosen to be consistent with \cite{OS1}, where the notation $\pi_2$ makes sense.)

Each holomorphic curve
connecting $\vx$ to $\vy$ represents an element of
$\pi_2(x,y)$.  For $A\in\pi_2(x,y)$, we denote by
$\cM^A$ the space of holomorphic curves connecting
$\vec{x}$ and $\vec{y}$ in the homology class $A$.  (We always
mod out by automorphism of the source $S$.)   Since we are considering cylindrical complex structures, $\bR$ acts on $\cM^A$ by translation.  Let $\hcM^A=\cM^A/\bR$.  We denote by
$\overline{\hcM^A}$ the compactification, as in \cite[Section 7]{Ya2},
of $\hcM^A$.

Given a homology class $A\in\pi_2(\vx,\vy)$, let $n_\Fz(A)$ denote
the intersection number of $A$ with $\{\Fz\}\times[0,1]\times\bR$.
Define $n_{\Fz_i}(A)$ similarly.  If $u$ is a curve in the homology
class $A$ we will sometimes write $n_\Fz(u)$ or $n_{\Fz_i}(u)$ for
$n_\Fz(A)$ or $n_{\Fz_i}(A)$.  We say that a homology class $A$
is \emph{positive} if $n_{\Fz_i}(A)\geq0$ for all $i$. Notice that if
$A$ has a holomorphic 
representative (with respect to any complex structure satisfying
(\textbf{J2})) then $A$ is positive; this is the
\emph{positivity of domains} mentioned twelve paragraphs above.  We shall let
$\hat{\pi}_2(\vx,\vy)=\{A\in\pi_2(\vx,\vy)|n_\Fz(A)=0\}.$  Elements of
$\hat{\pi}_2(\vx,\vx)$ are called \emph{periodic classes}.

\medskip
\textbf{Remark}\qua In fact, even without (\textbf{J2}), positivity of domains would still hold by Micallef and White \cite[Theorem 7.1]{White}. (See also \fullref{Lemma:White}.) On the other hand, by requiring (\textbf{J2}), which is easy to obtain, we can avoid invoking here this hard analytic result.

\medskip
Given a homology class $A$, we
define the \emph{domain of $A$} to be
the formal linear combination $\sum n_{\Fz_i}(A) D_i$.  If $u$
represents $A$ then we define the domain of $u$ to be the domain
of $A$.  The domains of periodic classes are called \emph{periodic domains}.

As in \cite{OS1}, concatenation makes $\pi_2(\vx,\vy)$ into a
$\pi_2(\vx,\vx)$--torseur.  We shall sometimes write concatenation
with a $+$ and sometimes with a $*$, depending on whether we are
thinking of domains or maps.

\section{Homotopy preliminaries}
\label{section:homotopy}
These issues are substantially simplified from \cite{OS1} because we
need only deal with homology, not homotopy.  This is reasonable:  by
analogy to the Dold--Thom theorem, the low--dimensional homotopy
theory of $\Sym^g(\Sigma)$ should agree with the homology theory of $\Sigma$.

Given an intersection point $\vx$, observe that projection from
$W$ gives rise to an isomorphism from $\pi_2(\vx,\vx)$ to
$H_2(\Sigma_g\times[0,1],\ba\times\{1\}\cup\bb\times\{0\})$.
Given intersection points $\vx$, $\vy$, either $\pi_2(\vx,\vy)$
is empty or $\pi_2(\vx,\vy)\cong
H_2(\Sigma_g\times[0,1],\ba\times\{1\}\cup\bb\times\{0\})$.  The
isomorphism is not canonical; it is given by fixing an element
of $\pi_2(\vx,\vy)$ and then subtracting the homology class it
represents from all other elements of $\pi_2(\vx,\vy)$.  We
calculate $ H_2(\Sigma_g\times[0,1],\ba\times\{1\}\cup\bb\times\{0\})$:
\begin{Lem}\label{homlem}{\rm(Compare \cite[Proposition 2.15]{OS1})}\qua
  There is a natural short exact sequence
\begin{displaymath}
0\to\bZ\to H_2(\Sigma\times[0,1],\ba\times\{1\}\cup\bb\times\{0\})
\to H_2(Y)\to 0.
\end{displaymath}
The choice of basepoint $\Fz$ gives a splitting
$n_\Fz\co H_2(\Sigma\times[0,1],\ba\times\{1\}\cup\bb\times\{0\})\to\bZ$ of
this sequence.
\end{Lem}

\proof
The long exact sequence for the pair
$(\Sigma\times[0,1],\ba\times\{1\}\cup\bb\times\{0\})$ gives
\begin{displaymath}
0\to H_2(\Sigma\times[0,1])\to H_2(\Sigma\times[0,1],
\ba\times\{1\} \cup \bb\times\{0\}) \to
H_1(\ba\times\{1\}\cup \bb\times\{0\}).
\end{displaymath}
The image of the last map is isomorphic to $H_1(\ba)\cap
H_1(\bb)$, viewed as a submodule of $H_1(\Sigma)$.

Let $Y=U_1\cup_\Sigma U_2$ be the Heegaard splitting.  The
Mayer--Vietoris sequence gives
\begin{displaymath}
H_2(U_1)\oplus H_2(U_2)\to H_2(Y)\to H_1(\Sigma) \to
H_1(U_1)\oplus H_1(U_2).
\end{displaymath}
Here, the kernel of the last map is $H_1(\ba)\cap H_1(\bb)$.
The groups $H_2(U_1)$ and $H_2(U_2)$ are both
trivial, so $H_2(Y)\cong H_1(\ba)\cap H_1(\bb)$.  Combining
this with the first sequence and using the fact that
$H_2(\Sigma\times[0,1])\cong \bZ$ gives the first part of the
claim.  With $n_\Fz$ defined as in \fullref{Section:Notation}
the second part of the claim is obvious.
\endproof

If we identify $\Sigma\times[0,1]$ with
$f^{-1}[3/2-\epsilon,3/2+\epsilon]$ for some self--indexing Morse
function $f$ on $Y$ then the map
$H_2(\Sigma\times[0,1],\ba\times\{1\}\cup\bb\times\{0\})\to H_2(Y)$ is simply
given by ``capping off'' a cycle with the ascending / descending disks
from the index 1 and 2 critical points of $f$.
Also, notice that a homology class $A\in\pi_2(\vx,\vy)$
specifies and is specified by its
domain.  (The domain need not, however, specify uniquely the intersection
points $\vx$ and $\vy$ which it connects.)

Following~\cite[Section 2.6]{OS1}, we observe that a choice of basepoint $\Fz$ and
intersection point $\vec{x}$ specify a $\Spin^\bC$--structure $\Ss$ 
on $Y$ as follows.  Choose a metric $\langle\cdot,\cdot\rangle$
and a self--indexing Morse
function $f$ which specify the Heegaard diagram $\cH$.  Then
$\vx$ specifies a $g$--tuple of flows of $\nabla f$ from the index $1$
critical points of $f$ to the index $2$ critical points of $f$.
The point $z$ lies on a flow from the index $0$ critical point of
$f$ to the index $3$ critical point of $f$.
Choose small ball neighborhoods of (the closure of) each of these flow
lines.  Call the union of these neighborhoods $B$.  Then, in the
complement of $B$, $\nabla f$ is nonvanishing.  One can extend $\nabla f$ to a nonvanishing vector
field $v$ on all of $Y$.  The vector field $v$ reduces the
structure group of $TY$ from $SO(3)$ to $SO(2)\oplus SO(1)\cong
U(1)\oplus 1\subset U(2)=\Spin^\bC(3)$, and thus determines a
$\Spin^\bC$--structure on $Y$.  We have, thus, defined a map
$\Ss_\Fz$ from the set of intersection points in $\cH$ to
the set of $\Spin^\bC$--structures on $Y$.

It is clear that the $\Spin^\bC$--structure $\Ss_\Fz(\vx)$ is
independent of the metric and particular Morse function used to define
it.

Given a $\Spin^\bC$--structure $\Ss$ on $Y$, we shall often
suppress $\Fz$ and write $\vx\in\Ss$ to mean $\Ss_\Fz(\vx)=\Ss$.

Note that by the previous construction, any nonvanishing vector field on
a 3--manifold $Y$ gives rise to a $\Spin^\bC$--structure.  It is
not hard to show that two nonvanishing vector fields give rise to the
same $\Spin^\bC$--structure if and only if they are
\emph{homologous}, ie, homotopic through
nonvanishing vector fields in the complement of some $3$--ball;
see~\cite{Turaev}.  We will use the analogous construction in the case
of $4$--manifolds in \fullref{Subsubsection:TriangleSpinc}.

Our reason for introducing $\Spin^\bC$--structures will become clear
in a moment.  First, one more definition.  Fix a pair of intersection
points $\vx$ and $\vy$, as well as a Morse function $f$ and
Riemannian metric $\langle\cdot,\cdot\rangle$ which realize the
Heegaard diagram.  This data specifies a homology class
$\epsilon(\vx,\vy)$ as follows.  Regard each of $\vx$ and $\vy$
as (the closure of) a $g$--tuple of gradient flow trajectories in
$Y$ from the $g$ index $1$ critical points to the $g$ index
$2$ critical points.  Then, $\vx-\vy$ is a $1$--cycle in $Y$.
We define $\epsilon(\vx,\vy)$ to be the homology class in $H_1(Y)$ of the
$1$--cycle $\vx-\vy$.

The element $\epsilon(\vx,\vy)$ can be calculated entirely in
$\cH$ by the following equivalent definition.  Let $\gamma_\alpha$
(respectively $\gamma_\beta$)
be a 1--cycle in $\ba$ (respectively $\bb$) such that
$\bdy\gamma_\alpha=\bdy\gamma_\beta=\vx-\vy$.  Then $\gamma_\alpha-\gamma_\beta$ is
a 1--cycle in $\Sigma$.  Define $\epsilon(\vx,\vy)$ to be the image of
$\gamma_\alpha-\gamma_\beta$ under the map 
$$H_1(\Sigma)\to
\frac{H_1(\Sigma)}{H_1(\ba)+H_1(\bb)}\cong H_1(Y).$$  

The equivalence of the two definitions is easy:  in the notation used
just above, $\vx+\gamma_\alpha-\gamma_\beta$ is homologous, rel
endpoints, to $\vy$.  It is obvious that the second definition is
independent of the choices of $\gamma_\alpha$ and $\gamma_\beta$.

The following lemma justifies our introduction of $\epsilon$ and of
$\Spin^\bC$--structures:
\begin{Lem} \label{Lemma:Epsilon}{\rm(Compare~\cite[Proposition 2.15, 
Lemma 2.19]{OS1})}\qua
Given a pointed Heegaard diagram $\cH$ and intersection
  points $\vx$ and $\vy$, the following are equivalent:
\begin{enumerate}
\item \label{item1}$\pi_2(\vec{x},\vec{y})$ is nonempty
\item \label{item2}$\epsilon(\vx,\vy)=0$
\item \label{item3}$\Ss_\Fz(\vx)=\Ss_\Fz(\vy).$
\end{enumerate}
\end{Lem}

\proof

$(1)\Rightarrow(2)$\quad
Let $A\in\pi_2(\vx,\vy)$.  View $A$ as a domain in $\Sigma$,
ie, a chain in $\Sigma$.  Then we can use $\bdy A$ to define
$\epsilon(\vx,\vy)$, which is thus zero in homology.

$(2)\Rightarrow(1)$\quad
Suppose that $\epsilon(\vx,\vy)=0$.  Then, using the same notation
as just before the lemma, for an appropriate choice of $\gamma_\alpha$ and
$\gamma_\beta$, $\gamma_\alpha-\gamma_\beta$ is null--homologous in
$H_1(\Sigma)$.  We can assume that $\gamma_\alpha$ and
$\gamma_\beta$ are cellular $1$--chains in the cellulation of
$\Sigma$ induced by the Heegaard diagram.  Then, there is a cellular
$2$--chain $A$ with boundary $\gamma_\alpha-\gamma_\beta$, and
$A$ is the domain of an element of $\pi_2(\vx,\vy)$.

$(2)\Leftrightarrow(3)$\quad
Let $v_\vx$ and $v_\vy$ denote the vector fields used to define
$\Ss_\Fz(\vx)$ and $\Ss_\Fz(\vy)$, respectively.  Let $v_\vy=Av_\vx$
where $A\co Y\to SO(3)$.  Let $Fr(v_\vx^\perp)$ and
$Fr(v_\vy^\perp)$ denote the principal $SO(2)=U(1)$--bundles of
frames of $v_\vx^\perp$ and $v_\vy^\perp$.  Then the principal
$\Spin^\bC$ bundles induced by $v_\vx$ and $v_\vy$ are
$\Ss_z(\vx)\co Fr(v_\vx^\perp)\times_{U(1)}U(2)\to Fr(TY)$ and
$\Ss_z(\vy)\co Fr(v_\vy^\perp)\times_{U(1)}U(2)\to Fr(TY)$.

Note that $\Ss_\Fz(\vx)$ and $\Ss_\Fz(\vy)$ are equivalent if and only
if $A$ is homotopic to a map $Y\to SO(2)$.  So, the two
$\Spin^\bC$--structures are equivalent if and only if the composition
$h\co Y\to SO(3)\to SO(3)/SO(2)=S^2$ is null homotopic.

Now, homotopy classes of maps from a $3$--manifold $Y$ to $S^2$
correspond to elements of $H^2(Y)$.  The Poincar\'e dual to such a map is the
homology class of the preimage of a regular value in $S^2$.

For a generic choice of the two Morse functions and metrics used to
define them, the flows through $\vx$
and $\vy$ glue together to a disjoint collection of circles
$\gamma$.  Let $\gamma'$ be a smoothing of $\gamma$.
The map $h$ is homotopic to a Thom collapse map of a neighborhood of
$\gamma'$.  It follows that the preimage of a regular value is
homologous to $\gamma$.

So, $\Ss_\Fz(\vx)=\Ss_\Fz(\vy)$ if and only if $\gamma$ is
null--homologous.  But $\gamma$ is a cycle defining
$\epsilon(\vx,\vy)$, so the result follows.
\endproof

Note that the previous proof in fact shows that the map $\Ss_\Fz$ from
intersection points to $\Spin^\bC$--structures is a
map of $H^2(Y)=H_1(Y)$--torseurs.

The following result (part of~\cite[Lemma 2.19]{OS1}) is nice to know,
but will not be used explicitly in this paper.  The reader can
imitate the proof of the previous proposition to prove it, or see~\cite{OS1}.
\begin{Lem}\label{OS:Lemma2.19}
Let $\vx$ be an intersection point of a Heegaard diagram
  $\cH$, and $\Fz_1$, $\Fz_2$ two different basepoints for $\cH$.
  Suppose that $\Fz_1$ can be joined to $\Fz_2$ by a path $\Fz_t$
  disjoint from the $\beta$ circles and such that $\#(\Fz_t\cap
  \alpha_i)=\delta_{i,j}$ (Kronecker delta).  Let $\gamma$ be a
  loop in $\Sigma$ such that $\gamma\cdot\alpha_i=\delta_{i,j}$.  Then,
  $\Ss_{\Fz_2}(\vx)-\Ss_{\Fz_1}(\vx)=PD(\gamma)$, the Poincar\'e
  dual to $\gamma$.
\end{Lem}

\section{Transversality}
\label{Section:Transversality}
We need to check that we can achieve transversality for the
generalized Cauchy--Riemann equations within the class of almost
complex structures satisfying (\textbf{J1})--(\textbf{J5}).  The
argument is relatively standard, and is almost
the same as the one found in \cite[Chapter 3]{MS2}.  This section is
somewhat technical, and the reader might want to skip most of it on a
first reading.

Before proving our transversality result we need a few lemmas about
the geometry of holomorphic curves in $W$.

\begin{Lem}\label{Lemma:White}Let $\pi\co E\to B$ be a smooth fiber bundle, with
  $\dim(E)=4$, $\dim(B)=2$.  Let $J$ be an almost complex
  structure on $E$ with respect to which the fibers are
  holomorphic.  Let $u\co S\to E$ be a $J$--holomorphic map, $S$
  connected, with $\pi\circ u$ not constant.  Let $p\in S$ be a
  critical point of $\pi\circ u$, $q=\pi\circ u(p)$.  Then there
  are neighborhoods $U\ni p$ and $V\ni q$, and $C^2$
  coordinate charts $z\co U\to\bC$, $w\co V\to \bC$ such that $w\circ
  (\pi\circ u)(z)=z^k$, for some $k>0$.
\end{Lem}
\proof
This follows immediately from~\cite[Theorem 7.1]{White} applied to the
intersection of $u$ with the fiber of $\pi$ over $q$.
\endproof

\begin{Cor}Let $\pi\co E\to B$ be a smooth fiber bundle, with
  $\dim(E)=4$, $\dim(B)=2$.  Let $J$ be an almost complex
  structure on $E$ with respect to which the fibers are
  holomorphic.  Let $u\co S\to E$ be a $J$--holomorphic map, $S$
  connected, with $\pi\circ u$ not constant.  Then the Riemann--Hurwitz
  formula applies to $\pi\circ u$.  That is, if $S$ is closed then $$\chi(S)=\chi(\pi\circ
  u(S))-\sum_{p\in S}\left(e_{\pi\circ u}(p)-1\right)$$ where
  $e_{\pi\circ u}(p)$ is the ramification index of $p$. If $S$ has boundary and punctures then the same formula holds with Euler measure in place of Euler characteristic.
\label{HurwitzHolds}
\end{Cor}
(The \emph{Euler measure} of a surface $S$ with boundary and punctures is $1/2\pi$ times the integral over $S$ of the curvature of a metric on $S$ for which $\bdy S$ is geodesic and the punctures of $S$ are right angles.  See \fullref{Section:Index} for further discussion of Euler measure.)

\begin{Lem} Let $\pi\co \Sigma\times[0,1]\times\bR\to
  \Sigma\times[0,1]$ denote projection.  Let $u$ be a
  holomorphic curve in $\Sigma\times[0,1]\times\bR$
  (with respect to some almost complex structure satisfying
  {\rm(\textbf{J1})--(\textbf{J5})}).  Let
  $S'$ be a component of $S$ on which $u$ is not a trivial disk
  and $\pi_\bD\circ u$ is not constant.  Then
  there is a nonempty, open subset $U$ of $S'$ on which $\pi\circ u$
  is injective and $\pi\circ u(U)\cap \pi\circ
  u(S\sm U)=\emptyset$.  Further, we can require that $u(U)$ be disjoint
  from $\overline{U_{\Fz_i}}$ and that $\pi_\Sigma\circ du$ and 
  $\pi_\bD\circ du$ be nonsingular on $U$.
\end{Lem}
(By a \emph{trivial disk} we mean  a component of $S$
  mapped diffeomorphically by $u$ to $\{x\}\times[0,1]\times\bR$ for some
  $x\in\Sigma$.)  

\proof
Let $x$ be such that $u|_{S'}$ is asymptotic to the Reeb chord
$\{x\}\times[0,1]$ at infinity.  Let $\overline{S}$ denote the
surface obtained by blowing--up $S$ at its punctures.  As discussed
earlier, we can extend
$u$ to a continuous map
$\overline{S}\to\overline{W}=\Sigma\times[0,1]\times[-1,1]$.  Let
$\overline{\pi}\co \overline{W}\to \Sigma\times[0,1]$ denote
projection.
Let $E$ denote the set of points $(x,s)\in
\Sigma\times[0,1]$ such that either $(\overline{\pi}\circ \overline{u})^{-1}(x,s)$ has
cardinality larger than 1 or contains the image of a critical point
 of $\pi_\Sigma\circ u$ or $\pi_\bD\circ u$.  Then $E$ is closed.

By the preceding corollary, there are only
finitely many critical points of $\pi_\Sigma\circ u$ or
$\pi_\bD\circ u$.  Further, ``positivity of intersections'' (eg, \cite[Theorem 7.1]{White}), applied to
$u$ and $\{x\}\times[0,1]\times\bR$, implies that
there are only finitely many points in $\pi\circ
u^{-1}(\{x\}\times[0,1])$.  So, there are only finitely many points
in $E\cap\{x\}\times[0,1]$.

However, $\{x\}\times[0,1]$ is contained in the image of
  $\overline{\pi}\circ \overline{u}$.  Choose $s\in[0,1]$ such that
  $(x,s)\in \{x\}\times[0,1]\setminus E$.  Let $V$ be an open
  neighborhood of $(x,s)$ disjoint from $E$.  Then $(\pi\circ
  u)^{-1}(V)$ has the desired properties.
\endproof

To prove transversality we need to specify precisely the spaces under
consideration.  Fix $p> 
2$, $k\geq 0$ ($k\in\bZ$) and $d>0.$ 
\begin{Def}For a Riemannian manifold $(M,\bdy M)$, a function $f:M\to\bR$ lies in $L^p_k(M)$ if $f$ has $k$ weak derivatives in $L^p$.  The $L^p_k$--norm of $f$ is 
$$\|f\|_{L^p_k}=\|f\|_{L^p}+\|f'\|_{L^p}+\cdots+\|f^{(k)}\|_{L^p}.$$  A function $f:M\to \bR^n$ lies in $L^p_k$ if each coordinate of $f$ does, and its $L^p_k$--norm is the sum of the $L^p_k$--norms of its coordinate functions.
\end{Def}

Fix a map $u\co (S,\bdy S)\to (W,C_\alpha\cup C_\beta)$.  Fix a Riemannian metric on $S$; the particular choice is unimportant.  Let $\{p_i^-\}$ denote the negative punctures of $S$ and $\{p_i^+\}$ the positive punctures.  Suppose $u$ is asymptotic to $x_i^\pm\times[0,1]$ at $p_i^\pm$.  Identify a neighborhood $U_i^-$ of each $p_i^-$ with $[0,1]\times (-\infty,0]]$ and a neighborhood $U_i^+$ of each $p_i^+$ with $[0,1]\times [0,\infty)$.  Let $(\sigma_i^\pm,\tau_i^\pm)$ denote the coordinates near $p_i^\pm$ induced by this identification.  Fix also a smooth embedding of $\Sigma$ in $\bR^{N-2}$ for some $N$.  This induces an embedding of $W=\Sigma\times[0,1]\times\bR$ in $\bR^N$ in an obvious way.  For the following definition we identify $W$ with its image in $\bR^N$.
\begin{Def}\label{Def:WPDK}We say that $u$ lies in $W^{p,d}_k\left((S,\bdy S);(W,C_\alpha\cup C_\beta)\right)$ if for some choice of constants $\{t_{0,i}^\pm\}\in\bR$,
\begin{itemize}

\item the restriction of $u$ to $S\setminus \left(U_1\cup\cdots\cup
U_g\cup V_1\cup\cdots\cup V_g\right)$ lies in $L^p_k$ (as a function
to $\bR^N$) and

\item on each $U_i^\pm$ the functions $e^{d|\tau_i^\pm|}\left(s\circ
u(\sigma_i^\pm,\tau_i^\pm)-\sigma_i^\pm\right)$,\newline
$e^{d|\tau_i^\pm|}\left(t\circ
u(\sigma_i^\pm,\tau_i^\pm)-\tau_i^\pm-t_{0,i}^\pm\right)$ and
$e^{d|\tau_i^\pm|}\left(u(\sigma_i^\pm,\tau_i^\pm)-x_i^\pm\right)$\newline
from $[0,\infty)\times\bR$ or $(-\infty,0]$ to $\bR$ lie in $L^p_k$.
\end{itemize}
\end{Def}
For $d$ small enough, all finite energy holomorphic curves (in the sense of~\cite[Section 5.3]{Ya2})  in $(W,C_\alpha\cup C_\beta)$ lie in $W^{p,d}_k$; see for instance~\cite[Chapter 3]{bourgeois}, particularly Propositions 3.5 and~3.6.  Conversely, all maps in $W^{p,d}_k$ have finite energy.

Choose a homology
class $A$ of maps to $W$ and a surface $S$.  Let  
$\cX^{p,d}_k$ denote the collection of maps $u\in W^{k,p}_\delta\left((S,\bdy
S); (W,C_\alpha\cup C_\beta)\right)$ in class
$A$.  

\begin{Def}
Let $E$ be a Riemannian vector bundle over $S$.  Let $f$ be a section of $E$.  Then the $L^{p,d}_k$--norm of $f$ is
$$
\|f\|_{L^{p,d}_k}=\|f|_{S\setminus\left(U_1^-\cup\cdots\cup U_g^+\right)}\|_{L^p_k} + \sum_{i=1}^g\left(\|e^{d|\tau_i^+|}f|_{U_i^+}\|_{L^p_k}+\|e^{d|\tau_i^-|}f|_{U_i^-}\|_{L^p_k}\right)
$$
Let $L^{p,d}_k(E)$ denote the Banach space of all sections of $E$ with finite $L^{p,d}_k$--norm.
\end{Def}

Note that the tangent space at $u$ to $\cX^{p,d}_k$ is
$\bR^{2g}\oplus L^{p,d}_k\left(u^*TW,\bdy\right)$ where\newline $L^{p,d}_k\left(u^*TW,\bdy\right)$ is the subspace of $L^{p,d}_k\left(u^*TW\right)$
of sections which lie in $u^*T(C_\alpha\cup C_\beta)$ over $\bdy S$.  The
$\bR^{2g}$ factor corresponds to varying the $2g$ constants
$t_{0,i}^\pm$ in \fullref{Def:WPDK}.  Choosing $2g$ smooth
vector fields $v_i^\pm$ given by $\frac{\partial}{\partial t}$
on a neighborhood of $p_i^\pm$ and zero near the other punctures
$p_j^\pm$, we can include the $\bR^{2g}$ into
$\Gamma\left(u^*TW\right)$ as $\mathrm{Span}\left(\{v_i^\pm\}\right)$.

  Let $\cJ^\ell$ denote the space of $C^\ell$ 
almost complex structures on 
$W$ which satisfy (\textbf{J1})--(\textbf{J5}).
  Let $\mathcal{J}^\ell(S)$ denote the space of 
$C^\ell$ almost complex structures on $S$.  Let 
$\cM^\ell=\{(u,j,J_s)\in \cX^{k,p}_\delta\times \cJ^\ell(S)\times 
\cJ^\ell | \dbar_{jJ_s}u=0\}$.

Let $\End(TS,j)$ denote the bundle whose fiber at $p\in S$ is the space 
of linear $Y\co T_pS\to T_pS$ such that $Yj+jY=0$.  Then the tangent space 
at $j$ to $\cJ^{\ell}(S)$ is the space of $C^\ell$ sections 
of $\End(TS,j)$.  Similarly, let $\End(TW,J_s)$ denote the space of 
$C^{\ell}$ paths $Y_s$ of linear maps $T\Sigma\to T\Sigma$ such that 
$Y_sJ_s+J_sY_s=0$.

By convention, if we omit the superscripts $k$, $\ell$, and $p$ then 
we are referring to smooth objects.

By an \emph{annoying curve} we mean a curve $u\co S\to W$ such that
there is a nonempty open subset of $S$ on which $\pi_\bD\circ u$ is
constant.

\begin{Prop} For $\ell\geq 1$ the space $\cM^\ell$ is a smooth Banach 
manifold away from annoying curves.
\end{Prop}
\proof
(This proof is a slight modification of~\cite[Proposition 3.4.1, page 
34]{MS2}.  The reader is referred there for a less terse exposition.)

Let $\cE^{p}_{k-1}$ be the bundle over $\cX^{p,d}_k\times 
\cJ^\ell(S)\times \cJ^\ell$ whose fiber over a point $(u,j,J_s)$ is 
$L^{p,d}_{k-1}(\Lambda^{0,1}T^*S\otimes_{j,J_s} u^*TW)$.

  We view 
$\dbar$ as a section of $\cE^{p}_{k-1}$, and want to show that it is 
transverse to the zero section.  At the zero section, the tangent space 
to $\cE$ splits as 
$$T(\cX^{p,d}_k\times \cJ^\ell(S)\times 
\cJ^\ell)\oplus 
L^{p,d}_{k-1}(\Lambda^{0,1}T^*S\otimes_{j,J_s}u^*TW).$$  Let 
{\setlength\arraycolsep{-2pt}
\begin{align*}
D\dbar (u,j,J_s)\co \bR^{2g}\times L^{p,d}_k
(u^*TW,\bdy)&\times C^\ell(\End(TS,j))\times
C^\ell(\End(TW,J_s)) \\ &\to\hspace{1ex}
L^{p,d}_{k-1}(\Lambda^{0,1}T^*S\otimes_{j,J_s}u^*TW)
\end{align*}}denote projection of the differential of $\dbar$ onto the vertical 
component of the tangent space to $\cE$ at a zero $(u,j,J_s)$ of
$\dbar$.  We must show that $D\dbar$ is surjective.  

The restriction of $D\dbar$ to any trivial disk is
surjective by~\cite[Theorem 2]{HLS}.  So, for the rest of the proof we
consider only the components of $S$ which are not trivial disks.
For these components we will, in
fact, show that the restriction of $D\dbar$ to $0\times
L^{p,d}_k(u^*TW,\bdy)\times C^\ell(\End(TS,j))\times
C^\ell(\End(TW,J_s))$ is surjective, and will focus on this
restriction from now on.

The differential $D\dbar$ is given by 
$$D\dbar(u,j,J_s)(\xi,Y,Y_s)=D_u\xi + \frac{1}{2}Y_s(u)\circ du\circ j + 
\frac{1}{2} J_s\circ du\circ Y$$
where $D_u\xi$ denotes the differential holding $j$ and $J_s$ 
fixed.  

The operator $D\dbar(u,j,J_s)$ has closed range since $D_u$ is 
Fredholm, and we only need to show that its range is dense.  First, take 
$k=1$.  If the range is not dense then there exists $\eta\in 
L^{q,d}(\Lambda^{0,1}T^*S\otimes_{j,J_s}u^*TW)$ (where $1/p+1/q=1$) 
which annihilates the range of $D\dbar$.  So, for any choice of 
$(\xi,Y,Y_s)$, we have
\begin{eqnarray}
\int_{S}\left\langle \eta,D_u\xi\right\rangle&=&0\\
\int_{S}\left\langle \eta, Y_s\circ du\circ j\right\rangle&=& 0\\
\int_{S}\left\langle \eta, J_s\circ du\circ Y\right\rangle &=& 0
\end{eqnarray}

The first equation says that $\eta$ is a weak solution of 
$D_u^*\eta=0$.  So, by elliptic regularity, $\eta\in 
L^{r,d}_{\ell+1}$ for any $r>0$.  Further, it suffices to show that 
$\eta$ vanishes on some open set to show that $\eta$ vanishes 
identically.

Let $U$ be as in the previous lemma and $z_0\in U$.  Choose 
coordinates $(x,y)$ on $S$ near $z_0$ with respect to which $j$ is 
represented by the matrix
$
\left(\begin{array}{cc}
0 & -1\\
1 & 0
\end{array}\right).
$
Choose coordinates $(x_1,y_1,x_2,y_2)$ near $u(z_0)$ preserving the 
splitting $TW=T\Sigma\oplus T\bD$ and with respect to which the complex 
structure $J_s$ on $W$ has the form
$$
\left(\begin{array}{cccc}
0 & -1 & 0 & 0\\
1 & 0 & 0 & 0\\
0 & 0 & 0 & -1\\
0 & 0 & 1 & 0\\
\end{array}\right).
$$
The map $\eta$ has the form
$
\left(\begin{array}{cc}
a & b\\
b & -a\\
c & d\\
d & -c
\end{array}\right).
$
If $Y=\left(\begin{array}{cc}
\alpha & \beta \\
\beta & -\alpha
\end{array}\right)
$
and 
$Y_s =\left(\begin{array}{cc}
\gamma & \delta\\
\delta & -\gamma
\end{array}\right).$
Let $u_1=x_1\circ u$, $u_2=x_2\circ u$.  Then, we have
$$J_s\circ du\circ Y =\left(\begin{array}{cc}
-\alpha\frac{\partial u_1}{\partial y} + \beta\frac{\partial u_1}{\partial 
x} &
  -\beta\frac{\partial u_1}{\partial y} - \alpha\frac{\partial 
u_1}{\partial x}\\
-\alpha\frac{\partial u_1}{\partial x} - \beta\frac{\partial u_1}{\partial 
y} &
  -\beta\frac{\partial u_1}{\partial x} + \alpha\frac{\partial 
u_1}{\partial y} \\
-\alpha\frac{\partial u_2}{\partial y} + \beta\frac{\partial u_2}{\partial 
x} &
  -\beta\frac{\partial u_2}{\partial y} - \alpha\frac{\partial 
u_2}{\partial x}\\
-\alpha\frac{\partial u_2}{\partial x} - \beta\frac{\partial u_2}{\partial 
y} &
  -\beta\frac{\partial u_2}{\partial x} + \alpha\frac{\partial 
u_2}{\partial y} \\
\end{array}\right)
$$
and
$$Y_s\circ du \circ j = \left(\begin{array}{cc}
-\gamma\frac{\partial u_1}{\partial y}-\delta\frac{\partial u_1}{\partial 
x} &
   \gamma\frac{\partial u_1}{\partial x} - \delta\frac{\partial 
u_1}{\partial y}\\
-\delta\frac{\partial u_1}{\partial y}+\gamma\frac{\partial u_1}{\partial 
x } &
   \delta\frac{\partial u_1}{\partial x} + \gamma\frac{\partial 
u_1}{\partial y}\\
   0 & 0\\
   0 & 0
\end{array}\right).$$
By choosing $\gamma$ and $\delta$ appropriately we can force $a=b=0$ 
near $z_0$.  (This uses the injectivity established in the lemma and the 
nonvanishing of $\pi_\Sigma\circ du(z_0)$.)  Then, choosing $\alpha$ 
and $\beta$ appropriately we can force $c=d=0$ near $z_0$.  (This 
uses the nonvanishing of $\pi_\bD\circ du(z_0)$.)  This establishes the 
surjectivity of $D\dbar$ and hence the $k=1$ case.

For general $k$, suppose $\eta\in
L^{p,d}_{k-1}(\Lambda^{0,1}T^*S\otimes_{j,J_s}u^*TW).$  From the
$k=1$ case, choose a triple $\xi\in L^{p}_1(u^*TW)$, $Y\in
C^\ell(\End(TS,j))$, and $Y_s\in C^\ell(\End(TW,J_s))$
    such that $D\dbar(u,j,J_s)(\xi,Y,Y_s)=\eta$.  Then, elliptic
    regularity implies $\xi\in L^{p,d}_k$, so $D\dbar$ is
    surjective.  Since $D\dbar$ is Fredholm, it follows from the
    infinite--dimensional implicit function theorem that $\cM^\ell$
    is a Banach manifold.
\endproof

\begin{Prop} For a dense set $J_{reg}$ of $C^\infty$ paths 
of smooth complex structures on $\Sigma$, the moduli space of holomorphic 
curves satisfying {\rm(\textbf{M1})}, {\rm(\textbf{M2})}, {\rm(\textbf{M4})} and {\rm(\textbf{M5})}, and without multiply covered components, is a smooth manifold.
\end{Prop}
\proof
Observing that (\textbf{M2}) implies the absence of annoying curve components, this follows easily from the previous result.  The set $J_{reg}$ is 
exactly the set of regular values for the projection of $\cM$ onto 
$\cJ$.  For $\cJ^\ell$ it is immediate from Smale's 
infinite--dimensional version of Sard's theorem that $J^\ell_{reg}$ is 
dense.  For the $C^\infty$ statement a short approximation 
argument is required.  We refer the reader to~\cite[page 36]{MS2}; our 
case is just the same as theirs.
\endproof

\medskip\textbf{Remark}\qua Note that (\textbf{M6}) implies that $u$ has no multiply covered components.

We will often say a complex structure $J$ \emph{achieves
  transversality} to mean $J\in J_{reg}$.

There is a second way that we can sometimes achieve transversality,
which is more convenient for computations:  by keeping the complex
structure on $W$ split and perturbing the
$\alpha$ and $\beta$ circles.  Specifically we have:
\begin{Prop}\label{Prop:BdyInjectivity} Suppose that a
homology class $A\in\pi_2(\vx,\vy)$ with $\ind(A)=1$ is such
that any $j_\Sigma\times j_\bD$--holomorphic curve $u\co S\to W$ in
the homology class $A$ must have
$\pi_\Sigma\circ u|_{\bdy S}$ somewhere injective.
Then
for a generic perturbation of the $\alpha$-- and $\beta$--circles,
for any $u$ in the homology class $A$ the linearization
$D\dbar$, computed with respect to the complex structure
$j_\Sigma\times j_\bD$ on $W$, is surjective.
\end{Prop}
The proof of this proposition is analogous to the argument
in~\cite{Oh}.  It is also a corollary of~\cite[Proposition 3.9]{OS1}, so we omit the proof.  

We shall refer to the condition in the preceding proposition as
\emph{boundary injectivity.}  One obvious time when boundary injectivity
holds is the following:
\begin{Lem}\label{Lem:BdyInjectivity}Suppose the homology class $A\in\pi_2(\vx,\vy)$ is
  represented by a domain
  $D=\sum_i n_i D_i$ such that for some $i$ and $j$, $n_i=1$,
  $n_j=0$, and $\bdy D_i\cap \bdy D_j\neq\emptyset$.  Then $A$
  satisfies the boundary injectivity hypothesis.
\end{Lem}

For computing the homologies defined in \fullref{Section:ChainComplexes},
if the boundary injectivity criterion is met by every domain with
index $1$ it will suffice to take a generic perturbation of the
boundary conditions and the split complex structure $j_\Sigma\times
j_\bD$ rather than a generic path $J_s$ of complex structures.
(The only time this is relevant in this paper is
\fullref{Section:Handleslides}, but in practice it is necessary
for most direct computations.)

\section{Index}
\label{Section:Index}
In this section we compute the index of the linearized
$\dbar$--operator $D\dbar$ at a holomorphic map $u\co (S,j)\to (W,J)$.  
We start by reducing to a result discussed in~\cite{bourgeois} via a doubling
argument similar to the one found in \cite{HLS}.  We then reinterpret
this index several times, obtaining the Chern class formula for the index of
periodic domains (\fullref{Index:PeriodicCor}, which is
\cite[Theorem 4.9]{OS1}), J Rasmussen's formula (\cite[Theorem
9.1]{Rasmussen}, proved here in
\fullref{Prop:RasIndex}) and a combinatorial formula for the
index near an embedded curve (\fullref{Cor:Index}) and, consequently, the
Maslov index in traditional Heegaard Floer homology (\fullref{Cor:HFIndex}).

\subsection{First formulas for the index}\label{Subsection:FirstIndex}
We may assume that $J$ is split, since deformations of
$J$ will not change the index.  Also, we assume that the $\alpha$
and $\beta$ curves meet in right angles.

Let $a_1,\ldots,a_g$ be the components of the boundary of $S$ (in
the complement of the punctures) mapped to $\alpha$--cylinders
 and $b_1,\ldots,b_g$ the components
mapped to $\beta$--cylinders.  We define the quadruple of
$S$, denoted $4\ltimes S$, by gluing four copies of $S$, denoted
$S_1$, $S_2$, $S_3$, and $S_4$, as follows.  Glue each $a_i$
in $S_1$ to $a_i$ in $S_2$ and each $a_i$ in $S_3$ to
$a_i$ in $S_4$.  Similarly, glue each $b_i$ in $S_1$ to
$b_i$ in $S_3$ and each $b_i$ in $S_2$ to $b_i$ in $S_4$.
Define a complex structure on $4\ltimes S$ by taking the complex structure
$j$ on $S_1$ and $S_4$ and its conjugate $\overline{j}$ on
$S_2$ and $S_3$.  Notice that these complex structures glue
together correctly.

The complex vector bundle $(u^*TW,u^*J)$ extends to a vector bundle
over $4\ltimes S$, which we
denote $(\overline{u^*TW},\overline{J})$, in an obvious
way, and $D\dbar$ extends to an operator $4\ltimes D\dbar$ on the sections of 
this
vector bundle.  We restrict $4\ltimes D\dbar$ to the space of sections which
approach zero near each puncture, as we require fixed asymptotics.

By~\cite[Corollary 5.4, page 53]{bourgeois}, the index of $4\ltimes D\dbar$
is
\begin{equation}
-\chi(4\ltimes S)+2c_1(A).\label{Index:Formula1}
\end{equation}
Here, $c_1(A)$ is defined as follows.  Choose a small disk
 near
each point in $\ba\cap\bb$.  Trivialize
$(T\Sigma,J)$ over these disks.  This gives a trivialization
of $u^*TW$ in a neighborhood of the punctures
in $4\ltimes S$ which extends to a trivialization of $\overline{u^*TW}$
over the surface $\overline{4\ltimes S}$ obtained by filling in the punctures.
  Then,
$c_1(A)$ is the pairing of the first Chern class of
$\overline{u^*TW}$ with the fundamental class of $\overline{4\ltimes S}$.
  Note that
since $T([0,1]\times\bR)$ is trivial, to compute $c_1$ we need only
look at the $\Sigma$ factor.  
Also, it is necessary to observe that
F. Bourgeois's calculations in~\cite[Section 5]{bourgeois} are all done in the pullback bundle, so the
fact that our index problem does not correspond to a genuine map is a
nonissue.

We convert Formula~\eqref{Index:Formula1} into one not involving the
quadruple of $S$.  First we compute that $\chi(4\ltimes S)=4\chi(S)-4g$.
Indeed, after doubling along the $\alpha$--arcs the Euler characteristic is
$2\chi(S)-g$.  Doubling again we obtain
$\chi(4\ltimes S)=2(2\chi(S)-g)-2g$.  (The last summand of $-2g$ comes
from the $2g$ punctures in $4\ltimes S$.)

Second,
$c_1(A)$ can be computed from
Maslov--type indices as follows.  Choose the trivializations of $TW$ over
the neighborhoods $\mathcal{V}$ of
$\ba\cap\bb$ above so
that for $p\in\ba\cap\bb$, 
$T_p\beta = \bR\subset \bC$ and $T_p\alpha = i\bR\subset \bC$.
Trivialize $(\pi_\Sigma\circ u)^*T\Sigma$ over $S$ so that this
trivialization agrees with the specified trivialization of $TW$ over the
neighborhoods $\mathcal{V}$.  Then, each boundary arc
$a_i$ or $b_i$ gives a loop of lines in $\bC$, and so has a
well--defined Maslov index $\mu(a_i)$ or $\mu(b_i)$.  It is not
hard to see that
$c_1(A)=2\left(\sum_{i=1}^g\mu(a_i)-\mu(b_i)\right)$.  This is
independent of the choice of trivialization subject to the specified criteria.

Fix a map $u\co (S,\bdy S)\to (W, C_\alpha\cup C_\beta)$.
An argument almost exactly like the one in \cite{HLS} shows that
\begin{equation}
\ind(D\dbar)=\frac{1}{4}\ind(4\ltimes D\dbar) = g - \chi(S) + \sum_{i=1}^g \mu(a_i) -
\sum_{i=1}^g \mu(b_i).\label{Index:Formula2}
\end{equation}
The factor of $1/4$ comes from the ``matching conditions'' on the
boundary of $S$.  

We again reinterpret the Maslov indices.  Given a domain $D$, we
define the Euler measure of $D$ as follows. Suppose first that $D$ is a surface with boundary and corners.
Choose a metric on $D$ such
that $\bdy D$ is geodesic and such that the corners of $D$ are
right angles.  Then the Euler measure $e(D)$ is defined to be 
$\frac{1}{2\pi}$ times the integral over $D$ of the curvature of the metric.
(This is normalized so that the Euler measure of a sphere 
is $2$, agreeing with its Euler characteristic.)
From this definition it is clear that the Euler measure is additive
under disjoint unions and gluing of components along boundaries, and so the definition extends naturally to domains (linear combinations of regions in $\Sigma$).

\begin{figure}
\centering
\includegraphics[scale=.8]{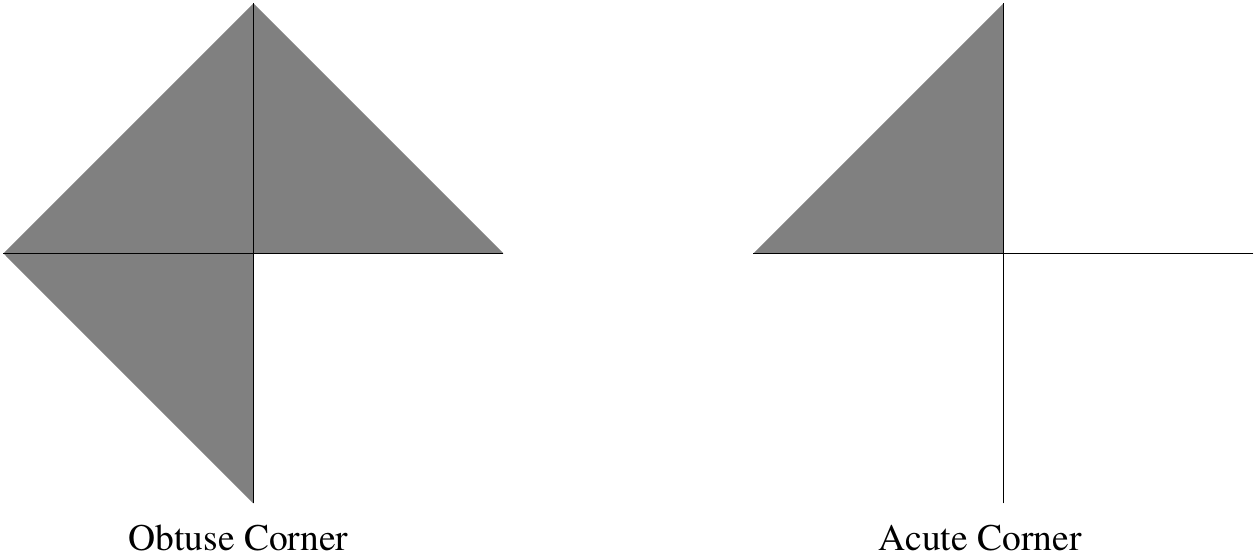}
\caption{$S$ is the shaded region.}
\label{Figure:Corners}
\end{figure}

It follows from the Gauss--Bonnet theorem that the Euler measure
of a surface $S$ with $k$ acute right--angled corners (see
\fullref{Figure:Corners}) and $\ell$ 
obtuse right--angled
corners is $\chi(S)-k/4+\ell/4$.  As with the previous formulation of Euler measure, this formula is 
additive, so the Euler measure of a domain $D=\sum_i D_i$ is 
$e(D)=\sum_i e(D_i)$.

From the Gauss--Bonnet theorem, we also know that if we endow $D$ with 
a flat metric such that all corners are right angles then the Euler
measure $e(D)$ is $\frac{1}{2\pi}$ times 
the geodesic curvature of $\bdy D$.  It is then clear that for $D$ the domain corresponding to $u$, 
$\sum_{i=1}^g
\mu(a_i)-\mu(b_i)=2e(D)$.
So, we can recast the index formula as
\begin{equation}
\label{Index:Formula3}
\ind(D\dbar)=g-\chi(S)+2e(D).
\end{equation}

\subsection{Determining $S$ from $A$}\label{Subsection:Determining}
Note that the formulas for the index derived so far depend not only on
the homology class but also on the topological type of the source.
This is as it should be.  However, as we will show presently, for
\emph{embedded} holomorphic curves the Euler characteristic of the source is
determined by the homology class.  (Actually, we prove this more generally for
  any curves satisfying certain hypotheses described
  in \fullref{Lemma:Represent} below, not just holomorphic ones.)  In fact, we can give an explicit
formula for the Euler characteristic, allowing us to give a
combinatorial formula for the index.  We will
see in \fullref{Subsection:IndexComparison} that this formula calculates
the Maslov index in the setup of~\cite{OS1} as well.  Before proving
this claim we introduce some more terminology and notation.

\begin{figure}
\centering
\includegraphics[scale=.75]{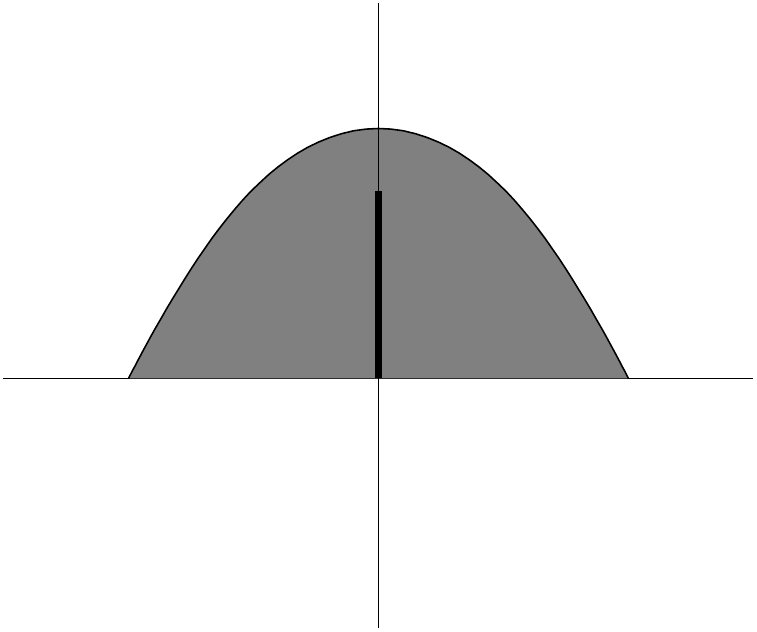}
\caption{A degenerate corner}
\end{figure}

Given an intersection point $\vx$ we call each $x_i\in\vx$ a
\emph{corner} of $\vx$.  Following Rasmussen, we define a corner
$x_i$ of $\vx$ to 
be \emph{degenerate} for a homology class $A\in\pi_2(\vx,\vy)$ if 
$x_i=y_j$ for some $y_j\in\vy$.  This definition will be convenient
presently.

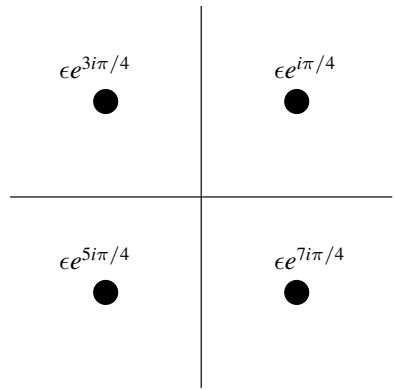
\begin{figure}
\centering
\begin{picture}(0,0)%
\includegraphics{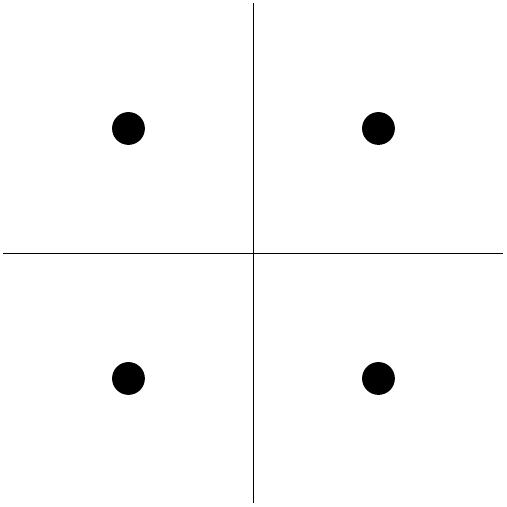}%
\end{picture}%
\setlength{\unitlength}{3947sp}%
\begingroup\makeatletter\ifx\SetFigFont\undefined%
\gdef\SetFigFont#1#2#3#4#5{%
  \reset@font\fontsize{#1}{#2pt}%
  \fontfamily{#3}\fontseries{#4}\fontshape{#5}%
  \selectfont}%
\fi\endgroup%
\begin{picture}(2424,2424)(589,-2173)
\put(901,-211){\makebox(0,0)[lb]{\smash{\SetFigFont{12}{14.4}{\rmdefault}{\mddefault}{\updefault}{\color[rgb]{0,0,0}\(\epsilon e^{3i\pi/4}\)}%
}}}
\put(901,-1411){\makebox(0,0)[lb]{\smash{\SetFigFont{12}{14.4}{\rmdefault}{\mddefault}{\updefault}{\color[rgb]{0,0,0}\(\epsilon e^{5i\pi/4}\)}%
}}}
\put(2251,-1411){\makebox(0,0)[lb]{\smash{\SetFigFont{12}{14.4}{\rmdefault}{\mddefault}{\updefault}{\color[rgb]{0,0,0}\(\epsilon e^{7i\pi/4}\)}%
}}}
\put(2251,-211){\makebox(0,0)[lb]{\smash{\SetFigFont{12}{14.4}{\rmdefault}{\mddefault}{\updefault}{\color[rgb]{0,0,0}\(\epsilon e^{i\pi/4}\)}%
}}}
\end{picture}
\caption{$n_p(A)=\frac{1}{4}\left( 
n_{\epsilon e^{i\pi/4}}(A)+n_{\epsilon e^{3i\pi/4}}(A)+n_{\epsilon e^{5i\pi/4}}(A)
+n_{\epsilon e^{7i\pi/4}}(A)\right)$
}
\label{Index:Figure0}
\end{figure}

Let $p\in\alpha_i\cap\beta_j$.  For a homology class 
$A\in\pi_2(\vx,\vy)$, define $n_p(A)$ to be the average of 
the coefficients of $A$ of the four cells with corners at $p$.  More 
precisely, choose coordinates identifying a neighborhood of $p$ in 
$\Sigma$ with the unit disk in $\bC$, $\alpha_i$ with the real axis, 
and $\beta_j$ with the imaginary axis.  Then $n_p(A)=\frac{1}{4}\big( 
n_{\epsilon e^{i\pi/4}}(A)+n_{\epsilon e^{3i\pi/4}}(A)+n_{\epsilon e^{5i\pi/4}}(A)
+n_{\epsilon e^{7i\pi/4}}(A)\big)$ for some $\epsilon<1$.  See \fullref{Index:Figure0}.
Define $n_\vx(A)$ to be $\sum_{x_i\in\vx} n_{x_i}(A)$, and 
$n_\vy(A)=\sum_{y_i\in\vy}n_{y_i}(A)$.  (See~\cite[page 1202]{OS2}.)

We need a lemma about representability of homology classes with
positive coefficients:
\begin{Lem}\label{Lemma:Represent}
Suppose $A\in\pi_2(\vx,\vy)$ is a positive homology class.
  Then there is a Riemann surface
  with boundary and corners $\overline{S}$ and smooth map $u\co S\to
  W$ (where $S$ denotes the complement in $\overline{S}$ of the
  corners of $\overline{S}$) in the homology class $A$ such that: 
\begin{enumerate}
\item $u^{-1}(C_\alpha\cup C_\beta)=\bdy S$.
\item \label{item:1} For each $i$, $u^{-1}(\alpha_i\times\{1\}\times\bR)$ and
  $u^{-1}(\beta_i\times\{0\}\times\bR)$ each consists of one arc in
  $\bdy S$.
\item The map $u$ is $J$--holomorphic in a
  neighborhood of $(\pi_\Sigma\circ u)^{-1}(\ba\cup\bb)$ for some
  $J$ satisfying {\rm(\textbf{J1})--(\textbf{J5})} (in fact, for
  $j_\Sigma\times j_\bD$).
\item For each component of $S$, either
\begin{itemize}
\item The component is a disk with two boundary punctures and the map is
  a diffeomorphism to $\{x_i\}\times[0,1]\times\bR$ for some
  $x_i\in\ba\cap\bb$ (such a component is a \emph{degenerate disk}) or
\item The map $\pi_\Sigma\circ u$ extends to a branched covering map
  $\overline{\pi_\Sigma\circ u}$, none
  of whose branch points map to points in $\ba\cap\bb$.
\end{itemize}
\item All the corners of $S$ are acute (see
  \fullref{Figure:Corners}).
\item \label{item:last}The map $u$ is an embedding.
\end{enumerate}
\end{Lem}
(Note that it follows from the conditions in the lemma that the map to
$\Sigma$ is orientation--preserving.  Also, observe that a generic
holomorphic representative of the homology class satisfies all of
the properties of the lemma.)

\proof
(For a similar construction of a $u$ with slightly different
properties, see Rasmussen,~\cite[Lemma 9.3]{Rasmussen}.  His construction is slightly more
subtle than we need.  Another inspiring construction can be found
in~\cite[Lemma 2.17]{OS1}.)

Let $D_1,\cdots,D_N$ denote the closures of the components of
$\Sigma\setminus (\ba\cup\bb)$, enumerated so that $\Fz_i\in
D_i$.  Form a surface $\ol{S}_0$ by gluing together $n_{\Fz_i}(A)$
copies of $D_i$ ($i=1,\cdots,N$) pairwise along common boundaries
maximally.  There is then an obvious orientation--preserving map
$p_{\Sigma,0}\co \ol{S}_0\to \Sigma$ which covers $\Fz_i$
$n_{\Fz_i}(A)$ times.

\begin{figure}
\centering
\begin{picture}(0,0)%
\includegraphics[scale=.8]{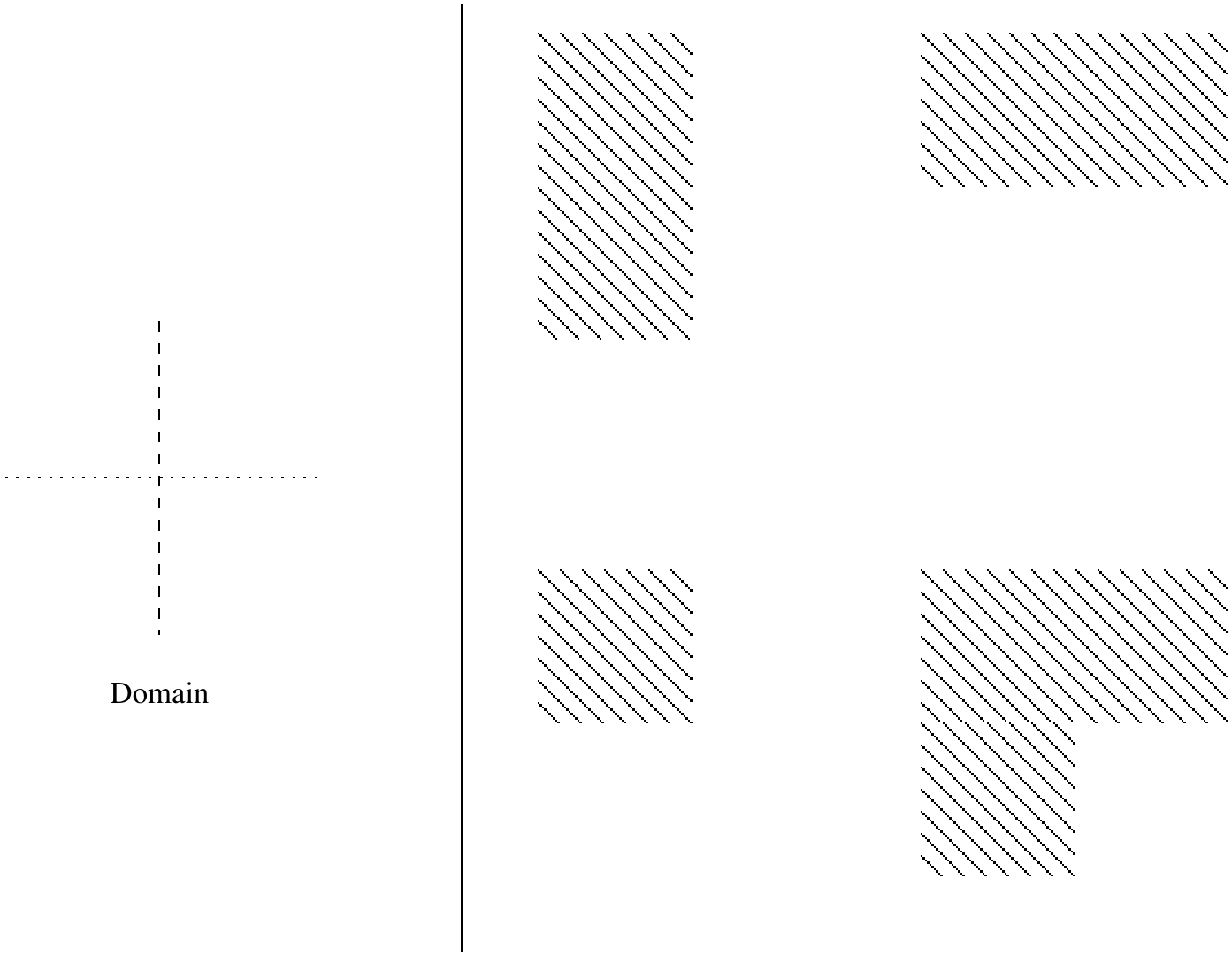}%
\end{picture}%
\setlength{\unitlength}{3158sp}%
\begingroup\makeatletter\ifx\SetFigFont\undefined%
\gdef\SetFigFont#1#2#3#4#5{%
  \reset@font\fontsize{#1}{#2pt}%
  \fontfamily{#3}\fontseries{#4}\fontshape{#5}%
  \selectfont}%
\fi\endgroup%
\begin{picture}(6677,5170)(589,-5271)
\put(4708,-3752){\makebox(0,0)[lb]{\smash{{\SetFigFont{9}{10.8}{\rmdefault}{\mddefault}{\updefault}{\color[rgb]{0,0,0}\(\bigcup\)}%
}}}}
\put(1508,-3587){\makebox(0,0)[lb]{\smash{{\SetFigFont{9}{10.8}{\rmdefault}{\mddefault}{\updefault}{\color[rgb]{0,0,0}\(\alpha\)}%
}}}}
\put(2360,-2788){\makebox(0,0)[lb]{\smash{{\SetFigFont{9}{10.8}{\rmdefault}{\mddefault}{\updefault}{\color[rgb]{0,0,0}\(\beta\)}%
}}}}
\put(921,-3109){\makebox(0,0)[lb]{\smash{{\SetFigFont{9}{10.8}{\rmdefault}{\mddefault}{\updefault}{\color[rgb]{0,0,0}\(1\)}%
}}}}
\put(975,-2308){\makebox(0,0)[lb]{\smash{{\SetFigFont{9}{10.8}{\rmdefault}{\mddefault}{\updefault}{\color[rgb]{0,0,0}\(2\)}%
}}}}
\put(1774,-2308){\makebox(0,0)[lb]{\smash{{\SetFigFont{9}{10.8}{\rmdefault}{\mddefault}{\updefault}{\color[rgb]{0,0,0}\(1\)}%
}}}}
\put(1828,-3161){\makebox(0,0)[lb]{\smash{{\SetFigFont{9}{10.8}{\rmdefault}{\mddefault}{\updefault}{\color[rgb]{0,0,0}\(0\)}%
}}}}
\put(4863,-1205){\makebox(0,0)[lb]{\smash{{\SetFigFont{9}{10.8}{\rmdefault}{\mddefault}{\updefault}{\color[rgb]{0,0,0}\(\bigcup\)}%
}}}}
\put(4708,-2192){\makebox(0,0)[lb]{\smash{{\SetFigFont{9}{10.8}{\rmdefault}{\mddefault}{\updefault}{\color[rgb]{0,0,0}One possible \(S\)}%
}}}}
\put(4708,-5051){\makebox(0,0)[lb]{\smash{{\SetFigFont{9}{10.8}{\rmdefault}{\mddefault}{\updefault}{\color[rgb]{0,0,0}Another possible \(S\)}%
}}}}
\end{picture}%
\caption{Two possible gluings $\ol{S}_0$ for a specified domain.
  The latter leads to extra corners of $\ol{S}_0$.}
\label{Figure:IndexNew1}
\end{figure}
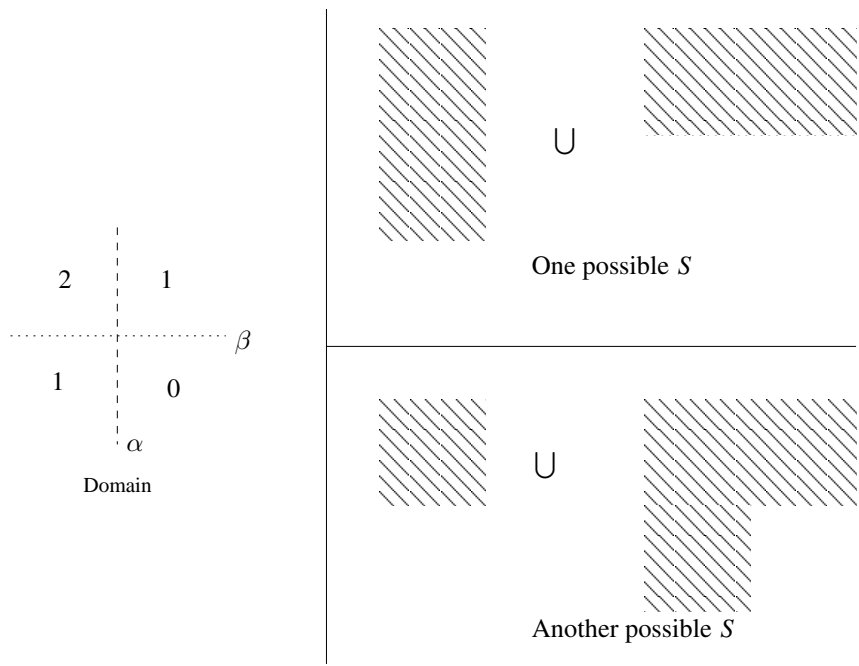

It is possible to perform the specified gluing so that the only
corners of $\ol{S}_0$ correspond one--to--one with non--degenerate corners of the domain.
This is not automatic; see \fullref{Figure:IndexNew1}.  One way to
achieve this is to first glue maximally along
$\alpha$--arcs.  Then glue along the $\beta$--arcs as much as
possible without introducing any obtuse corners.  After doing so, any
remaining corners must correspond to non--degenerate corners of the
domain.

The surface $\ol{S}_0$ lacks corners at degenerate corners of the
domain.  Let $\bD^2_2$ denote a disk with two punctures on the
boundary.  Let $\ol{S}_1$ denote the disjoint union of $\ol{S}_0$
with a copy of $\bD^2_2$ for each degenerate corner of the domain.
Extend $p_{\Sigma,0}$ to a map $p_{\Sigma,1}$ from $\bD^2_2$ by
mapping one copy of $\bD^2_2$ to each degenerate corner of the
domain.

Now, $S_0$ inherits a complex structure from $\Sigma$.  Extend
  this complex structure arbitrarily over the new disks in $S_1$ to
  obtain a complex structure on $S_1$.
  It is easy
to choose a map $p_{\bD,1}\co S_1\to[0,1]\times\bR$
such that the map $(p_{\Sigma,1}\times p_{\bD,1})\co S_1\to W$ satisfies
all of the properties specified in the statement of the lemma except perhaps
numbers~\eqref{item:1} and~\eqref{item:last}.  Perturbing
  $p_{\bD,1}$ we may assume that
$(p_{\Sigma,1}\times p_{\bD,1})$ is an embedding except for a
collection of transverse double points.  

For each $\alpha_i$ (respectively $\beta_j$), there will be exactly
one arc in $\bdy\Sigma_1$ mapped by $p_{\Sigma,1}$ to
$\alpha_i$ (respectively $\beta_j$), and possibly some circles in
$\bdy\Sigma_1$ mapped by $p_{\Sigma,1}$ to $\alpha_i$
(respectively $\beta_j$).  The map $p_{\bD,1}$ must be constant
near each closed component of $\bdy S_1$, and the image of the arc under
$(p_{\Sigma,1}\times p_{\bD,1})$ must intersect the image of each
closed component of $\bdy S_1$ mapped to $\alpha_i$ (respectively
$\beta_j$) exactly once.

Modifying $S_1$ and $p_{\Sigma,1}\times p_{\bD,1}$ near the double
points of $p_{\Sigma,1}\times p_{\bD,1}$
we can obtain a new map $u\co S\to W$ satisfying all of
the stated properties:  in the process of deforming away the double
points, we necessarily achieve property~\eqref{item:1} as well.
\endproof

\begin{Prop} Let $u\co S\to W$ be a map satisfying the conditions
  enumerated in the previous lemma,
representing a homology class $A$.  Then the Euler characteristic
$\chi(S)$ is given by  
\begin{equation}\label{Index:Formula4}
\chi(S)=g-n_\vx(A)-n_\vy(A)+e(A).\end{equation}
\label{Index:ChiCalcProp}
\end{Prop}
\proof
Applying the Riemann--Hurwitz formula to $\pi_\Sigma\circ u$, we 
only need to
calculate the degree of branching of $\pi_\Sigma\circ u$.

To calculate the number of branch points of
$\pi_\Sigma\circ u$ we reinterpret this number as a
self--intersection number.  We will assume all branch
points of $\pi_\Sigma\circ u$ have order $2$; we can clearly
arrange this.  Observe that since $S$
has no obtuse corners, by the Riemann--Hurwitz formula,
\begin{align*}
\chi(S)=e(S)+g/2=e(A)&-(\textrm{number of
  branch points})\\&+{\textstyle\frac{1}{2}}(\textrm{number of trivial disks})+g/2. 
\end{align*} 
(Branch points on $\bdy S$ should each be counted as half of a
branch point.)

Assume for the time being that $u$ contains no trivial disks, and in
fact has no degenerate corners.

Notice that the number of branch points of $\pi_\Sigma\circ u$ is
equal to the number of times the vector field  $\frac{\partial}{\partial t}$
is tangent to $u$.  (Tangencies on $\bdy S$ should each be counted
as half of a tangency.)  Let $u'$ denote the curve obtained from $u$
by translating a distance $R$ in the $\bR$--direction.  Then, for
small $R$, the number of branch points of $\pi_\Sigma\circ u$ is
equal to the intersection number of $u$ and $u'$.  (Intersections
on $\bdy S$ should each be counted as half of an intersection.)

This intersection number is invariant under isotopies of $u'$ such that 
all intersection points of $u$ and $u'$ remain in a compact subset of 
of $W$.  
(The only thing to check is that
when an intersection point in the interior of $W$ 
hits the boundary it gives rise to a pair of intersection points on the 
boundary.  It is not hard to check this using a doubling 
argument in a neighborhood of the boundary.)
We will calculate the intersection number by 
translating $u'$ far in the $\bR$--direction of $W$. 

Translate $u'$ by some $R\gg0$ in the $\bR$--factor of $W$.  
All intersection points between $u$ and $u'$ stay in a compact subset 
of $W$, so the intersection number $\#u\cap u'$ is unchanged.

We can modify $u'$ so that near each negative puncture (corresponding to 
some $x_i$) $u'$ agrees with the trivial disk 
$x_i\times[0,1]\times\bR$.  Further, we can do this modification so that 
all intersection points between $u$ and $u'$ stay within some compact 
subset of $W$.  (This follows from the simple asymptotic behavior of 
$u'$ near $-\infty$.)

Similarly, we can modify $u$ so that near the 
positive punctures of $S$, $u$ agrees with the trivial disks 
$\{y_i\}\times[0,1]\times\bR$ ensuring in the process that all intersection points 
between $u$ and $u'$ stay within some compact subset of $W$.

Finally, for $R$ large enough, we can assume that after the two 
modifications every intersection point between $u$ and $u'$ 
corresponds to an intersection point between $u$ and 
$\{x_i\}\times[0,1]\times\bR$ or between $\{y_j\}\times[0,1]\times\bR$ 
and $u'$.

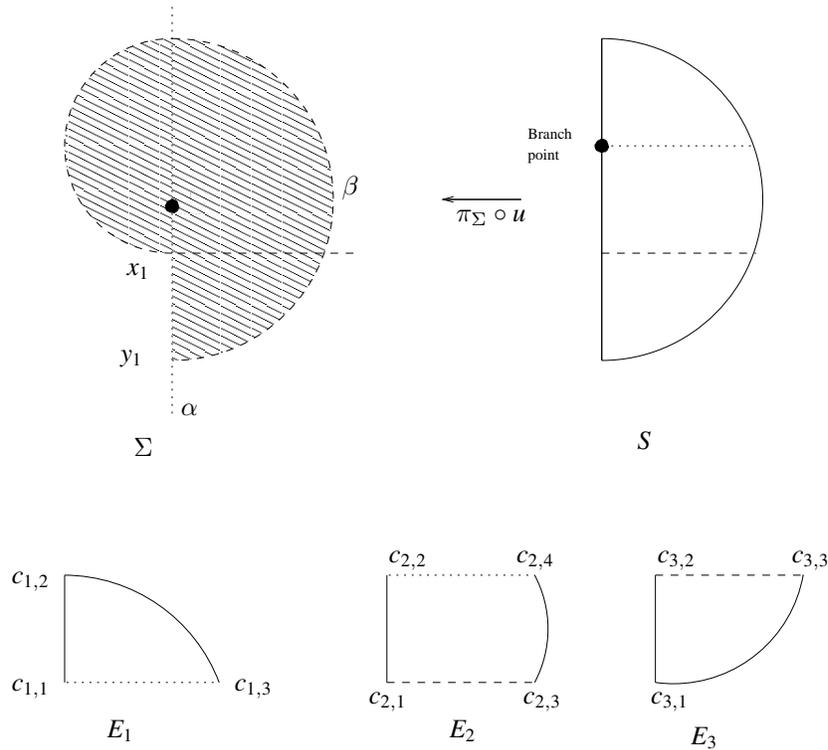
\begin{figure}
\centering
\begin{picture}(0,0)%
\includegraphics[scale=.8]{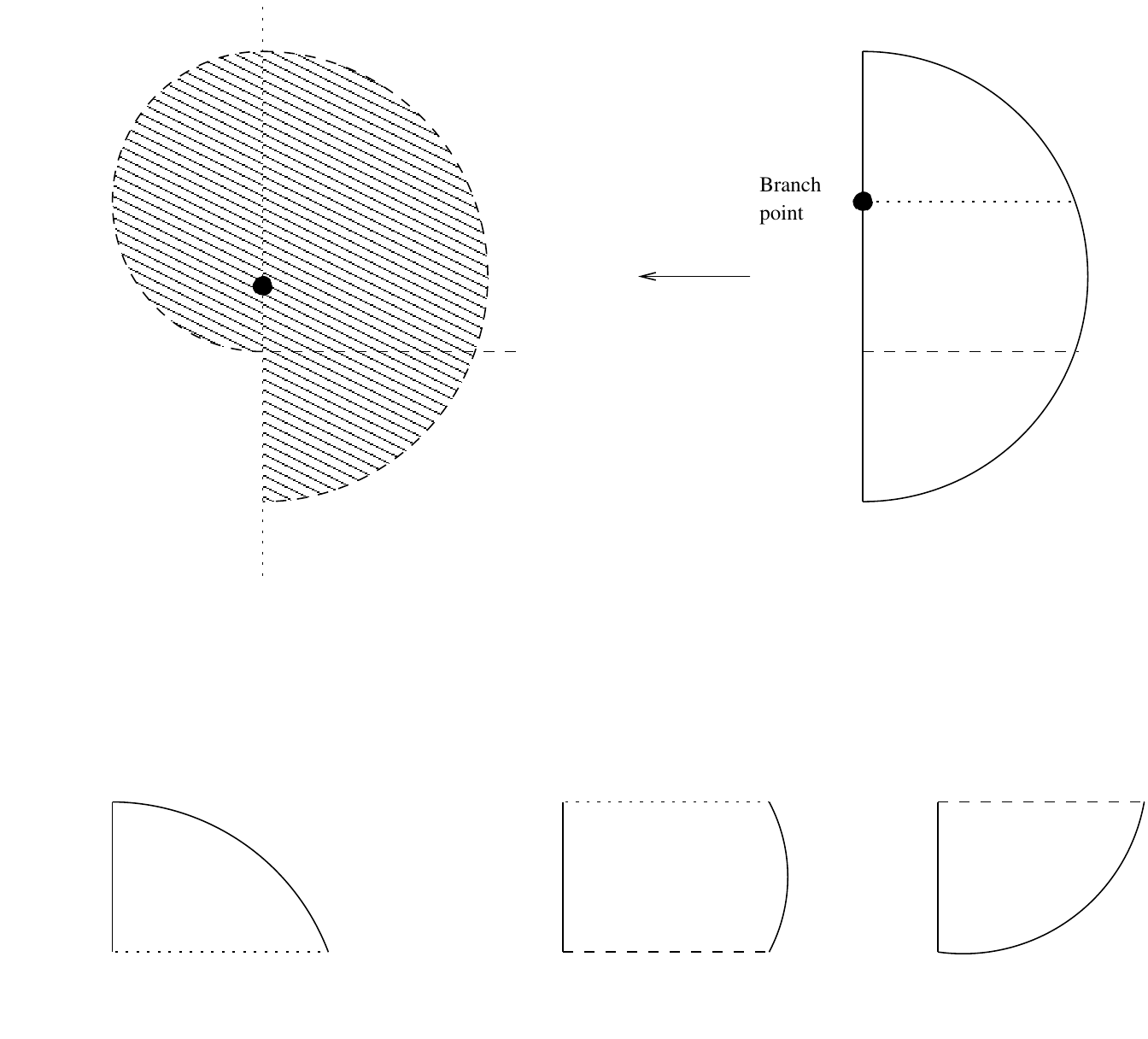}%
\end{picture}%
\setlength{\unitlength}{3158sp}%
\begingroup\makeatletter\ifx\SetFigFont\undefined%
\gdef\SetFigFont#1#2#3#4#5{%
  \reset@font\fontsize{#1}{#2pt}%
  \fontfamily{#3}\fontseries{#4}\fontshape{#5}%
  \selectfont}%
\fi\endgroup%
\begin{picture}(6442,5809)(301,-5183)
\put(1460,-2917){\makebox(0,0)[lb]{\smash{\SetFigFont{8}{9.6}{\rmdefault}{\mddefault}{\updefault}{\color[rgb]{0,0,0}\(\Sigma\)}%
}}}
\put(3990,-1072){\makebox(0,0)[lb]{\smash{\SetFigFont{8}{9.6}{\rmdefault}{\mddefault}{\updefault}{\color[rgb]{0,0,0}\(\pi_\Sigma\circ u\)}%
}}}
\put(501,-4762){\makebox(0,0)[lb]{\smash{\SetFigFont{8}{9.6}{\rmdefault}{\mddefault}{\updefault}{\color[rgb]{0,0,0}\(c_{1,1}\)}%
}}}
\put(501,-3918){\makebox(0,0)[lb]{\smash{\SetFigFont{8}{9.6}{\rmdefault}{\mddefault}{\updefault}{\color[rgb]{0,0,0}\(c_{1,2}\)}%
}}}
\put(2251,-4762){\makebox(0,0)[lb]{\smash{\SetFigFont{8}{9.6}{\rmdefault}{\mddefault}{\updefault}{\color[rgb]{0,0,0}\(c_{1,3}\)}%
}}}
\put(4517,-4867){\makebox(0,0)[lb]{\smash{\SetFigFont{8}{9.6}{\rmdefault}{\mddefault}{\updefault}{\color[rgb]{0,0,0}\(c_{2,3}\)}%
}}}
\put(3937,-5131){\makebox(0,0)[lb]{\smash{\SetFigFont{8}{9.6}{\rmdefault}{\mddefault}{\updefault}{\color[rgb]{0,0,0}\(E_2\)}%
}}}
\put(3463,-3753){\makebox(0,0)[lb]{\smash{\SetFigFont{8}{9.6}{\rmdefault}{\mddefault}{\updefault}{\color[rgb]{0,0,0}\(c_{2,2}\)}%
}}}
\put(5571,-3753){\makebox(0,0)[lb]{\smash{\SetFigFont{8}{9.6}{\rmdefault}{\mddefault}{\updefault}{\color[rgb]{0,0,0}\(c_{3,2}\)}%
}}}
\put(6625,-3753){\makebox(0,0)[lb]{\smash{\SetFigFont{8}{9.6}{\rmdefault}{\mddefault}{\updefault}{\color[rgb]{0,0,0}\(c_{3,3}\)}%
}}}
\put(5835,-5183){\makebox(0,0)[lb]{\smash{\SetFigFont{8}{9.6}{\rmdefault}{\mddefault}{\updefault}{\color[rgb]{0,0,0}\(E_3\)}%
}}}
\put(1250,-5131){\makebox(0,0)[lb]{\smash{\SetFigFont{8}{9.6}{\rmdefault}{\mddefault}{\updefault}{\color[rgb]{0,0,0}\(E_1\)}%
}}}
\put(5413,-2864){\makebox(0,0)[lb]{\smash{\SetFigFont{8}{9.6}{\rmdefault}{\mddefault}{\updefault}{\color[rgb]{0,0,0}\(S\)}%
}}}
\put(4464,-3753){\makebox(0,0)[lb]{\smash{\SetFigFont{8}{9.6}{\rmdefault}{\mddefault}{\updefault}{\color[rgb]{0,0,0}\(c_{2,4}\)}%
}}}
\put(3305,-4867){\makebox(0,0)[lb]{\smash{\SetFigFont{8}{9.6}{\rmdefault}{\mddefault}{\updefault}{\color[rgb]{0,0,0}\(c_{2,1}\)}%
}}}
\put(5519,-4867){\makebox(0,0)[lb]{\smash{\SetFigFont{8}{9.6}{\rmdefault}{\mddefault}{\updefault}{\color[rgb]{0,0,0}\(c_{3,1}\)}%
}}}
\put(1408,-1494){\makebox(0,0)[lb]{\smash{\SetFigFont{8}{9.6}{\rmdefault}{\mddefault}{\updefault}{\color[rgb]{0,0,0}\(x_1\)}%
}}}
\put(3094,-862){\makebox(0,0)[lb]{\smash{\SetFigFont{8}{9.6}{\rmdefault}{\mddefault}{\updefault}{\color[rgb]{0,0,0}\(\beta\)}%
}}}
\put(1829,-2601){\makebox(0,0)[lb]{\smash{\SetFigFont{8}{9.6}{\rmdefault}{\mddefault}{\updefault}{\color[rgb]{0,0,0}\(\alpha\)}%
}}}
\put(1355,-2179){\makebox(0,0)[lb]{\smash{\SetFigFont{8}{9.6}{\rmdefault}{\mddefault}{\updefault}{\color[rgb]{0,0,0}\(y_1\)}%
}}}
\end{picture}
\caption{With respect to $x_1$, corner $c_{1,2}$ has type (2),
  corners $c_{2,1}$ and $c_{3,2}$ have type (3), and all others have
  type (1).}
\label{Index:Figure0.5}
\end{figure}

Now, for each corner $c_{k,\ell}$ of each component $E_k$ of
$S\setminus (\pi_\Sigma\circ
u)^{-1}(\ba\cup\bb)$, one of the following four phenomena 
occurs:
\begin{enumerate}
\item The corner $c_{k,\ell}$ is mapped by $\pi_\Sigma\circ u$
  somewhere other than
  $x_i$.  That is, $\lim_{p\to c_{k,\ell}}\pi_\Sigma\circ
  u(p)\neq x_i$.
\item 
The corner $c_{k,\ell}$ is mapped by $\pi_\Sigma\circ u$  to
  $x_i$, but at
  $-\infty$.  
That is, $\lim_{p\to c_{k,\ell}}\pi_\Sigma\circ
  u(p)=x_i$ but $\lim_{p\to c_{k,\ell}}\pi_\bR\circ u(p)=-\infty$.
\item 
The corner $c_{k,\ell}$ is mapped by $\pi_\Sigma\circ u$ to
$x_i$, and is mapped by $u$ to the boundary of $W$.  That is,
$\pi_\Sigma\circ u(c_{k,\ell})=x_i$ and $\pi_\bD\circ u(c_{k,\ell})\in\bdy
[0,1]\times\bR$.
\item 
The corner $c_{k,\ell}$ is mapped by $\pi_\Sigma\circ u$ to
$x_i$, and is mapped by $u$ to the interior of $W$.  That is,
$\pi_\Sigma\circ u(c_{k,\ell})=x_i$ and $\pi_\bD\circ u(c_{k,\ell})\in
(0,1)\times\bR$.
\end{enumerate}
(Compare \fullref{Index:Figure0.5}.)

In the first two cases, $E_k$ does not contribute to $\#(u\cap 
x_i\times[0,1]\times\bR)$.  In the last two, $E_k$ contributes $1/4$ 
to the intersection number.  Notice that the second case occurs
exactly once, since $S$ has no obtuse corners.  There are a total of
$4n_{\vx}$ corners satisfying one of conditions (2)--(4) for some
$x_i$.  Exactly $g$ of them satisfy condition (2).  So, $\#
\left(u\cap\{x_i\}\times[0,1]\times\bR\right) = n_\vx - g/4$.

A similar analysis
works for the intersection 
points between $u'$ and $\{y_j\}\times[0,1]\times\bR$.  So, it follows 
that the intersection number between $u$ and $u'$ is
$$\#(u\cap u')=n_\vx(A)+n_\vy(A) - g/2.$$
It follows that $\chi(S)=e(A)-n_\vx(A)-n_\vy(A)+g$.

In the proof so far we assumed that there were no 
trivial disks.  Suppose 
$u$ contains trivial disks corresponding to the intersection points
$x_{i_1},\cdots,x_{i_k}$.   Since we are
considering only embedded curves, $n_{x_{i_j}}(A)=0$ for
$j=1,\cdots,k$.  By the argument above, after ignoring the trivial
disks, we find $\#(u\cap u')=n_\vx(A)+n_\vy(A)-g/2 + k/2$.  So, we
have the same formula for $\chi(S)$ as before.

Finally, we deal with degenerate corners which are not trivial disks.
Since we are assuming $S$ has only acute corners, and acute
degenerate corners have exactly one shared boundary component under
$\pi_\Sigma\circ u$, it is easy to
see that after
translating in the $\bR$--direction, there will be one intersection
point along the boundary near the puncture, so the extra corner mapped
to $\pm\infty$
contributes $\frac{1}{2}$ to the intersection number, as one would
expect from our formula.  This concludes the proof.
\endproof

Note that if $u\co S\to W$ is an embedded holomorphic curve (with
respect to any complex structure $J$ on $W$ satisfying
(\textbf{J1})--(\textbf{J5})) then, after slitting $S$ near any
obtuse corners and perturbing $u$ slightly, $u\co S\to W$ satisfies
the conditions of \fullref{Lemma:Represent}.  It follows that
\fullref{Index:ChiCalcProp} calculates the Euler
characteristic of $S$.

\begin{Cor}\label{Cor:Index} For $A$ a positive homology class
  and $u\co S\to W$ a representative for $A$ satisfying the
  conditions of \fullref{Lemma:Represent}, the index of the
  $D\dbar$ operator near $u$ is given by
\begin{equation}\label{Equation:ComboIndex}
\ind(D\dbar)=e(A)+n_\vx(A)+n_\vy(A)
\end{equation}
\end{Cor}
\proof
This is immediate from formula~\eqref{Index:Formula3} and \fullref{Index:ChiCalcProp}.
\endproof

\begin{Def}Given a positive homology class $A$ define the index
  $\ind(A)$ of $A$ to be the index of the $D\dbar$ operator near
  any curve satisfying the conditions of \fullref{Lemma:Represent}.
\end{Def}

\begin{Cor}\label{Cor:IndexSigma}
If $A$ and $A+k[\Sigma]$ are both positive
  then $\ind(A+k[\Sigma])=\ind(A)+2k$.
\end{Cor}
\proof
By \fullref{Cor:Index}, 
{\setlength\arraycolsep{2pt}
\begin{eqnarray*}
\ind(A+k[\Sigma])&=&e(A+k[\Sigma])+n_\vx(A+k[\Sigma])+n_\vy(A+k[\Sigma])\\
&=&e(A)+(2-2g)k+n_\vx(A)+gk+n_\vy(A)+gk\\
&=&e(A)+n_\vx(A)+n_\vy(A)+2k\\
&=&\ind(A)+2k.
\end{eqnarray*}}
\endproof

\begin{Def}For any homology class $A$ define the index $\ind(A)$ of
  $A$ to be $\ind(A+k[\Sigma])-2k$ where $k$ is chosen large
  enough that $A+k[\Sigma]$ is positive.
\end{Def}

\begin{Cor}\label{Cor:IndexAdds} Suppose that $A\in\pi_2(\vx,\vy)$ and
  $A'\in\pi_2(\vy,\vz)$.  Then $\ind(A+A')=\ind(A)+\ind(A')$.
\end{Cor}
\proof
We may clearly assume that $A$ and $A'$ are both positive.  Let
$u\co S\to W$ and $u'\co S'\to W$ be maps satisfying the conditions
\fullref{Lemma:Represent} representing $A$ and $A'$
respectively.  Then we can glue $u$ and $u'$ to a map $u\natural
u'\co S\natural S'\to W$ representing $A+A'$.  It follows from general
gluing results for the index that $\ind(D\dbar)(u\natural
u')=\ind(D\dbar)(u)+\ind(D\dbar)(u')$.  (Alternately, it follows from
the additivity of Formula~\eqref{Index:Formula3} under gluing.)
\endproof

\medskip\textbf{Remark}\qua Formula~\eqref{Equation:ComboIndex} was suggested to
me by Z Szab\'o.  Specifically, he suggested that it seemed the Maslov index in~\cite{OS1} can be calculated by this formula.  In particular, in the special case when $A\in\pi_2(\vx,\vx)$, Ozsv\'ath and Szab\'o proved (\cite[Theorem 4.9]{OS1} and~\cite[Proposition 7.5]{OS2}) that Formula~\eqref{Equation:ComboIndex} does computes the Maslov index.
Note that it is not even clear \emph{a priori} that Formula~\eqref{Equation:ComboIndex} is additive.  In fact, I do not know a more direct proof than the one we used to obtain \fullref{Cor:IndexAdds}.

\subsection{Comparison with classical Heegaard Floer homology}\label{Subsection:IndexComparison}
In this subsection we assume familiarity with~\cite{OS1}.

By considering domains, for example, there is a natural identification
of our $\pi_2(\vx,\vy)$ with $\pi_2(\vx,\vy)$ as defined
in~\cite[Section 2.4]{OS1}.  For $A\in\pi_2(\vx,\vy)$, let $\mu(A)$ denote the
Maslov index of $A$, viewed as a homotopy class of maps disks in
$(\Sym^g(\Sigma),T_{\alpha}\cup T_{\beta})$.
The goal of this subsection is to prove the following

\begin{Prop}\label{Prop:RasIndex}For $A\in\pi_2(\vx,\vy)$ we have $\ind(A)=\mu(A)$.
\end{Prop}
\proof
It is possible to give a direct proof (see~\cite[proof of Theorem
9.1]{Rasmussen}), but instead of doing so we will show our formula
agrees with the one given by Rasmussen in~\cite[Theorem 9.1]{Rasmussen}.
He proves that at a disk
$\phi\co (\bD,\bdy\bD)\to(\Sym^g(\Sigma),T_\alpha\cup T_\beta)$, 
\begin{equation}
\mu(\phi)=\Delta\cdot\phi + 2e(\phi).\end{equation}
Here, $e$ is the Euler measure defined in \fullref{Subsection:FirstIndex} and $\Delta\cdot\phi$ is the
algebraic intersection number of $\phi$
with the diagonal in $\Sym^g(\Sigma)$.  (The diagonal is an algebraic
subvariety of $\Sym^g(\Sigma)$ of real codimension 2 so the intersection
number is well--defined.)  

To compare his result with ours, we need a slight strengthening
of \fullref{Lemma:Represent}:
\begin{Lem}Suppose $A$ is a positive homology class.  Then we can
  represent $A+[\Sigma]$ by a map $u\co S\to W$ satisfying all the
  conditions of \fullref{Lemma:Represent} and such that,
  additionally:
\begin{itemize}
\item The map $\pi_\bD\circ u$ is a $g$--fold branched covering
  map with
  all its branch points of order $2$.
\item The map $u$ is holomorphic near the preimages of the branch
  points of $\pi_\bD\circ u$.
\end{itemize}
\end{Lem}

\proof
Construct a map $u_1\co S_1\to W$ representing $A$ as in
\fullref{Lemma:Represent}.  We would like to say that we can then
choose a branched cover $p_{\bD,1}\co S_1\to [0,1]\times\bR$ (mapping
arcs on the boundary appropriately).  This may not, however, be the
case:  suppose, for instance, that $g=2$ and $S_1$ were the
disjoint union of a disk and a surface of genus one with one boundary component.

However, note that $[\Sigma]\in\pi_2(\vy,\vy)$ can be represented by
a map with connected source.  Specifically, let $S_\Sigma$ be
obtained by making small slits in $\Sigma$ along $\alpha_i$ and
$\beta_i$ starting at $y_i\in\vy$ for $i=1,\cdots,g$.  There is
an obvious map $S_\Sigma\to \Sigma$.

Gluing the negative corners of $S_\Sigma$ to the positive corners of
$S_1$ we obtain a connected surface
$S_2$ and map $p_{\Sigma,2}\co S_2\to \Sigma$.  Since $S_2$ is
connected it is
possible to choose a branched covering map $p_{\bD,2}\co S_2\to \bD$ with
appropriate boundary behavior.  Perturbing this map we can assume all
of its branch points have order $2$.  Finally, deforming away the
double points of $p_{\Sigma,2}\times p_{\bD,2}$ and perturbing it to
be holomorphic in appropriate places, we obtain an
embedding satisfying the specified conditions.
\endproof

Now, fix a positive homology class $A$ in $\pi_2(\vx,\vy)$ and a
map $u$ representing $A$ as in the previous lemma.  The map $u$
induces a map $\phi\co \bD\to \Sym^g(\Sigma)$ as follows.
For $a\in \bD$ let $(\pi_\bD\circ u)^{-1}(a) = \{a_1,\ldots,a_g\}$.
Then define $\phi(a)=\{\pi_\Sigma\circ u(a_1),\ldots,\pi_\Sigma\circ
u(a_g)\}$.

There is a one--to--one
correspondence between order $2$ branch 
points of $\pi_\bD\circ u$ and transverse intersections of
$\phi$ with the top--dimensional stratum of the diagonal.  By the
Riemann--Hurwitz formula, $\chi(S)=g\chi(\bD^2) - \phi\cdot\Delta =
g-\phi\cdot\Delta$, so
$\ind(D\dbar)(u)=g-\chi(S)+2e(A+[\Sigma])=\phi\cdot\Delta+2e(A+[\Sigma])$.

This is exactly Rasmussen's formula for the Maslov index.  Thus, we
have shown that $\ind(A+[\Sigma])=\mu(A+[\Sigma])$ for $A$
positive.  But both $\ind$ and $\mu$ are additive, and have
$\ind([\Sigma])=\mu([\Sigma])=2$.  Thus, it follows that
$\ind(A)=\mu(A)$ for all $A$.
\endproof

\begin{Cor}\label{Cor:HFIndex}
In Heegaard Floer homology, the Maslov index of a domain
  $D$ is given by
$$
\mu(D)=n_\vx(D)+n_\vy(D)+e(D).
$$
\end{Cor}
\medskip\textbf{Remark}\qua There are no assumptions on the domain.

\subsection{Index for $A\in\pi_2(\vx,\vx)$}
The following result, which we will use below, is proved by Ozsv\'ath and Szab\'o in~\cite[Proposition 7.5]{OS2} by direct geometrical argument.
\begin{Lem}\label{OSTheorem4.9} If $\vx\in\Ss$ and $A\in\hat{\pi}_2(\vx,\vx)$
  then $\left\langle c_1(\Ss),A\right\rangle=e(A)+2n_\vx(A)$.  Here
  $\langle\cdot,\cdot\rangle$ denotes the
  natural pairing between homology and cohomology, $c_1(\Ss)$ the
  first Chern class of $\Ss$, and $A$ is viewed as an element of
  $H_2(Y)$.
\end{Lem}

The following is completely analogous to~\cite[Theorem 4.9]{OS1}.
\begin{Cor}\label{Index:PeriodicCor}
 Let $P$ be a homology class in $\pi_2(\vx,\vx)$.  Then
 $$
 \ind(P)=\left\langle c_1(\Ss), P\right\rangle + 2n_{\Fz}(P).
 $$
\end{Cor}

\section{Admissibility criteria}\label{Section:Admissibility}
In order to define the differential in our chain complexes it will be
important that for any intersection points $\vx$ and $\vy$, only finitely many homology
classes $A\in\pi_2(\vx,\vy)$ with $\ind(A)=1$ support holomorphic
curves.  For a rational homology sphere, this is automatic: there
are only finitely many homology classes in $\pi_2(\vx,\vy)$.  In
general, following~\cite{OS1}, we use special Heegaard diagrams and
``positivity of domains'' to ensure that only finitely many homology
classes support holomorphic curves.  Our definitions are the same as
theirs.  For the reader's amusement, we provide slightly different
proofs of two of the fundamental lemmas about admissibility.

\begin{Def}(Compare~\cite[Definition 4.10]{OS1})\qua The pointed Heegaard diagram $(\Sigma,\va,\vb,\Fz)$ is
  called \emph{weakly admissible} for the $\Spin^\bC$--structure
  $\Ss$ if every nontrivial periodic domain $P$ with
  $\left\langle c_1(\Ss),P\right\rangle =0$ has both positive and
  negative coefficients.  
\end{Def}
\begin{Def}(Compare~\cite[Definition 4.10]{OS1})\qua
The pointed Heegaard diagram $(\Sigma,\va,\vb,\Fz)$ is
  called \emph{strongly admissible} for the $\Spin^\bC$--structure
  $\Ss$ if every nontrivial periodic domain $P$ with $\left\langle
  c_1(\Ss),P\right\rangle=2n>0$ has $n_{\Fz_i}(P)>n$ for some
  $\Fz_i$.
  
\end{Def}

\medskip\textbf{Remark}\qua  Notice that for any $\Spin^\bC$--structure,
    $c_1(\Ss)$ is an even cohomology class:  $c_1(\Ss)$ is the
    first Chern class of $v^\perp$ for some nonvanishing vector
    field $v$.  Then $c_1(\Ss)\equiv w_2(v^\perp)$ $\textrm{mod
    }2$.  Since $TM$ is trivial and the line field determined by
    $v$ is obviously trivial, $1=(1+w_1(v^\perp)+w_2(v^\perp))$,
    so $w_2(v^\perp)=0$.

We now need two kinds of result.  The first is the finiteness
mentioned just above in the case of weak / strong admissibility.  The
second is that the admissibility criteria can be achieved, and that
any two admissible Heegaard diagrams can be connected by a sequence of
Heegaard moves through admissible diagrams.

First, a few simple observations.  A $\Spin^\bC$--structure is called
``torsion'' if $c_1(\Ss)$ is a torsion homology class.  For a
torsion $\Spin^\bC$--structure, $\langle c_1(\Ss), P\rangle=0$ for
any periodic class $P$.  So, the two definitions of admissibility
agree.  Further, if a Heegaard diagram is weakly (or equivalently strongly)
admissible for some torsion $\Spin^\bC$--structure then it is weakly
admissible for every $\Spin^\bC$--structure.  This point is useful
for computations.

Both admissibility criteria are, obviously, vacuous for a rational
homology sphere.

It will be useful to have equivalent definitions of weak / strong
admissibility:
\begin{Lem} \label{Lemma:AltAdmiss}
Fix a pointed Heegaard diagram $\cH=(\Sigma,\va,\vb,\Fz)$
  and a $\Spin^\bC$--structure $\Ss$.
\begin{itemize}
\item {\rm(Compare~\cite[Lemma 4.12]{OS1})}\qua The diagram $\cH$ is weakly admissible for $\Ss$ if and only
  if there is an area form on $\Sigma$ with respect to which every periodic
  domain $P$ with $\langle c_1(\Ss),P\rangle=0$ has zero signed
  area.  
\item The diagram $\cH$ is strongly admissible for $\Ss$ if there
  is an area form on $\Sigma$ with respect to which every periodic
  domain $P$ with $\langle c_1(\Ss),P\rangle = 2n$ has signed area
  equal to $n$, and with respect to which $\Sigma$ has area $1$.
\end{itemize}
\end{Lem}
\proof
The proofs of the two statements are very similar, and the proof of
the first statement is in~\cite[Lemma 4.12]{OS1}.  We
give here only the proof of the second statement.

Let $\{D_i\}$, $i=1,\cdots,N$ denote the components of
$\Sigma\setminus(\ba\cup\bb)$.  We can view the space of periodic
domains as a linear subspace $V$ of $\bZ^N\subset\bR^N$.
Suppose an area form assigns the area $a_i$ to $D_i$.  Then the
area assigned to $P$ is $P\cdot(a_i)$, the dot product of the
vector $P\in\bR^N$ and the vector $(a_i)$.  Since this is the only
way the area form enters the discussion, we will refer to the vector
$(a_i)$ as the area form.

Suppose there  is an area form $(a_i)$ on $\Sigma$ with respect to which
every periodic domain $P$ with $\langle c_1(\Ss),P\rangle = 2n>0$
has signed area equal to $n$ and $\Sigma$ has area $1$.  Suppose $\langle
c_1(\Ss),P\rangle=2n$.  Then by assumption $P\cdot (a_i)=n$.  So,
$(P-n[\Sigma])\cdot (a_i)=0$.  Hence, $P-n[\Sigma]$ must have some
positive coefficient.  Hence, $P$ must have some coefficient greater
than $n$.

The converse is slightly more involved.  Note that since
$area(-P)=-area(P)$, it suffices to construct an area form with the
desired property for periodic domains with $\langle c_1(\Ss),P\rangle
\geq 0$.

By \fullref{OSTheorem4.9}, the function which assigns to a periodic
domain $P$ the number $\langle c_1(\Ss),P\rangle$ extends to an
$\bR$--linear functional $\ell$ on $V$.  
The map $v\mapsto v-\left(\ell(v)/2\right)[\Sigma]$ gives a linear projection map
$p\co V\to \ker(\ell)$.  Let $V'=p(V).$

Now, we want to choose $a=(a_i)$ orthogonal to $V'$ so that $a_i>0$
for all $i$.  We will show that one can choose such an $a$
presently; for now, assume that such an $a$ has been chosen.  Multiplying
$a$ by some positive real number, we can assume that
$a\cdot[\Sigma]=1$.  Now, for $v\in V$,
$$
a\cdot v = a\cdot p(v) +\left( \ell(v)/2\right) a\cdot [\Sigma]=\ell(v)/2=\left\langle c_1(\Ss),P\right\rangle/2
$$
as desired.

Finally, we need to show such an $a$ exists.  The linear space $V$
is spanned by the periodic domains $P$ with $\langle
c_1(\Ss),P\rangle\geq 0$, so $V'$ is spanned by their
images under $p$.  Every periodic domain $P$ with $\langle
c_1(\Ss), P\rangle=2n\geq0$ has a coefficient bigger than $n=\ell(P)/2$,
so every $p(P)$ has a positive coefficient.  It is also
true that every $p(P)$ has a negative coefficient:  if $\ell(P)=0$
this follows by applying the hypothesis to $-P$; if $\ell(P)>0$
this follows from the fact that $n_\Fz(P)=0$.

Now, we are reduced to showing the following:  let $V'$ be a subspace
of $\bR^N$ such that every nonzero vector in $V'$ has both
positive and negative coefficients.  Then there is a vector orthogonal
to $V'$ with all its entries positive.  The proof of this claim is
a linear algebra exercise; see~\cite[Lemma 4.12]{OS1}.
\endproof

Now we get to the two lemmas justifying the introduction of our
admissibility criteria. 
The following is~\cite[Lemma 4.13]{OS1}.
\begin{Lem}\label{Lemma:WeakAdmis}If $(\Sigma,\va,\vb,\Fz)$ is weakly admissible for $\Ss$
  then for each $\vx,\vy\in\Ss$ and $j,k\in\bZ$ there are only
  finitely many positive homology classes $A\in\pi_2(\vx,\vy)$ with
  $\ind(A)=j$ and $n_\Fz(A)=k$.
\end{Lem}
\proof
If $A,B\in\pi_2(\vx,\vy)$ and $\ind(A)=\ind(B)=j$,
$n_\Fz(A)=n_\Fz(B)=k$ then $P=A-B$ is a periodic domain with $\langle
c_1(\Ss),P\rangle=0.$  So, we must show that there are only finitely
many periodic domains $P$ with $\langle c_1(\Ss),P\rangle=0$ such
that $A+P$ is positive.

Choose an area form on $\Sigma$ so that the signed area of any
periodic domain $P$ with $\langle c_1(\Ss),P\rangle=0$ is zero.
The condition that $A+P$ be positive obviously gives a lower
bound for every coefficient of $P$.  This and the condition that the
signed area of $P$ is zero gives an upper bound for every
coefficient of $P$.  The coefficients are all integers, so the
result is immediate.
\endproof

The following is~\cite[Lemma 4.14]{OS1}.
\begin{Lem}\label{Lemma:StrongAdmis}If $(\Sigma,\va,\vb,\Fz)$ is strongly admissible for
  $\Ss$ then for each $\vx,\vy\in\Ss$ and $j\in\bZ$ there are
  only finitely many positive homology classes $A\in\pi_2(\vx,\vy)$ with
  $\ind(A)=j$.
\end{Lem} 
\proof
Fix a homology class $A\in\pi_2(\vx,\vy)$ with $n_\Fz(A)=0$ and
$\ind(A)=j_0$.  Then any other homology class $B\in\pi_2(\vx,\vy)$
can be written as $A+P+k[\Sigma]$ for some integer $k$ and
periodic domain $P$.  We have $\ind(B)=j_0+k+\langle
c_1(\Ss),P\rangle$.  If we assume that $\ind(B)=j$ then $\langle
c_1(\Ss),P\rangle = j-j_0-2k$.

Fix an area form such that the area of $\Sigma$ is $1$ and the
area of any periodic domain $P$ is $\frac{1}{2}\langle
c_1(\Ss),P\rangle$.  Then, the area of $P$ is $\frac{j-j_0}{2}-k$.

If we impose the condition that $B$ be positive
then we automatically get lower bounds for every coefficient of
$P$ (which are independent of $k$).  Note that $k\geq 0$, since $k=n_{\Fz_i}(B)$ for some $i$.
The condition that the area of $P$ be $\frac{j-j_0}{2}-k\leq
\frac{j-j_0}{2}$ and the lower bound for the coefficients of $P$
gives an upper bound for the coefficients of $P$, independent of $k$.  This completes
the proof.
\endproof

The following is~\cite[Lemma 5.8 and Proposition 7.2]{OS1}.  We refer
the reader there for its (somewhat involved but essentially
elementary) proof.
\begin{Prop}\label{Prop:MaintainAdmis} Fix a 3--manifold $Y$ and $\Spin^\bC$--structure
  $\Ss$ on $Y$.
\begin{enumerate}
\item There is a weakly (respectively strongly) admissible Heegaard
  diagram for $\Ss$.
\item Suppose that $\mathcal{H}_1=(\Sigma,\va,\vb,\Fz)$ and
  $\mathcal{H}_2=(\Sigma',\va',\vb',\Fz')$ are weakly (respectively
  strongly) admissible Heegaard diagrams for $\Ss$.  Then there is a
  sequence of pointed Heegaard moves (ie, Heegaard moves
    supported in the complement of $\Fz$) connecting $\mathcal{H}_1$
  to $\mathcal{H}_2$ such that each intermediate Heegaard diagram is
  weakly (respectively strongly) admissible for $\Ss$.
\end{enumerate}
\end{Prop}

\section{Orientations}
\label{Section:Orientations}
In order to be able to work with $\bZ$ coefficients, we need to be
able to choose orientations for the moduli spaces $\hcM^A$ in a
coherent way.  First we need to know that each $\hcM^A$ (or
equivalently, each $\cM^A$) is orientable.  
Then we will discuss what we mean by a coherent orientation and why
such orientations exist.  This is all somewhat technical, and we will
sometimes supply references rather than details. 

Suppose that we have chosen an almost complex structure $J$ that achieves transversality for the
moduli space $\cM^A$.  Then for $u\in\cM^A$ the tangent space
$T_u\cM^A$ is naturally identified with the kernel
$\ker(D_u\dbar)$ of the linearized $\dbar$ operator at $u$.  In
fact, the
spaces $\ker(D_u\dbar)$ fit together to form a vector bundle over
$\cM^A$ naturally isomorphic to $T\cM^A$.  So, orienting $\cM^A$
is the same as trivializing the top exterior power of the vector
bundle $\ker(D\dbar)$ over $\cM^A$.

Rather than working with $\ker(D\dbar)$ it is better to work with the
line bundle $\cL=\det(D\dbar)$ which is defined to be the tensor product
of the top exterior power of $\ker(D\dbar)$ with the dual of the top
exterior power of $\coker(D\dbar)$.  (This is the ``determinant line
bundle of the virtual index bundle of the $\dbar$--operator.'')
Note that if $J$ achieves transversality at a curve $u$ then $\cL_u$ is just the
top exterior power of $\ker(D_u\dbar)$.

To keep the exposition clean we will assume that our sources are
stable, ie, have no components which are twice--punctured disks.
Let $\cT$ denote the \emph{Teichm\"uller configuration space} of pairs
$(j,u)$ where $j$ is a complex structure on $S$ (ie, an
element of the Teichm\"uller space of Riemann surfaces) and $u:S\to
W$ is a map satisfying (\textbf{M1}), (\textbf{M3}) and (\textbf{M4}) and which is asymptotic to the planes
$\{x_i\}\times[0,1]\times\bR$ or $\{y_i\}\times[0,1]\times\bR$ at
the appropriate punctures.  Say $(j,u)$ and $(j',u')$ are
equivalent if there is an isomorphism of Riemann surfaces
$\phi:(S,j)\to (S',j')$ so that 
$$
\xymatrix{S'\ar[dr]_{u}\ar[rr]^\phi & & S\ar[dl]^u\\
& W &
}
$$
commutes.  Let $\cB$ be the quotient of $\cT$ by this equivalence
relation, so $\cB$ is the \emph{moduli configuration space of maps to $W$.}
The reason to work with $\cL$ is that the bundle $\cL$ is defined
over all of $\cB$.

Different topological types of $S$
correspond to different components of $\cB$.  So from now on we
restrict attention to the subspace $\cB^S$ corresponding to maps
from a single topological type of source $S$.

Note in particular that we can talk about choosing an orientation over the
homology class $A$ even if $\cM^A$ is empty.  This will be useful
when we discuss coherence.

The determinant line bundle $\cL$ is, in our case, always
trivial.  To prove this we combine constructions from~\cite{EES}
and~\cite{Liu}.

Under our stability assumption, we have a fiber bundle
$$
\xymatrix{
\Map(S,W)\ar[r] & \cB^S\ar[d] \\ & \cM_S
}
$$
where $\Map(S,W)$ consists of maps $S\to W$ satisfying
(\textbf{M1}), (\textbf{M3}) and (\textbf{M4}) with appropriate asymptotics, and
$\cM_S$ corresponds to the moduli space of conformal structures on $S$. 

Call a space $X$ \emph{homotopy discrete} if every connected
component of $X$ is contractible.  The following proposition is
somewhat stronger than we need.  It is, however, of some independent
interest, and will be mentioned again in \fullref{Section:OtherRemarks}.

\begin{Prop}\label{Prop:Contractible}The space  $\Map(S,W)$ is
  homotopy discrete.
\end{Prop}
\proof
The space $\Map(S,W)$ is the product $\Map(S,\Sigma)\times
\Map(S,[0,1]\times\bR)$.  The space $\Map(S,[0,1]\times\bR)$ is
convex, so it suffices to prove that $\Map(S,\Sigma)$ is homotopy
discrete.  (Here, both $\Map(S,[0,1]\times\bR)$ and $\Map(S,\Sigma)$
refer to spaces of maps with certain obvious boundary and asymptotic conditions.)

There is a fibration 
$$\Map(S,\Sigma)\to \Omega(\alpha_1)\times\cdots\times\Omega(\alpha_g)\times\Omega(\beta_1)\times\cdots\times\Omega(\beta_g) $$
given by restricting a map to the boundary and identifying the space of
paths in $\alpha_i$ (or $\beta_j$) with endpoints $x_i$ and
$y_i$ with $\Omega(\alpha_i)$ (or $\Omega(\beta_j)$).  (Here,
$\Omega$ denotes the based loop space.)  Since
$\Omega(\alpha_1)\times\cdots\times\Omega(\beta_g)$ is homotopy discrete
it suffices to prove that each fiber of the
fibration is homotopy discrete.  Let $\Map(S,\Sigma;\bdy)$ denote a fiber of
the fibration.

Let $\Map'(S,\Sigma)$ denote the space of all maps $S\to\Sigma$ in
the homotopy class of $\Map(S,\Sigma;\bdy)$ with no boundary
conditions.  There is a fibration
$$
\xymatrix{ \Map(S,\Sigma;\bdy) \ar[r] & \Map'(S,\Sigma)\ar[d]\\
& \Map(\bdy S,\Sigma)
}.
$$
Different fibers of this fibration are homotopy equivalent.  So, we
can replace\break $\Map(S,\Sigma;\bdy)$ with the space
$\Map'(S,\Sigma;\bdy)$ of maps with all $2g$ punctures of $S$
mapped to a single point $p\in \Sigma$ and the boundary arcs mapped
to a fixed list $C_1,\cdots,C_{2g}$ of circles.  

There is a fibration
$$
\xymatrix{ \Map'(S,\Sigma;\bdy)\ar[r] & \Map'_*(S,\Sigma) \ar[d]\\
& \left(\Omega\Sigma\right)^{2g}}
$$
where $\Map'_*(S,\Sigma)$ denotes the component of the space of based
maps $S\to \Sigma$ containing $\Map'(S,\Sigma;\bdy)$.  Since
$\Omega\Sigma$ is homotopy discrete, it
suffices to prove that $\Map'_*(S,\Sigma)$ is homotopy discrete.
This follows from the following lemma. 

\begin{Lem}Let $K$ be a $K(\pi,1)$ and $X$ any finite, connected CW
  complex.  Fix basepoints in $K$ and $X$.  Then the space of based
  maps from $X$ to $K$ is homotopy discrete.
\end{Lem}
\proof
We may assume the zero--skeleton $X^{(0)}$ of $X$ consists of just the
basepoint, so the one--skeleton $X^{(1)}$ consists of a bouquet of circles.
Let $\Map_*(X,K)$ denote based maps from $X$ to $K$.  There is a
fibration $\Map_*(X,K)\to \Map_*(X^{(1)},K)=(\Omega K)^N$ (where $N$
is the number of $1$--cells in $X$).  The base is homotopy discrete since $K$ is a $K(\pi,1)$.  So, it
suffices to prove that the nonempty fibers are homotopy discrete.  Each fiber $F_1$ consists of maps of $X$ to $K$ fixed
on $X^{(1)}$.

The proof now proceeds by induction over the skeleta of $X$.  There
is a fibration 
$$F_1\to \prod_{\textrm{2--cells}} \Map_*(\bD^2,K;\bdy)
$$
where $\Map_*(\bD^2,K;\bdy)$ denotes (based) maps of a disk with the boundary
mapped as specified by the map of the one--skeleton.  The base is
homotopy discrete, as one deduces from the
fibration
$$
\xymatrix{
\Map_*(\bD^2,K;\bdy)\ar[r] & \Map_*(\bD^2,K)\ar[d]\\
& \Omega(K)
}
$$
in which the total space is contractible and the base space is
homotopy discrete.  Thus, it suffices to show
that any fiber $F_2$ of $F_1\to \prod_{\textrm{2--cells}}
\Map_*(\bD^2,K;\bdy)$ is contractible.  Notice that each $F_2$ consists of
maps of $X$ specified on the two--skeleton.

Proceeding as before, we have a fibration $F_3\to
\prod_{\textrm{3--cells}} \Map_*(\bD^3,K;\bdy)$, now with contractible base.  It suffices to prove the
fiber $F_4$ is contractible, and so on.
Since $X$ is finite, the process terminates at $F_n$ for
$n=\dim(X)$.  
\endproof
Since $\Sigma$ is a $K(\pi,1)$, this completes the proof of the proposition.
\endproof

If the space $\cM_S$ is contractible, or even just has trivial
$H^1$, then we are finished:  we then have $H^1(\cB^S)$ trivial
and hence all line bundles over $\cB^S$ are orientable.  This
is, in fact, the case if $S$ is a union of disks with boundary
punctures.  In general, however, $\cM_S$ can have interesting topology.

\medskip\textbf{Remark}\qua  It is not hard to show, by an argument similar to but
simpler than the proof of \fullref{Prop:Contractible}, that
the configuration spaces of disks in the original construction of
Heegaard Floer homology \cite{OS1} have trivial
$\pi_1$, proving orientability of the moduli spaces in that setting.

\begin{Prop}\label{Prop:Orientable}The determinant line bundle $\cL$ over $\cB^S$ is trivial.
\end{Prop}
\proof
By \fullref{Prop:Contractible}
for a given source $S$ and homology class $A$, the
space of possible maps of $\bdy S$ is contractible.
Then, it follows from the argument of~\cite[Lemma 6.37]{Liu}
or~\cite[Lemma 3.8]{EES}, which we
sketch briefly, that the line bundle $\cL$ is trivial.  
Fix a map $u_0$ of a collar neighborhood $C$ of $\bdy S$.
(By a collar neighborhood, we mean that $\bdy C$ is the union of
$\bdy S$ and a collection of circles, and each component of $C$ is
an annulus with punctures on one of its boundary circles.)  We
restrict to maps $u:S\to W$ agreeing with $u_0$ on $C$.  (That
it suffices to consider such maps is the only place we use anything
special about our situation.  If we choose $u_0$ to be an embedding,
say, it makes sense to say $u$ agrees with $u_0$ on a collar
neighborhood of the boundary, even though $\cB^S$ is a fiber bundle
and not a product of $\cM_S$ and $\Map(S,W)$.)

Let
$\bdy_{in}(C)$ denote the ``interior boundary'' of $C$, ie, the
circles without punctures.

As in \fullref{Section:Transversality}, we are interested in the kernel and cokernel of 
$$
D\dbar: L^{p,d}_1(u^*TW,\bdy)\oplus\bR^{2g} \to L^{p,d}\left(\Lambda^{0,1}u^*TW\right).$$
Let $S'$ denote the surface obtained from $S$ by collapsing $\bdy_{in}
C$.  The operator $D\dbar$ induces an operator $D'\dbar$ between
bundles over $S'$, and it is not hard to see that the determinant
lines of $D\dbar$ and $D'\dbar$ are naturally identified.  Let
$C'$ be the image of $C$ in $S'$, so $S'$ consists of $C'$
and some closed components.  The determinant line of $D'\dbar$ is
the tensor product of the determinant lines of the restriction of
$D'\dbar$ to each component of $S'$.

Over the closed components of $S'$ the determinant line has a
natural ``complex'' orientation.  Fix once and for all trivializations
of the determinant lines of $D'\dbar$ over the components of
$C'$.  (Since $u|_C= u_0$, the restriction of $D'\dbar$ to
$C'$ does not depend on $u$.)  Then, these trivializations induce
an orientation of the determinant of $D'\dbar$, and hence of
$D\dbar$, independent of $u$.  This proves orientability of $\cL$.
\endproof

\medskip\textbf{Remark}\qua  For comparison, in~\cite{OS1}, Ozsv\'ath and Szab\'o only use the fact
that the tangent bundles to the Lagrangian submanifolds in question are trivial.
In~\cite{EES}, although their setup is not quite the same as
standard symplectic field theory, T. Ekholm, J. Etnyre and M. Sullivan prove orientability under the
assumption that the Lagrangian manifolds in question are $\Spin$.

The coherence we want for our orientation system is the following.
Suppose that we have maps $u\co S\to W$ and $u'\co S'\to W$ in homology
classes $A\in\pi_2(\vx,\vx')$ and $A'\in\pi_2(\vx',\vx'')$, respectively.
Then by gluing the positive corners of $u$ to the negative corners of
$u'$ (in some fixed, concrete way) one can construct a
$1$--parameter family of curves $u\natural_r u'\co S\natural_r S'\to W$.  One obtains, thus, an
inclusion $\Map_A(T,W)\times \Map_A'(T,W)\times(R,\infty)\to
\Map_{A+A'}(T,W)$.  This inclusion is covered by a map of determinant
line bundles.  Coherence means that this map is
orientation--preserving.  A few more details can be found in the proof of \fullref{Prop:OrId}, and a complete discussion in any of~\cite[Section 1.8]{Ya1},~\cite{BM}
or~\cite{EES}.

There are strong general results about the existence of coherent
orientations (see~\cite{Ya1},~\cite{BM} or~\cite{EES}).  In our case,
however, we can construct them quite concretely:

As in \cite[Section 3.6]{OS1}, given a Heegaard diagram,
a $\Spin^\bC$--structure $\Ss$, and an intersection point
$\vx_0\in\Ss$, a \emph{complete set of paths for 
  $\Ss$} is a choice of homology class
$T_i\in\pi_2(\vx_0,\vx_i)$ for every other $\vx_i\in\Ss$.
Fix $\vx_0\in\Ss$, a representative $P_j$ of each homology class in
some basis for $\hat{\pi}_2(\vx_0,\vx_0)=H_2(Y)$, and a
complete set of paths for $\Ss$.
By the gluing properties of coherent orientations described two
paragraphs earlier, specifying
a particular coherent orientation of the moduli space is equivalent to
specifying an orientation over each $T_i$,
each $P_j$, and over $[\Sigma]\in\pi_2(\vx_0,\vx_0)$.  

Over $[\Sigma]$ there is a
canonical orientation, given by viewing the map
$\Sigma\coprod\left(\amalg_{k=1}^g\bD^2\right)\to
\Sigma\times\bD$ given by $id\times\{0\}$ on $\Sigma$ and
$\{x_{0,i}\}\times id$ on the $i^{\textrm{th}}$ $\bD^2$ as
positively oriented.  We shall always use this orientation over
$[\Sigma]$.  Over the $T_i$ and $P_j$ orientations can be
chosen arbitrarily.  This is exactly as in \cite{OS1}.

Following~\cite{OS1}, we shall let $\frako(A)$ denote the choice of
orientation over the homology class $A$.

A canonical choice of coherent orientation system is specified in
\cite{OS2}.  An analogous construction presumably works in our case as
well.  This is not, however, necessary for the current paper.

\section{Bubbling}
\label{Section:Bubbling}
By a \emph{holomorphic building} in $W$ we mean a list $v_1,v_2,\cdots,v_k$ of holomorphic curves in $W$, defined up to translation in $\bR$, such that the asymptotics of $v_i$ at $+\infty$ agree with the asymptotics of $v_{i+1}$ at $-\infty$.  We call $k$ the \emph{height} of the holomorphic building, the $v_i$ the \emph{stories} of the building, and $i$ the \emph{level} of $v_i$. It is proved in~\cite[Theorem 10.1]{Ya2} that any sequence $\{u_i\}$ of holomorphic curves in
$W$ (possibly disconnected) in a given homology class has a subsequence converging to some (possibly nodal)
holomorphic building.  

The meaning of convergence is, as usual for holomorphic curves, somewhat involved.  Roughly, convergence to a several story holomorphic building means that some parts of $u_i$ go to infinity (in the $\bR$--factor) with respect to other parts.  In the source, such \emph{level splitting} corresponds to the degeneration of the complex structure along a collection of disjoint arcs or, in principle, circles.  It is also possible for the source to degenerate along a collection of circles or arcs without level splitting occuring; this corresponds to formation of nodes as in traditional Deligne--Mumford theory.  In principle, disks or spheres can also bubble off.  See~\cite[Section 7]{Ya2} for a precise definition of convergence taking into account all of these possibilities.

For the definition of the Floer homology in subsequent sections to work, we need the following:

\begin{Prop}\label{CompactProp}  Fix an almost complex structure $J$ on $W$ satisfying {\rm(\textbf{J1})--(\textbf{J5})} and achieving transversality.
 Let $u_i\co S\to W$ be a sequence of
  holomorphic curves in a homology
  class $A$ satisfying {\rm(\textbf{M0})--(\textbf{M6})}
  and converging to some holomorphic building
  $v$.  Assume that
  $\ind(A)\leq 2$ and that $A\neq[\Sigma]$.  Then, each story
  $v_j$ of $v$ satisfies
  {\rm(\textbf{M0})--(\textbf{M6})}. 
\end{Prop}
Here, by $[\Sigma]$ we mean the homology class with $n_{\Fz_i}=1$
for every $i$.

\proof
First we check that Deligne--Mumford type degenerations which do not
correspond to level splitting are impossible.  Bubbling of spheres or
disks is impossible because $\pi_2(W)$ and $\pi_2(W,C_\alpha\cup
C_\beta)$ vanish.  This leaves us to rule out Deligne--Mumford type
degenerations along non--contractible curves and arcs.

Given transversality, Deligne--Mumford type degenerations in the
interior of $S$ are
prohibited for a generic choice of $J$ because they have
codimension 2 in the space of all holomorphic curves.  That is,
  these degenerations have codimension 2 after quotienting by
  translation.  Hence, they only occur if $\ind(A)\geq 3$.  This
rules out all but one kind of degeneration along the interior:  bubbling of
a copy of $\Sigma\times(s_0,t_0)$.  (Recall that our complex
structure does not achieve transversality for maps $u$ with
$\pi_\bD\circ u$ constant.)  However,
\fullref{Cor:Index} shows that if a copy of $\Sigma$ were
to bubble off then the remainder of $v$ would have to be a
collection of trivial disks, and so $n_{\Fz_i}(A)=1$ for every $i$. Bubbling more than one copy of $\Sigma$ or a multiply covered copy of $\Sigma$ is prohibited by \fullref{Cor:Index}.

So, we need to check that cusp degenerations (ie, degenerations
along arcs with boundary on $\bdy S$) are impossible.  

Suppose that a cusp degeneration occurs. Let $A$ denote an arc in $S$ which collapses. Since different components of $\bdy S$ are mapped by $u$ to different cylinders in $W$, both endpoints of $A$ must lie on the same boundary component $C$ of $S$.
Without loss of generality, suppose that $C$ is
mapped by the $u_i$ to 
$\alpha_1\times\{1\}\times\bR$.  Let $S'$ denote the nodal surface
obtained from $S$ by collapsing the arcs along which the complex
structure degenerates.  In the limit curve $v$ there
will be more than one boundary component mapped to
$\alpha_1\times\{1\}\times\bR$.  Let $C'$ denote all the boundary
components mapped to $\alpha_1\times\{1\}\times\bR$.  Consider the
map $\pi_\bD\circ
v|_{C'}\co C'\to i\bR$.  One of the boundary components in $C'$ will
have a local (in fact, global) maximum.  But this and the open mapping
theorem imply that the map $\pi_\bD\circ v$ must be constant near that 
boundary
component, hence constant on a connected component of $S'$.  In
particular, this implies that all of the boundary of $S'$ is mapped
to the $\alpha$--circles.

However, the $\alpha$--circles are non--separating, so the component
of $S'$ on which $\pi_\bD\circ v$ is constant must be mapped
diffeomorphically onto $\Sigma$.  As above, it then follows from our
index computation that the rest of $v$ must consist of $g$ trivial
disks, and $A=[\Sigma]$.

Thus, we have proved the each $u_i$ satisfies (\textbf{M0}).
Condition (\textbf{M1}) is automatically satisfied.

The only way that $\pi_\bD\circ v$
could be constant on some component of $S$ would be as a result of
bubbling, which we showed is prohibited.
This implies (\textbf{M2}) for the limit curves.

Condition (\textbf{M4}) and the condition that there are exactly
$2g$ punctures on each story also follow from the maximum modulus,
or open mapping, theorem:  applying the open mapping theorem to
$\pi_\bD\circ u_j$, one sees that the restriction of $\pi_\bR\circ u_j$ to
any component of $\bdy S\setminus\{p_1,\cdots,p_g,q_1,\cdots,q_g\}$
must be monotone, so no new Reeb chords could form as $j\to\infty$.

Since projection onto $[0,1]\times\bR$ is holomorphic, the open
mapping theorem prohibits any new boundary components from forming.
The maximum modulus theorem prohibits boundary components from
disappearing, so (\textbf{M3}) is satisfied by the limit curves.

That condition (\textbf{M5}) is satisfied by the limit curves is part
of the statement of the compactness theorem~\cite[Theorem 10.1]{Ya2}.

Let $S_j$ denote the
source of $v_j$.  Suppose that $v$ has height $\ell$.  Then
since all degenerations in $S$ correspond to level splitting,
$\chi(S)=(1-\ell)g + \sum_{j=1}^\ell\chi(S_j).$  Let $D_j$
be the domain in $\Sigma$ corresponding to $v_j$.  Then we have
$e(D)=\sum_{j=1}^\ell e(D_j)$.  So, $\ind(u_i)=g-\chi(S)+2e(D)=\ell
g - \sum_{j=1}^\ell\left(\chi(S_j) + 2e(D_j)\right) = \sum_{j=1}^\ell
\ind(v_j).$  That is, the index formula adds over levels.

Finally, near any immersed curve with the equivalent of $k$ double
points and with respect to which the complex structure achieves
transversality there is a $2k$--dimensional family of
embedded holomorphic curves.  This shows that the $v_j$ must all be
embedded.
\endproof

For an almost complex structure $J$ which achieves transversality, 
it follows from the previous proposition, \fullref{Gluing:Cylindrical}
and the compactness Theorem 10.1 from \cite{Ya2} that:
\begin{Cor}\label{CompactCor}
The space of holomorphic curves satisfying {\rm(\textbf{M0})--(\textbf{M6})} in a given
homology class $A$ (modulo translation in the $\bR$--factor of
$W$) with $\ind(A)\leq 2$ and $A\neq[\Sigma]$
forms a smooth manifold $\hcM^A$ of dimension $\ind(A)-1$ which is
the interior of a
compact manifold with boundary $\overline{\hcM^A}$.  The boundary of
that manifold $\bdy\overline{\hcM^A}$
consists of all multi--story holomorphic buildings each component of which satisfies {\rm(\textbf{M0})--(\textbf{M6})} in the homology class $A$. 
\end{Cor}

\section{Chain complexes}\label{Section:ChainComplexes}
Here we define the four basic chain complexes used by
Ozsv\'ath--Szab\'o -- $\widehat{CF}$, $CF^\infty$, $CF^-$ and
$CF^+$.  Once these are defined $CF^\pm_{\textrm{red}}$ are defined
in exactly the same way as in~\cite{OS1}.  Generalizing our
definitions to include the twisted theories of~\cite[Section 8]{OS2} is
straightforward.

Fix: 
a pointed Heegaard diagram $(\Sigma,\va,\vb,\Fz)$, an almost complex structure $J$ satisfying (\textbf{J1})--(\textbf{J5}) and achieving transversality, and a
coherent orientation system.  We shall assume that the Heegaard
diagram satisfies the weak or strong admissibility criterion as necessary.  

Fix a $\Spin^\bC$--structure $\Ss$.

First we define the chain complex $\widehat{CF}$.  Suppose that our Heegaard diagram
is weakly admissible for $\Ss$.
Recall that, by \fullref{Lemma:WeakAdmis}, for $\vec{x}, \vec{y}\in\Ss$,
there are only finitely many
$A\in\hat{\pi}_2(\vx,\vy)$ and
$n_{z_i}(A)\geq 0$ for all $i$ and $\ind(A)=1$.  Let
$\widehat{CF}(Y,\Ss)$ be 
the free Abelian group generated by the intersection points
$\vec{x}$.  

Define
$$\hcM(\vx,\vy)=\bigcup_{A\in \hat{\pi}_2(\vx,\vy),
  \textrm{ }\ind(A)=1}\hcM^A$$

We define a boundary operator on $\widehat{CF}$
by
$$
\bdy \vx = 
\sum_{\vy} \left(\#\hcM(\vx,\vy)\right)\vy.
$$
This sum is finite by the admissibility criterion and
\fullref{CompactCor}.
\begin{Lem}\label{Lemma:HatDefined}
With this boundary operator, $\widehat{CF}$ is a chain complex.
\end{Lem}
\proof
Fix $\vx$, $\vz$.  We will show that
the coefficient of $\vz$ in $\bdy^2\vx$ is zero.  Consider the space
$$\hcM_2(\vx,\vz) = \bigcup_{A\in \hat{\pi}_2(\vx,\vz),
  \textrm{ }\ind(A)=2}\hcM^A$$ with compactification
  $\overline{\hcM_2(\vx,\vz)}$.  (Note that by \fullref{Lemma:WeakAdmis} there are only finitely
  many classes $A$ with $\ind(A)=2$ and $\hcM^A\neq\emptyset$,
  so $\overline{\hcM_2(\vx,\vz)}$ is
  compact.)  From \fullref{CompactCor},
  $\overline{\hcM_2(\vx,\vz)}$ is a 1--manifold with boundary and
  $\bdy\left(\overline{\hcM_2(\vx,\vz)}\right)$ consists of broken trajectories
  connecting $\vx$ to $\vz$.   
Thus, $0=\#\left(\bdy\overline{\hcM(\vx,\vz)}\right)$ is the
  coefficient of $\vz$ in $\bdy^2\vx$.  Note that the definition of a coherent orientation system is chosen
  exactly to make this argument work.
\endproof

Next we define $CF^\infty$.  Assume the strong admissibility criterion is
satisfied.  The chain group of $CF^\infty(Y;\Ss)$ is freely
generated by pairs $[\vx,n]$ where $\vx\in\Ss$ and $n\in\bZ$.  

The boundary operator is given by
$$
\bdy [\vx,n]=\sum_{\vy}\sum_{ \substack{
      A\in\pi_2(\vx,\vy)\\\ind(A)=1}}
\left(\#\hcM^A\right)[\vy,n-n_\Fz(A)]
$$
The coefficient of each $[\vy,m]$ is finite by
\fullref{Lemma:StrongAdmis} and \fullref{CompactCor}.

\begin{Lem}\label{Lemma:InfinityDefined}Suppose that $J|_{s=0}$ and
  $J|_{s=1}$ have been chosen appropriately (in a sense to be
  specified in the proof) and that $J$ achieves transversality for
  all holomorphic curves of index $\leq1$.  Then
  $CF^\infty$ is a chain complex.
\end{Lem}
\proof
The proof is almost exactly the same as for $\widehat{CF}$.  The only
nuance is that the homology class $[\Sigma]$ was an exception to
\fullref{CompactCor}.  One resolution of this difficulty is the
following:

Recall that an annoying curve is a curve is a curve $u\co S\to W$ with a component on which $\pi_\bD\circ u$ is constant.  The difficulty with \fullref{CompactCor} for the homology class $[\Sigma]$ is the possibility of a curve degenerating to a collection of trivial disks together with a copy of $\Sigma$ mapped to a constant by $\pi_\bD\circ u$.

Suppose that instead of considering almost complex structures
satisfying (\textbf{J1})--(\textbf{J5}) we considered the broader
class of almost complex structures satisfying
(\textbf{J1})--(\textbf{J4}) and (\textbf{J5$'$}).  All the results
proved so far would
still hold.  (The only potential issue is the proof of
\fullref{CompactProp}, but the requirement that $\xi$ be horizontal
near $\ba\cup\bb$ is sufficient for that proof.)  The transversality
result would then hold for all
holomorphic curves except for annoying curves mapped entirely into
$\bdy W$.  Thus, annoying curves mapped into the interior of $W$, which are
non--generic, would cease to exist.  

To eliminate annoying curves mapped into $\bdy W$ we adopt an idea
from~\cite{OS1}.  Suppose that a sequence of holomorphic curves
$\{u_j\}$ converges to a holomorphic curve $u$ containing an annoying
curve, which is mapped by $u$ to
$p\in\ba\times\{1\}\times\bR$, say.  Then, rescaling near $p$ we obtain
from $\lim_{j\to\infty}\pi_\bD\circ u_j$ a $g$--fold branched
covering map $(S_0,\bdy S_0)\to(\bH,\bR)$
(where $\bH$ denotes the upper half plane).  Here, $S_0$ is a
surface obtained by cutting $\Sigma$ along the $\alpha$--circles.

Suppose that the complex structure $J$ is given by
$j'_\Sigma\times j_\bD$ at $s=1$.  We will show that for
appropriate choice of $j'_\Sigma$ there is no map $(S_0,\bdy
S_0)\to (\bH,\bR)$.  Choose curves on $\Sigma$ not intersecting the
$\alpha$--circles whose complement in
$\Sigma$ is a disjoint union of punctured tori (with one
$\alpha$--circle contained in each).  Let $\{j'_{\Sigma,n}\}$ be a
sequence of complex structures on $\Sigma$ obtained by stretching
$\Sigma$ along the chosen curves.  So, as $n\to\infty$, $\Sigma$
degenerates to a wedge sum of tori.

Suppose that for all $n$ there were sequences of holomorphic curves
converging to annoying curves.  Then for each $n$ we obtain
a map $\pi_n\co (S_{0,n},\bdy S_{0,n})\to (\bH,\bR)$, where $S_{0,n}$
is obtained from $(\Sigma,j'_{\Sigma,n})$ by cutting along the
$\alpha$--circles.  Choosing a convergent subsequence of the
$\pi_n$, in the limit we obtain a $g$--fold covering map from a
disjoint union of $g$ punctured tori to $\bH$.  Such a map clearly
does not exist.

So, for large enough $n$, if $J$ agrees with
$j'_{\Sigma,n}\times j_\bD$ for $s=1$ then there are no annoying
curves mapped by $\pi_\bD$ to $\{1\}\times\bR$.

A similar argument shows that if we choose $J|_{s=0}$ appropriately then
there are no annoying curves mapped by $\pi_\bD$ to
$\{0\}\times\bR$.

It follows that, with respect to a generic complex
structure extending the specified $J|_{s=0}$ and $J|_{s=1}$ and
satisfying (\textbf{J1})--(\textbf{J4}) and (\textbf{J5$'$}),
$\bdy^2=0$, by the same argument as for $\widehat{CF}$ above.

Now, let $J$ be an almost complex structure on $W$ satisfying
(\textbf{J1})--(\textbf{J5}), extending the specified $J|_{s=0}$ and
$J|_{s=1}$ and achieving transversality for
holomorphic curves of index $1$.  Let $\{J_n\}$ be a sequence of
almost complex structures satisfying (\textbf{J1})--(\textbf{J4}) and
(\textbf{J5$'$}) and achieving transversality which converges to $J$.

Let $\bdy_J$ denote the boundary map in $CF^\infty$ computed with
respect to $J$, $\bdy_{J_n}$ the boundary map computed with
respect to $J_n$.  Given a finite collection of
homology classes $\{A_j\in\pi_2(\vx_j,\vy_j)\}$ such that $\ind(A_j)=1$
for all $j$, there is some $N$
such that for $n>N$, $\hcM^{A_j}_{J_n}\cong\hcM^{A_j}_J$.  So,
since $\bdy_{J_n}^2=0$ for all $n$, $\bdy_J^2=0$.
\endproof

Note that in \fullref{Section:Isotopy} we will show that for any
two choices of (generic) almost complex structure on $W$
the pairs $(CF^\infty,\bdy)$ are chain homotopy equivalent.  This
implies in
particular that the restriction on $J|_{s=0}$ and
$J|_{s=1}$ are unnecessary.  Until then, we shall assume that
$J|_{s=0}$ and $J|_{s=1}$ have been chosen so that the preceding proof works.

Since $n_\Fz(A)\geq 0$ for any $A$ supporting a holomorphic disk,
$CF^\infty$ has a subcomplex $CF^-$ generated by the $[\vx,n]$
with $n<0$.  The quotient complex we denote $CF^+$.  The
homologies of $CF^\infty$, $CF^+$, $CF^-$ and $\widehat{CF}$ we
denote by $HF^\infty$, $HF^+$, $HF^-$ and $\widehat{HF}$
respectively.

For computations, it is helpful to observe that $CF^+$ is defined
even if one only assumes weak admissibility:  the sum in the
definition of $\bdy[\vx,n]$ only involves the $[\vx,n-j]$ for
$0\leq j\leq n$.  The weak admissibility criterion implies that for
each of these finite collection of $j$'s, there are at most finitely
many homology classes for which the moduli space is nonempty.

There is a natural action $U\co CF^\infty(Y;\Ss)\to
CF^\infty(Y;\Ss)$ given by $U[\vx,i]=[\vx,i-1]$.  This action
obviously descends to $HF^\infty$, making $HF^\infty(Y;\Ss)$ into
a module over $\bZ[U,U^{-1}]$.  Further, the action of $U$
preserves $CF^-$, so $HF^-$ and $HF^+$ are modules over
$\bZ[U]$.

There are relative gradings on all four homology theories.  On
$\widehat{CF}$, define $gr(\vec{x},\vec{y})=\ind(A)$ for any
$A\in\hat{\pi}_2(\vx,\vy)$.
It follows from \fullref{Index:PeriodicCor} that this
relative grading is
defined $\mod n$, where $n=\gcd\{ \langle c_1(\Ss),A\rangle \}$ for
$A\in H_2(Y)$. Obviously the boundary map lowers the relative
grading by $1$, so the relative grading descends to $\widehat{HF}(Y;\Ss)$. 

Similarly for $CF^\infty(Y,\Ss)$, $CF^-(Y,\Ss)$ and
$CF^+(Y,\Ss)$ there are relative mod $n$ gradings, $n=\gcd\{\langle c_1(\Ss),A\rangle\}$,
given by $gr([\vec{x},i],[\vec{y},j])=\ind(A)+2(i-j)$,
where $A\in\hat{\pi}_2(\vx,\vy)$.  As
before, this relative grading descends to $HF^\infty(Y,\Ss)$, $HF^+(Y,\Ss)$
and $HF^-(Y,\Ss)$.  Note that the action of $U$
lowers the relative grading by $2$.

The chain complexes we have defined depend on the choice of coherent
orientation system (see \fullref{Section:Orientations}).  It turns out, however, that some orientation
systems give isomorphic chain complexes.  Recall that the orientation
system was given by specifying orientations over a complete set of
paths $T_i$ and a basis $P_j$ for
$\hat{\pi}_2(\vx_0,\vx_0)$ (for some choice of $\vx_0\in\Ss$).
\begin{Prop}\label{Prop:OrientationIndependence}Different choices of
  orientation over the $T_i$ yield isomorphic chain complexes.
\end{Prop}
\proof
Let $\mathfrak{o}$ and $\mathfrak{o}'$ be two choices of
orientation system which agree over $P_j$ for all $j$.  Define
$\sigma_i$ to be $1$ if $\mathfrak{o}$ agrees with
$\mathfrak{o}'$ over $T_i$ and $-1$ otherwise.  (Define
$\sigma_0$ to be $1$.)  It is easy
to check that the map sending $[\vx_i,k]$ to $\sigma_i[\vx_i,k]$
(respectively $\vx_i$ to $\sigma_i\vx_i$) induces isomorphisms on
$CF^\infty$, $CF^+$, and $CF^-$ (respectively an isomorphism on
$\widehat{CF}$).
\endproof

It follows that there are ``only'' $2^{b_2(Y)}$ genuinely different choices of
coherent orientation system.

Although the homologies depend on the choice of orientation system, we
shall usually suppress it from the notation.  Similarly, we shall
often suppress the $\Spin^\bC$--structure $\Ss$ from the notation.

\subsection{Action of $H_1(Y,\bZ)/Tors$}
We now define an action of $H_1(Y,\bZ)/Tors$ on the Floer homologies;
cf~\cite[Section 4.2.5]{OS1}.  We will give the details
only for $HF^\infty$.  The corresponding results for $HF^\pm$ are
immediate consequences, and the results for $\widehat{HF}$ require only slight
modifications of the proofs.

 Recall
that for any $\vx_i$, $\pi_2(\vx_i,\vx_i)$ is identified
with $H_2(\Sigma\times[0,1],\ba\times\{1\}\cup\bb\times\{0\})$.
 From \fullref{homlem}, choosing a basepoint
$\Fz$ gives an isomorphism 
$\Ag\co H_2(\Sigma\times[0,1],\ba\times\{1\}\cup\bb\times\{0\}),\bZ\stackrel{\sim}{\to}
\bZ\oplus H_2(Y,\bZ)$ and hence
$Hom(H_2(\Sigma\times[0,1],\ba\times\{1\}\cup\bb\times\{0\}),\bZ)\cong
\bZ\oplus Hom(H_2(Y),\bZ)\cong \bZ\oplus
H_1(Y)/Tors$.

Choose a
complete set of paths $\{T_i\in\pi_2(\vx_0,\vx_i)\}$.  For each pair
of intersection points $\vx_i$, $\vx_j$ this
gives an isomorphism $\Ah_{\{T_i\}}\co \pi_2(\vx_i,\vx_j)\stackrel{\sim}{\to}
H_2(\Sigma\times[0,1],\ba\times\{1\}\cup\bb\times\{0\})$ via
$\Ah_{\{T_i\}}(A)=\Ag(T_i+A-T_j)$.  

Let $\zeta\in Hom(\bZ\oplus H_2(Y),\bZ)$.  Define
$A_{\zeta,\{T_i\}}\co CF^\infty(Y,\Ss)\to CF^\infty(Y,\Ss)$ by 
$$
A_{\zeta,\{T_i\}}([\vx,i])=
\sum_{\vy\in\Ss}\sum_{\substack{A\in\pi_2(\vx,\vy)\\\ind(A)=1}}
  \zeta(\Ah_{\{T_i\}}(A))\cdot\left(\#\hcM^A\right)[\vy,i-n_z(A)].
$$
Notice that our definition is superficially different from the one used
in~\cite[Section 4.2.5]{OS1}, although their definition makes sense in our language, too.

\begin{Lem}\label{Lemma:H1ChainMap} $A_{\zeta,\{T_i\}}$ is a chain map.
\end{Lem}
\proof
The proof is the same as for~\cite[Lemma 4.18]{OS1}.  Notice that
$\zeta(\Ah_{\{T_i\}}(A+B))=\zeta(\Ah_{\{T_i\}}(A)) +
\zeta(\Ah_{\{T_i\}}(B))$.  Suppose $C\in\pi_2(\vx,\vy)$,
$\ind(C)=2$.  Then $0=\left(\#\bdy\overline{\hcM^C}\right)$, so
\begin{eqnarray}
0&=&\sum_{\substack{ C\in\pi_2(\vx,\vz)\\\ind(C)=2\\n_z(C)=k}}\zeta(\Ah_{\{T_i\}}(C))\left(\#\left(\bdy\overline{\hcM^C}\right)\right)\nonumber\\
 &=&\sum_{\substack{C\in\pi_2(\vx,\vz)\\\ind(C)=2\\n_z(C)=k}}\sum_{\substack{A+B=C\\\ind(A)=\ind(B)=1}}\left(
\zeta(\Ah_{\{T_i\}}(A))\right)
 + \left(\zeta(\Ah_{\{T_i\}}(B))\right)
\left(\#\hcM^A\right)\left(\#\hcM^B\right).\nonumber
\end{eqnarray}
This is the coefficient of $[\vz,i-k]$ in $\left(\bdy\circ
A_{\zeta,\{T_i\}} + A_{\zeta,\{T_i\}}\circ \bdy\right)\left([\vx,i]\right),$
proving the result.
\endproof

The following lemma is analogous to~\cite[Lemma 4.19]{OS1}.
\begin{Lem} If $T'_i$ is another complete set of paths then
  $A_{\zeta,\{T_i\}}$ and $A_{\zeta,\{T'_i\}}$ are chain
  homotopic.
\end{Lem}
\proof
Let $P_i=T_i-T'_i$.  Consider the map
$H\co CF^\infty(Y,\Ss)\to CF^\infty(Y,\Ss)$ given by
$H([\vx_i,j])=\zeta(\Ag(P_i))[\vx_i,j]$.  Then
{\setlength\arraycolsep{2pt}
\begin{eqnarray}
(A_{\zeta,\{T_i\}}&-&A_{\zeta,\{T'_i\}})([\vx_i,j])
\nonumber\\
&=&
\sum_{\vx_k\in\Ss}\sum_{\substack{A\in\pi_2(\vx_i,\vx_k)\\\ind(A)=1}}
  \zeta(\Ah_{\{T_i\}}(A)-\Ah_{\{T'_i\}}(A))
    \cdot\left(\#\hcM^A\right)[\vx_k,j-n_\Fz(A)]\nonumber\\ 
&=&\vrule width 0pt height 12.5pt
\smash{\sum_{\vx_k\in\Ss}\sum_{\substack{A\in\pi_2(\vx_i,\vx_k)\\\ind(A)=1}}
  \zeta(\Ag(A+T_i-T_k)-\Ag(A+T'_i-T'_k))}\nonumber\\
&&\hspace{2.7in}    \cdot\left(\#\hcM^A\right)[\vx_k,j-n_\Fz(A)]\nonumber\\ 
&=&
\sum_{\vx_k\in\Ss}\sum_{\substack{A\in\pi_2(\vx_i,\vx_k)\\\ind(A)=1}}
  \zeta(\Ag(T_i-T'_i)-\Ag(T_k-T'_k))
    \cdot\left(\#\hcM^A\right)[\vx_k,j-n_\Fz(A)]\nonumber\\ 
&=&
\sum_{\vx_k\in\Ss}\sum_{\substack{A\in\pi_2(\vx_i,\vx_k)\\\ind(A)=1}}
  \zeta(\Ag(P_i)-\Ag(P_k))
    \cdot\left(\#\hcM^A\right)[\vx_k,j-n_\Fz(A)]\nonumber\\
&=&
\bdy\circ H +H\circ\bdy\nonumber
\end{eqnarray} }
\endproof

We are now justified in denoting the map on $HF^\infty$ induced by
$A_{\zeta,\{T_i\}}$ (for any complete set of paths
$\{T_i\}$) by simply $A_\zeta$.

The following is~\cite[Proposition 4.17]{OS1}.
\begin{Prop} 
There is a natural action of the group $\bZ\oplus \Hom(H_2(Y),\bZ)$ on
  $HF^\infty(Y,\Ss)$, $HF^+(Y,\Ss)$, $HF^-(Y,\Ss)$ and
  $\widehat{HF}(Y,\Ss)$ lowering degree by one.  This induces an action of
  the exterior algebra $\bigwedge^*(H_1(Y)/Tors)$ on each group.
\label{ActionProp}
\end{Prop}
\proof
The action of $\zeta$ is given by $A_\zeta$.  Obviously
$A_\zeta$ lowers the grading by one and
$A_{\zeta+\zeta'}=A_\zeta+A_{\zeta'}$.  We need to check that
$A_\zeta\circ A_\zeta=0$.

Choose a curve $K$ in
$f^{-1}(3/2-\epsilon,3/2+\epsilon)=\Sigma\times[0,1]$ representing
the homology class $\zeta$ in $H_1(Y)/Tors$.  That is, choose
$K$ so that $A_\zeta(A)$ is the intersection number of $A$ with
$K$ for $A\in\pi_2(\vx,\vx)$.  Perturb $K$ so that
$K\times\bR$ meets transversely every holomorphic curve $u$ with
$\ind(u)=1$, and is transverse to families of holomorphic curves with $\ind=2$.

Let $\cM_1$ denote the moduli space of Riemann surfaces with one
marked point $p$.  Let
$\cM_K^A$ denote the space of holomorphic curves $u\co (S,p)\to
(W,K\times\bR)$ with $(S,p)\in\cM_1$ which, after forgetting the
marked point, represent the homology class $A$.  Let
$\hcM_K^A=\cM_K^A/\bR$.  Then, assuming the appropriate
transversality result, which is left to the reader, we have
$$A_\zeta([\vx,i])=\sum_{\substack{A\in\pi_2(\vx,\vy)\\\ind(A)=1}}
\left(\#\hcM_K(A)\right)[\vy,i-n_\Fz(A)].$$

To prove the proposition, consider the space $\cM_2$ of Riemann surfaces
with $2$ marked points $\{p_1,p_2\}$.  Let $\cM_{K,2}^A$
denote the space of holomorphic maps
$u\co (S,\{p_1,p_2\})\to(W,K\times\bR)$, for
$(S,\{p_1,p_2\})\in\cM_2$.  Let
$\hcM_{K,2}/\bR=\cM_{K,2}^A/\bR$.  

For a generic choice of
$J_s$, and a homology class $A$ with $\ind(A)=2$,
$\hcM_{K,2}(A)$ is a smooth 1--manifold.  The manifold
$\hcM_{K,2}(A)$ has four kinds of ends:
\begin{enumerate}
\item Ends corresponding to $\pi_\bR\circ u(p_1)-\pi_\bR\circ
  u(p_2)\to\infty$.  These correspond to $A_\zeta\circ A_\zeta$.
\item Ends corresponding to
$\pi_\bR\circ u(p_1)-\pi_\bR\circ u(p_2)\to-\infty$.  These also
correspond to $A_\zeta\circ A_\zeta$.
\item Ends where there is a $q$ in $S$ with $\pi_\bR\circ u(q) -
  \pi_\bR\circ u(p_1)\to\infty$ but $\pi_\bR\circ u(p_1) -
  \pi_\bR\circ u(p_2)$ stays bounded.  These correspond to $\bdy
  \circ A_\zeta$.
\item Ends where there is a $q$ in $S$ with $\pi_\bR\circ u(q) -
  \pi_\bR\circ u(p_1)\to-\infty$ but $\pi_\bR\circ u(p_1) -
  \pi_\bR\circ u(p_2)$ stays bounded.  These correspond to $A_\zeta
  \circ \bdy$.
\end{enumerate}
There is also a free action of $\bZ/2\bZ$ on $\hcM_{K,2}^A$
  exchanging the labeling of $p_1$ and $p_2$.  Counting the ends
  of $\hcM_{K,2}^A/(\bZ/2\bZ)$ and summing over $A$ then shows
  that $A_\zeta\circ A_\zeta([\vx,i])$ is chain homotopic to $0$.
 
The result for $HF^\infty$ is immediate.  The constructions of the
lemmas and this proof preserve $CF^-$, so the results for $HF^\pm$
follow.  The proofs for $\widehat{HF}$ are analogous.
\endproof

\section{Isotopy invariance}
\label{Section:Isotopy}
In this section we prove that the homologies defined in \fullref{Section:ChainComplexes} are independent of
deformations of the almost complex structure and isotopies of the
$\alpha$-- and $\beta$--circles not crossing the basepoint
$\Fz$.
Two parts of this story are slightly nonstandard.  One is extending the coherent
orientation system through isotopies which introduce new intersection
points.  The other is that we want to allow seemingly non--Hamiltonian
isotopies
of the $\alpha$ and $\beta$ curves.  The rest of the proof is
analogous to the discussion in~\cite[Section  1.9]{Ya1}.  

We discuss how to extend the orientation system first.  On a first reading the
reader might want to skip the next three paragraphs.

Suppose that we have two pointed Heegaard diagrams
$D_1=(\Sigma,j,\vec{\alpha},\vec{\beta},\Fz)$ and
$D_2=(\Sigma,j',\vec{\alpha}',\vec{\beta}',\Fz)$ which differ by pointed
isotopies (ie, isotopies during which none of the curves
  cross $\Fz$)
and deformations of the complex structure.  The only interesting case
is when the $\Spin^\bC$--structure 
$\Ss$, viewed as a collection of intersection points, is not empty
in either Heegaard diagram.  (If in one of the diagrams $\Ss$
  contains no
intersection points then any choice of orientation system for the
other will be fine.)  By choosing an
appropriate isotopy and then breaking it into a sequences of
isotopies, we can assume that some intersection point $\vx_0\in\Ss$
exists in both Heegaard diagrams.  
For convenience, we fix a parametrization of $\Sigma$ such that
the $\alpha$--circles stay fixed during the isotopy.  Then, with this
parametrization, the $\beta$--circles move and the complex structure
deforms during the isotopy.

We can identify $\pi_2(\vx_0,\vx_0)$ in $D_1$ with
$\pi_2(\vx_0,\vx_0)$ in $D_2$ as follows.  Denote the isotopy by
$I\co \bb\times[0,1]\to\Sigma$.   Let
$\phi\co (C,\bdy C)\to (W,\ba\cup\bb)$ be any singular
2--chain in $\pi_2^{D_1}(\vx_0,\vx_0)$, $C$ a simplicial complex.
Let $\bdy_\beta C=\bdy C\cap I^{-1}(\bb)$.  Define $\psi\co \bdy_\beta
C\times[0,1]\to W$ by $\psi(x,t)=(I(\pi_\Sigma\circ
\phi(x),t),\pi_\bD\circ\phi(x))$.
Then identify $\phi$ with
the chain in $\pi_2^{D_2}(\vx_0,\vx_0)$ given by $\phi+\psi$.

Thus, an orientation system over a
complete set of paths in $D_2$ and an orientation system for $D_1$
determine an orientation system for $D_2$.  We assume that we are
computing the homologies of $D_2$ with respect to an orientation system
determined from the orientation system of $D_1$ in this way; as
observed in the section on chain complexes, different orientation systems
over a complete set of paths lead to isomorphic homologies.  So, which
particular one we choose is unimportant.  Note
also that our choices determine, for $\vx_i^k$ an intersection point
in $D_k$, an isomorphism $\pi_2(\vx_i^1,\vx_j^2)\cong \bZ\oplus
H_2(Y)$.

By a \emph{basic isotopy} we mean an isotopy $(\ba_t,\bb_t)$ with one
of the following two properties.
\begin{enumerate}
\item For all times $t$, $\ba_t$ intersects $\bb_t$
  transversally.  (These are \emph{basic isotopies of the first type})
\item The isotopy introduces one pair of transverse intersections
  between $\ba_t$ and $\bb_t$ by a Hamiltonian deformation (``finger move'')
  of the $\alpha$--circles.  (These are \emph{basic isotopies of the
  second type.})
\end{enumerate}
We only consider isotopies which are sequences of basic isotopies.
Call such an isotopy \emph{strongly admissible} (respectively
\emph{weakly admissible}) if before
and after each basic isotopy the Heegaard diagram is strongly
(respectively weakly) admissible (for $\Ss$).

It is clear that if two Heegaard diagrams are isotopic then they are
isotopic through a sequence of basic isotopies.  In fact, by
\fullref{Prop:MaintainAdmis}, any two isotopic strongly admissible
Heegaard diagrams are isotopic through a sequence of strongly
admissible basic isotopies,
and any two isotopic weakly admissible Heegaard diagrams are isotopic through a
sequence of weakly admissible basic isotopy.

Following~\cite[Section 7.3]{OS1}, we use the fact that basic isotopies of the first type
are equivalent to deformations of the complex structures on
$\Sigma$ and $W$.  
We make this precise.
 Suppose that $D_2$ differs from $D_1$ by a basic isotopy of the
 first type.
Then, there is an orientation--preserving diffeomorphism $\psi$ of
$\Sigma$ taking $D_1$ to 
$D_2$ and mapping $U^{D_1}_{z_i}$ onto $U^{D_2}_{z_i}$.  It
follows that computing the homologies of $D_2$ with
respect to the complex structures $j_\Sigma$ on $\Sigma$ and $J_s$ on
$W$ is the same as computing the homologies of $D_1$ with respect
to $\psi^*j_\Sigma$ on $\Sigma$ and $(\psi\times Id)^*J_s$ on $W$.
Note that if $J_s$ satisfies (\textbf{J1})--(\textbf{J5}) with
respect to $j_\Sigma$ then so does $(\psi\times Id)^*J_s$ with respect to 
$\psi^*j_\Sigma$.
Consequently, independence of the homologies with respect to isotopies
preserving transversality of $\ba\cap\bb$ will follow from
independence with respect to complex structure.

The other point we need to check is that basic isotopies of the second type
do not change
the homologies, either.  Suppose that $D_2$ is obtained from $D_1$
by a basic isotopy of the second type.  Then we can find a
collection of Lagrangian cylinders in $\Sigma\times[0,1]\times\bR$
which agrees
with $\left(\ba\times\{1\}\times\bR\right)\cup\left(\bb\times\{0\}\times\bR\right)$ near
$t=-\infty$ and with
$\left(\ba'\times\{1\}\times\bR\right)\cup\left(\bb'\times\{0\}\times\bR\right)$ near
$t=\infty$.  Call the collection of these Lagrangian cylinders $C$.

We combine invariance under both types of basic isotopy into one:
\begin{Prop}\label{Prop:IsotopyInvariance}  Suppose that either $D_1$ differs from $D_2$ only by
  an isotopy preserving transversality of the $\ba$ and $\bb$
  circles, or that they differ by such an isotopy and a pair creation,
  and that both $D_1$ and $D_2$ are strongly (or, in the case of
  $CF^+$ and $\widehat{CF}$, weakly) admissible Heegaard diagrams.
  Suppose $J_1$ (respectively $J_2$) satisfies
  {\rm(\textbf{J1})--(\textbf{J5})} and achieves transversality
  for $D_1$ (respectively $D_2$), with respect to $j_1$
  (respectively $j_2$) on $\Sigma$.
  Then the chain complexes  $CF^\infty_{D_1}$ and
  $CF^\infty_{D_2}$ (respectively $CF^-_{D_1}$ and $CF^-_{D_2}$;
  $CF^+_{D_1}$ and $CF^+_{D_2}$; and $\widehat{CF}_{D_1}$ and
  $\widehat{CF}_{D_2}$) are chain homotopy equivalent.
  Further, the isomorphisms induced on homologies respect the
  $H_1(Y)/Tors$--module and, where appropriate, the
  $\bZ[U,U^{-1}]$-- or $\bZ[U]$--module
  structures.\label{Isotopy:Independence}
\end{Prop}

By the discussion preceding the proposition, proving this proposition proves the
independence of the Floer homologies under isotopies.  Note that it
also implies that the restrictions we imposed in
\fullref{Section:ChainComplexes} on $J|_{s=0}$ and $J|_{s=1}$
for $CF^\infty$ and $CF^\pm$ are not needed for $\bdy^2$ to be $0$.

Fix $T>0$.  Choose an almost complex structure $J$ on
$\Sigma\times[0,1]\times\bR$ which satisfies (\textbf{J1}), (\textbf{J2})
and (\textbf{J4}), agrees
with $J_1$ on $(-\infty,-T]$ and with $J_2$
  on $[T,\infty)$, and achieves transversality for
holomorphic curves of the form $u\co (S,\bdy S)\to
(\Sigma\times[0,1]\times\bR,C)$.
This is possible by essentially the same argument as in
\fullref{Section:Transversality}. 

We define a chain map from each chain complex defined on $D_1$ to
the corresponding chain complex on $D_2$ by counting
$J$--holomorphic curves in
$W$.  We carry out the details for $CF^\infty$; the results for
$CF^+$ and $CF^-$ will follow, and the proof for $\widehat{CF}$ is
similar.  Let $(CF^\infty,\bdy_i)$ be the chain complex defined on
$D_i$, $i=1,2$. 

For $\vx^1$ (respectively $\vx^2$) an intersection point in
$D_1$ (respectively $D_2$), let
$\cM_0(\vx^1,\vx^2)$ (respectively $\cM_1(\vx^1,\vx^2)$) denote the space of all holomorphic curves
$u$ in $(W,J)$ connecting $\vx^1$ to $\vx^2$ in homology
classes $A$ with $\ind(A)=0$ (respectively $\ind(A)=1$), and satisfying
(\textbf{M0})--(\textbf{M6}) (with respect to the new Lagrangian cylinders $C$).

Consider the map
$\Phi\co CF_1^\infty\to CF_2^\infty$ defined by
$$
\Phi([\vx^1,i])=\sum_{\vx^2\in\Ss}\sum_{u\in\cM_0(\vx^1,\vx^2)}[\vx^2,i-n_z(u)].
$$
We need to check that the coefficient of $[\vx^2,j]$ is a finite sum
for each $\vx^2$ and
$j$:
\begin{Lem}Given $j\in\bZ$ there are at most finitely many homology
  classes $A$
  connecting $\vx^1$ to $\vx^2$ with $\ind(A)=j$ which admit
  a holomorphic curve.
\end{Lem}
\proof 
This follows from the strong admissibility criterion; our proof is
essentially the same as~\cite[Lemma 7.4]{OS1}.  

Choose a point $\Fz_i$ in each component of
$\Sigma\setminus(\ba^1\cup\bb^1)$ in such a way that none of the
$\alpha$-- or $\beta$--circles cross any of the $\Fz_i$ during the
isotopy.  If $A$ supports a $J$--holomorphic curve then
$n_{\Fz_i}(A)\geq 0$ for all $i$.  Choose any homology class in
$B\in\pi_2^{D_2}(\vx^2,\vx^1)$.  We can view $A+B$ as an element
of $\pi_2^{D_1}(\vx^1,\vx^1)$, and $n_{\Fz_i}(A+B)\geq
n_{\Fz_i}(B)$ for all $i$.  Now the argument used in
\fullref{Lemma:StrongAdmis} gives bounds for the coefficients of
$A+B$, and hence bounds for the coefficients of $A$.  This proves
the result.
\endproof

We shall sometimes denote $\Phi$ by $\Phi_{12}$ to
emphasize that $\Phi$ is induced from a bordism from $D_1$ to $D_2$.

\begin{Lem}$\Phi$ is a chain map.\label{Isotopy:MapsDefined}
\end{Lem}
\proof
Let $\bdy_i$ denote the boundary map on the chain complex for $D_i$.

We consider the compactified moduli space
$\overline{\hcM_1(\vx^1,\vy^2,k)}$ of index $1$ holomorphic curves
in homology classes $A$ with $n_\Fz(A)=k$.  There is still no bubbling, so
this is a compact one--manifold whose boundary consists
of height two holomorphic buildings one story of
which lies in $(W,J)$ and the other of which lies in either $(W,J_1)$ or
$(W,J_2)$.  Hence,
$$0=\# \left(\bdy \overline{\hcM_1(\vx^1,\vy^2,k)}\right).$$
But this is the coefficient of $[\vy^2,i-k]$ in
$\Phi\circ\bdy_1([\vx^1,i]) + \bdy_2\circ\Phi([\vx^1,i])$.
\endproof

\begin{Lem}Given two different choices of complex structure $J$
  and $J'$ connecting $J_1$ to $J_2$ and two different choices
  $C$ and $C'$ for the cylinders connecting
 $\ba\times\{1\}\cup\bb\times\{0\}$ to
 $\ba'\times\{1\}\cup\bb'\times\{0\}$ the
  maps $\Phi$ and $\Phi'$ are chain
  homotopic.\label{Isotopy:MapsIndependent}
\end{Lem}
\proof
We outline the proof; further details are left to the reader.
Choose a generic path $J_t$ from $J$ to $J'$ and a Lagrangian
isotopy $C_t$ from $C$ to $C'$.  Let $\cM_{-1,t}(\vx^1,\vx^2)$
denote the moduli space of holomorphic curves $u$ connecting
$\vx^1$ to $\vx^2$ in homology classes of index $-1$.
For a finite collection of $t_i$'s, $0<t_1<\ldots<t_k<1$,
$J_{t_i}$ is degenerate in such a way that $\cM_{-1,t}(\vx^1,\vx^2)$
is nonempty.  (As usual, the finiteness uses the admissibility
hypothesis.)   Then we define a chain
homotopy $\Delta$ by
$$
\Delta([\vx^1,i])=\sum_t\sum_{\vx^2}\sum_{u\in\cM_{-1,t}}[\vx^2,i-n_\Fz(u)].
$$
Similarly, define $\cM_{0,t}(\vx^1,\vy^2)$ to be the
moduli space of holomorphic curves $u$ connecting
$\vx^1$ to $\vy^2$ in homology classes of index $0$.
Then $\cup_{t}\cM_{0,t}(\vx^1,\vy^2)$ is a 1--manifold with boundary
with four types of ends:
\begin{enumerate}
\item Ends corresponding to height two holomorphic buildings, with a
  $J_{t_k}$--holomor\-phic curve of index
  $-1$ and a $J_1$--holomorphic curve of index $1$.  These
  correspond to $\Delta\circ\bdy_1$.
\item Ends corresponding to height two holomorphic buildings, with a
  $J_{t_k}$--holomor\-phic curve of index
  $-1$ and a $J_2$--holomorphic curve of index $1$.  These
  correspond to $\bdy_2\circ\Delta$.
\item Ends corresponding to $t=0$ and a $J$--holomorphic curve.
  These correspond to $-\Phi$.
\item Ends corresponding to $t=1$ and a $J'$--holomorphic curve.
  These correspond to $\Phi'$.
\end{enumerate}
So, counting the ends gives the result.
\endproof

\begin{Lem}If we have a third diagram $D_3$ and a bordism 
  from $D_2$ to $D_3$ then on homology
  $\Phi_{23}\circ\Phi_{12}=\Phi_{13}$.\label{Isotopy:MapsCompose}
\end{Lem}
Note that since different choices of bordism give chain homotopic
$\Phi_{ij}$ the exact choices of $\Phi_{ij}$ are unimportant here.

\proof
This is immediate from compactness and \fullref{Gluing:Cylindrical}.
\endproof

\begin{Lem}The maps on the homologies induced by $\Phi$ preserve the
  $H_1(Y)/Tors$--module structure.\label{Isotopy:PreserveStructure}
\end{Lem}
\proof
Fix $\zeta\in H_1(Y)/Tors$.  Let $K$ be a knot in
$\Sigma\times[0,1]$ representing $\zeta$, as in the proof of
\fullref{ActionProp}.

Let $\cM_1$ denote the moduli space of Riemann surfaces with one
marked point $p$.  For $i=1,2$,
 let
$\cM_{K,i}^A$ denote the space of $J_i$--holomorphic curves $u\co (S,p)\to
(W,K\times\bR)$ with $(S,p)\in\cM_1$ which, after forgetting the
marked point, represent the homology class $A$.  Let
$\hcM_{K,i}^A=\cM_{K,i}^A/\bR$.  Assume $K$ (or $J_1$ and
$J_2$) is chosen so that we have transversality.  Then, on the chain
level the action of $H_1(Y)/Tors$ on $HF^\infty_{D_i}$ is given by
$$A_{\zeta,i}([\vx^i,j])=\sum_{\substack{A\in\pi_2(\vx^i,\vy^i)\\\ind(A)=1}}
\left(\#\hcM_{K,i}^A\right)[\vy^i,j-n_z(A)].$$

Let
$\cM_{K}^A$ denote the space of $J$--holomorphic curves $u\co (S,p)\to
(W,K\times\bR)$ with $(S,p)\in\cM_1$ which, after forgetting the
marked point, represent the homology class $A$.
Assume $K$, $J_1$, $J_2$ and $J$ are
chosen so that we have transversality.
Consider the ends of 
$$\cM_{K}(\vx^1,\vy^2,k)=\bigcup_{\substack{A\in\pi_2(\vx^1,\vy^2)\\\ind(A)=1\\n_\Fz(A)=k}}\cM_K(A).$$
There are two kinds of ends:
\begin{enumerate}
\item Ends where $\pi_\bD\circ u(p)\to\infty$.  These ends correspond to the
  coefficient of $[\vy,i-k]$ in
  $A_{\zeta,2}\circ \Phi([\vx,i])$
\item Ends where $\pi_\bR\circ u(p)\to-\infty$.  These ends correspond to the
  coefficient of $[\vy,i-k]$ in
  $\Phi\circ A_{\zeta,1}([\vx,i])$.
  \item Ends where $\pi_\bR\circ u(p)$ stays bounded and an index $1$ curve splits--off at $+\infty$.  These correspond to the coefficient of $[\vy,i-k]$ in $\bdy\circ \Phi([\vx^1,i-k])$.
  \item Ends where $\pi_\bR\circ u(p)$ stays bounded and an index $1$ curve splits--off at $-\infty$.  These correspond to the coefficient of $[\vy,i-k]$ in $\Phi\circ \bdy([\vx^1,i-k])$.
\end{enumerate}
Counting the ends, and using the fact that $\Phi\circ \bdy + \bdy\circ\Phi=0$, gives the result.
\endproof

\textbf{Proof of \fullref{Prop:IsotopyInvariance}.}

The proposition follows immediately from these four lemmas.  From
\fullref{Isotopy:MapsDefined} we have
chain maps $\Phi_{12}$ and $\Phi_{21}$ whose induced maps preserve the
$H_1(Y)/Tors$--structure by \fullref{Isotopy:PreserveStructure}.
From \fullref{Isotopy:MapsIndependent} and \fullref{Isotopy:MapsCompose},
 $\Phi_{12}\circ\Phi_{21}$ is chain homotopic to the identity map, as
is $\Phi_{21}\circ\Phi_{12}$.  All the maps are obviously maps of
$\bZ[U,U^{-1}]$--modules.  This
proves \fullref{Isotopy:Independence} for $HF^\infty$.  To
conclude the
result for $HF^+$ and $HF^-$ it is only necessary to make the
trivial observation that all of the maps used preserve $HF^-$.

The
proof for $\widehat{HF}$ is completely analogous -- one simply restricts
in each case to holomorphic curves with $n_\Fz=0$.
\endproof

\section{Triangles}
\label{Section:Triangles}
By a pointed Heegaard triple--diagram we mean a Riemann surface $\Sigma$ together with three $g$--tuples of pairwise disjoint, homologically linearly independent simple closed curves
 $\va=\{\alpha_1,\cdots,\alpha_g\}$, $\vb=\{\beta_1,\cdots,\beta_g\}$ and $\vc=\{\gamma_1,\cdots,\gamma_g\}$, together with a basepoint
  $\Fz\in\Sigma\setminus(\ba\cup\bb\cup\bc)$.  A Heegaard triple--diagram specifies three 3--manifolds $\Yab$, $\Ybc$ and $\Yac$.  It also specifies a $4$--manifold $X_{\alpha,\beta,\gamma}$ with
  boundary $\Yab\cup\Ybc\cup-\Yac$ as follows.  The curves $\va$ 
(respectively $\vb$, $\vc$) specify a
handlebody $U_\alpha$ (respectively $U_\beta$, $U_\gamma$) with
boundary $\Sigma$.  Let $\ol{T}$ be a triangle with edges $e_1$,
$e_2$, and $e_3$.  Define $X_{\alpha,\beta,\gamma}$ to be the
manifold obtained by gluing 
$U_\alpha\times[0,1]$, $U_\beta\times[0,1]$, and
$U_\gamma\times[0,1]$ to $\Sigma\times \ol{T}$ along $\Sigma\times e_1$,
$\Sigma\times e_2$, and $\Sigma\times e_3$ by identifying
$\ba\in \bdy U_\alpha\times\{p\}$ with $\ba\in \Sigma\times\{p\}$
($p\in e_1=[0,1]$) and similarly for $\beta$ and $\gamma$.

In this (rather long) section we associate to a Heegaard triple--diagram maps between the Floer homologies of $\Yab$, $\Ybc$ and $\Yac$.  Specifically

\begin{Construct}\label{Construct:TriangleMaps}
To an admissible (see \fullref{Subsection:TriangleAdmissibility}) pointed Heegaard triple--diagram $(\Sigma,\va,\vb,\vc,\Fz)$ and a $\Spin^\bC$--structure $\St$ on $X_{\alpha,\beta,\gamma}$, as well as a coherent orientation system, as described in \fullref{Subsection:TriangleOrientations}, we associate $U$--equivariant homomorphisms
\begin{eqnarray}
\hat{F}_{\alpha,\beta,\gamma}&\co &\widehat{HF}(\Yab;\St|_{\Yab})\otimes_{\bZ}\widehat{HF}(\Ybc;\St|_{\Ybc}) \to \widehat{HF}(\Yac;\St|_{\Yac})\nonumber\\
F^\infty_{\alpha,\beta,\gamma}&\co &HF^\infty(\Yab;\St|_{\Yab})\otimes_{\bZ[U]}HF^\infty(Ybc;\St|_{\Ybc})\to HF^\infty(\Ybc;\St|_{\Yac})\nonumber\\
F^+_{\alpha,\beta,\gamma}&\co &HF^+(\Yab;\St|_{\Yab})\otimes_{\bZ[U]}HF^{\leq0}(\Ybc;\St|_{\Ybc})\to HF^+(\Yac;\St|_{\Yac})\nonumber\\
F^-_{\alpha,\beta,\gamma}&\co &HF^{\leq0}(\Yab;\St|_{\Yab})\otimes_{\bZ[U]}HF^{\leq0}(\Ybc;\St|_{\Ybc})\to HF^{\leq0}(\Yac;\St|_{\Yac}).\nonumber
\end{eqnarray}
These maps satisfy an associativity property stated in
\fullref{AssocProp}.
\end{Construct}

Of course, the zero map satisfies these conditions, but our maps are
usually more interesting.  In particular, an instance of the
construction will be used to prove handleslide invariance in the next
section.  (They are also a key computational tool, have been used to
define $4$--manifold invariants, contact invariants, and so on.)
They will be produced by counting holomorphic curves in the product of
$\Sigma$ and a disk with three boundary punctures, ie, a
triangle.

The outline of this section is as follows.  In
\fullref{Subsection:TriangleBasics} we discuss basic
topological properties of maps to $\Sigma\times T$, $T$ a
``triangle''.  In particular, we discuss
when such maps exist, how many homology classes of them there are, and
how such homology classes specify $\Spin^\bC$--structures.  In
\fullref{Subsection:TriangleModuli1} we discuss some basic
prerequisites for construction of the triangle maps (the complex
structures we consider, the index) and then define the triangle maps,
conditional on certain technicalities to be addressed in the following
two sections.  In \fullref{Subsection:TriangleOrientations} we
address the first of these technicalities:  coherent orientations of the moduli
spaces of triangles.  In
\fullref{Subsection:TriangleAdmissibility} we deal with the
second technicality:  admissibility criteria for Heegaard
triple--diagrams necessary for the triangle maps to be defined.  In
\fullref{Subsection:TriangleModuli2} we return to the
definition of the triangle maps, proving that the maps are chain maps,
and are independent of the complex structure on $\Sigma\times T$ and
isotopies of the $\alpha$--, $\beta$-- and $\gamma$--circles.
Finally, in \fullref{Subsection:TriangleAssoc} we prove an
associativity property of triangle maps.

\subsection{Topological preliminaries on triangles}\label{Subsection:TriangleBasics}
Fix a pointed Heegaard triple--diagram $\cH^3=(\Sigma,\va,\vb,\vc,\Fz)$.  Let $\cH_{\alpha,\beta}$, $\cH_{\beta,\gamma}$ and $\cH_{\alpha,\gamma}$ denote the pointed Heegaard diagrams for $Y_{\alpha,\beta}$, $Y_{\beta,\gamma}$ and $Y_{\alpha,\gamma}$ specified by $\cH^3$.

\subsubsection{Homological preliminaries on triangles}
\begin{figure}
\centering
\begin{picture}(0,0)%
\includegraphics[scale=0.7]{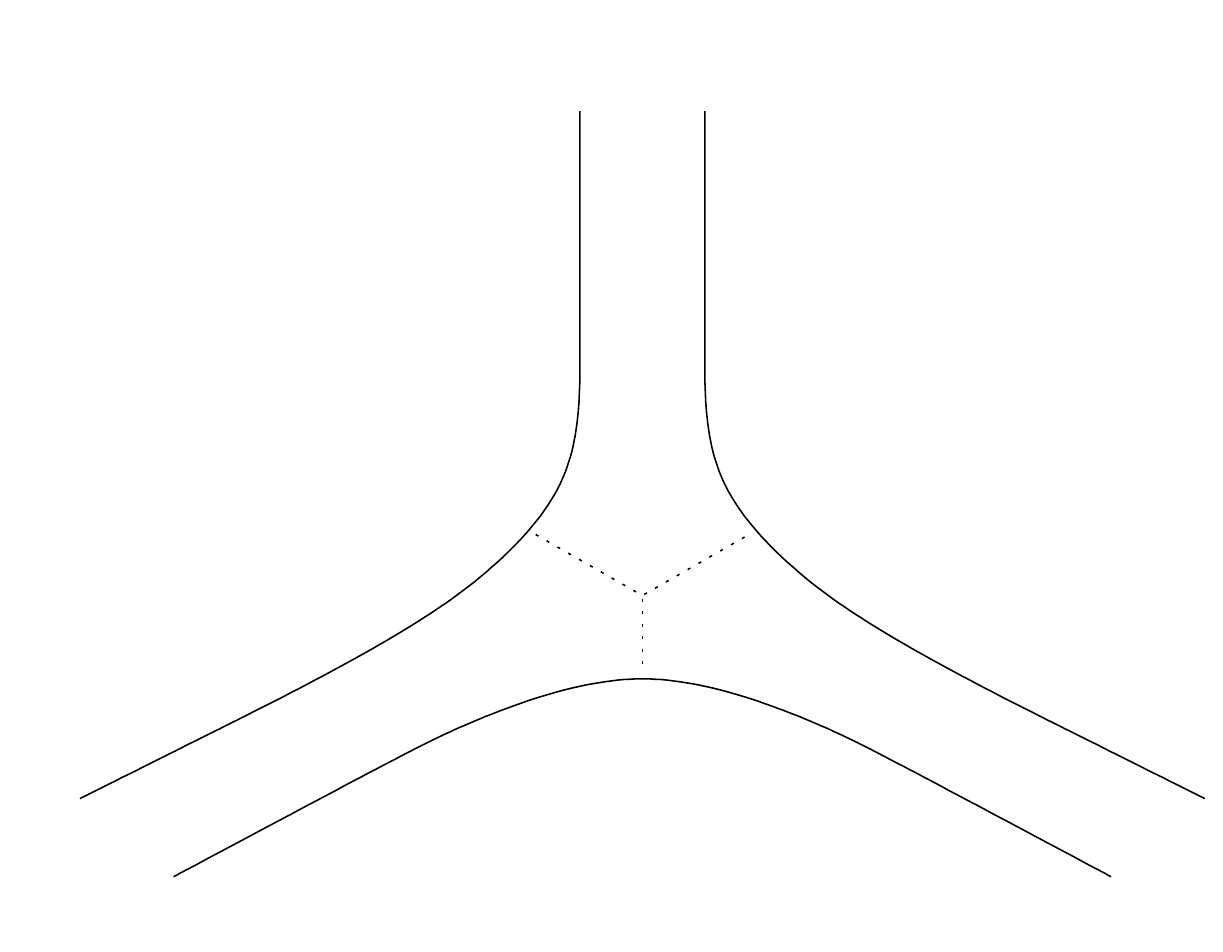}%
\end{picture}%
\setlength{\unitlength}{2763sp}%
\begingroup\makeatletter\ifx\SetFigFont\undefined%
\gdef\SetFigFont#1#2#3#4#5{%
  \reset@font\fontsize{#1}{#2pt}%
  \fontfamily{#3}\fontseries{#4}\fontshape{#5}%
  \selectfont}%
\fi\endgroup%
\begin{picture}(5850,4548)(226,-4390)
\put(2401,-1486){\makebox(0,0)[lb]{\smash{\SetFigFont{12}{14.4}{\rmdefault}{\mddefault}{\updefault}{\color[rgb]{0,0,0}\(e_2\)}%
}}}
\put(3226,-3286){\makebox(0,0)[lb]{\smash{\SetFigFont{12}{14.4}{\rmdefault}{\mddefault}{\updefault}{\color[rgb]{0,0,0}\(v_1\)}%
}}}
\put(2401,-2311){\makebox(0,0)[lb]{\smash{\SetFigFont{12}{14.4}{\rmdefault}{\mddefault}{\updefault}{\color[rgb]{0,0,0}\(v_2\)}%
}}}
\put(3901,-2311){\makebox(0,0)[lb]{\smash{\SetFigFont{12}{14.4}{\rmdefault}{\mddefault}{\updefault}{\color[rgb]{0,0,0}\(v_3\)}%
}}}
\put(3451,-2761){\makebox(0,0)[lb]{\smash{\SetFigFont{12}{14.4}{\rmdefault}{\mddefault}{\updefault}{\color[rgb]{0,0,0}\(G\)}%
}}}
\put(601,-3961){\makebox(0,0)[lb]{\smash{\SetFigFont{12}{14.4}{\rmdefault}{\mddefault}{\updefault}{\color[rgb]{0,0,0}\(v_{12}\)}%
}}}
\put(6076,-4336){\makebox(0,0)[lb]{\smash{\SetFigFont{12}{14.4}{\rmdefault}{\mddefault}{\updefault}{\color[rgb]{0,0,0}\(\alpha,\gamma\)}%
}}}
\put(226,-4336){\makebox(0,0)[lb]{\smash{\SetFigFont{12}{14.4}{\rmdefault}{\mddefault}{\updefault}{\color[rgb]{0,0,0}\(\alpha,\beta\)}%
}}}
\put(5851,-4036){\makebox(0,0)[lb]{\smash{\SetFigFont{12}{14.4}{\rmdefault}{\mddefault}{\updefault}{\color[rgb]{0,0,0}\(v_{13}\)}%
}}}
\put(2776, 14){\makebox(0,0)[lb]{\smash{\SetFigFont{12}{14.4}{\rmdefault}{\mddefault}{\updefault}{\color[rgb]{0,0,0}\(\beta,\gamma\)}%
}}}
\put(3151,-286){\makebox(0,0)[lb]{\smash{\SetFigFont{12}{14.4}{\rmdefault}{\mddefault}{\updefault}{\color[rgb]{0,0,0}\(v_{23}\)}%
}}}
\put(2176,-3661){\makebox(0,0)[lb]{\smash{\SetFigFont{12}{14.4}{\rmdefault}{\mddefault}{\updefault}{\color[rgb]{0,0,0}\(e_1\)}%
}}}
\put(3676,-1486){\makebox(0,0)[lb]{\smash{\SetFigFont{12}{14.4}{\rmdefault}{\mddefault}{\updefault}{\color[rgb]{0,0,0}\(e_3\)}%
}}}
\end{picture}
\caption{The triangle $T$}
\label{Triangle1}
\end{figure}
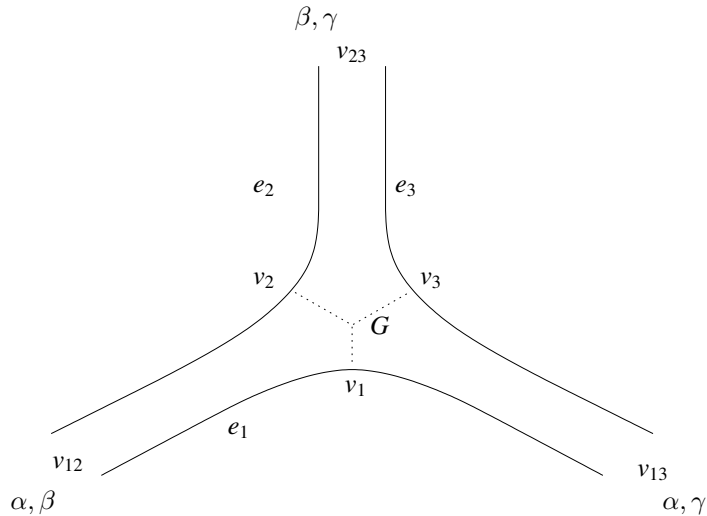

From now on, by the triangle $T$ we mean a Y--shaped
region in $\bC$ with three cylindrical ends, as shown in
\fullref{Triangle1}.  Note that $T$ is conformally equivalent to
a (in fact, any) triangle with punctures at the corners.  (We
  will occasionally use a closed triangle, which we
  continue to denote
  $\ol{T}$ and think of as the closure of $T$.)
Let $e_1$, $e_2$,
and $e_3$ denote the three boundary components of $T$, ordered
clockwise, and $v_{12}$ ,$v_{23}$ and $v_{13}$ the ends between $e_1$ and $e_2$, $e_2$ and $e_3$, and $e_1$ and $e_3$ respectively.  Let $W_{\alpha,\beta,\gamma}=\Sigma\times T$.
For I--chord collections $\vx$ of $\cH^1$, $\vy$ of $\cH^2$ and
$\vz$ of $\cH^3$, let $\pi_2(\vx,\vy,\vz)$ denote the collection
of homology classes of maps $(S,\bdy S)\to (W_{\alpha,\beta,\gamma},
\ba\times e_1\cup\bb\times e_2\cup\bc\times e_3)$ ($S$ a Riemann
surface with boundary and punctures on the  boundary) which are
asymptotic to $\vx$ at $v_{12}$, $\vy$ at $v_{23}$ and $\vz$
at $v_{13}$.  As before, there is a map
$n_\Fz\co \pi_2(\vx,\vy,\vz)\to\bZ$.  Let
$\hat{\pi}_2(\vx,\vy,\vz)=\{A\in\pi_2(\vx,\vy,\vz)|n_\Fz(A)=0\}$. 

Let $W_{\alpha,\beta}$, $W_{\beta,\gamma}$ and $W_{\alpha,\gamma}$ denote the three cylindrical manifolds which are the ends of $W_{\alpha,\beta,\gamma}$.

Note that there are concatenation maps
$*\co \pi_2(\vx,\vx')\times\pi_2(\vx,\vy,\vz)\to\pi_2(\vx',\vy,\vz)$
(and similarly for $\vy$ and $\vz$).  

Suppose that $\pi_2(\vx,\vy,\vz)$ is nonempty.  Fix an element
$A\in\pi_2(\vx,\vy,\vz)$.  Let $G$ denote the tree with one
vertex of valence three and three vertices of valence 1 (so $G$
looks like a figure Y).  Let $v_1$, $v_2$, and $v_3$ denote the
three valence--one vertices in $G$.  Then $A$ determines a map
$$\pi_2(\vx,\vy,\vz)\to H_2(\Sigma\times
G,\ba\times\{v_1\}\cup\bb\times\{v_2\}\cup\bc\times\{v_3\})$$
by subtracting from each element of $\pi_2(\vx,\vy,\vz)$ a
representative for $A$ and pushing forward via the retract suggested by \fullref{Triangle1}.  It is easy to see that this map is bijective.

Recall that $\pi_2(\vx,\vx)$ is canonically identified with
$H_2(\Sigma\times[0,1],\ba\times\{1\}\cup\bb\times\{0\})$.  Viewing
$[0,1]$ as the path from $v_2$ through the trivalent vertex of $G$ to
$v_1$, we obtain an inclusion
$H_2(\Sigma\times[0,1],\ba\times\{1\}\cup\bb\times\{0\})\to
H_2(\Sigma\times G,\ba\times
\{v_1\}\cup\bb\times\{v_2\}\cup\bc\times\{v_3\})$.  Under these
identifications, concatenation
$\pi_2(\vx,\vx)\times\pi_2(\vx,\vy,\vz)\to\pi_2(\vx,\vy,\vz)$
corresponds to addition in $H_2(\Sigma\times G,\ba\times
\{v_1\}\cup\bb\times\{v_2\}\cup\bc\times\{v_3\})$.  Similar remarks
apply to $\vy$ and $\vz$, of course.

We have:

\begin{Lem}{\rm(Compare~\cite[Proposition 8.2 and Proposition 8.3]{OS1})}\label{Lemma:TriPi}\qua
There is a natural short exact sequence
$$
0\to\bZ\to 
H_2(\Sigma\times G,\ba\times\{v_1\}\cup\bb\times\{v_2\}\cup\bc\times\{v_3\})
\to H_2(X_{\alpha,\beta,\gamma})\to 0.
$$
The basepoint $\Fz\in\Sigma\setminus(\ba\cup\bb\cup\bc)$
determines a splitting 
$$n_\Fz\co H_2(\Sigma\times
G,\ba\times\{v_1\}\cup\bb\times\{v_2\}\cup\bc\times\{v_3\})\to\bZ$$
of this sequence.
\end{Lem}
\proof
From the long exact sequence for the pair $(X_{\alpha,\beta,\gamma},\Sigma\times \ol{T})$ we have:
$$
H_2(\Sigma\times T)\to H_2(X_{\alpha,\beta,\gamma})\to
H_2(X_{\alpha,\beta,\gamma},\Sigma\times \ol{T})\to
H_1(\Sigma\times \ol{T}).
$$
The first map is trivial since $\Sigma$ is null--homologous in $X_{\alpha,\beta,\gamma}$
(it bounds in $U_\alpha$, for example).  Since the boundary map in
the long exact sequence for the pair $(U_\alpha,\Sigma)$ takes
$H_2(U_\alpha,\Sigma)$ one--to--one onto $H_1(\ba)\subset H_1(\Sigma)$,
and similarly for $\beta$ and $\gamma$, the kernel of the last map
is the same as 
$\ker\left( H_1(\ba)\oplus H_1(\bb)\oplus H_1(\bc)\to
H_1(\Sigma)\right).$
  (Here, the map is induced from including $\ba$,
$\bb$, and $\bc$ in $\Sigma$.) Thus,
$H_2(X_{\alpha,\beta,\gamma})$ is isomorphic to this kernel.

From the long exact sequence for the pair $(\Sigma\times
G,\ba\times\{v_1\}\cup\bb\times\{v_2\}\cup\bc\times\{v_3\})$ we have
{\setlength\arraycolsep{0pt} \begin{eqnarray*}
0\to H_2(\Sigma\times G)\to 
H_2(\Sigma\times 
G,\ba\times\{v_1\}&\cup&\bb\times\{v_2\}\cup\bc\times\{v_3\})\\
&\to& H_1(\ba)\oplus H_1(\bb)\oplus H_1(\bc) \to H_1(\Sigma\times G)
\end{eqnarray*}
The kernel of the last map is
$\ker\left( H_1(\ba)\oplus H_1(\bb)\oplus H_1(\bc)\to
H_1(\Sigma)\right)\cong H_2(X_{\alpha,\beta,\gamma}).$ }
$H_2(\Sigma\times G)\cong \bZ$.
The statement about the splitting is clear.
\endproof

The previous lemma tells us what $\pi_2(\vx,\vy,\vz)$ is if
nonempty.  It is worth knowing when $\pi_2(\vx,\vy,\vz)$ is in fact
empty.  Define $\epsilon(\vx,\vy,\vz)$ as follows.  Choose a chain
$p_\alpha$ (respectively $p_\beta$, $p_\gamma$) in $\ba$
(respectively $\bb$, $\bc$) with $\bdy
p_\alpha=\vx-\vy$ (respectively $\bdy p_\beta=\vy-\vz$, $\bdy
p_\gamma=\vz-\vx$).  Then $\epsilon(\vx,\vy,\vz)$ is the image of
$p_\alpha+p_\beta+p_\gamma$ under the map 
$$
H_1(\Sigma)\to\frac{H_1(\Sigma)}{H_1(\ba)+H_1(\bb)+H_1(\bc)}\cong H_1(X_{\alpha,\beta,\gamma}).
$$
\begin{Lem}{\rm(Compare~\cite[Proposition 8.3]{OS1})}\label{Lemma:TriEpsilon}\qua
The set $\pi_2(\vx,\vy,\vz)$ is nonempty if and only if
  $\epsilon(\vx,\vy,\vz)=0$.  
\end{Lem}
\proof
If $\pi_2(\vx,\vy,\vz)$ is nonempty, choose an element
$A\in\pi_2(\vx,\vy,\vz)$.  View $A$ as a chain in $\Sigma$.
Then the boundary of $A$ is a chain which defines
$\epsilon(\vx,\vy,\vz)$, and is obviously zero in homology.

Conversely, if $\epsilon(\vx,\vy,\vz)$ is zero then we can choose $p_\alpha$, $p_\beta$ and $p_\gamma$ to be cellular chains (with
respect to the cellulation of $\Sigma$ induced by $\ba$, $\bb$
and $\bc$) so that $p_\alpha+p_\beta+p_\gamma$ is null--homologous.
Any chain with boundary $p_\alpha+p_\beta+p_\gamma$ is an element of
$\pi_2(\vx,\vy,\vz)$.
\endproof

\subsubsection{$\Spin^\bC$--structures and triangles}\label{Subsubsection:TriangleSpinc}
While each intersection point in a Heegaard diagram specifies a
$\Spin^\bC$--structure on the underlying $3$--manifold, for a Heegaard triple--diagram it is the
elements of $\pi_2(\vx,\vy,\vz)$ that specify
$\Spin^\bC$--structures on the corresponding $4$--manifold.  This correspondence is somewhat more
complicated than in the 3--manifold case.  Our exposition will be
very close to that in~\cite[Section 8.1.4]{OS1}, but the reader may find some points clearer in one treatment or the other.

Recall (\fullref{section:homotopy}) that the definition of a $\Spin^\bC$--structure on a
$3$--manifold which we have used is a ``homology class of
nonvanishing vector fields.''  On a $4$--manifold, the analogous
definition is:
\begin{Def}Fix a connected $4$--manifold $M$.  Suppose that
  $J_1$ and $J_2$ are almost complex structures on $M$, defined
  in the complement of some $4$--ball in $M$.  (We will say that
  $J_1$ and $J_2$ are \emph{almost defined}.)  We say that $J_1$
  and $J_2$ are \emph{homologous} if $J_1$ and $J_2$ are
  isotopic in the complement of some (larger) $4$--ball in $M$.
  We define a $\Spin^\bC$--structure on $M$ to be a homology class
  of almost defined almost complex structures.
\end{Def}

We sketch the identification with the standard definition of
$\Spin^\bC$--structures.  Suppose we are given an almost complex $J$
structure defined in the complement of some $4$--ball $B$.  On the
complement of $B$, the almost complex structure $J$ determines
canonically a
$\Spin^\bC$ lifting of the bundle of frames.  The obstruction to
extending the $\Spin^\bC$--structure over $B$ lies in $H^3(B,\bdy
B)=0$ and the collection of distinct extensions correspond to
$H^2(B,\bdy B)=0$.

Conversely, a $\Spin^\bC$ lifting of the bundle of frames determines
complex positive and negative spinor bundles.  Choosing a section
$s$ of the positive spinor bundle vanishing at a finite number of
points, Clifford multiplication by $s$ gives an isomorphism of the
negative spinor bundle with $TM$ away from a finite number of
points.  Choose a ball $B$ containing these points.  Then the
identification of $TM$ with the negative spinor bundle determines an
almost complex structure on $M\setminus B$.

Now, fix a homology class $A\in\pi_2(\vx,\vy,\vz)$.  Between here and \fullref{Lemma:TriangleSpinDependence} we associate a $\Spin^\bC$--structure on
$X_{\alpha,\beta,\gamma}$ to $A$ and some extra data.  (It will turn out that the $\Spin^\bC$--structure does not depend on the extra data.)

Fix a height function
$f_\alpha$ (respectively $f_\beta$, $f_\gamma$) on the
handlebody $U_\alpha$ (respectively $U_\beta$, $U_\gamma$) with
one index $0$ critical point and $g$ index $1$ critical points,
such that $f_{\alpha}|_{\bdy U_{\alpha}}$ (respectively
$f_{\beta}|_{\bdy U_{\beta}}$, $f_{\gamma}|_{\bdy U_{\gamma}}$) is
constant.

Choose a smooth immersion
$\phi\co S\to W_{\alpha,\beta,\gamma}$ representing $A$.  We will
place some requirements on $\phi$ presently.

Let $F$ denote the surface obtained by capping off $\phi(S)\cup
\{\Fz\}\times T$ with the downward gradient flows of
$f_{\alpha}$, $f_\beta$ and $f_\gamma$.  So, $F$ is an
immersed surface with boundary on the $3g+3$ critical lines in $X_{\alpha,\beta,\gamma}$.  (By a critical line we mean a line of the form $(\textrm{critical point of $f_i$})\times e_i$.)

Let $\cL$ be the $2$--plane field on $X_{\alpha,\beta,\gamma}\setminus
F$ given by 
\begin{itemize}
\item $\cL(p)=T(\{p_1\}\times T)\subset T_{(p_1,p_2)}\Sigma\times T$
for $p=(p_1,p_2)\in\Sigma\times T$
\item $\cL(p)=\ker df_\alpha(p)$ (respectively $\ker df_\beta(p)$,
  $\ker df_\gamma(p)$) for $p\in U_{\alpha}\times e_1$
  (respectively $p\in U_{\beta}\times e_2$, $p\in U_{\gamma}\times
  e_3$).
\end{itemize}

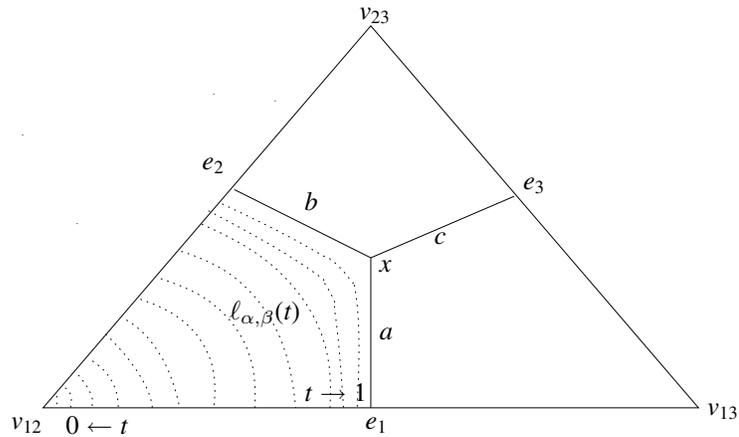
\begin{figure}
\centering
\begin{picture}(0,0)%
\includegraphics[scale=0.7]{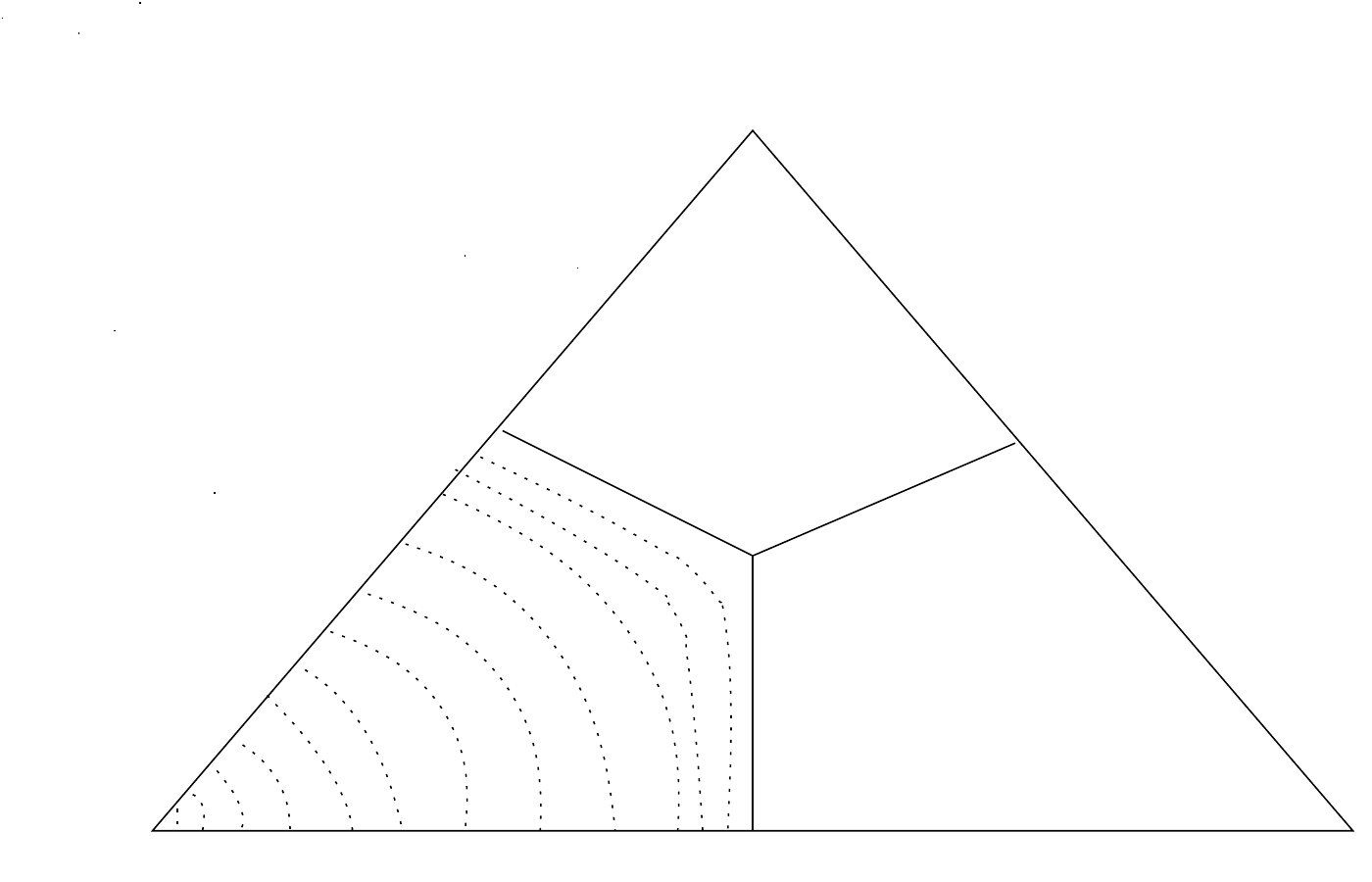}%
\end{picture}%
\setlength{\unitlength}{2763sp}%
\begingroup\makeatletter\ifx\SetFigFont\undefined%
\gdef\SetFigFont#1#2#3#4#5{%
  \reset@font\fontsize{#1}{#2pt}%
  \fontfamily{#3}\fontseries{#4}\fontshape{#5}%
  \selectfont}%
\fi\endgroup%
\begin{picture}(6703,4315)(-761,-3914)
\put(2987,-2442){\makebox(0,0)[lb]{\smash{{\SetFigFont{10}{12.0}{\rmdefault}{\mddefault}{\updefault}{\color[rgb]{0,0,0}\(x\)}%
}}}}
\put(2987,-3055){\makebox(0,0)[lb]{\smash{{\SetFigFont{10}{12.0}{\rmdefault}{\mddefault}{\updefault}{\color[rgb]{0,0,0}\(a\)}%
}}}}
\put(2314,-1891){\makebox(0,0)[lb]{\smash{{\SetFigFont{10}{12.0}{\rmdefault}{\mddefault}{\updefault}{\color[rgb]{0,0,0}\(b\)}%
}}}}
\put(3477,-2197){\makebox(0,0)[lb]{\smash{{\SetFigFont{10}{12.0}{\rmdefault}{\mddefault}{\updefault}{\color[rgb]{0,0,0}\(c\)}%
}}}}
\put(1640,-2871){\makebox(0,0)[lb]{\smash{{\SetFigFont{10}{12.0}{\rmdefault}{\mddefault}{\updefault}{\color[rgb]{0,0,0}\(\ell_{\alpha,\beta}(t)\)}%
}}}}
\put(170,-3900){\makebox(0,0)[lb]{\smash{{\SetFigFont{10}{12.0}{\rmdefault}{\mddefault}{\updefault}{\color[rgb]{0,0,0}\(0\gets t\)}%
}}}}
\put(2859,-3851){\makebox(0,0)[lb]{\smash{{\SetFigFont{10}{12.0}{\rmdefault}{\mddefault}{\updefault}{\color[rgb]{0,0,0}\(e_1\)}%
}}}}
\put(1395,-1523){\makebox(0,0)[lb]{\smash{{\SetFigFont{10}{12.0}{\rmdefault}{\mddefault}{\updefault}{\color[rgb]{0,0,0}\(e_2\)}%
}}}}
\put(4274,-1707){\makebox(0,0)[lb]{\smash{{\SetFigFont{10}{12.0}{\rmdefault}{\mddefault}{\updefault}{\color[rgb]{0,0,0}\(e_3\)}%
}}}}
\put(2314,-3606){\makebox(0,0)[lb]{\smash{{\SetFigFont{10}{12.0}{\rmdefault}{\mddefault}{\updefault}{\color[rgb]{0,0,0}\(t\to1\)}%
}}}}
\put(-320,-3858){\makebox(0,0)[lb]{\smash{{\SetFigFont{10}{12.0}{\rmdefault}{\mddefault}{\updefault}{\color[rgb]{0,0,0}\(v_{12}\)}%
}}}}
\put(2804,-176){\makebox(0,0)[lb]{\smash{{\SetFigFont{10}{12.0}{\rmdefault}{\mddefault}{\updefault}{\color[rgb]{0,0,0}\(v_{23}\)}%
}}}}
\put(5927,-3728){\makebox(0,0)[lb]{\smash{{\SetFigFont{10}{12.0}{\rmdefault}{\mddefault}{\updefault}{\color[rgb]{0,0,0}\(v_{13}\)}%
}}}}
\end{picture}%
\caption{Foliation of $T$}
\label{TriangleSpin}
\end{figure}

To use $\cL$ to define a $\Spin^\bC$--structure, we need to extend
it further.  Fix a point $x\in T$ and line segments $a$, $b$ and
$c$ from $x$ to the edges $e_1$, $e_2$ and $e_3$ of $\cL$,
respectively.  Let $v_{ij}=e_1\cap e_j$.  Let $\ell_{\alpha,\beta}(t)$, $\ell_{\beta,\gamma}(t)$,
and $\ell_{\alpha,\gamma}(t)$ denote the foliations of $T\setminus
(a\cup b\cup c)$, parametrized by $(0,1)$, so that as $t\to 0$,
$\ell_{\alpha,\beta}(t)$ degenerates to the corner $v_{12}$ and as
$t\to 1$ $\ell_{\alpha,\beta}(t)$ degenerates to $a\cup b$.  See
\fullref{TriangleSpin}.  
The map $\pi_T$ extends in an obvious way to a
map $\pi_T\co X_{\alpha,\beta,\gamma}\to T$.  Let
$\tilde{\ell}_{\alpha,\beta}(t)=\pi_{T}^{-1}(\ell_{\alpha,\beta}(t))\subset
X_{\alpha,\beta,\gamma}$,
$\tilde{\ell}_{\beta,\gamma}(t)=\pi_{T}^{-1}(\ell_{\beta,\gamma}(t))\subset
X_{\alpha,\beta,\gamma}$,
and $\tilde{\ell}_{\alpha,\gamma}(t)=\pi_T^{-1}(\ell_{\alpha,\gamma}(t))\subset
X_{\alpha,\beta,\gamma}$.

Choose $\phi$ so that:
\begin{enumerate}
\item The intersection of $F$ with each $\tilde{\ell}_{\alpha,\beta}(t)$, $\tilde{\ell}_{\beta,\gamma}(t)$ or $\tilde{\ell}_{\alpha,\gamma}(t)$ is a finite disjoint union of contractible $1$--complexes.
\item For all but finitely many $t$, $F\cap\tilde{\ell}_{\alpha,\beta}(t)$ consists of $g+1$ disjoint embedded arcs.
\item In some small neighborhood of the corner $v_{12}$
  (respectively $v_{23}$, $v_{13}$) of
  $T$, $\phi$ agrees with $\vx\times T\subset \Sigma\times T$
  (respectively $\vy\times T\subset \Sigma\times T$, $\vz\times
  T\subset \Sigma\times T$).
\item The intersections $F\cap \pi_T^{-1}(a)$, $F\cap\pi_T^{-1}(b)$ and $F\cap\pi_T^{-1}(c)$ each consist of $g+1$ disjoint embedded arcs.
\item The preimage under $\phi$ of
  $\alpha_i\times e_1$, $\beta_i\times e_2$ and $\gamma_i\times
  e_3$ is a connected arc for each $i$.
\end{enumerate}
Such $\phi$ exist.

Fix a $t$ such that $F\cap\tilde{\ell}_{\alpha,\beta}(t)$ consists of $g+1$ disjoint embedded arcs.  Observe that $\tilde{\ell}_{\alpha,\beta}(t)$ is
diffeomorphic to $Y_{\alpha,\beta}$.  For an appropriate choice of
the height functions $f_\alpha$ and $f_\beta$, they and a
parameter for the interval $\ell_{\alpha,\beta}(t)$ determine a
Morse function $f_t$ on $Y_{\alpha,\beta}$.  For an appropriate
choice of metric, the $2$--plane field
$\cL$ is the orthogonal complement of $\nabla f_t$.  As when we associated $\Spin^\bC$--structures on $3$--manifolds to intersection points, one can then use the
$(g+1)$--tuple of paths $F\cap\tilde{\ell}_{\alpha,\beta}(t)$ to replace $\cL$ with a $2$--plane field
defined on all of $\tilde{\ell}_{\alpha,\beta}(t)$.

Doing this construction uniformly in $t$, we can extend $\cL$ over
all of $\tilde{\ell}_{\alpha,\beta}(t)$ for
all $t$ such that $F\cap\tilde{\ell}_{\alpha,\beta}(t)$ consists
of $g+1$ disjoint arcs.  The same construction obviously works for
$\tilde{\ell}_{\beta,\gamma}(t)$ and $\tilde{\ell}_{\alpha,\gamma}(t)$.

Note that we can also extend $\cL$ over the boundary
$Y_{\alpha,\beta}\cup Y_{\beta,\gamma}\cup Y_{\alpha,\gamma}$, by
exactly the same method.

Now we have defined $\cL$ except on the intersection of $F$ with
\begin{itemize}
\item $\pi_T^{-1}(a\cup b\cup c)$ and
\item $\tilde{\ell}_{\alpha,\beta}(t_i)$,
  $\tilde{\ell}_{\beta,\gamma}(t_i')$ and
  $\tilde{\ell}_{\alpha,\gamma}(t_i^{\prime\prime})$ for some finite
  collection of $t_i$, $t_i'$ and $t_i^{\prime\prime}$.
\end{itemize}

Thus the region to which we have not extended $\cL$ consists of a
collection of disjoint contractible $1$--complexes.   So, we can
find an open ball $B$ in $X_{\alpha,\beta,\gamma}$ such that
$\cL$ is defined on $X_{\alpha,\beta,\gamma}\setminus B$.

Choose a metric on $X_{\alpha,\beta,\gamma}$.  Then the metric,
orientation, and $2$--plane field $\cL$ define an almost complex
structure on $X_{\alpha,\beta,\gamma}\setminus B$, and hence a
$\Spin^\bC$--structure on $X_{\alpha,\beta,\gamma}$.

The first question we address is how this construction
depends on $\phi$.
\begin{Lem}\label{Lemma:TriangleSpinDependence}The $\Spin^\bC$--structure just constructed depends only
  on the restriction of $\phi$ to the boundary of $S$.
\end{Lem}
\proof
Observe that in the construction, the restriction of $\phi$ to the
boundary of $S$ determined the $\Spin^\bC$--structure on
$X_{\alpha,\beta,\gamma}\setminus \Sigma\times \ol{T}$.  For a manifold
$M$, let $\Spin^\bC(M)$ denote the collection of
$\Spin^\bC$--structures on $M$.  We check that the restriction map
$\Spin^\bC(X_{\alpha,\beta,\gamma})\to
\Spin^\bC(X_{\alpha,\beta,\gamma}\setminus \Sigma\times \ol{T})$ is
injective; this then proves the result.

Recall that $\Spin^\bC(M)$ is an affine copy of $H^2(M;\bZ)$.  Further,
if $N\subset M$ then the restriction map from $\Spin^\bC(M)$ to
$\Spin^\bC(N)$ commutes with the $H^2$--action.  That is, let
$\xi\in \Spin^\bC(M)$, $\xi|_N$ its restriction to $\Spin^\bC(N)$,
$a\in H^2(M;\bZ)$ and $a|_N$ its restriction to $H^2(N;\bZ)$.  Then,
$(a\cdot \xi)|_N=a|_N\cdot \xi|_N$.  This can be proved, for
example, by thinking about \v Cech cocycles.

So, it suffices to prove that the restriction map
$H^2(X_{\alpha,\beta,\gamma};\bZ)\to H^2(X_{\alpha,\beta,\gamma}\setminus
\Sigma\times \ol{T};\bZ)$ is injective.  Consider the following commutative
diagram:
$$
\xymatrix@=22pt{
H^2(X_{\alpha,\beta,\gamma},X_{\alpha,\beta,\gamma}\setminus
\Sigma\times \ol{T};\bZ)\ar[r]\ar[d]^{\cong} &
H^2(X_{\alpha,\beta,\gamma};\bZ) \ar[r] &
H^2(X_{\alpha,\beta,\gamma}\setminus\Sigma\times \ol{T};\bZ)\\
H^2(\Sigma\times\overline{\bD},\Sigma\times\bdy\overline{\bD};\bZ)
\ar[r]^0 &
H^2(\Sigma\times\overline{\bD};\bZ) \ar@{^{(}->}[r]\ar[u] &
H^2(\Sigma\times\bdy\overline{\bD};\bZ)
}.
$$
(The top row is from the long exact sequence for the pair
$(X_{\alpha,\beta,\gamma}, X_{\alpha,\beta,\gamma}\setminus
\Sigma\times \ol{T})$.  The bottom is from the long exact sequence for the
pair $(\Sigma\times \ol{T}, \Sigma\times T)$.)  It follows from the
diagram that the map $H^2(X_{\alpha,\beta,\gamma};\bZ)\to
H^2(X_{\alpha,\beta,\gamma}\setminus\Sigma\times \ol{T};\bZ)$ is an injection.
\endproof

It would be nice to have a better way of presenting the construction
of a $\Spin^\bC$--structure on $X_{\alpha,\beta,\gamma}$ from an
element of $\pi_2(\vx,\vy,\vz)$.  Unfortunately I do not know one.

The following is~\cite[Proposition 8.4]{OS1}:
\begin{Lem}\label{Lemma:SzWellDefined}
The assignment described above induces a well--defined map
  $$\Ss_\Fz\co \pi_2(\vx,\vy,\vz)\to \Spin^\bC(X_{\alpha,\beta,\gamma}).$$
\end{Lem}
\proof
By the previous lemma, the construction depends only on the
restriction of $\phi$ to $\bdy S$.  This restriction is
defined up to isotopy by the element in $\pi_2(\vx,\vy,\vz)$.  It is
clear from the construction that an isotopy of $\phi$ does not
change the $\Spin^\bC$--structure constructed.
\endproof

For $\Ss\in \Spin^\bC(X_{\alpha,\beta,\gamma})$ we shall often write
$A\in\Ss$ to mean that $\Ss_{\Fz}(A)=\Ss$.

\begin{Def}
Given two triples of intersection points
$(\vx,\vy,\vz)$ and $(\vx',\vy',\vz')$, and 
$A\in\pi_2(\vx,\vy,\vz)$,
$A'\in\pi_2(\vx',\vy',\vz')$ define $\psi$ and $\psi'$ to be
\emph{$\Spin^\bC$--equivalent} if there exist elements
$B_{\alpha,\beta}\in\pi_2(\vx,\vx')$, $B_{\beta,\gamma}\in\pi_2(\vy,\vy')$, and
$B_{\alpha,\gamma}\in\pi_2(\vz,\vz')$ such that
$A+B_{\alpha,\beta}+B_{\beta,\gamma}+B_{\alpha,\gamma}=A'$.  We let $S_{\alpha,\beta,\gamma}$
denote the set of $\Spin^\bC$--equivalence classes of triangles.
\end{Def}

The following is~\cite[Proposition 8.5]{OS1}:
\begin{Lem}~\label{Lemma:SpinEquiv}
The assignment $\pi_2(\vx,\vy,\vz)\to
  \Spin^\bC(X_{\alpha,\beta,\gamma})$ defined above descends to a map
  $\Ss_{\Fz}\co S_{\alpha,\beta,\gamma}\to \Spin^\bC(X_{\alpha,\beta,\gamma})$.
  This new map is injective.  The image of $\Ss_\Fz$ consists of all
  $\Spin^\bC$--structures whose restrictions to $\bdy
  X_{\alpha,\beta,\gamma}$ are represented by intersection points.
\end{Lem}
\proof
Our proof is the same as in~\cite{OS1}.

The fact that $\Ss_\Fz$ descends to $S_{\alpha,\beta,\gamma}$
follows from the fact that the restriction of the 2--plane field
used to define $\Ss_\Fz(A)$ to $Y_{\alpha,\beta}$ is homologous to
the $2$--plane field used to define $\Ss_\Fz(\vx)$ (and similarly
for the restrictions to $Y_{\beta,\gamma}$ and
$Y_{\alpha,\gamma}$).  That is, suppose $A$ and $A'$
are $\Spin^\bC$--equivalent elements of $\pi_2(\vx,\vy,\vz)$,
$A'=A+B_{\alpha,\beta}+B_{\beta,\gamma}+B_{\alpha,\gamma}$.  Let
$\xi_A$ and $\xi_{A'}$ denote the $2$--plane fields constructed
above to define $\Ss_\Fz(A)$ and $\Ss_\Fz(A')$ respectively.  We
can assume that $\xi_A$ and $\xi_{A'}$ agree outside some collar
neighborhood $U$ of $\bdy W_{\alpha,\beta,\gamma}$.  Further, on
$Y_{\alpha,\beta}\times [0,1)\subset U$ (respectively
$Y_{\beta,\gamma}\times [0,1)\subset U$,
$Y_{\alpha,\gamma}\times[0,1)\subset U$) we can assume that
$\xi_A$ and $\xi_{A'}$ are given by $\nabla
f_{\alpha,\beta}^\perp$ (respectively $\nabla
f_{\beta,\gamma}^\perp$, $\nabla f_{\alpha,\gamma}^\perp$) for some
Morse function $f_{\alpha,\beta}$ on $Y_{\alpha,\beta}$
(respectively $f_{\beta,\gamma}$ on $Y_{\beta,\gamma}$,
$f_{\alpha,\gamma}$ on $Y_{\alpha,\gamma}$) outside some ball
neighborhoods of $(\vx\cup\Fz)\times[0,1)$ (respectively $(\vy\cup\Fz)\times[0,1)$,
$(\vz\cup\Fz)\times[0,1)$).  But it is then immediate from the
definition that
$\xi_A$ and $\xi_{A'}$ define the same $\Spin^\bC$--structure on
$X_{\alpha,\beta,\gamma}$.

Recall that the restriction map $\Spin^\bC(X_{\alpha,\beta,\gamma})\to
\Spin^\bC(\bdy X_{\alpha,\beta,\gamma})$ commutes with the
(transitive) actions of $H^2(X_{\alpha,\beta,\gamma})$ and
$H^2(\bdy X_{\alpha,\beta,\gamma})$ respectively.  It follows that
the cokernel $\Spin^\bC(\bdy
X_{\alpha,\beta,\gamma})/\Spin^\bC(X_{\alpha,\beta,\gamma})$ is
naturally identified with the cokernel of the restriction map
$H^2(X_{\alpha,\beta,\gamma};\bZ)\to H^2(\bdy
X_{\alpha,\beta,\gamma};\bZ)$.  This in turn is identified with the
image of the connecting homomorphism $\delta\co  H^2(\bdy
X_{\alpha,\beta,\gamma};\bZ)\to H^3(X_{\alpha,\beta,\gamma},\bdy
X_{\alpha,\beta,\gamma};\bZ)$.

Summarizing, we have a map $\epsilon'\co \Spin^\bC(\bdy
X_{\alpha,\beta,\gamma})\to H^3(X_{\alpha,\beta,\gamma},\bdy
X_{\alpha,\beta,\gamma};\bZ)$ given by the composition of the
coboundary map 
$$\delta\co H^2(\bdy
X_{\alpha,\beta,\gamma};\bZ)/H^2(X_{\alpha,\beta,\gamma};\bZ)\to
H^3(X_{\alpha,\beta,\gamma},\bdy X_{\alpha,\beta,\gamma};\bZ)$$
 with the projection 
$$\Spin^\bC(\bdy X_{\alpha,\beta,\gamma})\to \Spin^\bC(\bdy
X_{\alpha,\beta,\gamma})/\Spin^\bC(X_{\alpha,\beta,\gamma})=H^2(\bdy
X_{\alpha,\beta,\gamma})/H^2(X_{\alpha,\beta,\gamma}).$$ The element
$\epsilon'(\Ss)$ vanishes if and only if $\Ss$ can be extended to
all of $X_{\alpha,\beta,\gamma}$.

Recall that we defined earlier in this section the
obstruction $\epsilon(\vx,\vy,\vz)\in
H_1(X_{\alpha,\beta,\gamma})$ to the existence of elements in
$\pi_2(\vx,\vy,\vz)$.  Let $\Ss_\Fz(\vx,\vy,\vz)$ denote the
$\Spin^\bC$--structure induced by $\vx$, $\vy$ and $\vz$ on
$\bdy X_{\alpha,\beta,\gamma}$.  We next check that
$\epsilon(\vx,\vy,\vz)=PD(\epsilon'(\Ss_\Fz(\vx,\vy,\vz)))$, where
$PD$ denotes the Poincar\'e dual.

To show this, isotope the $\alpha$--, $\beta$-- and
$\gamma$--circles so that there are intersection points $\vx'$,
$\vy'$ and $\vz'$ with $\pi_2(\vx',\vy',\vz')\neq\emptyset$
(this is easy).  Then, $\epsilon(\vx',\vy',\vz')=0$.  We already showed
 that the $\Spin^\bC$--structure
$\Ss_\Fz(\vx',\vy',\vz')$ extends to all of
$X_{\alpha,\beta,\gamma}$, so
$\epsilon'(\Ss_\Fz(\vx',\vy',\vz'))=0$.  It is obvious from the
definitions that, up to a universal sign,
$\epsilon(\vx,\vy,\vz)-\epsilon(\vx',\vy',\vz')=i(\epsilon(\vx,\vx'))+i(\epsilon(\vy,\vy'))+i(\epsilon(\vz,\vz')).$
(Here, $i\co H_1(\bdy X_{\alpha,\beta,\gamma})$ $\to
H_1(X_{\alpha,\beta,\gamma})$ is the map induced by inclusion.)
Naturality of Poincar\'e duality gives the commutative diagram
\begin{displaymath}
\xymatrix{ 
H_1(\bdy X_{\alpha,\beta,\gamma})\ar[r]^i\ar@{<->}[d]^{PD} & 
H_1(X_{\alpha,\beta,\gamma})\ar@{<->}[d]^{PD}\\
H^2(\bdy X_{\alpha,\beta,\gamma})\ar[r]^(.43)\delta &
H^3(X_{\alpha,\beta,\gamma},\bdy X_{\alpha,\beta,\gamma})
}
\end{displaymath}
which implies that
$$PD\left(\epsilon(\vx,\vy,\vz)\right)-PD\left(\epsilon(\vx',\vy',\vz')\right)=\delta\left(PD(\epsilon(\vx,\vx')\oplus\epsilon(\vy,\vy')\oplus\epsilon(\vz,\vz'))\right).$$ 
From its definition, 
\begin{align*}
\epsilon'(\Ss_{\Fz}(\vx,\vy,\vz))-\epsilon'(\Ss_\Fz(\vx',\vy',\vz'))=\delta\left(
\left(\Ss_{\Fz}(\vx)-\Ss_{\Fz}(\vx')\right)\right.&\oplus
\left(\Ss_{\Fz}(\vy)-\Ss_{\Fz}(\vy')\right)\\&\left.\oplus
\left(\Ss_{\Fz}(\vz)-\Ss_{\Fz}(\vz')\right)\right).
\end{align*}
Now,
$PD\left(\epsilon(\vx,\vx')\right)=\Ss_{\Fz}(\vx)-\Ss_{\Fz}(\vx')$
(and similarly for $\vy$ and $\vz$).  So, since
$$PD\left(\epsilon(\vx',\vy',\vz')\right)=\epsilon'\left(\Ss_{\Fz}(\vx',\vy',\vz')\right),$$
it follows that
$PD\left(\epsilon(\vx,\vy,\vz)\right)=\epsilon'(\vx,\vy,\vz)$.

Now suppose we had a $\Spin^\bC$--structure $\Ss$ on
$X_{\alpha,\beta,\gamma}$ whose restriction to the boundary is
realized by intersection points $\vx$, $\vy$, and $\vz$.  Then
$0=\epsilon'(\vx,\vy,\vz)=\epsilon(\vx,\vy,\vz)$.  It follows that
$\Ss$ is in the image of $\Ss_\Fz$.  The converse is obvious.

Finally, it remains to show injectivity.  Fix a homology class
$A\in\pi_2(\vx,\vy,\vz)$.  The following diagram commutes:

\hbox to \hsize{\hss\small$
\xymatrix@C-10pt{
\hat{\pi}_2(\vx,\vx)\oplus\hat{\pi}_2(\vy,\vy)\oplus\hat{\pi}_2(\vz,\vz)
\ar[dd]^\cong\ar[r]^(.67){+A}
& \hat{\pi}_2(\vx,\vy,\vz)\ar[dd]^\cong_{-A}\ar[r]
& S_{\alpha,\beta,\gamma}\ar[r]\ar[dr]^{\Ss_{\Fz}}\ar@{^{(}->}[dd]_{[-A]}
& 0\\
& & & \Spin^\bC(X_{\alpha,\beta,\gamma})\ar[dl]^\cong_{-\Ss_{\Fz}(A)} \\
H_2(\bdy X_{\alpha,\beta,\gamma})\ar[r] 
& H_2(X_{\alpha,\beta,\gamma})\ar[r]
& H_2(X_{\alpha,\beta,\gamma},\bdy X_{\alpha,\beta,\gamma}) 
& .
}
$\hss}%

Injectivity of $\Ss_\Fz$ is immediate.
\endproof

\subsection{Definitions of the moduli spaces and maps}\label{Subsection:TriangleModuli1}

We now deal with the analysis involved in defining the triangle maps. We start by discussing the almost complex structures with which we will work. Fix a point
$\Fz_i$ in each component of
$\Sigma\setminus(\ba\cup\bb\cup\bc)$.  Fix an almost complex
structure $J_{\alpha,\beta}$ (respectively $J_{\beta,\gamma}$,
$J_{\alpha,\gamma}$) on $W_{\alpha,\beta}$ (respectively
$W_{\beta,\gamma}$, $W_{\alpha,\gamma}$) satisfying
(\textbf{J1})--(\textbf{J5}) and which achieves transversality for
holomorphic curves of index $\leq1$.   We will work with complex structures $J$ on $W_{\alpha,\beta,\gamma}$ such that:
\begin{itemize}
\item[(\textbf{J$'$1})] $J$ is tamed by $\eta$, the split symplectic form on
  $\Sigma\times T$.
\item[(\textbf{J$'$2})] In a neighborhood $U_{\{\Fz_i\}}$ of $\{\Fz_i\}\times T$,
 $J = j_\Sigma\times j_T$ is split.  (Here, $U_{\{\Fz_i\}}$ is small
 enough that its closure does not intersect
 $(\ba\cup\bb\cup\bc)\times[0,1]\times T$).
\item[(\textbf{J$'$3})] Near $\Sigma\times\{v_1\}$, $J$ agrees with
  $J_{\alpha,\beta}$
  Similarly, $J$ agrees with
  $J_{\beta,\gamma}$ near $\Sigma\times \{v_2\}$ and with 
$J_{\alpha,\gamma}$
  near $\Sigma\times \{v_3\}$.  
\item[(\textbf{J$'$4})] Projection $\pi_T$ onto $T$ is holomorphic and each fiber of $\pi_\Sigma$ is holomorphic.
\end{itemize}

For now, fix an almost complex structure $J$ satisfying (\textbf{J$'$1})--(\textbf{J$'$4}). For $A\in\pi_2(\vx,\vy,\vz)$, let $\cM^A$ denote the moduli space
of embedded holomorphic curves $u\co (S,\bdy S)\to
(W_{\alpha,\beta,\gamma}, \ba\times e_1\cup\bb\times e_2\cup\bc\times
e_3)$ asymptotic to $\vx$, $\vy$ and $\vz$ at the three ends of
$\Sigma\times T$ and in the homology class $A$.  We require that $u$
map exactly one component of $\bdy S$ to each of the $3g$
cylinders $\alpha_i\times e_1$, $\beta_i\times e_2$ and
$\gamma_i\times e_3$.  We also require that there be no components
of $S$ on which $\pi_T\circ u$ is constant.

We digress briefly to discuss the index of the
$\dbar$--operator for triangles.  Fix a homology class
$A\in\pi_2(\vx,\vy,\vz)$.  Suppose that $u\co S\to
W_{\alpha,\beta,\gamma}$ is a map in the homology class $A$.  Then
an argument similar to the one given in the beginning of
\fullref{Section:Index} shows that
$$
\ind(D\dbar)(u)=\frac{1}{2}g-\chi(S)+2e(A).$$
Here, $e(A)$ denotes the Euler measure of $A$, as described in
\fullref{Section:Index}.  The $\frac{1}{2}g$ looks strange, but
for $g$ odd it is easy to see that $2e(A)$ is a half--integer. When deriving this formula, one should keep in mind an extra $-g$, not appearing in \fullref{Section:Index}, coming from the Maslov index of $\pi_T\circ u$.

Analogously to \fullref{Subsection:Determining} we have:
\begin{Prop}\label{Prop:TriangleIndex}
The Euler characteristic of an embedded holomorphic curve
  $u\co S\to W_{\alpha,\beta,\gamma}$ is determined by the homology
  class of $u$.
\end{Prop}
\proof
The only part of the index which does not depend \emph{a priori} only
on $A$ is $\chi(S)$.
As in the proof of \fullref{Index:ChiCalcProp} we will
re--interpret the Euler characteristic as an intersection number.  

Let
$u\co S\to W_{\alpha,\beta,\gamma}$ be a holomorphic map in the
homology class $A$.  By the
Riemann--Hurwitz formula, the degree of branching of $\pi_\Sigma\circ
u$ determines $\chi(S)$, and vice versa.

Choose a diffeomorphism $T\to \ol{\bD^2}\setminus\{z^3=1\}$ from
$T$ to the unit disk with three boundary punctures.   Let $\xi_T$
be the vector field on $T$ induced by the vector field
$r\frac{\partial}{\partial\theta}$ on $\bD^2$.  Here, we choose
the diffeomorphism so that the preimage of $0$ is not a branch point
of $\pi_T\circ u$.
Let $\xi$ be the vector field $0\times \xi_T$ on
$W_{\alpha,\beta,\gamma}$.  Then, for $\epsilon$ small, the degree
of branching of $\pi_\Sigma\circ u$ is given by $\#(u\cap \exp_{\epsilon
  \xi}(u))-g$.  

Fix embedded holomorphic curves $u\co S\to W_{\alpha,\beta,\gamma}$, $u'\co S'\to W_{\alpha,\beta,\gamma}$ in the homology class $A$.
Because
we are in a low dimension (two), we can choose a (proper) bordism $v\co R\to
W_{\alpha,\beta,\gamma}$ from $u$ to $u'$.  Choose a Morse
function $f\co R\to\bR$ with $\bdy R=f^{-1}(\{0,1\})$, and so that
$v$ restricts to $u$ on $f^{-1}(0)$ and to $u'$ on
$f^{-1}(1)$.  We will think of $f$ as a coordinate on $R$, and
write $R_a$ for $f^{-1}(a)$ and $v_a=v|_{R_a}$.

For an appropriate choice of $R$, $v$ and $f$ there is a
partition $0=t_0<s_1<t_1<s_2<t_2<\cdots<s_k=1$ of $[0,1]$ such that
\begin{itemize}
\item The function $f$ has no critical values in $[s_i,t_i]$,
  $i=1,\cdots k$ (and so $f^{-1}([s_i,t_i])$ is a product).
\item For each point $s_i$ (respectively $t_i$),
  $v_{s_i}$ (respectively $v_{t_i}$) is an embedding in the homology class
  $A$ which projects as a branched cover to $T$.
\item The map $v|_{f^{-1}((t_i,s_{i+1}))}$ is constant near infinity
  (as a function of $f(p)$).
\end{itemize}
In words, we have chosen $R$, $v$ and $f$ so that we can
subdivide $[0,1]$ into subintervals over which either the map
doesn't change near the punctures or the topology of $S$ doesn't
change near the punctures (and so that each interval starts and ends
with maps for which the degree of branching of $\pi_\Sigma\circ v_a$
determines $\chi(R_a)$).

For each $i$, $\chi(R_{s_i})=\chi(R_{t_i})$, by the first property.

Also, if $a,a',a''\in(t_i,s_{i+1})$ are regular values then $\#(v_a\cap
v_{a'})=\#(v_a\cap v_{a''})$.  It follows from the Riemann--Hurwitz
formula that $\chi(R_{t_i})=\chi(R_{t_{i+1}}).$

But this implies that
$\chi(S)=\chi(R_{t_0})=\chi(R_{s_k})=\chi(S')$, completing the
proof.
\endproof

For $A\in\pi_2(\vx,\vy,\vz)$, let $\ind(A)$ denote the index of the
$D\dbar$ problem for embedded curves in the homology class $A$, if
such an embedded curve exists.
Note that the index is additive in the sense that for
$B_{\alpha,\beta}\in\pi_2(\vx,\vx')$,
$\ind(A+B_{\alpha,\beta})=\ind(A)+\ind(B_{\alpha,\beta})$ (and similarly
for the other ends), if both sides are defined.  We shall
omit the words ``if defined'' from all subsequent discussion of the
index -- that we are only discussing homology classes representable by
holomorphic curves will be implicit throughout.

\medskip\textbf{Remark}\qua  We could have proved 
\fullref{Prop:TriangleIndex} in more generality by
using an analog of \fullref{Lemma:Represent}.  Then, $\ind(A)$
would have a natural meaning for any homology class
$A\in\pi_2(\vx,\vy,\vz)$.  The only cases in which we are interested in the
index, however, are when there is an embedded holomorphic curve.

We now must check that, for generic $J$ satisfying (\textbf{J$'$1})--(\textbf{J$'$4}), the moduli spaces of $J$--holomorphic curves $\cM^A$ are reasonably well behaved.  Again, we achieve transversality by varying $J$.  
The argument to show that one can achieve transversality among
$J$ satisfying (\textbf{J$'$1})--(\textbf{J$'$4}) is analogous to the one given in
\fullref{Section:Transversality}.

Note that bubbling and Deligne--Mumford type degenerations in moduli spaces with $\ind\leq 1$ are
prohibited by the argument used in
\fullref{Section:Bubbling}.  So, by the compactness theorem~\cite[Theorem
10.2]{Ya2}, the compactification of $\cM^A$ consists of
holomorphic buildings with one story in $W_{\alpha,\beta,\gamma}$
and all of their other stories in $W_{\alpha,\beta}$,
$W_{\beta,\gamma}$ or $W_{\alpha,\gamma}$.  (As is standard in
symplectic field theory, the stories in the cylindrical bordisms
$W_{\alpha,\beta}$, $W_{\beta,\gamma}$, and $W_{\alpha,\gamma}$
are only defined up to translation.)

For the rest of this section, fix a $\Spin^\bC$ equivalence class of
triangles.  Denote it $\Ss_{\alpha,\beta,\gamma}$.  
Fix a complex structure
satisfying (\textbf{J$'$1})--(\textbf{J$'$4}) and achieving transversality.

There are still some technical details to address before we can define
the triangle maps.  However, we will give the definitions now, asking
the reader to trust that all the symbols make sense and sums are
finite.  We will justify this trust presently.

Given the choice of $\Ss_{\alpha,\beta,\gamma}\in
\Spin^\bC(X_{\alpha,\beta,\gamma})$ we will define a
map
$$f^\infty_{\alpha,\beta,\gamma}\co CF^\infty(Y_{\alpha,\beta})\otimes_{\bZ[U]}
CF^\infty(Y_{\beta,\gamma})\to
CF^\infty(Y_{\alpha,\gamma})$$
by
$$
f^\infty_{\alpha,\beta,\gamma}([\vx,i]\otimes[\vy,j]) =
\sum_{\vz}\sum_{\substack{A\in\pi_2(\vx,\vy,\vz)\cap\Ss_{\alpha,\beta,\gamma}\\\ind(A)=0}}
\left(\#\cM^A\right)[\vz,i+j-n_z(A)].
$$

\medskip\textbf{Remark}\qua  
For the complexes $CF^\infty$, $CF^-$, $CF^+$ and
$CF^{\leq0}$ our tensor products shall always be over $\bZ[U]$.
For $\widehat{CF}$ they shall be over $\bZ$.  In the case of
$CF^\infty$ it would be equivalent to take the tensor product over
$\bZ[U,U^{-1}]$.  It is not equivalent, and quickly leads to
nonsense, to take all the tensor products over $\bZ$.  The corresponding
remark also applies to the $Hom$ functor, if one wanted to obtain
cohomology theories; cf~\cite[Section 2]{OS2}.

There are two obvious issues that need to be addressed.  Firstly,
since we have been working with $\bZ$--coefficients, the symbol
``$\#$'' implies that we have chosen orientations for the $\cM^A$,
which should presumably be consistent with the orientations for the
moduli spaces for $\cH_{\alpha,\beta}$, $\cH_{\beta,\gamma}$ and
$\cH_{\alpha,\gamma}$.  We will address this issue in
\fullref{Subsection:TriangleOrientations}.  Secondly, we need
to know that the coefficient of $[\vz,k]$ in
$f^\infty_{\alpha,\beta,\gamma}([\vx,i],[\vy,j])$ is a finite sum.
This will require that we impose an admissibility condition on the
Heegaard triple--diagram, as we will discuss in
\fullref{Subsection:TriangleAdmissibility}.  Note, however,
that if we work with $\bZ/2$--coefficients and if
$H_2(X_{\alpha,\beta,\gamma})$ is finite then the formula defining
$f^\infty_{\alpha,\beta,\gamma}$ already makes perfect sense.

Before addressing the issues of orientations and admissibility, we define the rest of the maps that will appear.

We will define a map 
$$\hat{f}_{\alpha,\beta,\gamma}\co \widehat{CF}(Y_{\alpha,\beta})\otimes_{\bZ}
\widehat{CF}(Y_{\beta,\gamma})\to
\widehat{CF}(Y_{\alpha,\gamma})$$
by
$$
\hat{f}_{\alpha,\beta,\gamma}(\vx\otimes\vy) =
\sum_{\vz}\sum_{\substack{A\in\hat{\pi}_2(\vx,\vy,\vz)\cap\Ss_{\alpha,\beta,\gamma}\\\ind(A)=0}}
\left(\#\cM^A\right)\vz.
$$

There is a subcomplex $CF^{\leq0}$ of $CF^\infty$ generated by all
$[\vx,i]$ with $i\leq 0$.  The homology $HF^{\leq0}$ of
$CF^{\leq0}$ is naturally isomorphic to $HF^-$, but the results on
triangles are phrased most simply in terms of $CF^{\leq0}$.  The map
$f^\infty_{\alpha,\beta,\gamma}$ restricts to a map
$f^{\leq0}_{\alpha,\beta,\gamma}\co CF^{\leq0}(Y_{\alpha,\beta})\otimes
CF^{\leq0}(Y_{\beta,\gamma})\to
CF^{\leq0}(Y_{\alpha,\gamma}).$  Hence it also induces a map
$$f^+_{\alpha,\beta,\gamma}\co CF^+(Y_{\alpha,\beta})\otimes
CF^{\leq0}(Y_{\beta,\gamma};M_{\beta,\gamma})\to
CF^+(Y_{\alpha,\gamma}).$$

Note that $f^{\infty}_{\alpha,\beta,\gamma}$ is a map of
$\bZ[U]$--modules (in fact, $\bZ[U,U^{-1}]$--modules), so
$f^{\leq}_{\alpha,\beta,\gamma}$ and $f^+_{\alpha,\beta,\gamma}$
are also maps of $\bZ[U]$--modules.  Further, since $CF^{\leq0}(Y_{\beta,\gamma})$ is a subcomplex of
$CF^\infty(Y_{\beta,\gamma})$, the map $f^\infty_{\alpha,\beta,\gamma}$ restricts to a map
$CF^\infty(Y_{\alpha,\beta})\otimes_{\bZ[U]}CF^{\leq0}(Y_{\beta,\gamma})$ $\to
CF^\infty(Y_{\alpha,\gamma})$, which we also denote
$f^\infty_{\alpha,\beta,\gamma}$. 

The next step is to show that the maps just defined are chain maps.  The proof is completely standard, but before giving it we digress to deal with orientations and admissibility.

\subsection{Orientations}\label{Subsection:TriangleOrientations}
By the same arguments as in \fullref{Section:Orientations}, it
follows that the moduli space $\cM^A(X_{\alpha,\beta,\gamma})$ are
orientable, and by general arguments we can find orientations for
all the $\cM^A(X_{\alpha,\beta,\gamma})$,
$\cM^{B_{\alpha,\beta}}(Y_{\alpha,\beta})$,
$\cM^{B_{\beta,\gamma}}(Y_{\beta,\gamma})$ and
$\cM^{B_{\alpha,\gamma}}(Y_{\alpha,\gamma})$ (or rather, for the
determinant line bundles over the corresponding configuration spaces)
which are consistent with gluings.  However, we
would like somewhat more.  Specifically: 

\begin{Lem}{\rm(Compare~\cite[Lemma 8.7]{OS1})}%
\label{Lemma:TriangleOrientations}\qua
Given coherent orientation systems
$\frako_{\alpha,\beta}(B_{\alpha,\beta})$ and
$\frako_{\beta,\gamma}(B_{\beta,\gamma})$ (for all $B_{\alpha,\beta}$
and $B_{\beta,\gamma}$) 
there are orientation systems $\frako_{\alpha,\gamma}$
and $\frako_{\alpha,\beta,\gamma}$ consistent with
$\frako_{\alpha,\beta}$ and $\frako_{\beta,\gamma}$.
\end{Lem}

Note that we have not claimed that $\frako_{\alpha,\gamma}$ and $\frako_{\alpha,\beta,\gamma}$ are unique.  The indeterminacy will be clear from the proof.

\proof[Proof of \fullref{Lemma:TriangleOrientations}]
Our proof is the same as in~\cite[Section 8.2]{OS1}.

Fix a $\Spin^\bC$--structure $\Ss_{\alpha,\beta,\gamma}$ on
  $X_{\alpha,\beta,\gamma}$ and intersection points $\vx_0\in
  \Ss_{\alpha,\beta,\gamma}|_{\Yab}$, $\vy_0\in\Ss_{\alpha,\beta,\gamma}|_{\Ybc}$ and $\vz_0\in\Ss_{\alpha,\beta,\gamma}|_{\Yac}$.
  Fix $A_0\in\pi_2(\vx_0,\vy_0,\vz_0)$.  Choose any orientation over
  $A_0$.

Let 
{\setlength\arraycolsep{3pt}
\begin{eqnarray*}
K=\{B_{\alpha,\gamma}\in\hat{\pi}_2(\vz_0,\vz_0)&|&\exists
B_{\alpha,\beta}\in\hat{\pi}_2(\vx_0,\vx_0),B_{\beta,\gamma}\in\hat\pi_2(\vy_0,\vy_0)\\ && \textrm{
  such that }
A_0+B_{\alpha,\gamma}=A_0+B_{\beta,\gamma}+B_{\alpha,\gamma}\}.
\end{eqnarray*}}

\begin{Sublem}\label{Sublemma:1}
\begin{enumerate}
\item $\hat{\pi}_2(\vz_0,\vz_0)\cong K\oplus\bZ^N$ for some $N$.
\item Given $B_{\alpha,\gamma}\in K$ there is only one pair $B_{\beta,\gamma}$, $B_{\alpha,\gamma}$ such that $A_0+B_{\alpha,\gamma}=A_0+B_{\beta,\gamma}+B_{\alpha,\gamma}$.
\end{enumerate}
\end{Sublem}

Assuming the sublemma, the lemma is almost immediate.  By Part 2 of
the sublemma, $\frako_{\alpha,\beta,\gamma}(A_0)$,
$\frako_{\alpha,\beta}$, and $\frako_{\beta,\gamma}$ determine
$\frako_{\alpha,\gamma}$ over $K$.  Choosing $\frako_{\alpha,\gamma}$
arbitrarily over a basis of $\bZ^N$ determines
$\frako_{\alpha,\gamma}$ over $\hat{\pi}_2(\vz_0,\vz_0)$.  The orientation over $[\Sigma]$ is determined as in \fullref{Section:Orientations}.  Choosing a homology class
$B_{\alpha,\gamma,\vz}\in\pi_2(\vz_0,\vz)$ for each intersection
point $\vz$ (ie, a complete set of paths for $Y_{\alpha,\gamma}$
in the sense of \fullref{Section:Orientations}) and then choosing
orientations arbitrarily over the $B_{\alpha,\gamma,\vz}$ determines
$\frako_{\alpha,\gamma}$ over all of $\Ss_{\alpha,\beta,\gamma}|_{\Yac}$.  The orientation
over $A_0$ and $\frako_{\alpha,\beta}$, $\frako_{\beta,\gamma}$ and
$\frako_{\alpha,\gamma}$ together determine
$\frako_{\alpha,\beta,\gamma}$ over all of $\Ss_{\alpha,\beta,\gamma}$.  This completes the
proof, except for the

\proof[Proof of \fullref{Sublemma:1}]
The subgroup $K$ is canonically identified with the intersection of $H_2(\Yab\cup\Ybc)$ and $H_2(\Yac)$ in $H_2(X_{\alpha,\beta,\gamma})$.  From the fragment
$$
H_3(X_{\alpha,\beta,\gamma},Y_{\alpha,\beta}\cup Y_{\beta,\gamma})\to H_2(\Yab\cup\Ybc)\stackrel{j}{\to} H_2(X_{\alpha,\beta,\gamma})\stackrel{p}{\to}H^2(X_{\alpha,\beta,\gamma},\Yab\cup\Ybc)
$$
of the long exact sequence for the pair $(X_{\alpha,\beta,\gamma},Y_{\alpha,\beta}\cup Y_{\beta,\gamma})$, one sees
$\hat{\pi}_2(\vz_0,\vz_0)/K\cong p\circ i(H_2(\Yac)).$  By
excision, $H_*(X_{\alpha,\beta,\gamma},\Yab\cup\Ybc)\cong
H_*(U_\beta\times[0,1],U_\beta\times\bdy[0,1]).$  So,
$H_3(X_{\alpha,\beta,\gamma},\Yab\cup\Ybc)=0$, so $j$ is
injective, implying part (2) of the sublemma.  Further, 
\begin{eqnarray}
H_2(X_{\alpha,\beta,\gamma},\Yab\cup\Ybc)&\cong&
H_2(U_\beta\times[0,1],U_\beta\times\bdy[0,1]) \nonumber\\
&\cong&
H_2\left((\bigvee_{i=1}^gS^1)\times S^1,
  (\bigvee_{i=1}^gS^1)\times\{\textrm{pt}\}\right)\nonumber\\
&\cong& \bZ^g\nonumber
\end{eqnarray}
(from the long exact sequence for the pair), so $\hat{\pi}_2(\vz_0,\vz_0)/K$ is free Abelian.  It follows that the sequence 
$$
0\to K\to \hat{\pi}_2(\vz_0,\vz_0)\to\hat{\pi}_2(\vz_0,\vz_0)/K\to 0
$$
splits, yielding the result.
\endproof
\endproof

In the rest of this section, we shall always assume that coherent orientations have been chosen for the moduli spaces under consideration, but shall suppress the orientation systems from the notation.
\subsection{Admissibility}\label{Subsection:TriangleAdmissibility}
As when we defined the chain complexes, we will need the Heegaard
triple--diagram to satisfy certain admissibility criteria in order to
ensure finiteness when we define maps between Floer homologies.

\begin{Def} The pointed Heegaard triple--diagram $\cH^3$ is \emph{weakly
    admissible} if the following condition is met.  For
  any $B_{\alpha,\beta}\in\hat{\pi}_2(\vx,\vx)$,
  $B_{\beta,\gamma}\in\hat{\pi}_2(\vy,\vy)$, and
  $B_{\alpha,\gamma}\in\hat{\pi}_2(\vz,\vz)$ we require that 
  $B_{\alpha,\beta}+B_{\beta,\gamma}+B_{\alpha,\gamma}$ have both positive and negative coefficients (or be identically zero).
\end{Def}
Note that the definition given in \fullref{Section:Admissibility} of weak admissibility here corresponds to the definition of
weak admissibility for all $\Spin^\bC$--structures.
\begin{Def}Fix a $\Spin^\bC$--structure $\Ss$ on
  $X_{\alpha,\beta,\gamma}$ and let $\Ss_{\alpha,\beta}$,
  $\Ss_{\beta,\gamma}$ and $\Ss_{\alpha,\gamma}$ be the
  restrictions of $\Ss_{\alpha,\beta,\gamma}$ to
  $Y_{\alpha,\beta}$, $Y_{\beta,\gamma}$ and $Y_{\alpha,\gamma}$
  respectively.  We say that $\cH^3$ is
  \emph{strongly admissible for $\Ss_{\alpha,\beta,\gamma}$} if for
  any $\vx\in\Ss_{\alpha,\beta}$, $\vy\in\Ss_{\beta,\gamma}$,
  $\vz\in\Ss_{\alpha,\gamma}$,
  $B_{\alpha,\beta}\in\hat{\pi}_2(\vx,\vx)$,
  $B_{\beta,\gamma}\in\hat{\pi}_2(\vy,\vy)$
  and $B_{\alpha,\gamma}\in\hat{\pi}_2(\vz,\vz)$ with
$$
\left\langle c_1(\Ss_{\alpha,\beta}),B_{\alpha,\beta}\right\rangle + \left\langle
c_1(\Ss_{\beta,\gamma}),B_{\beta,\gamma}\right\rangle + \left\langle
c_1(\Ss_{\alpha,\gamma}),B_{\alpha,\gamma}\right\rangle = 2n\geq 0
$$
and $B_{\alpha,\beta}+B_{\beta,\gamma}+B_{\alpha,\gamma}$ not identically zero
there is some coefficient of $A+B+C$ strictly greater than $n$.  
\end{Def}
Note that weak (respectively strong) admissibility for $\cH^3$ implies weak
(strong) admissibility for each of $\cH_{\alpha,\beta}$,
$\cH_{\beta,\gamma}$ and $\cH_{\alpha,\gamma}$.

The proof of the following alternate characterization of weak
admissibility is the same as the proof of  \fullref{Lemma:AltAdmiss}.
\begin{Lem}
With notation as above:
\begin{itemize}
\item The diagram $\cH^3$ is weakly admissible if and only
  if there is an area form on $\Sigma$ with respect to which for any
  $B_{\alpha,\beta}$, $B_{\beta,\gamma}$ and $B_{\alpha,\gamma}$
  as in the definition of weak admissibility, the domain
  $B_{\alpha,\beta}+B_{\beta,\gamma}+B_{\alpha,\gamma}$ has zero
  signed area.  
\item The diagram $\cH^3$ is strongly admissible for $\Ss$ if there
  is an area form on $\Sigma$ with respect to which for any
  $B_{\alpha,\beta}$,
  $B_{\beta,\gamma}$ and $B_{\alpha,\gamma}$ as in the definition
  of strong admissibility,
  $B_{\alpha,\beta}+B_{\beta,\gamma}+B_{\alpha,\gamma}$ has signed
  area
  equal to $n$, and with respect to which $\Sigma$ has area $1$.
\end{itemize}
\end{Lem}

Recall that we call a homology class $A$ positive if the corresponding domain has no negative coefficients.

The following is~\cite[Lemma 8.9]{OS1}
\begin{Lem}\label{Lemma:WeakTriangleAdmisImplies}Suppose $\cH^3$ is weakly admissible.  Fix intersection
  points $\vx$, $\vy$, and $\vz$ and a $\Spin^\bC$--structure
  $\Ss$ on $X_{\alpha,\beta,\gamma}$.  Then for each $j,k\in\bZ$
  there are only finitely many positive $A\in\Ss_{\alpha,\beta,\gamma}\cap\pi_2(\vx,\vy,\vz)$ such that
\begin{itemize}
\item $\ind(A)=j$,
\item $n_\Fz(A)=k$,
\end{itemize}
\end{Lem}
\proof
Suppose that $A,A'\in\pi_2(\vx,\vy,\vz)$, $n_\Fz(A)=n_\Fz(A')=k$, and
$A,A'\in\Ss_{\alpha,\beta,\gamma}$.  Then by \fullref{Lemma:SpinEquiv},
$A$ and $A'$ are $\Spin^\bC$--equivalent, so
$A'=A+B_{\alpha,\beta}+B_{\beta,\gamma}+B_{\alpha,\gamma}$ where
$B_{\alpha,\beta}\in\hat{\pi}_2(\vx,\vx)$,
$B_{\beta,\gamma}\in\hat{\pi}_2(\vy,\vy)$, and
$B_{\alpha,\gamma}\in\hat{\pi}_2(\vz,\vz)$.
By the previous lemma, we can choose an area form on $\Sigma$ so
that $B_{\alpha,\beta}+B_{\beta,\gamma}+B_{\alpha,\gamma}$ has zero
signed area.  The result then follows as in \fullref{Lemma:WeakAdmis}.
\endproof
It follows from this lemma and compactness that if $\cH^3$ satisfies the weak admissibility criterion then the sums  defining $\hat{f}_{\alpha,\beta,\gamma}$ and $f^+_{\alpha,\beta,\gamma}$ are finite.

The following is~\cite[Lemma 8.10]{OS1}.
\begin{Lem}\label{Lemma:StrongTriangleAdmisImplies} Fix $j\in\bZ$, intersection points $\vx$, $\vy$, $\vz$,
  and a $\Spin^\bC$--structure $\Ss_{\alpha,\beta,\gamma}$. Suppose $\cH^3$
  is strongly admissible for $\Ss$.  Then there are only finitely
  many positive
  $A\in\pi_2(\vx,\vy,\vz)\cap\Ss_{\alpha,\beta,\gamma}$ such that
  $\ind(A)=j$.
\end{Lem}
\proof
As in the previous proof, given $A,A'\in\pi_2(\vx,\vy,\vz)$
satisfying the hypotheses,
$A-A'=B_{\alpha,\beta}+B_{\beta,\gamma}+B_{\alpha,\gamma}$ for some $B_{\alpha,\beta}\in\pi_2(\vx,\vx)$, $B_{\beta,\gamma}\in\pi_2(\vy,\vy)$ and $B_{\alpha,\gamma}\in\pi_2(\vz,\vz)$.
The proof then follows as in \fullref{Lemma:StrongAdmis}.
\endproof
It follows from this lemma and compactness that if $\cH^3$ satisfies the strong admissibility criterion then the sums  defining $f^\infty_{\alpha,\beta,\gamma}$ and $f^{\leq0}_{\alpha,\beta,\gamma}$ are finite.

The following is~\cite[Lemma 8.11]{OS1}.  We refer the reader there for its proof.
\begin{Prop} \label{Prop:TriangleAdmisExists}
Given any pointed Heegaard triple--diagram
  $(\Sigma,\va,\vb,\vc,\Fz)$ there is an isotopic weakly admissible
  Heegaard triple--diagram.  Given any pointed Heegaard
  triple--diagram $(\Sigma,\va,\vb,\vc,\Fz)$ and a
  $\Spin^\bC$--structure $\Ss_{\alpha,\beta,\gamma}$ on
  $X_{\alpha,\beta,\gamma}$ there is an isotopic pointed Heegaard triple--diagram
  which is strongly admissible for $\Ss_{\alpha,\beta,\gamma}$.
\end{Prop}

\subsection{Moduli spaces and maps, part 2}\label{Subsection:TriangleModuli2}
If we wish to make a statement about all of $f^\infty_{\alpha,\beta,\gamma}$,
$f^{\leq0}_{\alpha,\beta,\gamma}$, $f^+_{\alpha,\beta,\gamma}$, or
$\hat{f}_{\alpha,\beta,\gamma}$ at once we will simply write
$f_{\alpha,\beta,\gamma}$.  For example
\begin{Lem}\label{Lemma:TriangleMapsChain}The maps $f_{\alpha,\beta,\gamma}$ are chain maps.
\end{Lem}
\proof
This follows by considering the 1--dimensional moduli spaces $\cM^{A_{\alpha,\beta,\gamma}}$
where $A_{\alpha,\beta,\gamma}\in\pi_2(\vx,\vy,\vz)$, $\ind(A_{\alpha,\beta,\gamma})=1$. 
The proof of \fullref{CompactProp} still works,
so the boundary of $\cM^{A_{\alpha,\beta,\gamma}}$ consists of
height two holomorphic buildings in the homology class $A_{\alpha,\beta,\gamma}$.
 Each of of these height two
holomorphic buildings consists of
\begin{enumerate}
\item a curve of index 0 in $\Sigma\times T$ and
\item a curve of index $1$ (defined up to translation) in one of
  $W_{\alpha,\beta}$, $W_{\beta,\gamma}$ or $W_{\gamma,\alpha}$,
\end{enumerate}
and every such building is in $\bdy\cM^{A_{\alpha,\beta,\gamma}}$
for some $A_{\alpha,\beta,\gamma}$ of index $1$.
This follows from~\cite[Theorem 10.2]{Ya2} and
\fullref{Gluing:Cylindrical}.
Hence,
\begin{eqnarray}
0=\#\left(\bdy\cM^C\right)&=& 
\sum_{\substack{B_{\alpha,\beta}\in\pi_2(\vx,\vx')\\\ind(B_{\alpha,\beta})=1}}
\left(\#\hcM^{B_{\alpha,\beta}}\right)\left(\#\cM^{A_{\alpha,\beta,\gamma}-B_{\alpha,\beta}}\right)
\nonumber\\
&+&
\sum_{\substack{B_{\beta,\gamma}\in\pi_2(\vy,\vy')\\\ind(B_{\beta,\gamma})=1}}
\left(\#\hcM^{B_{\beta,\gamma}}\right)\left(\#\cM^{A_{\alpha,\beta,\gamma}-B_{\beta,\gamma}}\right)\nonumber\\
&+& \sum_{\substack{B_{\alpha,\gamma}\in\pi_2(\vz',\vz)\\\ind(B_{\alpha,\gamma})=1}}
\left(\#\hcM^B_{\alpha,\gamma}\right)\left(\#\cM^{A_{\alpha,\beta,\gamma}-B_{\alpha,\gamma}}\right).\nonumber
\end{eqnarray}
But summing this over $A_{\alpha,\beta,\gamma}$ with $n_\Fz(A_{\alpha,\beta,\gamma})=k$ gives the coefficient
of $[\vz,i+j-k]$ in $\bdy\circ 
f^\infty_{\alpha,\beta,\gamma}+f^\infty_{\alpha,\beta,\gamma}\circ
\bdy$, proving the result for $f^{\infty}_{\alpha,\beta,\gamma}.$

The results for $f^{\leq0}_{\alpha,\beta,\gamma}$ and
$f^+_{\alpha,\beta,\gamma}$ follow immediately.  The proof for
$\hat{f}_{\alpha,\beta,\gamma}$ is analogous, restricting to curves
with $n_\Fz=0$; we leave the details of this case to the reader.
\endproof

We use $F_{\alpha,\beta,\gamma}$ (appropriately decorated) to denote
the maps on homology induced by $f_{\alpha,\beta,\gamma}$.

\begin{Lem}\label{Lemma:TriangleCXIndep}The maps
  $F_{\alpha,\beta,\gamma}$ just defined are independent of the
  choice of complex structure $J$ on $W_{\alpha,\beta,\gamma}$
  satisfying {\rm(\textbf{J$'$1})--(\textbf{J$'$4})}.
\end{Lem}
\proof
Suppose $J$ and $J'$ are two complex structures on
$W_{\alpha,\beta,\gamma}$ satisfying (\textbf{J$'$1})--(\textbf{J$'$4}).  Note in particular that $J$ and $J'$ agree on the ends of
$W_{\alpha,\beta,\gamma}$.  Let $f^\infty_{\alpha,\beta,\gamma,J}$
and $f^\infty_{\alpha,\beta,\gamma,J'}$ denote the maps defined
above, computed with respect to $J$ and $J'$ respectively.

Choose a generic path $J_t$ connecting
$J$ to $J'$, which is fixed on the ends of
$W_{\alpha,\beta,\gamma}$.  Then for any $k\in\bZ$ there are a
finite collection of 
$t\in(0,1)$ such that $\cM^{A_{\alpha,\beta,\gamma}}\neq\emptyset$ for some $A_{\alpha,\beta,\gamma}$ with
$\ind(A_{\alpha,\beta,\gamma})=-1$ and
$n_\Fz(A_{\alpha,\beta,\gamma})\leq k$.  (This uses the
admissibility
hypothesis.)  Define a map
$\Phi\co CF^\infty(Y_{\alpha,\beta})\otimes CF^\infty(Y_{\beta,\gamma})\to
CF^\infty(Y_{\alpha,\gamma})$ by
$$
\Phi([\vx,i],[\vy,j])=\sum_{\vz}\sum_{\substack{(A_{\alpha,\beta,\gamma},t)\in\pi_2(\vx,\vy,\vz)\times(0,1)\\\ind(A_{\alpha,\beta,\gamma})=-1}}
\#\cM^{A_{\alpha,\beta,\gamma}}[\vz,i+j-n_\Fz(A_{\alpha,\beta,\gamma})].
$$
The coefficient of each $[\vz,k]$ is a finite sum.  By exactly the
same argument as used in \fullref{Isotopy:MapsIndependent}, $\Phi$ is a
chain homotopy from $f^\infty_{\alpha,\beta,\gamma,J}$ to
$f^\infty_{\alpha,\beta,\gamma,J'}$.

The results for $F^\infty_{\alpha,\beta,\gamma}$,
$F^{\leq0}_{\alpha,\beta,\gamma}$ and 
$F^+_{\alpha,\beta,\gamma}$ follow.  The result for
$\hat{F}_{\alpha,\beta,\gamma}$ is proved in an analogous way; as
has become our habit we leave the details of this case to the reader.
\endproof

\begin{Lem}\label{TrianglesIndependent}
The maps
  $F_{\alpha,\beta,\gamma}$ are independent of the choices of
  complex structures $J_{\alpha,\beta}$, $J_{\beta,\gamma}$ and
  $J_{\alpha,\gamma}$ satisfying {\rm(\textbf{J1})--(\textbf{J5})} and
  achieving transversality, and
  are independent of isotopies of the
  $\alpha$, $\beta$ and $\gamma$ preserving the admissibility
  hypotheses and not crossing $\Fz$.
\end{Lem}
\proof
Let $\alpha'$, $\beta'$ and $\gamma'$ be isotopic to
$\alpha$, $\beta$, and $\gamma$.  As in \fullref{Section:Isotopy}, we may assume the isotopy
introduces or cancels only one pair of intersection points, and can
be realized by Lagrangian cylinders.
Let $J_{\alpha',\beta'}$, $J_{\beta',\gamma'}$ and 
$J_{\alpha',\gamma'}$ be a complex
structures on $W_{\alpha,\beta}$, $W_{\beta,\gamma}$ and
$W_{\alpha,\gamma}$ respectively satisfying
(\textbf{J1})--(\textbf{J5}) and achieving transversality.  Choose a path of complex
structures $J_{t,\alpha,\beta}$ (respectively $J_{t,\beta,\gamma}$, 
$J_{t,\alpha,\gamma}$), $t\in[0,1]$, interpolating between
$J_{\alpha,\beta}$
and $J_{\alpha',\beta'}$ (respectively $J_{\beta,\gamma}$ and
$J_{\beta',\gamma'}$, $J_{\alpha,\gamma}$ and
$J_{\alpha',\gamma'}$), as in \fullref{Section:Isotopy}.

\begin{figure}
\centering
\begin{picture}(0,0)%
\includegraphics[scale=0.8]{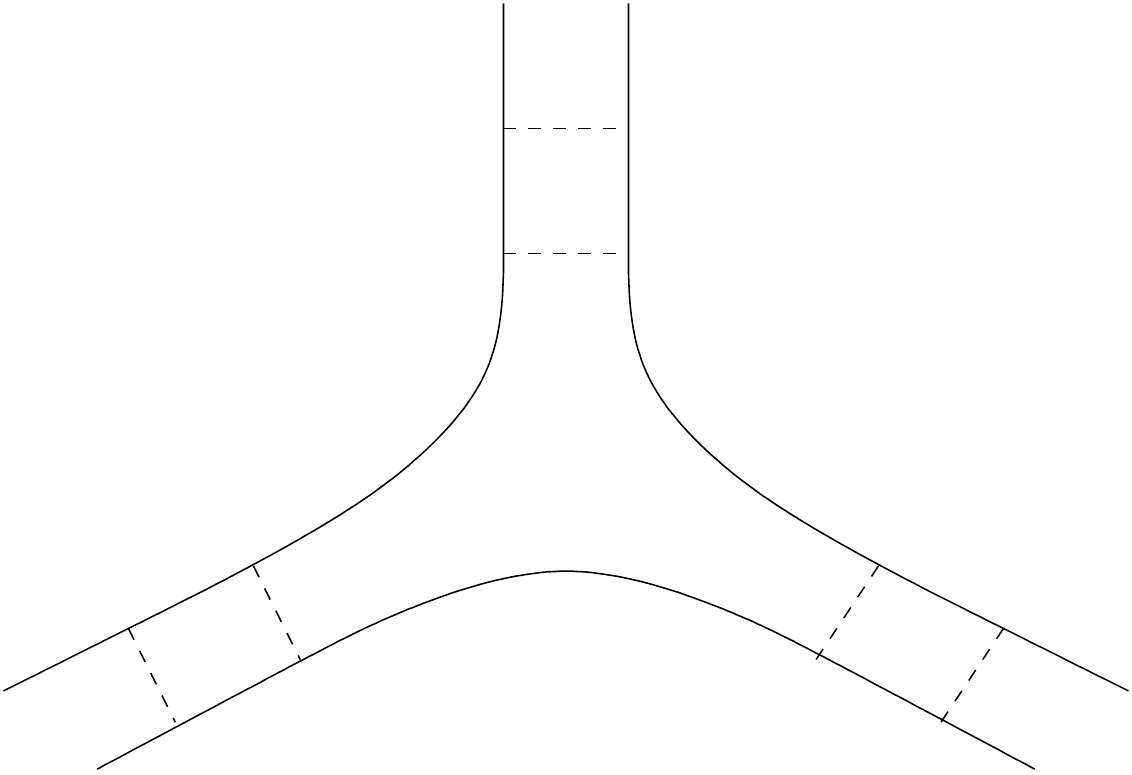}%
\end{picture}%
\setlength{\unitlength}{3157sp}%
\begingroup\makeatletter\ifx\SetFigFont\undefined%
\gdef\SetFigFont#1#2#3#4#5{%
  \reset@font\fontsize{#1}{#2pt}%
  \fontfamily{#3}\fontseries{#4}\fontshape{#5}%
  \selectfont}%
\fi\endgroup%
\begin{picture}(5424,3699)(589,-4048)
\put(676,-3811){\makebox(0,0)[lb]{\smash{\SetFigFont{12}{14.4}{\rmdefault}{\mddefault}{\updefault}{\color[rgb]{0,0,0}\(E_{\alpha,\beta}\)}%
}}}
\put(1351,-3511){\makebox(0,0)[lb]{\smash{\SetFigFont{12}{14.4}{\rmdefault}{\mddefault}{\updefault}{\color[rgb]{0,0,0}\(I_{\alpha,\beta}\)}%
}}}
\put(3076,-1336){\makebox(0,0)[lb]{\smash{\SetFigFont{12}{14.4}{\rmdefault}{\mddefault}{\updefault}{\color[rgb]{0,0,0}\(I_{\beta,\gamma}\)}%
}}}
\put(3076,-736){\makebox(0,0)[lb]{\smash{\SetFigFont{12}{14.4}{\rmdefault}{\mddefault}{\updefault}{\color[rgb]{0,0,0}\(E_{\beta,\gamma}\)}%
}}}
\put(4726,-3436){\makebox(0,0)[lb]{\smash{\SetFigFont{12}{14.4}{\rmdefault}{\mddefault}{\updefault}{\color[rgb]{0,0,0}\(I_{\alpha,\gamma}\)}%
}}}
\put(5326,-3811){\makebox(0,0)[lb]{\smash{\SetFigFont{12}{14.4}{\rmdefault}{\mddefault}{\updefault}{\color[rgb]{0,0,0}\(E_{\alpha,\gamma}\)}%
}}}
\end{picture}
\caption{}
\label{Triangle2}
\end{figure}

Divide $W_{\alpha',\beta',\gamma'}$ into seven regions:  the three ends
$E_{\alpha,\beta}$, 
$E_{\beta,\gamma}$ and $E_{\alpha,\gamma}$; three interpolation regions
$I_{\alpha,\beta}$, $I_{\beta,\gamma}$, and $I_{\gamma,\alpha}$
just below the ends;
 and the rest of $W_{\alpha,\beta,\gamma}$,
which we call the heart, as in \fullref{Triangle2}.  Let $J$ be a complex
structure on $W_{\alpha,\beta,\gamma}$ which
\begin{enumerate}
\item satisfies (\textbf{J$'$1}) and (\textbf{J$'$4}).
\item agrees with $J_{\alpha',\beta'}$ over $E_{\alpha,\beta}$,
  with ${J}_{\beta',\gamma'}$ over $E_{\beta,\gamma}$, and with
  $J_{\alpha',\gamma'}$ over $E_{\alpha,\gamma}$;
\item agrees with 
${J}_{t,\alpha,\beta}$ over $I_{\alpha,\beta}$, with
${J}_{t,\beta,\gamma}$ over $I_{\beta,\gamma}$, and with
$J_{t,\alpha,\gamma}$ over $I_{\alpha,\gamma}$;
\item is standard in a neighborhood of $\{\Fz_i\}\times T$ for all
  but at most one $i$;
\item achieves transversality.
\end{enumerate}
Let $J_s$ be a the complex structure obtained from $J$ by
inserting a neck of length $s$ just between each interpolation
region and the heart.  (See \fullref{Subsection:Splitting} for a precise definition of this process.)

Let $C_{\alpha}$ (respectively $C_\beta$, $C_\gamma$) be
Lagrangian cylinders which interpolate between the $\alpha$ and
$\alpha'$ (respectively $\beta$ and $\beta'$, $\gamma$ and
$\gamma'$) curves in the regions $I_{\alpha,\beta}$ and
$I_{\alpha,\gamma}$ (respectively $I_{\alpha,\beta}$ and
$I_{\beta,\gamma}$, $I_{\alpha,\gamma}$ and $I_{\beta,\gamma}$).

Counting $J_s$--holomorphic curves with boundary on the
$C_{\alpha}$ we obtain maps $f_{s,\alpha',\beta',\gamma'}$.  By the
same proof as the previous proposition, the maps
$f_{s,\alpha',\beta',\gamma'}$ are all chain homotopic to the map
$f_{\alpha',\beta',\gamma'}$.

Taking $s\to\infty$ and using the compactness result~\cite[Theorem 10.3]{Ya2}
and \fullref{Gluing:Split} we find that
$f_{\alpha',\beta',\gamma'}$ is chain 
homotopic to
$\Phi_{\alpha,\gamma}\circ
f_{\alpha,\beta,\gamma}\circ\left(\Phi_{\alpha,\beta}\otimes\Phi_{\beta,\gamma}\right).$
Here, $\Phi_{\alpha,\beta}$ (respectively $\Phi_{\beta,\gamma}$,
$\Phi_{\alpha,\gamma}$) is the chain map defined
in \fullref{Section:Isotopy}, for the isotopy between
$(\va,\vb)$ and $(\va',\vb')$ (respectively between
$(\vb,\vc)$ and $(\vb',\vc')$, between
$(\va,\vc)$ and $(\va',\vc')$).

Orientation systems, which have been implicit in the discussion,
are extended as discussed in \fullref{Section:Isotopy}.
\endproof
\subsection{Associativity of triangle maps}\label{Subsection:TriangleAssoc}
Next we show that the maps $F_{\alpha,\beta,\gamma}$ satisfy an associativity property.  Before stating it, however, we need some basic properties of Heegaard quadruple--diagrams.

\subsubsection{Heegaard quadruple--diagrams}
As the reader can presumably guess, a\break pointed Heegaard
quadruple--diagram consists of a genus $g$ surface $\Sigma$, four
$g$--tuples of pairwise disjoint homologically linearly independent
simple closed curves $\va,$ $\vb$, $\vc$ and $\vd$, and a
distinguished point
$\Fz\in\Sigma\setminus(\ba\cup\bb\cup\bc\cup\bd)$.   In this section
we will state the analogs for Heegaard quadruple--diagrams of the
basic definitions and lemmas stated earlier for Heegaard triple
diagrams.  Except as noted, the proofs are the same as for Heegaard
triple--diagrams, and hence are omitted.

Fix a pointed Heegaard quadruple--diagram
$\cH^4=(\Sigma,\va,\vb,\vc,\vd,\Fz)$.  The diagram $\cH^4$
specifies a $4$--manifold $X_{\alpha,\beta,\gamma,\delta}$ with
boundary $\bdy X_{\alpha,\beta,\gamma,\delta}=Y_{\alpha,\beta}\cup
Y_{\beta,\gamma}\cup Y_{\gamma,\delta}\cup (-Y_{\alpha,\delta})$ by
gluing the elongated handlebodies $U_\alpha\times[0,1]$,
$U_\beta\times[0,1]$, $U_\gamma\times[0,1]$ and
$U_\delta\times[0,1]$ to the product of $\Sigma$ and a square.

Let $W_{\alpha,\beta,\gamma,\delta}$ be the product of $\Sigma$
and a disk $R$ with four boundary punctures, thought of as a
topological space without complex structure for the moment.  Let
$e_1,\cdots,e_4$ denote the four boundary arcs of $R$, enumerated
clockwise, $v_{12}$, $v_{23}$, $v_{34}$ and $v_{41}$ the four
punctures, with $v_{ij}$ between $e_i$ and $e_j$; see
\fullref{Triangle3}.  In $W_{\alpha,\beta,\gamma,\delta}$ we have
$4g$ cylinders:  $\ba\times e_1$, $\bb\times e_2$, $\bc\times
e_3$ and $\bd\times e_4$.  Let $\pi_2(\vx,\vy,\vz,\vw)$ denote
the collection of homology classes of maps $(S,\bdy S)\to
(W_{\alpha,\beta,\gamma,\delta},\ba\times e_1\cup\cdots\cup\bd\times
e_4)$ asymptotic to the I--chord collection $\vx$ for
$\cH_{\alpha,\beta}$ at $v_1$, $\vy$ for $\cH_{\beta,\gamma}$
at $v_2$, and so on.

Given $\vx,\cdots,\vw$ there is an element
$\epsilon(\vx,\vy,\vz,\vw)\in H_1(X_{\alpha,\beta,\gamma,\delta})$
defined as we defined $\epsilon(\vx,\vy,\vz)$ for triangles in
\fullref{Subsection:TriangleBasics}.
\begin{Lem}
\begin{enumerate}
\item $\pi_2(\vx,\vy,\vz,\vw)$ is nonempty if and only if
  $\epsilon(\vx,\vy,\vz,\vw)=0$. 
\item If $\epsilon(\vx,\vy,\vz,\vw)=0$ then
  $\pi_2(\vx,\vy,\vz,\vw)\cong\bZ\oplus
  H_2(X_{\alpha,\beta,\gamma,\delta})$.  The map\break
  $\pi_2(\vx,\vy,\vz,\vw)\to\bZ$ is given by $n_\Fz$.  The
  identification of $\pi_2(\vx,\vy,\vz,\vw)/\bZ$ with
  $H_2(X_{\alpha,\beta,\gamma,\delta})$ is affine (but canonical up
  to translation).
\end{enumerate}
\end{Lem}
\proof
See \fullref{Lemma:TriPi} and \fullref{Lemma:TriEpsilon}.
\endproof

As with triangles, each element $A$ of $\pi_2(\vx,\vy,\vz,\vw)$
specifies a $\Spin^\bC$--structure $\Ss_\Fz(A)$ on
$X_{\alpha,\beta,\gamma,\delta}$; the construction is the same as in
\fullref{Subsubsection:TriangleSpinc}.  Also as before,
there are obvious concatenation maps
$\pi_2(\vx',\vx)\times\pi_2(\vx,\vy,\vz,\vw)\to\pi_2(\vx',\vy,\vz,\vw)$,
and similarly for $\vy$, $\vz$ and $\vw$.  Again we say
$A\in\pi_2(\vx,\vy,\vz,\vw)$ is $\Spin^\bC$--equivalent to
$A+B_{\alpha,\beta}+B_{\beta,\gamma}+B_{\gamma,\delta}+B_{\alpha,\delta}$ for $B_{\alpha,\beta}\in\pi_2(\vx,\vx),\cdots,B_{\alpha,\delta}\in\pi_2(\vw,\vw)$,
and let $S_{\alpha,\beta,\gamma,\delta}$ denote the collection of
$\Spin^\bC$--equivalence classes.  Again we have:
\begin{Lem}The map $\pi_2(\vx,\vy,\vz,\vw)\to
  \Spin^\bC(X_{\alpha,\beta,\gamma,\delta})$ descends to an injective
  map $S_{\alpha,\beta,\gamma,\delta}\mapsinto
  \Spin^\bC(X_{\alpha,\beta,\gamma,\delta})$ whose image consists of
  all those $\Spin^\bC$--structures whose restrictions to $\bdy
  X_{\alpha,\beta,\gamma,\delta}$ are realized to intersection
  points. 
\end{Lem}
\proof
See \fullref{Lemma:SpinEquiv}.
\endproof

Now, however, there is somewhat more structure.  The manifold
$X_{\alpha,\beta,\gamma,\delta}$ decomposes as
$X_{\alpha,\beta,\gamma}\cup_{Y_{\alpha,\gamma}}X_{\alpha,\gamma,\delta}$
and as
$X_{\alpha,\beta,\delta}\cup_{Y_{\beta,\delta}}X_{\beta,\gamma,\delta}$.
Let $\delta_{\alpha,\gamma}$ (respectively $\delta_{\beta,\gamma}$) be the coboundary map for the Mayer--Vietoris sequence for the former (respectively the latter) decomposition.
Working for the moment with the former decomposition, we have
restriction maps
$\Spin^\bC(X_{\alpha,\beta,\gamma,\delta})\stackrel{r}{\to}\Spin^\bC(X_{\alpha,\beta,\gamma})\times
\Spin^\bC(X_{\alpha,\gamma,\delta})$.  These maps commute with the
$H^2$--actions, and so by the Mayer--Vietoris theorem, the fibers of
$r$ are the orbits of the action of $\delta_{\alpha,\gamma}
H^1(X_{\alpha,\gamma})\subset H^2(X_{\alpha,\beta,\gamma,\delta})$ on
$\Spin^\bC(X_{\alpha,\beta,\gamma,\delta})$.
Corresponding statements hold for the decomposition
$X_{\alpha,\beta,\gamma,\delta}\cong
X_{\alpha,\beta,\delta}\cup_{Y_{\beta,\delta}}X_{\beta,\gamma,\delta}.$

Rather than fixing a single $\Spin^\bC$--structure over $X_{\alpha,\beta,\gamma,\delta}$ we shall fix a $\delta_{\alpha,\gamma} H^1(Y_{\alpha,\gamma})+\delta_{\beta,\delta} H^1(Y_{\beta,\delta})$--orbit of $\Spin^\bC$--structures.  The reason is easier to see in terms of domains.  Concatenation (addition) gives a well--defined map 
$$\pi_2^{\alpha,\beta,\gamma}(\vx,\vy,\vec{a})\times \pi_2^{\alpha,\gamma,\delta}(\vec{a},\vz,\vw)\to \pi_2^{\alpha,\beta,\gamma,\delta}(\vx,\vy,\vz,\vw).$$
This map, however, does \emph{not} descend to a map
$S_{\alpha,\beta,\gamma}\times S_{\alpha,\gamma,\delta}\to
S_{\alpha,\beta,\gamma,\delta}$:  for
$A_{\alpha,\beta,\gamma}\in\pi_2(\vx,\vy,\vec{a})$,
$A_{\alpha,\gamma,\delta}\in\pi_2(\vec{a},\vz,\vw)$ and
$B_{\alpha,\gamma}\in\pi_2(\vec{a},\vec{a})$ the domains
  $A_{\alpha,\beta,\gamma}$ and
  $A_{\alpha,\beta,\gamma}+B_{\alpha,\gamma}$ define the same
  element of $S_{\alpha,\beta,\gamma}$, but
  $A_{\alpha,\beta,\gamma}+A_{\alpha,\gamma,\delta}$ and
  $A_{\alpha,\beta,\gamma}+B_{\alpha,\gamma}+A_{\alpha,\gamma,\delta}$ may not define the same element of $S_{\alpha,\beta,\gamma,\delta}$.  So, the composition $F_{\alpha,\beta,\gamma}\left(F_{\alpha,\gamma,\delta}(\cdot,\cdot),\cdot\right)$ involves domains belonging to an entire $\pi_2^{\alpha,\gamma}(\vec{a},\vec{a})$--orbit of elements of $S_{\alpha,\beta,\gamma,\delta}$.
Since $\delta_{\alpha,\gamma} H^1(Y_{\alpha,\gamma})$ and $\delta_{\beta,\delta}
H^1(Y_{\beta,\delta})$ may not coincide, the best we can expect to
prove is associativity for certain sums of triangle maps.  This is
what we will (eventually) prove.

The next issue to address is admissibility.
\begin{Def}The pointed Heegaard quadruple--diagram $\cH^4$ is
  \emph{weakly admissible} if given
  $B_{\alpha,\beta}\in\hat{\pi}_2(\vx,\vx)$,
  $B_{\beta,\gamma}\in\hat{\pi}_2(\vy,\vy)$,
  $B_{\gamma,\delta}\in\hat{\pi}_2(\vz,\vz)$,
  $B_{\alpha,\delta}\in\hat{\pi}_2(\vw,\vw)$ with
  $B_{\alpha,\beta}+B_{\beta,\gamma}+B_{\gamma,\delta}+B_{\alpha,\delta}\neq 0$, then $B_{\alpha,\beta}+B_{\beta,\gamma}+B_{\gamma,\delta}+B_{\alpha,\delta}$ has both positive and negative coefficients.
\end{Def}
\begin{Def}A pointed Heegaard quadruple--diagram $\cH^4$ is
  \emph{strongly admissible} for a $\delta_{\alpha,\gamma}
  H^1(Y_{\alpha,\gamma})+\delta_{\beta,\delta} H^1(Y_{\beta,\delta})$--orbit of
  $\Spin^\bC$--structures $\SS$ if for any
  $\Ss_{\alpha,\beta,\gamma,\delta}\in\SS$ and any six
  domains $B_{\xi,\eta}\in\hat{\pi}_2^{\xi,\eta},$
  $\{\xi,\eta\}\subset\{\alpha,\beta,\gamma,\delta\}$
  ($\xi\neq\eta$) such that
\begin{itemize}
\item $ \sum_{\{\xi,\eta\}}B_{\xi,\eta}\neq 0$ and
\item $\sum_{\{\xi,\eta\}}\left\langle
  c_1\left(\Ss_{\alpha,\beta,\gamma}|_{Y_{\xi,\eta}}\right),B_{\xi,\eta}\right\rangle=2n\geq0$
\end{itemize}
then some coefficient of $\sum_{\{\xi,\eta\}}B_{\xi,\eta}$ is greater than $n$.
\end{Def}

\begin{Lem}
\begin{enumerate}
\item Given any pointed Heegaard quadruple--diagram there is an
  isotopic weakly admissible pointed Heegaard quadruple--diagram. 
\item Suppose that the pointed Heegaard quadruple--diagram
  $(\Sigma,\va,\vb,\vc,\vd,\Fz)$ satisfies the conditions $\delta_{\beta,\delta}
  H^1(Y_{\beta,\delta})|_{Y_{\alpha,\gamma}}=0$ and $\delta_{\alpha,\gamma}
  H^1(Y_{\alpha,\gamma})|_{Y_{\beta,\delta}}=0$.  Fix a
  $\delta_{\alpha,\gamma} H^1(Y_{\alpha,\gamma})+\delta_{\beta,\delta}
  H^1(Y_{\beta,\delta})$--orbit of $\Spin^\bC$--structures $\SS$.  Then
  there is an isotopic pointed Heegaard quadruple--diagram which is strongly
  admissible for $\SS$.
\end{enumerate}
\end{Lem}
For the argument, see~\cite[Section 8.4.2]{OS1}.

\begin{Lem} Suppose $\cH^4$ is weakly admissible.  Fix
  $j,k\in\bZ$, intersection points $\vx,\vy,\vz,\vw$ and a
  $\delta H^1(Y_{\beta,\delta})+\delta
  H^1(Y_{\alpha,\gamma})$--orbit of $\Spin^\bC$--structures $\SS$
  on $X_{\alpha,\beta,\gamma,\delta}$.  Then there are only finitely
  many positive $A\in\pi_2(\vx,\vy,\vz,\vw)$ such that
\begin{itemize}
\item $\ind(A)=j$
\item $n_\Fz(A)=k$
\item $\Ss_\Fz(A)\in\SS$
\end{itemize}
\end{Lem}
\proof
See \fullref{Lemma:WeakTriangleAdmisImplies}.
\endproof

\begin{Lem}\label{Lemma:SquareStrongAdmis}  Suppose $\cH^4$ is
  strongly admissible for a $\delta_{\alpha,\gamma} H^1(Y_{\alpha,\gamma})+$\break $\delta_{\beta,\delta}
  H^1(Y_{\beta,\delta})$--orbit of $\Spin^\bC$--structures $\SS$ on $X_{\alpha,\beta,\gamma,\delta}$.
  Fix $j\in\bZ$ and intersection points $\vx,\vy,\vz,\vw$.  Then
  there are only finitely many positive $A\in\pi_2(\vx,\vy,\vz,\vw)$ such
  that $\Ss_\Fz(A)\in\SS$ and $\ind(A)=j$.
\end{Lem}
\proof
See \fullref{Lemma:StrongTriangleAdmisImplies}.
\endproof

\subsubsection{Moduli spaces of squares}
\begin{figure}
\centering
\begin{picture}(0,0)%
\includegraphics[scale=.8]{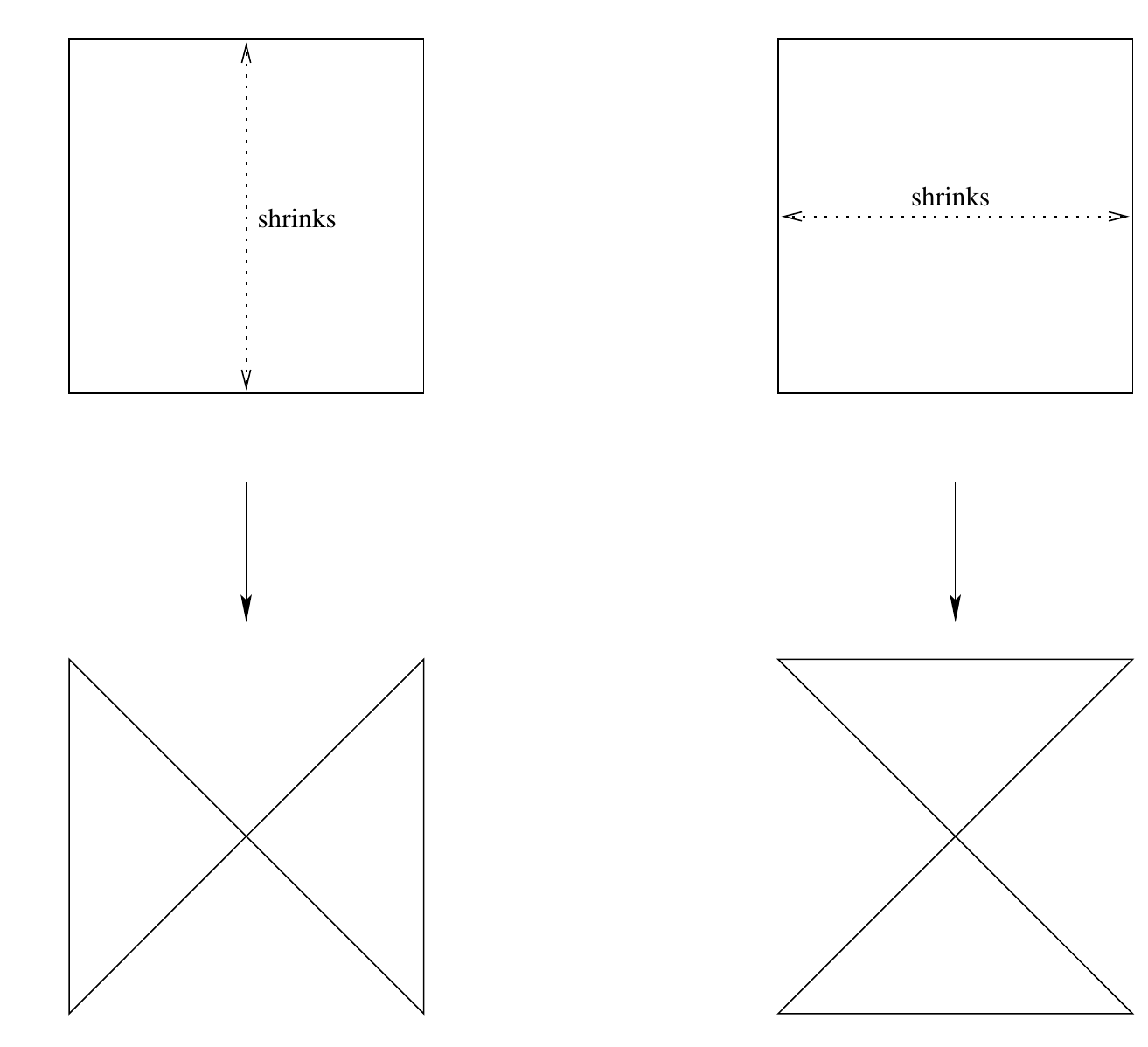}%
\end{picture}%
\setlength{\unitlength}{3157sp}%
\begingroup\makeatletter\ifx\SetFigFont\undefined%
\gdef\SetFigFont#1#2#3#4#5{%
  \reset@font\fontsize{#1}{#2pt}%
  \fontfamily{#3}\fontseries{#4}\fontshape{#5}%
  \selectfont}%
\fi\endgroup%
\begin{picture}(6294,5738)(736,-5253)
\put(5859,343){\makebox(0,0)[lb]{\smash{{\SetFigFont{9}{10.8}{\rmdefault}{\mddefault}{\updefault}{\color[rgb]{0,0,0}\(e_1\)}%
}}}}
\put(7015,-752){\makebox(0,0)[lb]{\smash{{\SetFigFont{9}{10.8}{\rmdefault}{\mddefault}{\updefault}{\color[rgb]{0,0,0}\(e_2\)}%
}}}}
\put(5799,-1847){\makebox(0,0)[lb]{\smash{{\SetFigFont{9}{10.8}{\rmdefault}{\mddefault}{\updefault}{\color[rgb]{0,0,0}\(e_3\)}%
}}}}
\put(4643,-752){\makebox(0,0)[lb]{\smash{{\SetFigFont{9}{10.8}{\rmdefault}{\mddefault}{\updefault}{\color[rgb]{0,0,0}\(e_4\)}%
}}}}
\put(3001,343){\makebox(0,0)[lb]{\smash{{\SetFigFont{9}{10.8}{\rmdefault}{\mddefault}{\updefault}{\color[rgb]{0,0,0}\(v_{12}\)}%
}}}}
\put(3001,-1847){\makebox(0,0)[lb]{\smash{{\SetFigFont{9}{10.8}{\rmdefault}{\mddefault}{\updefault}{\color[rgb]{0,0,0}\(v_{23}\)}%
}}}}
\put(933,-1847){\makebox(0,0)[lb]{\smash{{\SetFigFont{9}{10.8}{\rmdefault}{\mddefault}{\updefault}{\color[rgb]{0,0,0}\(v_{34}\)}%
}}}}
\put(2211,-2576){\makebox(0,0)[lb]{\smash{{\SetFigFont{9}{10.8}{\rmdefault}{\mddefault}{\updefault}{\color[rgb]{0,0,0}\(s\to-\infty\)}%
}}}}
\put(6041,-2576){\makebox(0,0)[lb]{\smash{{\SetFigFont{9}{10.8}{\rmdefault}{\mddefault}{\updefault}{\color[rgb]{0,0,0}\(s\to\infty\)}%
}}}}
\put(995,343){\makebox(0,0)[lb]{\smash{{\SetFigFont{9}{10.8}{\rmdefault}{\mddefault}{\updefault}{\color[rgb]{0,0,0}\(v_{41}\)}%
}}}}
\put(2576,343){\makebox(0,0)[lb]{\smash{{\SetFigFont{9}{10.8}{\rmdefault}{\mddefault}{\updefault}{\color[rgb]{0,0,0}\(\alpha\)}%
}}}}
\put(3123,-1299){\makebox(0,0)[lb]{\smash{{\SetFigFont{9}{10.8}{\rmdefault}{\mddefault}{\updefault}{\color[rgb]{0,0,0}\(\beta\)}%
}}}}
\put(2393,-1785){\makebox(0,0)[lb]{\smash{{\SetFigFont{9}{10.8}{\rmdefault}{\mddefault}{\updefault}{\color[rgb]{0,0,0}\(\gamma\)}%
}}}}
\put(872,-1299){\makebox(0,0)[lb]{\smash{{\SetFigFont{9}{10.8}{\rmdefault}{\mddefault}{\updefault}{\color[rgb]{0,0,0}\(\delta\)}%
}}}}
\put(2211,-3671){\makebox(0,0)[lb]{\smash{{\SetFigFont{9}{10.8}{\rmdefault}{\mddefault}{\updefault}{\color[rgb]{0,0,0}\(\alpha\)}%
}}}}
\put(3123,-3914){\makebox(0,0)[lb]{\smash{{\SetFigFont{9}{10.8}{\rmdefault}{\mddefault}{\updefault}{\color[rgb]{0,0,0}\(\beta\)}%
}}}}
\put(933,-4157){\makebox(0,0)[lb]{\smash{{\SetFigFont{9}{10.8}{\rmdefault}{\mddefault}{\updefault}{\color[rgb]{0,0,0}\(\delta\)}%
}}}}
\put(5434,-3854){\makebox(0,0)[lb]{\smash{{\SetFigFont{9}{10.8}{\rmdefault}{\mddefault}{\updefault}{\color[rgb]{0,0,0}\(\delta\)}%
}}}}
\put(5799,-3063){\makebox(0,0)[lb]{\smash{{\SetFigFont{9}{10.8}{\rmdefault}{\mddefault}{\updefault}{\color[rgb]{0,0,0}\(\alpha\)}%
}}}}
\put(5859,-5191){\makebox(0,0)[lb]{\smash{{\SetFigFont{9}{10.8}{\rmdefault}{\mddefault}{\updefault}{\color[rgb]{0,0,0}\(\gamma\)}%
}}}}
\put(6406,-3854){\makebox(0,0)[lb]{\smash{{\SetFigFont{9}{10.8}{\rmdefault}{\mddefault}{\updefault}{\color[rgb]{0,0,0}\(\beta\)}%
}}}}
\put(2088,-4583){\makebox(0,0)[lb]{\smash{{\SetFigFont{9}{10.8}{\rmdefault}{\mddefault}{\updefault}{\color[rgb]{0,0,0}\(\gamma\)}%
}}}}
\put(1967,343){\makebox(0,0)[lb]{\smash{{\SetFigFont{9}{10.8}{\rmdefault}{\mddefault}{\updefault}{\color[rgb]{0,0,0}\(e_1\)}%
}}}}
\put(3123,-752){\makebox(0,0)[lb]{\smash{{\SetFigFont{9}{10.8}{\rmdefault}{\mddefault}{\updefault}{\color[rgb]{0,0,0}\(e_2\)}%
}}}}
\put(1906,-1847){\makebox(0,0)[lb]{\smash{{\SetFigFont{9}{10.8}{\rmdefault}{\mddefault}{\updefault}{\color[rgb]{0,0,0}\(e_3\)}%
}}}}
\put(751,-752){\makebox(0,0)[lb]{\smash{{\SetFigFont{9}{10.8}{\rmdefault}{\mddefault}{\updefault}{\color[rgb]{0,0,0}\(e_4\)}%
}}}}
\end{picture}%
\caption{Degenerations of a rectangle}
\label{Triangle3}
\end{figure}
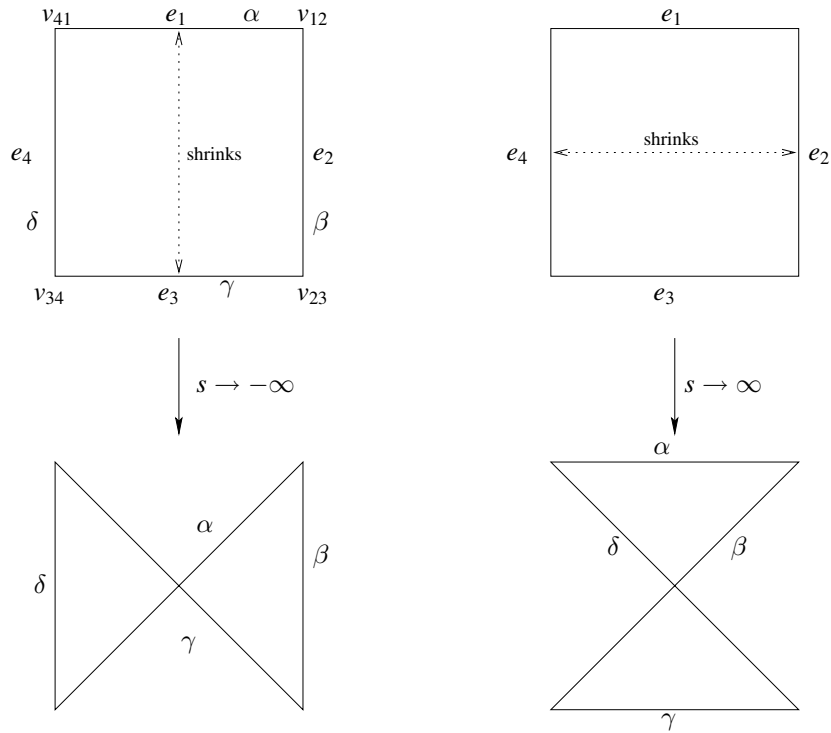

Fix a pointed Heegaard quadruple--diagram $\cH^4=(\Sigma,\va,\vb,\vc,\vd,\Fz)$.  Let $R$ denote the unit disk with four punctures on its boundary.  Let $e_1,\cdots,e_4$ denote the
four boundary components of $R$, enumerated clockwise, and
$v_{12},v_{23},v_{34}$ and $v_{4,1}$ the vertices of $R$, enumerated clockwise, with
$v_{12}$ between $e_1$ and $e_2$.  See \fullref{Triangle3}.  Let
$W_{\alpha,\beta,\gamma,\delta}=\Sigma\times R$.

The moduli space of conformal structures on a rectangle is
parameterized by $\bR$.  Let $j_a$, $a\in\bR$, sweep out this
space.  Do this in such a way that as $s\to-\infty$ an arc in $R$
connecting $e_1$ to $e_3$ collapses, while as
$s\to\infty$ an arc in $R$ connecting $e_2$ to $e_4$
collapses; see \fullref{Triangle3}.

Fix a point $\Fz_i$ in each component of
$\Sigma\setminus(\ba\cup\bb\cup\bc\cup\boldsymbol{\delta})$.  Choose
complex structures $J_{\xi,\eta}$,
$\{\xi,\eta\}\subset\{\alpha,\beta,\gamma,\delta\}$ which satisfy
(\textbf{J1})--(\textbf{J5}) and achieve
transversality.

Fix a path $J_a$ of complex structures on
$W_{\alpha,\beta,\gamma,\delta}$ such that
\begin{enumerate}
\item For every $a\in\bR$, $J_a$ is tamed by the split symplectic
  form on $\Sigma\times R$.
\item For every $a\in\bR$, projection $\pi_R$ onto $R$ is
  $(j_a,J_a)$--holomorphic.
\item In a neighborhood $U_{\{\Fz_i\}}$ of $\{\Fz_i\}\times R$,
 $J = j_\Sigma\times j_T$ is split.
\item 
Near $\Sigma\times\{v_{12}\}$, $J$ agrees with $J_{\alpha,\beta}$.
Near $\Sigma\times\{v_{23}\}$, $J$ agrees with $J_{\beta,\gamma}$.
Near $\Sigma\times\{v_{34}\}$, $J$ agrees with $J_{\gamma,\delta}$.
Near $\Sigma\times\{v_{41}\}$, $J$ agrees with $J_{\alpha,\delta}$.
\item As $a\to-\infty$, $J_a$ degenerates to complex structures
  which satisfy (\textbf{J$'$1})--(\textbf{J$'$4}) and achieve
  transversality for $W_{\alpha,\beta,\gamma}$ and
  $W_{\alpha,\gamma,\delta}$.
\item As $a\to\infty$, $J_a$ degenerates to complex structures
  which satisfy (\textbf{J$'$1})--(\textbf{J$'$4}) and achieve
  transversality for $W_{\delta,\alpha,\beta}$ and
  $W_{\gamma,\delta,\beta}$.
\item $J_a$ achieves transversality (as a path of almost complex
 structures) for holomorphic curves with index $\leq1$.
\end{enumerate}
Checking that such a $J_a$ exists is similar to the proof of transversality
in \fullref{Section:Transversality}.

Given I--chord collections $\vx,\vy,\vz,\vec{w}$ for
$Y_{\alpha,\beta}$, $Y_{\beta,\gamma}$, $Y_{\gamma,\delta}$ and
$Y_{\alpha,\delta}$ let $\pi_2(\vx,\vy,\vz,\vw)$ denote
homology classes of maps to $W_{\alpha,\beta,\gamma,\delta}$
connecting $\vx$, $\vy$, $\vz$, and $\vw$.  Given
$A\in\pi_2(\vx,\vy,\vz,\vw)$, let $\cM^A$ denote the union over
$a\in\bR$ of all embedded $J_a$--holomorphic curves $u\co (S,\bdy
S)\to (W_{\alpha,\beta,\gamma,\delta},\ba\times e_1\cup\bb\times
e_2\cup \bc\times e_3\cup \bd\times e_4)$ (without closed
components, and with one component of $\bdy S$ mapped to each
Lagrangian cylinder) in the
homology class $A$.  Let $\ind(A)$ denote the index of this
$\dbar$ problem in the homology class $A$, so $dim\cM^A=\ind(A)+1$, if $\cM^A$ is nonempty.  (The $+1$ appears because we
  consider a $1$--parameter family of almost complex structures.)

\subsubsection{Orienting squares}
\label{Subsubsection:OrientingSquares}
The moduli spaces $\cM^{A_{\alpha,\beta,\gamma,\delta}}$ are
orientable for the same reason all the other moduli spaces considered
so far have been.  Again we want to choose orientations for the
$\cM^{A_{\alpha,\beta,\gamma,\delta}}$,
$\cM^{A_{\kappa,\eta,\xi}}$
and $\cM^{B_{\eta,\xi}}$,
$\{\kappa,\eta,\xi\}\subset\{\alpha,\beta,\gamma,\delta\}$) consistent
with various gluings -- all possible gluings of a $2$--gon to a
rectangle or triangle, and the two gluings
$\cM^{A_{\alpha,\beta,\gamma}}\times\cM^{A_{\alpha,\gamma,\delta}}\times[R,\infty)\mapsinto\cM^{A_{\alpha,\beta,\gamma}+A_{\alpha,\gamma,\delta}}$
and
$\cM^{A_{\alpha,\beta,\delta}}\times\cM^{A_{\beta,\gamma,\delta}}\times[R,\infty)\mapsinto\cM^{A_{\alpha,\beta,\delta}+A_{\beta,\gamma,\delta}}.$
Again, it follows from standard arguments that there is some such
coherent orientation system.  And, again it is useful to have
something slightly stronger:

\begin{Lem}{\rm(Compare~\cite[Proposition 8.15]{OS1})}\label{Lemma:QuadOr}\qua
Suppose $\cH^4$ is a Heegaard quad\-ruple--diagram such that 
the image of $H_2(Y_{\beta,\delta})$ under the map $H_2(X_{\alpha,\beta,\gamma,\delta})\to H_2(X_{\alpha,\beta,\gamma,\delta},\bdy X_{\alpha,\beta,\gamma,\delta})$ is zero.
Fix a $\delta H^1(Y_{\alpha,\gamma})$--orbit of $\Spin^\bC$--structures $\SS$
  on $X_{\alpha,\beta,\gamma,\delta}$ and orientation systems
  $\frako_{\alpha,\beta,\gamma}$ and $\frako_{\alpha,\gamma,\delta}$ for
  $(\Sigma,\va,\vb,\vc,\Fz,\SS|_{X_{\alpha,\beta,\gamma}})$ and
  $(\Sigma,\va,\vc,\vd,\Fz,\SS|_{X_{\alpha,\gamma,\delta}})$
  inducing the same orientation over $Y_{\alpha,\gamma}$.  Then
  there is a coherent orientation system for $\cH^4$ extending
  $\frako_{\alpha,\beta,\gamma}$ and $\frako_{\alpha,\gamma,\delta}$.
\end{Lem}
\proof
Our proof is the same as the one given in~\cite{OS1}.

The orientation systems $\frako_{\alpha,\beta,\gamma}$ and $\frako_{\alpha,\gamma,\delta}$ determine orientations $\frako_{\alpha,\beta,\gamma,\delta}(A_{\alpha,\beta,\gamma,\delta})$ for all $A_{\alpha,\beta,\gamma,\delta}$ with $\Ss_\Fz(A_{\alpha,\beta,\gamma,\delta})\in\SS$, which is well--defined since $\frako_{\alpha,\beta,\gamma}$ and $\frako_{\alpha,\gamma,\delta}$ induce the same orientation system on $Y_{\alpha,\gamma}$.  The orientation systems $\frako_{\alpha,\beta,\gamma}$ and $\frako_{\alpha,\gamma,\delta}$ also determine orientation systems $\frako_{\alpha,\beta}$, $\frako_{\beta,\gamma}$, $\frako_{\gamma,\delta}$, $\frako_{\alpha,\delta}$ and $\frako_{\alpha,\gamma}$ over $Y_{\alpha,\beta}$, $Y_{\beta,\gamma}$, $Y_{\gamma,\delta}$, $Y_{\alpha,\delta}$ and $Y_{\alpha,\gamma}$.

Choose intersection points $\vx\in\SS|_{Y_{\alpha,\beta}}$,
$\vy\in\SS|_{\Ybc}$, $\vz\in\SS|_{Y_{\gamma,\delta}}$,
$\vw\in\SS|_{Y_{\alpha,\delta}}$,
$\vec{u}\in\SS|_{Y_{\alpha,\gamma}}$ and
$\vec{v}\in\SS|_{Y_{\beta,\delta}}$.  Choose any
$A_{\alpha,\beta,\gamma,0}\in\pi_2(\vx,\vy,\vec{u})$ and
$A_{\beta,\gamma,\delta,0}\in\pi_2(\vy,\vz,\vec{v})$.  Choose an
arbitrary orientation
$\frako_{\alpha,\beta,\delta}(A_{\alpha,\beta,\delta,0})$ over
$A_{\alpha,\beta,\delta,0}$.  Then
$\frako_{\alpha,\beta,\delta}(A_{\alpha,\beta,\delta,0})$ and
$\frako_{\alpha,\beta,\gamma,\delta}(A_{\alpha,\beta,\delta,0}+A_{\beta,\gamma,\delta,0})$
determine an orientation
$\frako_{\beta,\gamma,\delta}(A_{\beta,\gamma,\delta,0}).$

Now, we choose $\frako_{\beta,\delta}$ as follows.  For $B_{\beta,\delta}\in\pi_2(\vec{v},\vec{v})$, from the assumption that $H_2(Y_{\beta,\delta})$ is trivial inside $H_2(X_{\alpha,\beta,\gamma,\delta},\bdy X_{\alpha,\beta,\gamma,\delta})$, we can choose $B_{\alpha,\beta}$, $B_{\beta,\gamma}$, $B_{\gamma,\delta}$ and $B_{\alpha,\delta}$ so that
\begin{equation}\label{equation:AB}
A_{\alpha,\beta,\delta,0}+A_{\beta,\gamma,\delta,0}+B_{\beta,\delta}=A_{\alpha,\beta,\delta,0}+A_{\beta,\gamma,\delta,0}+B_{\alpha,\beta}+B_{\beta,\gamma}+B_{\gamma,\delta}+B_{\alpha,\delta}
\end{equation}
Choose $\frako_{\beta,\delta}(B_{\beta,\delta})$ so that the
orientations induced by the two decompositions in Equation~\eqref{equation:AB} agree.  We have constructed
$\frako_{\alpha,\beta}$, $\frako_{\beta,\gamma}$,
$\frako_{\gamma,\delta}$ and $\frako_{\alpha,\delta}$ so that this is
independent of the choice of
$B_{\alpha,\beta},\cdots,B_{\alpha,\delta}$.  Choose arbitrary
orientations over a complete set of paths for $Y_{\beta,\delta}$ to
finish defining $\frako_{\beta,\delta}$.

Finally, $\frako_{\alpha,\beta,\delta}(A_{\alpha,\beta,\delta,0})$,
$\frako_{\alpha,\beta}$, $\frako_{\beta,\delta}$ and
$\frako_{\alpha,\delta}$ together specify $\frako_{\alpha,\beta,\delta}$
completely, and a similar remark applies to
$\frako_{\beta,\gamma,\delta}$.  We have, thus, finished constructing
all the orientations.  It is clear that they are coherent.
\endproof

\subsubsection{Proof of associativity}
\begin{Prop}\label{AssocProp} Let $\cH^4=(\Sigma,\va,\vb,\vc,\vd,\Fz)$ be
  a pointed Heegaard quadruple--diagram.  Fix a $\delta_{\alpha,\gamma}
  H^1(Y_{\alpha,\gamma})+\delta_{\beta,\delta} H^1(Y_{\beta,\delta})$--orbit of
  $\Spin^\bC$--structures $\SS$ on\break
  $X_{\alpha,\beta,\gamma,\delta}$.  Assume that $\cH^4$ is
  strongly admissible for $\SS$.  Fix $*\in\{\infty,+,-\}$.
Then for any $\xi_{\alpha,\beta}\in HF^*(Y_{\alpha,\beta})$ and
  $\theta_{\beta,\gamma}\in HF^{\leq0}(Y_{\beta,\gamma})$,
  $\theta_{\gamma,\delta}\in HF^{\leq0}(Y_{\gamma,\delta})$ we have
{\setlength\arraycolsep{0pt}
\begin{eqnarray}
\sum_{\Ss\in\SS}
 &F^*_{\alpha,\gamma,\delta}&\left(
 F^*_{\alpha,\beta,\gamma}(\xi_{\alpha,\beta}\otimes\theta_{\beta,\gamma};\Ss|_{X_{\alpha,\beta,\gamma}})\otimes\theta_{\gamma,\delta};\Ss|_{X_{\alpha,\gamma,\delta}}\right)\nonumber\\
 &=& \sum_{\Ss\in\SS}
F^*_{\alpha,\beta,\delta}\left(\xi_{\alpha,\beta}\otimes
 F^{\leq0}_{\beta,\gamma,\delta}(\theta_{\beta,\gamma}\otimes\theta_{\gamma,\delta};\Ss|_{X_{\beta,\gamma,\delta}});\Ss|_{X_{\alpha,\beta,\delta}}\right).\nonumber
\end{eqnarray}}
An exactly similar statement holds with both $HF^*$ and
$HF^{\leq}$ replaced by $\widehat{HF}$.  For the case of $\widehat{HF}$,
weak admissibility of the Heegaard quadruple--diagram is sufficient.
\end{Prop}
(In the proposition, we assume the maps are computed with respect to a
coherent choice of orientations, as described in
\fullref{Subsubsection:OrientingSquares}.  Also note that
every $\Spin^\bC$--structure in $\SS$ restricts to the same
$\Spin^\bC$--structure on $Y_{\alpha,\beta}$ (respectively
$Y_{\beta,\gamma}$, $Y_{\gamma,\delta}$, $Y_{\alpha,\delta}$).
When we speak of the Floer homology groups of $Y_{\alpha,\beta}$
(respectively $Y_{\beta,\gamma}$, $Y_{\gamma,\delta}$,
$Y_{\alpha,\delta}$) we mean the Floer homology groups calculated
with respect to this $\Spin^\bC$--structure.)

\proof
To make notation clearer, we will give the proof for $*=\infty$.
The other cases are completely analogous.
Define $h\co CF^\infty(Y_{\alpha,\beta})\otimes
CF^{\leq0}(Y_{\beta,\gamma})\otimes CF^{\leq0}(Y_{\gamma,\delta})\to
CF^\infty(Y_{\alpha,\delta})$ by 
$$
h\left([\vx,i]\otimes[\vy,j]\otimes[\vz,k]\right)=\sum_{\vw}\sum_{\substack{
    A\in\pi_2(\vx,vy,\vz,\vw)\\\Ss_{\Fz}(A)\in\SS\\\ind(A)=-1}}\left(\#\cM^A\right)[\vw,i+j+k-n_{\Fz}(A)].
$$
The sum defining $h$ makes sense by compactness and
\fullref{Lemma:SquareStrongAdmis}.  Note that the map $h$ has
degree $-1$.

Counting the ends of the space
$$
\bigcup_{\substack{A\in\SS\\\ind(A)=0}}\cM^A
$$
we find that
{\setlength\arraycolsep{0pt}
\begin{eqnarray}
0=\sum_{\Ss\in\SS}
 &f^\infty_{\alpha,\gamma,\delta}&\left(
 f^\infty_{\alpha,\beta,\gamma}([\vx,i]\otimes[\vy,j];\Ss|_{X_{\alpha,\beta,\gamma}})\otimes[\vz,k];\Ss|_{X_{\alpha,\gamma,\delta}}\right)\nonumber\\
 &-& \sum_{\Ss\in\SS}
f^*_{\alpha,\beta,\delta}\left([\vx,i]\otimes
 f^{\leq0}_{\beta,\gamma,\delta}([\vy,j]\otimes[\vz,k];\Ss|_{X_{\beta,\gamma,\delta}});\Ss|_{X_{\alpha,\beta,\delta}}\right)\nonumber\\
&+&h\circ\bdy([\vx,i]\otimes[\vy,j]\otimes[\vz,k])+\bdy\circ h([\vx,i]\otimes[\vy,j]\otimes[\vz,k]).\nonumber
\end{eqnarray}}
(The first two terms correspond to contributions from the
degenerations shown in \fullref{Triangle3}.  The last two terms
correspond to the usual level splittings at the (four) cylindrical ends.)

The result is immediate.
\endproof

\section{Handleslides}
\label{Section:Handleslides}
In this section we will prove that the Floer homologies defined in \fullref{Section:ChainComplexes}
are unchanged by handleslides among the $\beta$-- (and,
symmetrically, the $\alpha$--) circles.  This follows from the
associativity of the triangle maps defined in
\fullref{Section:Triangles} and the calculation of a few
specific moduli spaces.

Our proof is essentially the same as the proof given by Ozsv\'ath and Szab\'o in~\cite[Section 9]{OS1}.
Indeed, their calculations are done in our language, and most of the
rest of the proof follows formally from the results of
\fullref{Section:Triangles}.  Our proof differs from theirs in the
following ways.  Firstly, we try to calculate the minimum amount
necessary in order to prove handleslide invariance.  We hope this will
provide a slightly better understanding of the formal properties
underlying handleslide invariance, and in turn lead to
generalizations.  Secondly, we give a more geometrical proof of
\fullref{Prop:HSB} (roughly~\cite[Proposition 9.8]{OS1}),
using $1$--gons.  Thirdly, some details
about the $H_1(Y)/Tors$--actions look slightly different in our
language. 

\begin{figure}
\centering
\begin{picture}(0,0)%
\includegraphics[scale=.9]{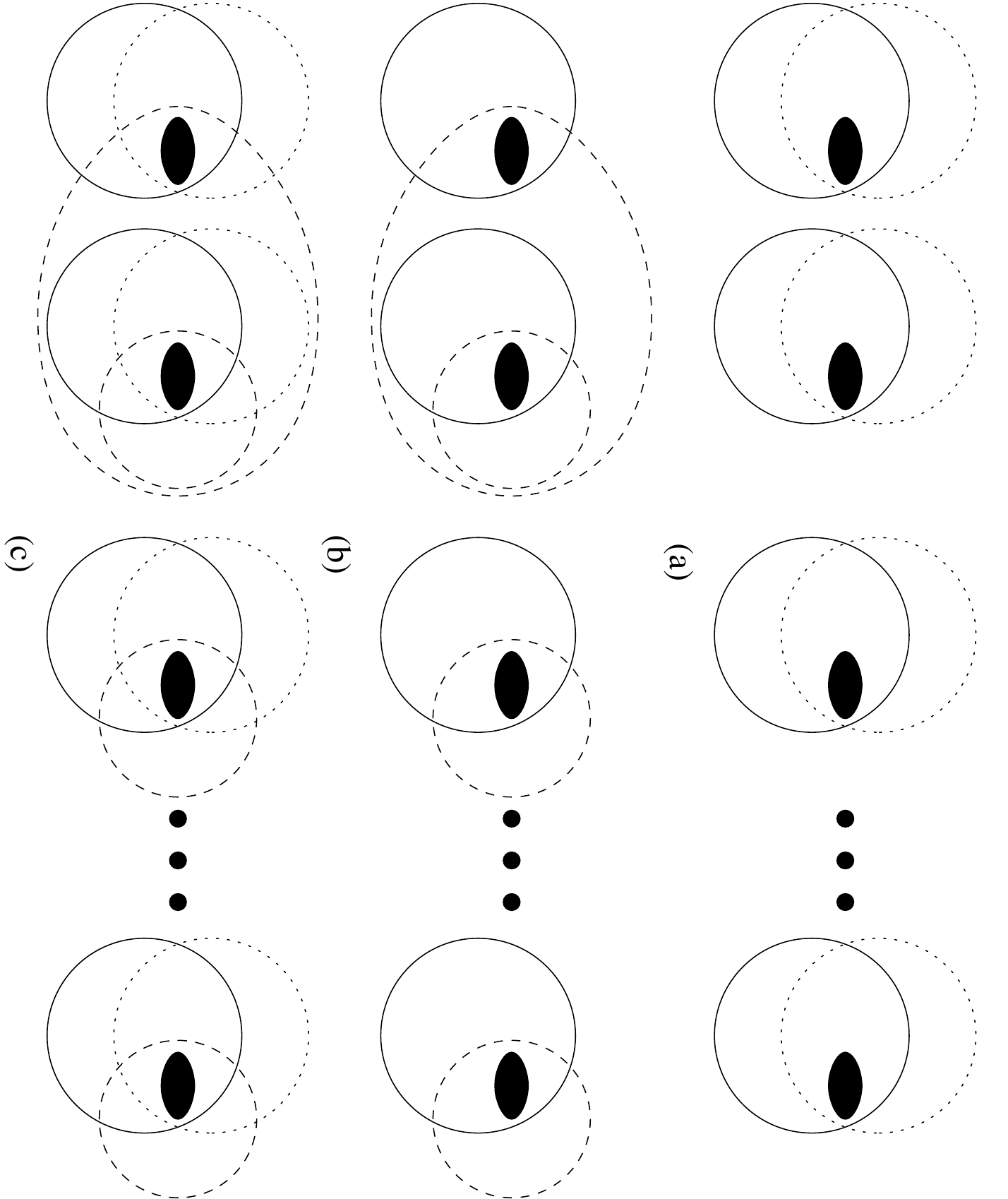}%
\end{picture}%
\setlength{\unitlength}{3552sp}%
\begingroup\makeatletter%
\gdef\SetFigFont#1#2#3#4#5{%
  \reset@font\fontsize{#1}{#2pt}%
  \fontfamily{#3}\fontseries{#4}\fontshape{#5}%
  \selectfont}%
\endgroup%
\begin{picture}(6847,8134)(133,-7831)
\put(4122,-3430){\rotatebox{270.0}{\makebox(0,0)[lb]{\smash{{\SetFigFont{9}{10.8}{\rmdefault}{\mddefault}{\updefault}{\color[rgb]{0,0,0}\(F^H\)}%
}}}}}
\put(4916, 30){\rotatebox{270.0}{\makebox(0,0)[lb]{\smash{{\SetFigFont{9}{10.8}{\rmdefault}{\mddefault}{\updefault}{\color[rgb]{0,0,0}\(\beta_1\)}%
}}}}}
\put(4916,-1501){\rotatebox{270.0}{\makebox(0,0)[lb]{\smash{{\SetFigFont{9}{10.8}{\rmdefault}{\mddefault}{\updefault}{\color[rgb]{0,0,0}\(\beta_2\)}%
}}}}}
\put(4916,-3714){\rotatebox{270.0}{\makebox(0,0)[lb]{\smash{{\SetFigFont{9}{10.8}{\rmdefault}{\mddefault}{\updefault}{\color[rgb]{0,0,0}\(\beta_3\)}%
}}}}}
\put(4859,-6323){\rotatebox{270.0}{\makebox(0,0)[lb]{\smash{{\SetFigFont{9}{10.8}{\rmdefault}{\mddefault}{\updefault}{\color[rgb]{0,0,0}\(\beta_g\)}%
}}}}}
\put(6788,-27){\rotatebox{270.0}{\makebox(0,0)[lb]{\smash{{\SetFigFont{9}{10.8}{\rmdefault}{\mddefault}{\updefault}{\color[rgb]{0,0,0}\(\beta_1'\)}%
}}}}}
\put(6788,-1558){\rotatebox{270.0}{\makebox(0,0)[lb]{\smash{{\SetFigFont{9}{10.8}{\rmdefault}{\mddefault}{\updefault}{\color[rgb]{0,0,0}\(\beta_2'\)}%
}}}}}
\put(6845,-3771){\rotatebox{270.0}{\makebox(0,0)[lb]{\smash{{\SetFigFont{9}{10.8}{\rmdefault}{\mddefault}{\updefault}{\color[rgb]{0,0,0}\(\beta_3'\)}%
}}}}}
\put(6845,-6494){\rotatebox{270.0}{\makebox(0,0)[lb]{\smash{{\SetFigFont{9}{10.8}{\rmdefault}{\mddefault}{\updefault}{\color[rgb]{0,0,0}\(\beta_g'\)}%
}}}}}
\put(5881,201){\rotatebox{270.0}{\makebox(0,0)[lb]{\smash{{\SetFigFont{9}{10.8}{\rmdefault}{\mddefault}{\updefault}{\color[rgb]{0,0,0}\(\theta_1\)}%
}}}}}
\put(5881,-1048){\rotatebox{270.0}{\makebox(0,0)[lb]{\smash{{\SetFigFont{9}{10.8}{\rmdefault}{\mddefault}{\updefault}{\color[rgb]{0,0,0}\(\eta_1\)}%
}}}}}
\put(5766,-1388){\rotatebox{270.0}{\makebox(0,0)[lb]{\smash{{\SetFigFont{9}{10.8}{\rmdefault}{\mddefault}{\updefault}{\color[rgb]{0,0,0}\(\theta_2\)}%
}}}}}
\put(5824,-3487){\rotatebox{270.0}{\makebox(0,0)[lb]{\smash{{\SetFigFont{9}{10.8}{\rmdefault}{\mddefault}{\updefault}{\color[rgb]{0,0,0}\(\theta_3\)}%
}}}}}
\put(5824,-6153){\rotatebox{270.0}{\makebox(0,0)[lb]{\smash{{\SetFigFont{9}{10.8}{\rmdefault}{\mddefault}{\updefault}{\color[rgb]{0,0,0}\(\theta_g\)}%
}}}}}
\put(5881,-2579){\rotatebox{270.0}{\makebox(0,0)[lb]{\smash{{\SetFigFont{9}{10.8}{\rmdefault}{\mddefault}{\updefault}{\color[rgb]{0,0,0}\(\eta_2\)}%
}}}}}
\put(5881,-4679){\rotatebox{270.0}{\makebox(0,0)[lb]{\smash{{\SetFigFont{9}{10.8}{\rmdefault}{\mddefault}{\updefault}{\color[rgb]{0,0,0}\(\eta_3\)}%
}}}}}
\put(5881,-7402){\rotatebox{270.0}{\makebox(0,0)[lb]{\smash{{\SetFigFont{9}{10.8}{\rmdefault}{\mddefault}{\updefault}{\color[rgb]{0,0,0}\(\eta_g\)}%
}}}}}
\put(5256,-83){\rotatebox{270.0}{\makebox(0,0)[lb]{\smash{{\SetFigFont{9}{10.8}{\rmdefault}{\mddefault}{\updefault}{\color[rgb]{0,0,0}\(D_1\)}%
}}}}}
\put(5256,-1728){\rotatebox{270.0}{\makebox(0,0)[lb]{\smash{{\SetFigFont{9}{10.8}{\rmdefault}{\mddefault}{\updefault}{\color[rgb]{0,0,0}\(D_2\)}%
}}}}}
\put(5256,-3827){\rotatebox{270.0}{\makebox(0,0)[lb]{\smash{{\SetFigFont{9}{10.8}{\rmdefault}{\mddefault}{\updefault}{\color[rgb]{0,0,0}\(D_3\)}%
}}}}}
\put(5256,-6494){\rotatebox{270.0}{\makebox(0,0)[lb]{\smash{{\SetFigFont{9}{10.8}{\rmdefault}{\mddefault}{\updefault}{\color[rgb]{0,0,0}\(D_g\)}%
}}}}}
\put(6504,-197){\rotatebox{270.0}{\makebox(0,0)[lb]{\smash{{\SetFigFont{9}{10.8}{\rmdefault}{\mddefault}{\updefault}{\color[rgb]{0,0,0}\(E_1\)}%
}}}}}
\put(6504,-1785){\rotatebox{270.0}{\makebox(0,0)[lb]{\smash{{\SetFigFont{9}{10.8}{\rmdefault}{\mddefault}{\updefault}{\color[rgb]{0,0,0}\(E_2\)}%
}}}}}
\put(6504,-3827){\rotatebox{270.0}{\makebox(0,0)[lb]{\smash{{\SetFigFont{9}{10.8}{\rmdefault}{\mddefault}{\updefault}{\color[rgb]{0,0,0}\(E_3\)}%
}}}}}
\put(6504,-6607){\rotatebox{270.0}{\makebox(0,0)[lb]{\smash{{\SetFigFont{9}{10.8}{\rmdefault}{\mddefault}{\updefault}{\color[rgb]{0,0,0}\(E_g\)}%
}}}}}
\put(2590, 30){\rotatebox{270.0}{\makebox(0,0)[lb]{\smash{{\SetFigFont{9}{10.8}{\rmdefault}{\mddefault}{\updefault}{\color[rgb]{0,0,0}\(\beta_1\)}%
}}}}}
\put(2987,-1332){\rotatebox{270.0}{\makebox(0,0)[lb]{\smash{{\SetFigFont{9}{10.8}{\rmdefault}{\mddefault}{\updefault}{\color[rgb]{0,0,0}\(\beta_2\)}%
}}}}}
\put(2590,-3714){\rotatebox{270.0}{\makebox(0,0)[lb]{\smash{{\SetFigFont{9}{10.8}{\rmdefault}{\mddefault}{\updefault}{\color[rgb]{0,0,0}\(\beta_3\)}%
}}}}}
\put(2590,-6380){\rotatebox{270.0}{\makebox(0,0)[lb]{\smash{{\SetFigFont{9}{10.8}{\rmdefault}{\mddefault}{\updefault}{\color[rgb]{0,0,0}\(\beta_g\)}%
}}}}}
\put(4122,-480){\rotatebox{270.0}{\makebox(0,0)[lb]{\smash{{\SetFigFont{9}{10.8}{\rmdefault}{\mddefault}{\updefault}{\color[rgb]{0,0,0}\(\theta_1^H\)}%
}}}}}
\put(4122,-1955){\rotatebox{270.0}{\makebox(0,0)[lb]{\smash{{\SetFigFont{9}{10.8}{\rmdefault}{\mddefault}{\updefault}{\color[rgb]{0,0,0}\(\theta_2^H\)}%
}}}}}
\put(4065,-4111){\rotatebox{270.0}{\makebox(0,0)[lb]{\smash{{\SetFigFont{9}{10.8}{\rmdefault}{\mddefault}{\updefault}{\color[rgb]{0,0,0}\(\theta_3^H\)}%
}}}}}
\put(4065,-6834){\rotatebox{270.0}{\makebox(0,0)[lb]{\smash{{\SetFigFont{9}{10.8}{\rmdefault}{\mddefault}{\updefault}{\color[rgb]{0,0,0}\(\theta_g^H\)}%
}}}}}
\put(4462,-1615){\rotatebox{270.0}{\makebox(0,0)[lb]{\smash{{\SetFigFont{9}{10.8}{\rmdefault}{\mddefault}{\updefault}{\color[rgb]{0,0,0}\(\beta_1^H\)}%
}}}}}
\put(4009,-2353){\rotatebox{270.0}{\makebox(0,0)[lb]{\smash{{\SetFigFont{9}{10.8}{\rmdefault}{\mddefault}{\updefault}{\color[rgb]{0,0,0}\(\beta_2^H\)}%
}}}}}
\put(3951,-5018){\rotatebox{270.0}{\makebox(0,0)[lb]{\smash{{\SetFigFont{9}{10.8}{\rmdefault}{\mddefault}{\updefault}{\color[rgb]{0,0,0}\(\beta_3^H\)}%
}}}}}
\put(3951,-7798){\rotatebox{270.0}{\makebox(0,0)[lb]{\smash{{\SetFigFont{9}{10.8}{\rmdefault}{\mddefault}{\updefault}{\color[rgb]{0,0,0}\(\beta_g^H\)}%
}}}}}
\put(2761,-764){\rotatebox{270.0}{\makebox(0,0)[lb]{\smash{{\SetFigFont{9}{10.8}{\rmdefault}{\mddefault}{\updefault}{\color[rgb]{0,0,0}\(\eta_1^H\)}%
}}}}}
\put(3158,-2353){\rotatebox{270.0}{\makebox(0,0)[lb]{\smash{{\SetFigFont{9}{10.8}{\rmdefault}{\mddefault}{\updefault}{\color[rgb]{0,0,0}\(\eta_2^H\)}%
}}}}}
\put(2930,-4622){\rotatebox{270.0}{\makebox(0,0)[lb]{\smash{{\SetFigFont{9}{10.8}{\rmdefault}{\mddefault}{\updefault}{\color[rgb]{0,0,0}\(\eta_3^H\)}%
}}}}}
\put(2930,-7344){\rotatebox{270.0}{\makebox(0,0)[lb]{\smash{{\SetFigFont{9}{10.8}{\rmdefault}{\mddefault}{\updefault}{\color[rgb]{0,0,0}\(\eta_g^H\)}%
}}}}}
\put(4179,-1445){\rotatebox{270.0}{\makebox(0,0)[lb]{\smash{{\SetFigFont{9}{10.8}{\rmdefault}{\mddefault}{\updefault}{\color[rgb]{0,0,0}\(E_1^H\)}%
}}}}}
\put(3271, 30){\rotatebox{270.0}{\makebox(0,0)[lb]{\smash{{\SetFigFont{9}{10.8}{\rmdefault}{\mddefault}{\updefault}{\color[rgb]{0,0,0}\(D_1^H\)}%
}}}}}
\put(3214,-3600){\rotatebox{270.0}{\makebox(0,0)[lb]{\smash{{\SetFigFont{9}{10.8}{\rmdefault}{\mddefault}{\updefault}{\color[rgb]{0,0,0}\(D_3^H\)}%
}}}}}
\put(3214,-6380){\rotatebox{270.0}{\makebox(0,0)[lb]{\smash{{\SetFigFont{9}{10.8}{\rmdefault}{\mddefault}{\updefault}{\color[rgb]{0,0,0}\(D_g^H\)}%
}}}}}
\put(3612,-2636){\rotatebox{270.0}{\makebox(0,0)[lb]{\smash{{\SetFigFont{9}{10.8}{\rmdefault}{\mddefault}{\updefault}{\color[rgb]{0,0,0}\(E_2^H\)}%
}}}}}
\put(3555,-4679){\rotatebox{270.0}{\makebox(0,0)[lb]{\smash{{\SetFigFont{9}{10.8}{\rmdefault}{\mddefault}{\updefault}{\color[rgb]{0,0,0}\(E_3^H\)}%
}}}}}
\put(3668,-7402){\rotatebox{270.0}{\makebox(0,0)[lb]{\smash{{\SetFigFont{9}{10.8}{\rmdefault}{\mddefault}{\updefault}{\color[rgb]{0,0,0}\(E_g^H\)}%
}}}}}
\put(321,201){\rotatebox{270.0}{\makebox(0,0)[lb]{\smash{{\SetFigFont{9}{10.8}{\rmdefault}{\mddefault}{\updefault}{\color[rgb]{0,0,0}\(\beta_1\)}%
}}}}}
\put(605,-1501){\rotatebox{270.0}{\makebox(0,0)[lb]{\smash{{\SetFigFont{9}{10.8}{\rmdefault}{\mddefault}{\updefault}{\color[rgb]{0,0,0}\(\beta_2\)}%
}}}}}
\put(321,-3771){\rotatebox{270.0}{\makebox(0,0)[lb]{\smash{{\SetFigFont{9}{10.8}{\rmdefault}{\mddefault}{\updefault}{\color[rgb]{0,0,0}\(\beta_3\)}%
}}}}}
\put(321,-6437){\rotatebox{270.0}{\makebox(0,0)[lb]{\smash{{\SetFigFont{9}{10.8}{\rmdefault}{\mddefault}{\updefault}{\color[rgb]{0,0,0}\(\beta_g\)}%
}}}}}
\put(2250,201){\rotatebox{270.0}{\makebox(0,0)[lb]{\smash{{\SetFigFont{9}{10.8}{\rmdefault}{\mddefault}{\updefault}{\color[rgb]{0,0,0}\(\beta_1'\)}%
}}}}}
\put(2079,-1615){\rotatebox{270.0}{\makebox(0,0)[lb]{\smash{{\SetFigFont{9}{10.8}{\rmdefault}{\mddefault}{\updefault}{\color[rgb]{0,0,0}\(\beta_2'\)}%
}}}}}
\put(2250,-3658){\rotatebox{270.0}{\makebox(0,0)[lb]{\smash{{\SetFigFont{9}{10.8}{\rmdefault}{\mddefault}{\updefault}{\color[rgb]{0,0,0}\(\beta_3'\)}%
}}}}}
\put(2250,-6323){\rotatebox{270.0}{\makebox(0,0)[lb]{\smash{{\SetFigFont{9}{10.8}{\rmdefault}{\mddefault}{\updefault}{\color[rgb]{0,0,0}\(\beta_g'\)}%
}}}}}
\put(1626,-7742){\rotatebox{270.0}{\makebox(0,0)[lb]{\smash{{\SetFigFont{9}{10.8}{\rmdefault}{\mddefault}{\updefault}{\color[rgb]{0,0,0}\(\beta_g^H\)}%
}}}}}
\put(1682,-5018){\rotatebox{270.0}{\makebox(0,0)[lb]{\smash{{\SetFigFont{9}{10.8}{\rmdefault}{\mddefault}{\updefault}{\color[rgb]{0,0,0}\(\beta_3^H\)}%
}}}}}
\put(1626,-2692){\rotatebox{270.0}{\makebox(0,0)[lb]{\smash{{\SetFigFont{9}{10.8}{\rmdefault}{\mddefault}{\updefault}{\color[rgb]{0,0,0}\(\beta_2^H\)}%
}}}}}
\put(2250,-2239){\rotatebox{270.0}{\makebox(0,0)[lb]{\smash{{\SetFigFont{9}{10.8}{\rmdefault}{\mddefault}{\updefault}{\color[rgb]{0,0,0}\(\beta_1^H\)}%
}}}}}
\put(492,-764){\rotatebox{270.0}{\makebox(0,0)[lb]{\smash{{\SetFigFont{9}{10.8}{\rmdefault}{\mddefault}{\updefault}{\color[rgb]{0,0,0}\(\eta_1^H\)}%
}}}}}
\put(889,-2522){\rotatebox{270.0}{\makebox(0,0)[lb]{\smash{{\SetFigFont{9}{10.8}{\rmdefault}{\mddefault}{\updefault}{\color[rgb]{0,0,0}\(\eta_2^H\)}%
}}}}}
\put(661,-4564){\rotatebox{270.0}{\makebox(0,0)[lb]{\smash{{\SetFigFont{9}{10.8}{\rmdefault}{\mddefault}{\updefault}{\color[rgb]{0,0,0}\(\eta_3^H\)}%
}}}}}
\put(661,-7288){\rotatebox{270.0}{\makebox(0,0)[lb]{\smash{{\SetFigFont{9}{10.8}{\rmdefault}{\mddefault}{\updefault}{\color[rgb]{0,0,0}\(\eta_g^H\)}%
}}}}}
\put(1286,-7344){\rotatebox{270.0}{\makebox(0,0)[lb]{\smash{{\SetFigFont{9}{10.8}{\rmdefault}{\mddefault}{\updefault}{\color[rgb]{0,0,0}\(\eta_g\)}%
}}}}}
\put(1286,-4622){\rotatebox{270.0}{\makebox(0,0)[lb]{\smash{{\SetFigFont{9}{10.8}{\rmdefault}{\mddefault}{\updefault}{\color[rgb]{0,0,0}\(\eta_3\)}%
}}}}}
\put(1286,-2522){\rotatebox{270.0}{\makebox(0,0)[lb]{\smash{{\SetFigFont{9}{10.8}{\rmdefault}{\mddefault}{\updefault}{\color[rgb]{0,0,0}\(\eta_2\)}%
}}}}}
\put(1286,-991){\rotatebox{270.0}{\makebox(0,0)[lb]{\smash{{\SetFigFont{9}{10.8}{\rmdefault}{\mddefault}{\updefault}{\color[rgb]{0,0,0}\(\eta_1\)}%
}}}}}
\put(1286,257){\rotatebox{270.0}{\makebox(0,0)[lb]{\smash{{\SetFigFont{9}{10.8}{\rmdefault}{\mddefault}{\updefault}{\color[rgb]{0,0,0}\(\theta_1\)}%
}}}}}
\put(1286,-1274){\rotatebox{270.0}{\makebox(0,0)[lb]{\smash{{\SetFigFont{9}{10.8}{\rmdefault}{\mddefault}{\updefault}{\color[rgb]{0,0,0}\(\theta_2\)}%
}}}}}
\put(1286,-3374){\rotatebox{270.0}{\makebox(0,0)[lb]{\smash{{\SetFigFont{9}{10.8}{\rmdefault}{\mddefault}{\updefault}{\color[rgb]{0,0,0}\(\theta_3\)}%
}}}}}
\put(1286,-6153){\rotatebox{270.0}{\makebox(0,0)[lb]{\smash{{\SetFigFont{9}{10.8}{\rmdefault}{\mddefault}{\updefault}{\color[rgb]{0,0,0}\(\theta_g\)}%
}}}}}
\put(1796,-480){\rotatebox{270.0}{\makebox(0,0)[lb]{\smash{{\SetFigFont{9}{10.8}{\rmdefault}{\mddefault}{\updefault}{\color[rgb]{0,0,0}\(\theta_1^H\)}%
}}}}}
\put(1796,-2012){\rotatebox{270.0}{\makebox(0,0)[lb]{\smash{{\SetFigFont{9}{10.8}{\rmdefault}{\mddefault}{\updefault}{\color[rgb]{0,0,0}\(\theta_2^H\)}%
}}}}}
\put(1796,-4054){\rotatebox{270.0}{\makebox(0,0)[lb]{\smash{{\SetFigFont{9}{10.8}{\rmdefault}{\mddefault}{\updefault}{\color[rgb]{0,0,0}\(\theta_3^H\)}%
}}}}}
\put(1796,-6777){\rotatebox{270.0}{\makebox(0,0)[lb]{\smash{{\SetFigFont{9}{10.8}{\rmdefault}{\mddefault}{\updefault}{\color[rgb]{0,0,0}\(\theta_g^H\)}%
}}}}}
\put(1456,-6437){\rotatebox{270.0}{\makebox(0,0)[lb]{\smash{{\SetFigFont{9}{10.8}{\rmdefault}{\mddefault}{\updefault}{\color[rgb]{0,0,0}\(T_g\)}%
}}}}}
\put(1456,-3714){\rotatebox{270.0}{\makebox(0,0)[lb]{\smash{{\SetFigFont{9}{10.8}{\rmdefault}{\mddefault}{\updefault}{\color[rgb]{0,0,0}\(T_3\)}%
}}}}}
\put(1456,-1615){\rotatebox{270.0}{\makebox(0,0)[lb]{\smash{{\SetFigFont{9}{10.8}{\rmdefault}{\mddefault}{\updefault}{\color[rgb]{0,0,0}\(T_2\)}%
}}}}}
\put(1456,-83){\rotatebox{270.0}{\makebox(0,0)[lb]{\smash{{\SetFigFont{9}{10.8}{\rmdefault}{\mddefault}{\updefault}{\color[rgb]{0,0,0}\(T_1\)}%
}}}}}
\put(775,-424){\rotatebox{270.0}{\makebox(0,0)[lb]{\smash{{\SetFigFont{9}{10.8}{\rmdefault}{\mddefault}{\updefault}{\color[rgb]{0,0,0}\(\theta_1'\)}%
}}}}}
\put(832,-2012){\rotatebox{270.0}{\makebox(0,0)[lb]{\smash{{\SetFigFont{9}{10.8}{\rmdefault}{\mddefault}{\updefault}{\color[rgb]{0,0,0}\(\theta_2'\)}%
}}}}}
\put(775,-4054){\rotatebox{270.0}{\makebox(0,0)[lb]{\smash{{\SetFigFont{9}{10.8}{\rmdefault}{\mddefault}{\updefault}{\color[rgb]{0,0,0}\(\theta_3'\)}%
}}}}}
\put(775,-6777){\rotatebox{270.0}{\makebox(0,0)[lb]{\smash{{\SetFigFont{9}{10.8}{\rmdefault}{\mddefault}{\updefault}{\color[rgb]{0,0,0}\(\theta_g'\)}%
}}}}}
\put(2079,-820){\rotatebox{270.0}{\makebox(0,0)[lb]{\smash{{\SetFigFont{9}{10.8}{\rmdefault}{\mddefault}{\updefault}{\color[rgb]{0,0,0}\(\eta_1'\)}%
}}}}}
\put(1910,-2466){\rotatebox{270.0}{\makebox(0,0)[lb]{\smash{{\SetFigFont{9}{10.8}{\rmdefault}{\mddefault}{\updefault}{\color[rgb]{0,0,0}\(\eta_2'\)}%
}}}}}
\put(1910,-4564){\rotatebox{270.0}{\makebox(0,0)[lb]{\smash{{\SetFigFont{9}{10.8}{\rmdefault}{\mddefault}{\updefault}{\color[rgb]{0,0,0}\(\eta_3'\)}%
}}}}}
\put(1910,-7288){\rotatebox{270.0}{\makebox(0,0)[lb]{\smash{{\SetFigFont{9}{10.8}{\rmdefault}{\mddefault}{\updefault}{\color[rgb]{0,0,0}\(\eta_g'\)}%
}}}}}
\put(3328,-1501){\rotatebox{270.0}{\makebox(0,0)[lb]{\smash{{\SetFigFont{9}{10.8}{\rmdefault}{\mddefault}{\updefault}{\color[rgb]{0,0,0}\(D_2^H\)}%
}}}}}
\end{picture}%
\caption{Heegaard diagrams for $\#^gS^1\times S^2$}
\label{Figure:S1S2}
\end{figure}
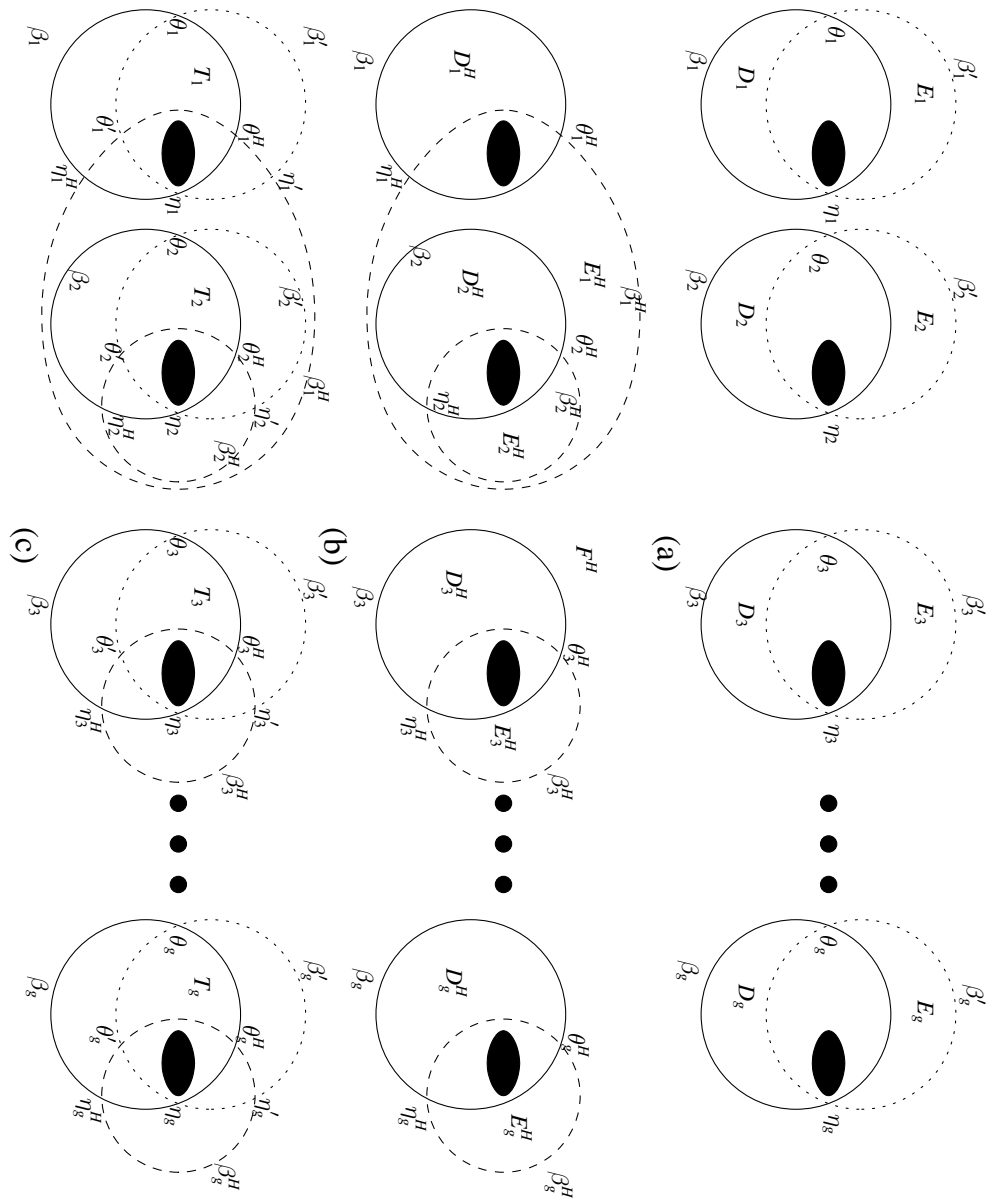

We now state the main steps, and then return later to the proofs.
First, notation.  Let $\cH=(\Sigma,\va,\vb,\Fz)$ be a pointed
Heegaard diagram.  Let $\beta_i'$ be a small Hamiltonian
perturbation of $\beta_i$ intersecting $\beta_i$ transversally in
two points and disjoint from $\beta_j$ for $i\neq j$.  Let
$\beta_1^H$ be a curve obtained by handlesliding $\beta_1$ over
$\beta_2$ in the complement of $\Fz$ and then doing a small
isotopy so that $\beta_1^H$ intersects each of $\beta_1$ and
$\beta_1'$ transversally in two points and is disjoint from
$\beta_i$ and $\beta_i'$ for $i>1$.  Let $\beta_i^H$, $i>1$,
be a small Hamiltonian perturbation of $\beta_i$, intersecting each
of $\beta_i$ and $\beta_i'$ transversally in two points and
disjoint from $\beta_j$ and $\beta_j'$ for $j\neq i$.  See
\fullref{Figure:S1S2}.  In the notation of~\cite[Section 9]{OS1},
our $\beta_i'$ is their $\delta_i$ and our $\beta_i^H$ is their
$\gamma_i$.

There are $2^g$ intersection points in $(\Sigma,\vb,\vb',\Fz)$.
Let $\vec{\theta}_{\beta,\beta'}=\{\theta_1,\cdots,\theta_g\}$
denote the unique intersection point of maximal grading.  Similarly, in
$(\Sigma,\vb,\vb^H,\Fz)$ (respectively $(\Sigma,\vb,\vb',\Fz)$)
there are exactly $2^g$ intersection points.  Let
$\vec{\theta}_{\beta,\beta^H}=\{\theta^H_1,\cdots,\theta_g^H\}$
(respectively
$\vec{\theta}_{\beta',\beta^H}=\{\theta'_1,\cdots,\theta'_g\}$)
denote the unique intersection point of maximal grading.  See
\fullref{Figure:S1S2}.

A few words about $\Spin^\bC$--structures.  We fix a
 $\Spin^\bC$--structure $\Ss$ on $Y_{\alpha,\beta}$ for the rest
 of the section; when we refer to Floer homology groups of
 $Y_{\alpha,\beta}$ we mean the groups computed with respect to
 $\Ss$.  Now, note that all of
the $3\cdot 2^g$ intersection points in $(\Sigma,\vb,\vb')$,
 $(\Sigma,\vb,\vb^H)$ and $(\Sigma,\vb',\vb^H)$ represent the unique torsion $\Spin^\bC$--structure on $\#^g S^1\times S^2$.  (This
follows, for instance, from \fullref{OSTheorem4.9}.)  Throughout
this section, when discussing Heegaard diagrams for $\#^gS^1\times
S^2$ we shall implicitly work with the torsion
$\Spin^\bC$--structure.  Later in this section, we shall consider the Heegaard
triple--diagrams $(\Sigma,\va,\vb,\vb^H)$,
$(\Sigma,\va,\vb,\vb')$, and $(\Sigma,\vb,\vb^H,\vb')$.  It is
easy to see that $X_{\alpha,\beta,\beta^H}$ and
$X_{\alpha,\beta,\beta'}$ are each diffeomorphic to the complement
of (a neighborhood of) a bouquet of $g$ circles in
$Y_{\alpha,\beta}\times[0,1]$.  (Compare~\cite[Example 8.1]{OS1}.)
By considering cohomology groups, say, it follows that given a
$\Spin^\bC$--structure $\Ss$ on $Y_{\alpha,\beta}$ there is a
unique $\Spin^\bC$--structure on $X_{\alpha,\beta,\beta^H}$
(respectively $X_{\alpha,\beta,\beta'}$) which restricts to $\Ss$
on $Y_{\alpha,\beta}$ and the torsion $\Spin^\bC$--structure on
$Y_{\beta,\beta^H}$ (respectively $Y_{\beta,\beta'}$).  It follows that
there is a unique $\Spin^\bC$--structure on
$X_{\beta,\beta^H,\beta'}$ which restricts to the torsion
$\Spin^\bC$--structure on the three boundary components.  When
discussing triangle maps we shall always assume that they are computed
with respect to these choices of $\Spin^\bC$--structures on
$X_{\alpha,\beta,\beta^H}$, $X_{\alpha,\beta,\beta'}$ and
$X_{\beta,\beta^H,\beta'}$.  We shall, however, tend to suppress
$\Spin^\bC$--structures from the notation.

\begin{Lem}\label{Lemma:HandleslideOrientations}
\begin{enumerate}
\item There is a coherent orientation system $\frako_{\beta,\beta'}$
  for $(\Sigma,\vb,\vb',\Fz)$ with respect to which
  $\vec{\theta}_{\beta,\beta'}$ is a cycle in $\widehat{CF}$ and
  $[\vec{\theta}_{\beta,\beta'},0]$ is a cycle in $CF^{\leq0}$. 
\item There is a coherent orientation system $\frako_{\beta,\beta^H}$
  (respectively $\frako_{\beta',\beta^H}$) for\break
  $(\Sigma,\vb,\vb^H,\Fz)$ (respectively
  $(\Sigma,\vb^H,\vb',\Fz)$) with respect to which
  $\vec{\theta}_{\beta,\beta^H}$ (respectively
  $\vec{\theta}_{\beta^H,\beta'}$) is a cycle in $\widehat{CF}$ and
  $[\vec{\theta}_{\beta,\beta^H},0]$ (respectively
  $[\vec{\theta}_{\beta^H,\beta'},0]$) is a cycle in
  $CF^{\leq0}$. 
\item The orientation systems $\frako_{\beta,\beta'}$,
  $\frako_{\beta,\beta^H}$ and $\frako_{\beta^H,\beta'}$ above can be
  chosen so that for some orientation system
  $\frako_{\beta,\beta^H,\beta'}$ for $(\Sigma,\vb,\vb^H,\vb',\Fz)$,
  the orientation systems $\frako_{\beta,\beta^H,\beta'}$,
  $\frako_{\beta,\beta'}$, $\frako_{\beta,\beta^H}$ and
  $\frako_{\beta^H,\beta'}$ are coherent.
\end{enumerate}
\end{Lem}
The proof will be given in \fullref{Subsection:PropHSC}.  We shall choose, once and for all, orientation systems as in the lemma.

The goal of this section is the following
\begin{Prop}{\rm(Compare~\cite[Theorem 9.5]{OS1})}\label{Prop:HandleslideInvariance}\qua
Fix a $\Spin^\bC$--structure $\Ss$ on $Y_{\alpha,\beta}=Y_{\alpha,\beta^H}$.  Then the map $$\hat{F}_{\alpha,\beta,\beta^H}(\cdot\otimes\vec{\theta}_{\beta,\beta^H})\co \widehat{HF}(\Sigma,\va,\vb,\Fz,\Ss)\to \widehat{HF}(\Sigma,\va,\vb^H,\Fz,\Ss)$$
is an isomorphism.  The map
$$
F^*_{\alpha,\beta,\beta^H}(\cdot\otimes[\vec{\theta}_{\beta,\beta^H},0])\co HF^*(\Sigma,\va,\vb,\Fz,\Ss)\to HF^*(\Sigma,\va,\vb,\Fz,\Ss)
$$
is an isomorphism for $*\in\{\infty,+,i\}$.  These isomorphisms commute with the long exact sequences and the $\bZ[U]\otimes_\bZ\Lambda^*H^1(Y)/Tors$--actions.  
\end{Prop}

{\bf Note}\qua  we have suppressed some discussion of orientation systems from
the statement of the proposition.  See \fullref{Lem:HSOrientations} below.

The essence of the proof of \fullref{Prop:HandleslideInvariance} is the following two propositions:
\begin{Prop}{\rm(Compare~\cite[Lemma 9.7]{OS1})}\label{Prop:HSC}\qua
\begin{eqnarray}
\hat{F}_{\beta,\beta^H,\beta'}(\vec{\theta}_{\beta,\beta^H}\otimes\vec{\theta}_{\beta^H,\beta'})&=&\vec{\theta}_{\beta,\beta'}.\nonumber\\
F^{\leq0}_{\beta,\beta^H,\beta'}([\vec{\theta}_{\beta,\beta^H},0]\otimes[\vec{\theta}_{\beta^H,\beta'},0])&=&[\vec{\theta}_{\beta,\beta'},0].\nonumber
\end{eqnarray}
\end{Prop}
(The proof is in \fullref{Subsection:PropHSC}.)

\begin{Prop}{\rm(Compare~\cite[Proposition 9.8]{OS1})}\label{Prop:HSB}\qua
The map
$\hat{F}_{\alpha,\beta,\beta'}(\cdot\otimes\vec{\theta}_{\beta,\beta'})$ (respectively the maps $F^*(\cdot\otimes[\vec{\theta}_{\beta,\beta'},0]$ for $*\in\{+,-,\infty\}$) is the isomorphism $\hat{\Phi}$ (respectively $\Phi^*$, $*\in\{+,-,\infty\}$) induced by the isotopy from $\vb$ to $\vb'$ as in \fullref{Section:Isotopy}.  
\end{Prop}
(The proof is in \fullref{Subsection:PropHSB}.)

We sketch the proof of \fullref{Prop:HandleslideInvariance} assuming \fullref{Prop:HSB} and \fullref{Prop:HSC}.

\proof[Proof of \fullref{Prop:HandleslideInvariance} (sketch)]
Observe that
\begin{eqnarray}
\hat{F}_{\alpha,\beta^H\!,\beta'}\left(\hat{F}_{\alpha,\beta,\beta^H}(\cdot\otimes\vec{\theta}_{\beta,\beta^H})\otimes\vec{\theta}_{\beta^H,\beta'}\right) & = & \hat{F}_{\alpha,\beta,\beta'}\!\left(\cdot\otimes \hat{F}_{\beta,\beta^H,\beta'}( \vec{\theta}_{\beta,\beta^H}\otimes\vec{\theta}_{\beta^H,\beta'})\right)\nonumber\\
&=& \hat{F}_{\alpha,\beta,\beta'}(\cdot\otimes\vec{\theta}_{\beta,\beta'})\nonumber\\
&=&\hat{\Phi}_{\beta,\beta'}(\cdot)\nonumber
\end{eqnarray}
where $\hat{\Phi}_{\beta,\beta'}$ is the isomorphism induced by the
isotopy from $\vb$ to $\vb'$.  (The first equality follows by the
associativity of triangle maps (\fullref{AssocProp}).  The
second follows from \fullref{Prop:HSC}.  The third follows
from \fullref{Prop:HSB}.) 

It follows that $\hat{F}_{\alpha,\beta^H,\beta'}$ is surjective and $\hat{F}_{\alpha,\beta,\beta^H}$ is injective.  The same argument with the roles of $\beta$ and $\beta^H$ exchanged shows that $\hat{F}_{\alpha,\beta,\beta^H}$ is surjective.  So, $\hat{F}_{\alpha,\beta,\beta^H}$ is an isomorphism, proving that $\widehat{HF}$ is invariant under handleslides.  The proofs for $HF^-$, $HF^+$ and $HF^\infty$ are just the same.
\endproof

\subsection{Proofs of Proposition \ref{Prop:HSC} and Lemma \ref{Lemma:HandleslideOrientations}}\label{Subsection:PropHSC}

\proof[Proof of \fullref{Prop:HSC}]
With the notation of \fullref{Figure:S1S2}, note that the domain $T_1+\cdots+T_g$ achieves transversality, and admits a unique holomorphic representative, so $\#\cM^{T_1+\cdots+T_g}=\pm1$.

For any intersection point $\xi\neq\vec{\theta}_{\beta,\beta'}$ of $(\Sigma,\vb,\vb',\Fz)$ and any $A\in\pi_2(\vec{\theta}_{\beta,\beta^H},\vec{\theta}_{\beta^H,\beta'},\xi)$, either $n_{\Fz}(A)<0$ or $\ind(A)>0$.

Any element of
$\pi_2(\vec{\theta}_{\beta,\beta^H},\vec{\theta}_{\beta^H,\beta'},\vec{\theta}_{\beta,\beta'})$
other than $T_1+\cdots+T_g$ either has $\ind>0$ or some negative
coefficient. 

It follows that for any choice of coherent orientation system, 
$$\eqalignbot
{\hat{f}_{\beta,\beta^H,\beta'}(\vec{\theta}_{\beta,\beta^H}\otimes\vec{\theta}_{\beta^H,\beta'})&=\pm\vec{\theta}_{\beta,\beta'}\cr 
f^{\leq 0}_{\beta,\beta^H,\beta'}([\vec{\theta}_{\beta,\beta^H},0]\otimes[\vec{\theta}_{\beta^H,\beta'},0)&=\pm[\vec{\theta}_{\beta,\beta'},0]}\proved$$

\proof[Proof of \fullref{Lemma:HandleslideOrientations}]
Firstly, note that Part (2), \fullref{Prop:HSC}, \fullref{Lemma:TriangleOrientations} and \fullref{Lemma:TriangleMapsChain} imply Parts (1) and (3).  So, it remains only to prove Part (2)

Label the components of $\Sigma\setminus(\bb\cup\bb^H)$ by $D_1^H,\cdots,D_g^H,E_1^H,\cdots,E_1^H$ and $F^H$, and the points in $\beta_i\cap\beta_i^H$ by $\theta_i^H$ and $\eta_i^H$ as in \fullref{Figure:S1S2}. 

Observe that $$gr\left(\vec{\theta}_{\beta,\beta^H},\{\theta_{i_1}^H,\cdots,\theta_{i_k}^H,\eta_{j_{k+1}}^H,\cdots,\eta_{j_g}^H\}\right) = g-k.$$
It follows from this and positivity of domains that the only homology classes which could contribute to $\hat{\bdy}\vec{\theta}_{\beta,\beta^H}$ or $\bdy^{\leq0}[\vec{\theta}_{\beta,\beta^H},0]$ are $D_1^H,\cdots,D_g^H,E_2^H,\cdots,E_g^H,E_1^H+E_2^H$ and $E_1^H+D_2^H$.  

By \fullref{Prop:BdyInjectivity} and \fullref{Lem:BdyInjectivity}, a generic perturbation of the $\alpha$-- and $\beta$--circles achieves transversality for all of the moduli spaces under consideration.  It follows from the Riemann mapping theorem that there is a unique equivalence class of holomorphic curves in each of $\hcM^{D_1^H},\cdots,\hcM^{D_g^H}$ and $\hcM^{E_2^H},\cdots,\hcM^{E_g^H}$.

We will show that for a generic perturbation of the $\alpha$-- and $\beta$--circles, one of $\hcM^{E_1^H+E_2^H}$ and $\hcM^{E_1^H+D_1^H}$ has one element and the other is empty.  To this end, we use the following
\begin{Sublem}\label{Sublem:Annuli}
Let $A$ be an annulus with boundary circles $C_1$
  and $C_2$ and marked points $x_i$, $y_i$ on $C_i$.  Define
  the \emph{conformal angle} between $x_1$ and $y_1$ (respectively
  $x_2$ and $y_2$) as follows.  The annulus $A$ is conformally
  equivalent to $\{1\leq |z|\leq R\}\subset\bC$ for some $R$.  We
  can arrange that the equivalence takes $x_1$ (respectively
  $x_2$) to $1$.  Then the conformal angle between $x_1$ and
  $y_1$ (respectively $x_2$ and $y_2$) as the image of $y_1$
  (respectively $y_2$) in $S^1$ under the equivalence.  Then:
\begin{enumerate}
\item The conformal angle is well--defined.
\item There is a holomorphic involution of $A$ exchanging $x_1$
  and $x_2$, and $y_1$ and $y_2$, if and only if the conformal
  angle between $x_1$ and $y_1$ equals the conformal angle between
  $x_2$ and $y_2$.
\end{enumerate}
\end{Sublem}
\proof
It is well--known that every conformal automorphism of the annulus
$A'=\{1\leq|z|\leq R\}\subset\bC$ is either of the form $z\mapsto
e^{i\theta}z$ or $z\mapsto\frac{Re^{i\theta}}{z}$.  (Consider the
universal cover.)  Both parts of the claim are immediate from this
observation.
\endproof
Note that the condition in Part (2) of the sublemma is the same as the
condition that there be a holomorphic $2$--fold branched cover
$c\co A\to\bD$ such that $c(x_1)=c(x_2)=-i$ and $c(y_1)=c(y_2)=i$.

Now, if $u\co S\to W$ is an element of $\hcM^{E_1^H+E_2^H}$ then
$\pi_\Sigma\circ u$ is an analytic isomorphism on the interior of $S$.  It
follows that $S$ can be obtained by cutting
$\overline{E_1^H+E_2^H}$ along $\beta_2^H$.  Let $x_1$ denote
the preimage of $\theta_1^H$, $y_1$ the preimage of $\eta_1^H$
and let $x_2$ and $y_2$ denote the two preimages of
$\theta_2^H$.  (Which preimage should be labeled $x_2$ and which $y_2$ is clear.)  Then elements of $\hcM^{E_1^H+E_2^H}$ correspond to
choices of
cuts of $\overline{E_1^H\cup E_2^H}$ along $\beta_2^H$ starting at
$\theta_2^H$ for which there are holomorphic involutions exchanging
$x_1$ and $x_2$, and $y_1$ and $y_2$.  A similar remark
applies to $\hcM^{E_1^H+D_2^H}$ with $\beta_2$ in place of
$\beta_2^H$. 

For $p\in\beta_2^H$, let $a^H(p)$ denote the ratio of the
conformal angle from $x_2$ to $y_2$ to the conformal angle from
$x_1$ to $y_1$, for the annulus obtained by cutting
$\overline{E_1^H+E_2^H}$ along $\beta_2^H$ to the point $p$.
Let $a(p)$ denote the same with $\overline{E_1^H+D_2^H}$ in place
of $\overline{E_1^H\cup E_2^H}$ and $\beta_2$ in place of
$\beta_2^H$.

Observe that:
\begin{enumerate}
\item $a(p)=1$ (respectively $a^H(p)=1$) if and only if cutting to $p$ gives an element of $\hcM^{E_1^H+D_2^H}$ (respectively $\hcM^{E_1^H+E_2^H}$).
\item $a(\theta_2^H)=a^H(\theta_2^H)$.
\item $a$ and $a^H$ are monotone as $p$ travels from $\theta_2^H$ to $\eta_2^H$.
\item $a\to0$ as $p\to\eta_2^H$.  $a^H\to\infty$ as $p\to\eta_2^H$.
\end{enumerate}
The first two claims are obvious.  The third is clear.  The fourth follows by considering the Gromov limit as $p\to\eta_2^H$.

For a generic choice of $\alpha$-- and $\beta$--circles, $a(\theta_2^H)\neq 1$.  It follows that one of $a$ or $a^H$ assumes the value $1$ once, and the other never assumes the value $1$.  Hence one of $\hcM^{E_1^H+D_2^H}$ or $\hcM^{E_1^H+E_2^H}$ has a unique element and the other is empty.

Now, specifying a coherent orientation system for
$(\Sigma,\vb,\vb^H,\Fz)$ is equivalent to specifying orientations
arbitrarily over $D_1^H,\cdots,D_g^H,E_2^H,\cdots,E_g^H,$
and either $E_1^H+D_2^H$ or $E_1^H+E_2^H$.  Choose the
orientations such that $\#\hcM^{D_i^H}=1$ ($i=1,\cdots,g$),
$\#\hcM^{E_i^H}=-1$ ($i=2,\cdots,g$) and
$\#\hcM^{E_1^H+D_2^H}+\#\hcM^{E_1^H+E_2^H}$=-1.  Then,
$\theta_{\beta,\beta^H}$ is a cycle in $\widehat{HF}$ and
$[\theta_{\beta,\beta^H},0]$ is a cycle in $HF^{\leq0}$.
\endproof

\subsection{Proof of Proposition \ref{Prop:HSB}}\label{Subsection:PropHSB}
\begin{figure}
\centering
\begin{picture}(0,0)%
\includegraphics[scale=0.7]{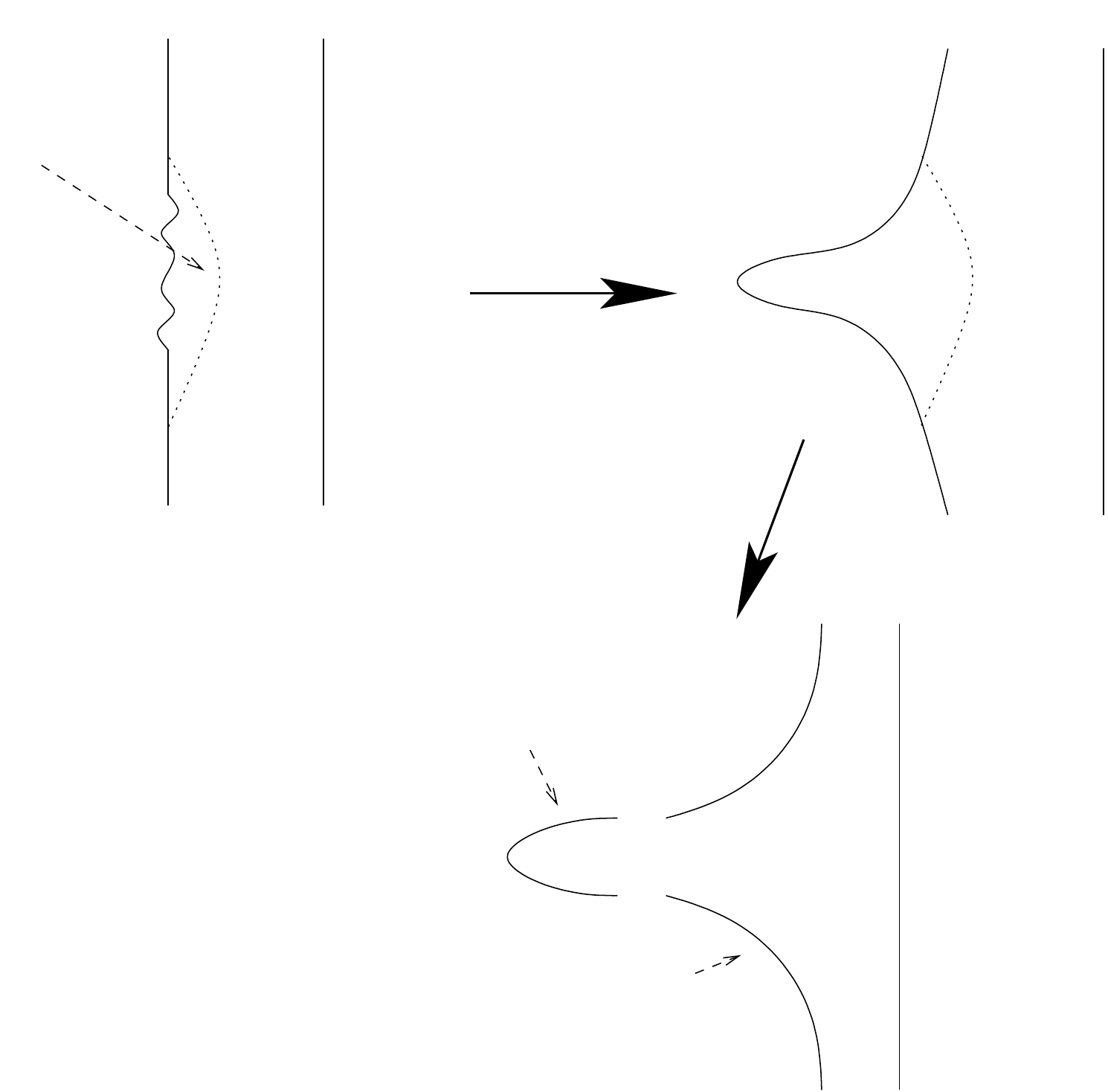}%
\end{picture}%
\setlength{\unitlength}{2763sp}%
\begingroup\makeatletter\ifx\SetFigFont\undefined%
\gdef\SetFigFont#1#2#3#4#5{%
  \reset@font\fontsize{#1}{#2pt}%
  \fontfamily{#3}\fontseries{#4}\fontshape{#5}%
  \selectfont}%
\fi\endgroup%
\begin{picture}(7135,7047)(587,-7573)
\put(3526,-5236){\makebox(0,0)[lb]{\smash{{\SetFigFont{10}{12.0}{\rmdefault}{\mddefault}{\updefault}{\color[rgb]{0,0,0}\(1\)--gons here.}%
}}}}
\put(1051,-4425){\makebox(0,0)[lb]{\smash{{\SetFigFont{10}{12.0}{\rmdefault}{\mddefault}{\updefault}{\color[rgb]{0,0,0}Noncylindrical cobordism}%
}}}}
\put(1051,-4613){\makebox(0,0)[lb]{\smash{{\SetFigFont{10}{12.0}{\rmdefault}{\mddefault}{\updefault}{\color[rgb]{0,0,0}induces \(\Phi\)}%
}}}}
\put(1295,-1339){\makebox(0,0)[rb]{\smash{{\SetFigFont{10}{12.0}{\rmdefault}{\mddefault}{\updefault}{\color[rgb]{0,0,0}Stretch the}%
}}}}
\put(1295,-1527){\makebox(0,0)[rb]{\smash{{\SetFigFont{10}{12.0}{\rmdefault}{\mddefault}{\updefault}{\color[rgb]{0,0,0}neck here}%
}}}}
\put(5251,-6961){\makebox(0,0)[rb]{\smash{{\SetFigFont{10}{12.0}{\rmdefault}{\mddefault}{\updefault}{\color[rgb]{0,0,0}Triangle map here.}%
}}}}
\put(2476,-661){\makebox(0,0)[rb]{\smash{{\SetFigFont{10}{12.0}{\rmdefault}{\mddefault}{\updefault}{\color[rgb]{0,0,0}\((\alpha,\beta')\)}%
}}}}
\put(2401,-4036){\makebox(0,0)[rb]{\smash{{\SetFigFont{10}{12.0}{\rmdefault}{\mddefault}{\updefault}{\color[rgb]{0,0,0}\((\alpha,\beta)\)}%
}}}}
\end{picture}%
\caption{Sketch of proof of \fullref{Prop:HSB}}
\label{HandleslideOutlineFig}
\end{figure}

Our proof, illustrated schematically in
\fullref{HandleslideOutlineFig}, proceeds by a neck--stretching
argument.  We will show that the moduli spaces used to define the chain map $\Phi$
corresponding to the isotopy from $\vb$ to $\vb'$ are the products of the moduli spaces used to define
$F_{\alpha,\beta,\beta'}$ and a moduli space of $1$--gons.  Before
showing this, we need to understand the moduli space of
$1$--gons.

Fix a Hamiltonian isotopy $\beta_t$ from $\beta$ to $\beta'$,
agreeing with $\beta$ for $t\ll0$ and with $\beta'$ for
$t\gg0$.  Let $\bH$ denote the upper half plane in $\bC$.  The
isotopy $\beta_t$ defines a $g$--tuple of Lagrangian cylinders
$C_\beta$ in $\Sigma\times\bdy\bH\subset\Sigma\times\bH=:W_\beta$.
We will consider maps $(S,\bdy S)\to(W_\beta,C_\beta)$, with one
component of $\bdy S$ mapped to each cylinder in $C_\beta$.  Such
maps will converge to some I--chord collection
$\vx=\{x_i\}\subset\bb\cap\bb'$ at $\infty\in\bdy\bH$.  Let $\pi_2(\vx)$ denote the
collection of homology classes of such maps.
\begin{Lem}The map $n_\Fz\co \pi_2(\vx)\to\bZ$ given by $n_\Fz(A)=\#\left(A\cap (\bH\times\{\Fz\})\right)$ is an isomorphism.
\end{Lem}
\proof
The $\beta$--circles are homologically linearly independent.
\endproof

For a given almost complex structure $J_{\beta,\beta'}$ on $W_{\beta,\beta'}$ we work with complex structures $J$ on $W_\beta$ satisfying the obvious analogs of (\textbf{J$'$1})--(\textbf{J$'$4}) from \fullref{Subsection:TriangleModuli1}.  For $A\in\pi_2(\vx)$, let $\cM^{A}$ denote the moduli space of embedded $J$--holomorphic curves in homology class $A$. 
\begin{Lem} For $A\in\pi_2(\vx)$, $\ind(A)=2n_{\Fz}(A)$.
\end{Lem}
\proof
Let $\gamma_1,\cdots,\gamma_g$ be curves in $\Sigma$ such that $(\Sigma,\vb,\vc)$ and $(\Sigma,\vb',\vc)$ are both the standard genus $g$ Heegaard diagram for $S^3$.  Let $\vy$ be the intersection point between $\vb$ and $\vc$, and $\vz$ the intersection point between $\vb'$ and $\vc$.  Let $T_i$ be a small triangle connecting $\vx$, $\vy$ and $\vz$, as in \fullref{Figure:1gons}.
\begin{centering}
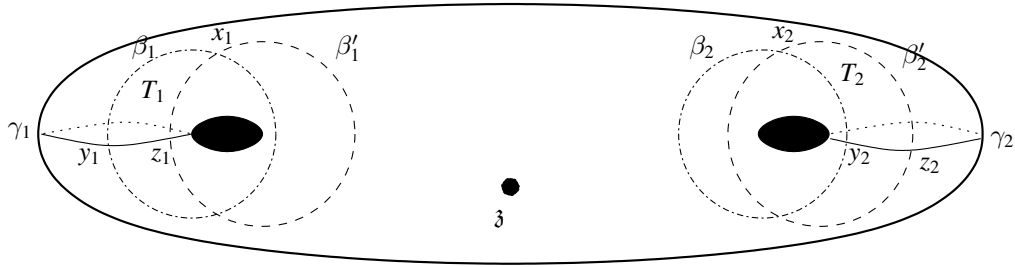
\begin{figure}
\begin{picture}(0,0)%
\includegraphics[scale=0.9]{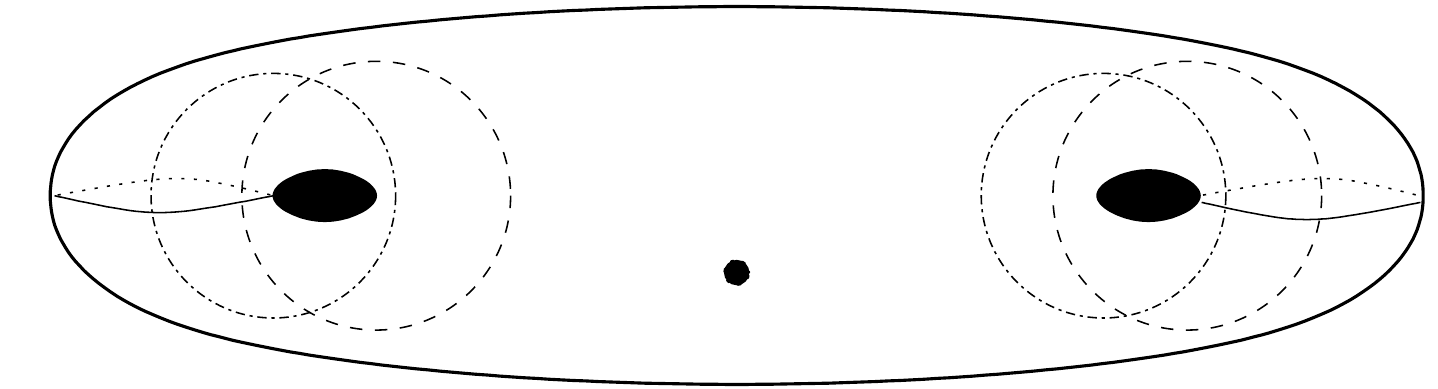}%
\end{picture}%
\setlength{\unitlength}{3552sp}%
\begingroup\makeatletter\gdef\SetFigFont#1#2#3#4#5{\small}%
\endgroup%
\begin{picture}(6887,1858)(511,-1365)
\put(1526,-621){\makebox(0,0)[lb]{\smash{{\SetFigFont{10}{12.0}{\rmdefault}{\mddefault}{\updefault}{\color[rgb]{0,0,0}\(z_1\)}%
}}}}
\put(1391,120){\makebox(0,0)[lb]{\smash{{\SetFigFont{10}{12.0}{\rmdefault}{\mddefault}{\updefault}{\color[rgb]{0,0,0}\(\beta_1\)}%
}}}}
\put(5283,120){\makebox(0,0)[lb]{\smash{{\SetFigFont{10}{12.0}{\rmdefault}{\mddefault}{\updefault}{\color[rgb]{0,0,0}\(\beta_2\)}%
}}}}
\put(6766, 59){\makebox(0,0)[lb]{\smash{{\SetFigFont{10}{12.0}{\rmdefault}{\mddefault}{\updefault}{\color[rgb]{0,0,0}\(\beta'_2\)}%
}}}}
\put(526,-437){\makebox(0,0)[lb]{\smash{{\SetFigFont{10}{12.0}{\rmdefault}{\mddefault}{\updefault}{\color[rgb]{0,0,0}\(\gamma_1\)}%
}}}}
\put(1950,220){\makebox(0,0)[lb]{\smash{{\SetFigFont{10}{12.0}{\rmdefault}{\mddefault}{\updefault}{\color[rgb]{0,0,0}\(x_1\)}%
}}}}
\put(5863,231){\makebox(0,0)[lb]{\smash{{\SetFigFont{10}{12.0}{\rmdefault}{\mddefault}{\updefault}{\color[rgb]{0,0,0}\(x_2\)}%
}}}}
\put(6395,-621){\makebox(0,0)[lb]{\smash{{\SetFigFont{10}{12.0}{\rmdefault}{\mddefault}{\updefault}{\color[rgb]{0,0,0}\(y_2\)}%
}}}}
\put(6889,-681){\makebox(0,0)[lb]{\smash{{\SetFigFont{10}{12.0}{\rmdefault}{\mddefault}{\updefault}{\color[rgb]{0,0,0}\(z_2\)}%
}}}}
\put(1452,-188){\makebox(0,0)[lb]{\smash{{\SetFigFont{10}{12.0}{\rmdefault}{\mddefault}{\updefault}{\color[rgb]{0,0,0}\(T_1\)}%
}}}}
\put(6334,-66){\makebox(0,0)[lb]{\smash{{\SetFigFont{10}{12.0}{\rmdefault}{\mddefault}{\updefault}{\color[rgb]{0,0,0}\(T_2\)}%
}}}}
\put(3924,-1053){\makebox(0,0)[lb]{\smash{{\SetFigFont{10}{12.0}{\rmdefault}{\mddefault}{\updefault}{\color[rgb]{0,0,0}\(\Fz\)}%
}}}}
\put(2811,120){\makebox(0,0)[lb]{\smash{{\SetFigFont{10}{12.0}{\rmdefault}{\mddefault}{\updefault}{\color[rgb]{0,0,0}\(\beta'_1\)}%
}}}}
\put(7383,-498){\makebox(0,0)[lb]{\smash{{\SetFigFont{10}{12.0}{\rmdefault}{\mddefault}{\updefault}{\color[rgb]{0,0,0}\(\gamma_2\)}%
}}}}
\put(1020,-621){\makebox(0,0)[lb]{\smash{{\SetFigFont{10}{12.0}{\rmdefault}{\mddefault}{\updefault}{\color[rgb]{0,0,0}\(y_1\)}%
}}}}
\end{picture}%
\caption{Case $g=2$}
\label{Figure:1gons}
\end{figure}
\end{centering}

Then, $T_1+\cdots+T_g\in\pi_2(\vx,\vy,\vz)$ has $\ind(T_1+\cdots+T_g)=0$.  Further, for $A\in\pi_2(\vx)$, $A+T_1+\cdots+T_g\in\pi_2(\vy,\vz)$ has $\ind(A+T_1+\cdots+T_g)=2n_{\Fz}(A)$.  The result follows by additivity of the index under gluings.
\endproof

\begin{Lem}\label{Lemma:1gons}
There is a choice of $C_\beta$ such that for any split complex structure on $W_\beta$, if $\vx\neq\vec{\theta}_{\beta,\beta'}$ then for $A\in\pi_2(\vx)$ with $\ind(A)=2n_{\Fz}(A)=0$, $\cM^A=\emptyset$.
\end{Lem}
\proof
Note that if one understands symplectic field theory in the Morse--Bott case this is easy to prove by taking $\beta'\to\beta$.  To avoid introducing this machinery we will instead give a somewhat perturbed argument.

\begin{figure}
\centering
\begin{picture}(0,0)%
\includegraphics[scale=0.7]{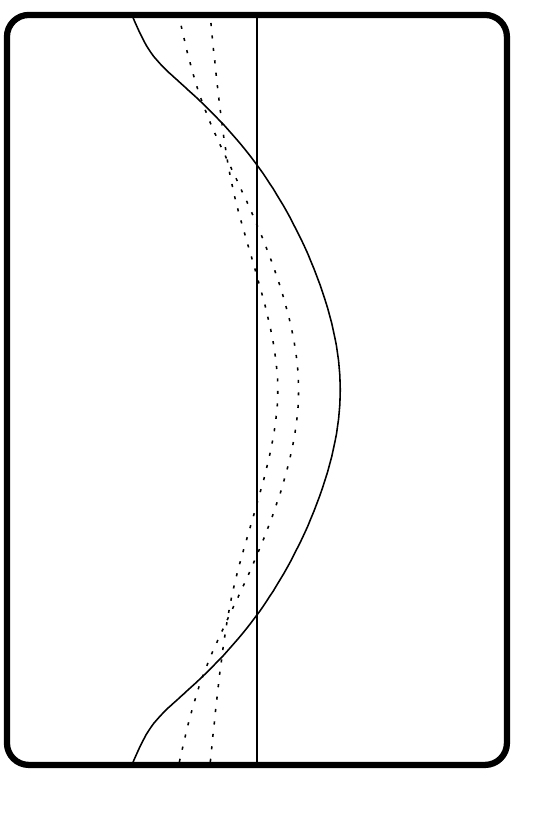}%
\end{picture}%
\setlength{\unitlength}{2763sp}%
\begingroup\makeatletter\ifx\SetFigFont\undefined%
\gdef\SetFigFont#1#2#3#4#5{%
  \reset@font\fontsize{#1}{#2pt}%
  \fontfamily{#3}\fontseries{#4}\fontshape{#5}%
  \selectfont}%
\fi\endgroup%
\begin{picture}(2583,3948)(1168,-4240)
\put(2876,-2236){\makebox(0,0)[lb]{\smash{{\SetFigFont{12}{14.4}{\rmdefault}{\mddefault}{\updefault}{\color[rgb]{0,0,0}\(\beta'_i\)}%
}}}}
\put(3751,-4036){\makebox(0,0)[lb]{\smash{{\SetFigFont{12}{14.4}{\rmdefault}{\mddefault}{\updefault}{\color[rgb]{0,0,0}\(\theta=0\)}%
}}}}
\put(3751,-436){\makebox(0,0)[lb]{\smash{{\SetFigFont{12}{14.4}{\rmdefault}{\mddefault}{\updefault}{\color[rgb]{0,0,0}\(\theta=2\pi\)}%
}}}}
\put(3751,-2236){\makebox(0,0)[lb]{\smash{{\SetFigFont{12}{14.4}{\rmdefault}{\mddefault}{\updefault}{\color[rgb]{0,0,0}\(\theta=\pi\)}%
}}}}
\put(2251,-4186){\makebox(0,0)[lb]{\smash{{\SetFigFont{12}{14.4}{\rmdefault}{\mddefault}{\updefault}{\color[rgb]{0,0,0}\(r=0\)}%
}}}}
\put(2026,-2236){\makebox(0,0)[lb]{\smash{{\SetFigFont{12}{14.4}{\rmdefault}{\mddefault}{\updefault}{\color[rgb]{0,0,0}\(\beta_i\)}%
}}}}
\end{picture}%
\caption{A Hamiltonian isotopy of $\beta_i$}
\label{IsotopyFig1v2}
\end{figure}
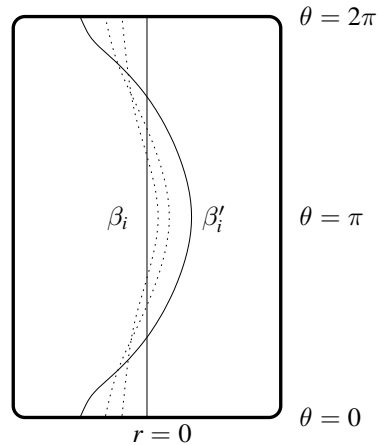

Let $t$ denote the first coordinate under the identification $\bH=(-\infty,\infty)\times[0,\infty)$.
We assume that for $t\in[-1/4,1/4]$,
$C_{\beta,\beta'}$ is the graph of an exact Hamiltonian isotopy of
the following form.  Fix  $\delta$ with $0<\delta<\pi$.  Identify a neighborhood in $\Sigma$ of each
$\beta_i$ with $S^1\times(-2\delta,2\delta)$.  Let
$\theta_i\in[0,2\pi)$ and $r_i\in(-\delta,\delta)$ be coordinates on
the $i^{th}$ neighborhood. 
Fix a bump function $b(\theta)$ on 
the circle which
is $\delta$ on $[\delta,\pi-\delta]$, $0$ on the interval
$[\pi+\delta,2\pi-\delta]$, and monotone on $[2\pi-\delta,\delta]$ and
$[\pi-\delta, \pi+\delta]$.
For some fixed
collection of constants $C_i$, $i=1,\cdots,g$,
consider the Hamiltonian $H$ given
by $C_i+\epsilon\left(\sin(\theta_i)+b(\theta)r_i\right)$ on the $i^{th}$
neighborhood, and extended arbitrarily outside the neighborhoods of
the $\beta_i$.  Here, we choose $\epsilon$ small enough that for
each $i$ the
graph with respect to $\beta_i$ of $H$ is contained in the chosen
neighborhood of
$\beta_i$ up to time $1$.  Then for each $i$,
$\beta_i'$ is the
time 1 graph with respect to $\beta_i$ of the Hamiltonian isotopy
specified by $H$.
  Thus, the Hamiltonian isotopy of the $\beta$--curves looks like
\fullref{IsotopyFig1v2}.

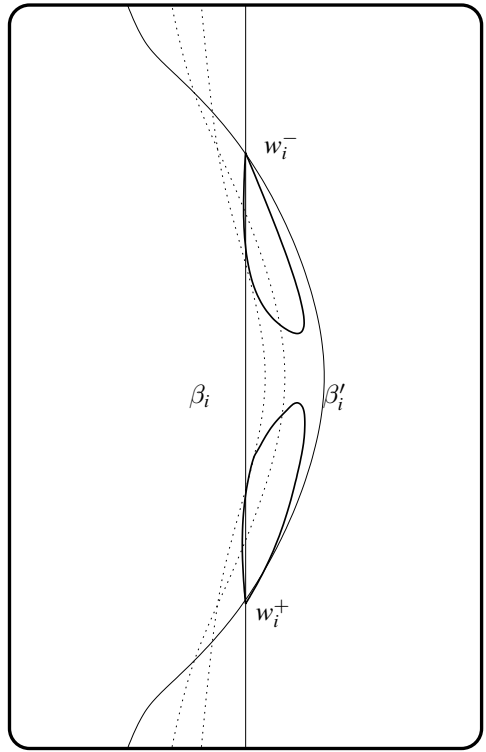
\begin{figure}
\centering
\begin{picture}(0,0)%
\includegraphics[scale=0.7]{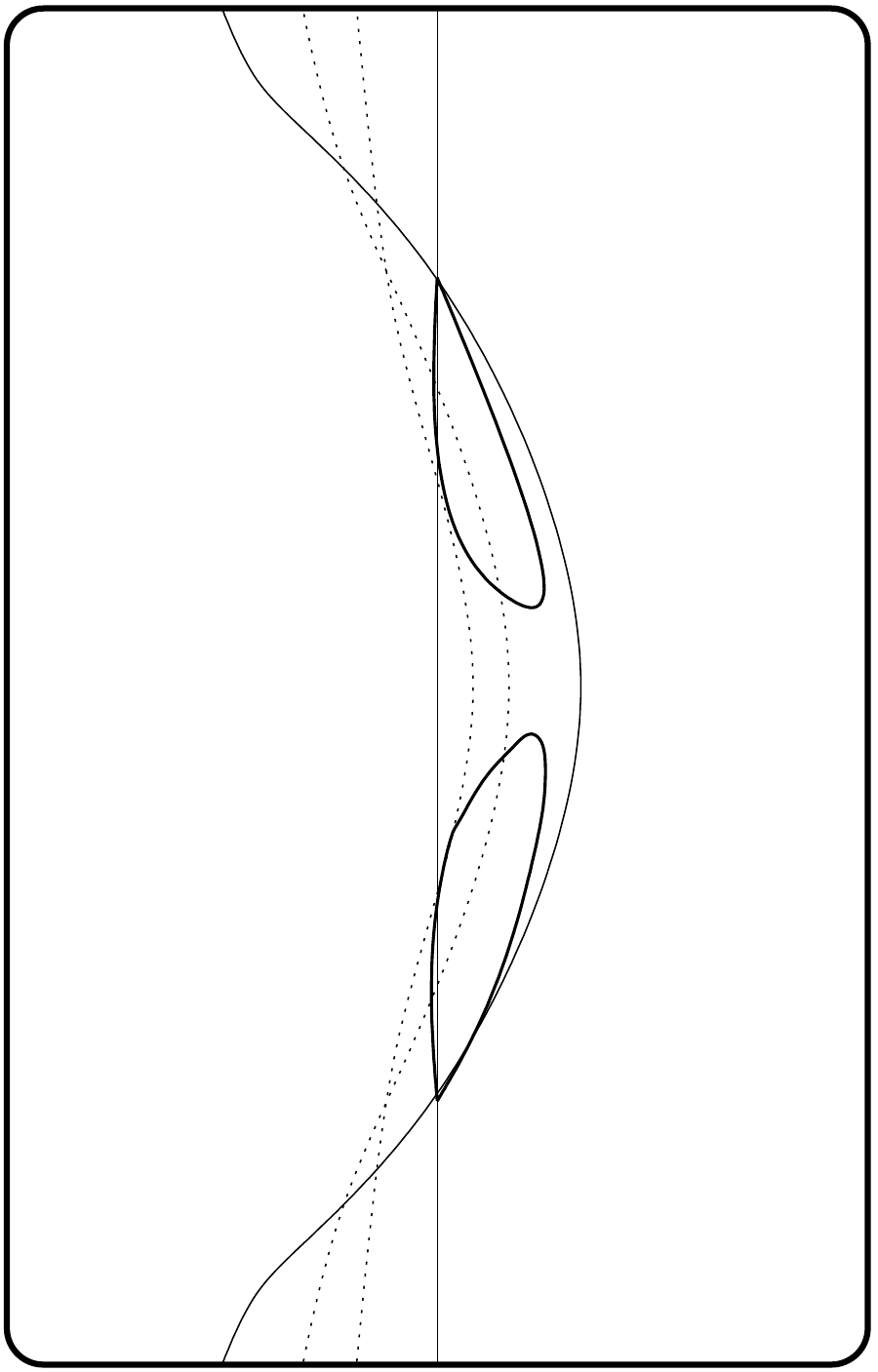}%
\end{picture}%
\setlength{\unitlength}{2763sp}%
\begingroup\makeatletter\ifx\SetFigFont\undefined%
\gdef\SetFigFont#1#2#3#4#5{%
  \reset@font\fontsize{#1}{#2pt}%
  \fontfamily{#3}\fontseries{#4}\fontshape{#5}%
  \selectfont}%
\fi\endgroup%
\begin{picture}(4302,6742)(1150,-7018)
\put(3976,-3886){\makebox(0,0)[lb]{\smash{{\SetFigFont{12}{14.4}{\rmdefault}{\mddefault}{\updefault}{\color[rgb]{0,0,0}\(\beta'_i\)}%
}}}}
\put(2776,-3886){\makebox(0,0)[lb]{\smash{{\SetFigFont{12}{14.4}{\rmdefault}{\mddefault}{\updefault}{\color[rgb]{0,0,0}\(\beta_i\)}%
}}}}
\put(3376,-5836){\makebox(0,0)[lb]{\smash{{\SetFigFont{12}{14.4}{\rmdefault}{\mddefault}{\updefault}{\color[rgb]{0,0,0}\(w_i^+\)}%
}}}}
\put(3451,-1636){\makebox(0,0)[lb]{\smash{{\SetFigFont{12}{14.4}{\rmdefault}{\mddefault}{\updefault}{\color[rgb]{0,0,0}\(w_i^-\)}%
}}}}
\end{picture}%
\caption{Two possible paths $\pi_\Sigma\circ u|_{\bdy S}$ (shown in bold)}
\label{IsotopyFig2v2}
\end{figure}

Let $A$  be as in the statement of the proposition and $u\co S\to
W_\beta$ be a holomorphic representative of $A$.  Note
that $S$ is a disjoint union of $g$ disks, and
$\pi_\Sigma\circ u(\bdy S)$ is a collection of $g$ pairwise disjoint
simple closed curves in $\bigcup_{\theta\in S^1}C_\theta$.  Now, by
considering orientations and the path taken by $\pi_\Sigma\circ
u|_{\bdy S}$, we see that $u$ can not exist unless
$\vx=\vec{\theta}_{\beta,\beta'}$.  See \fullref{IsotopyFig2v2}.
This proves the lemma.
\endproof

\proof[Proof of \fullref{Prop:HSB}]
We will write the proof in the notation for $HF^\infty$, but the proofs for the other theories are just the same.

As usual, we let $t$ denote the $\bR$--coordinate
on $W$.  Let
$C_\alpha=\ba\times\{1\}\times\bR\subset W$, and let $C_{\beta,\beta'}$ be
a $g$--tuple of Lagrangian cylinders which agree with
$\bb\times\{0\}\times\bR$ for $t<-\frac{1}{4}$ and with
$\bb'\times\{0\}\times\bR$ for $t>\frac{1}{4}$.

For a given complex structure $J$ on
$W$, let $J_R$ denote the complex structure obtained from $J$ by
inserting a neck of length $R$ along the hypersurface
$\Sigma\times\{|z|=1/2,\Re(z)>0\}\subset\Sigma\times[0,1]\times\bR$.  (See \fullref{Subsection:Splitting} for a discussion of the splitting process.)

Let $\Phi_R\co CF^\infty(\Sigma,\va,\vb,\Fz)\to
CF^\infty(\Sigma,\va,\vd,\Fz)$ be the map defined in
\fullref{Section:Isotopy},
with respect to the complex structure $J_R$.  We showed in
\fullref{Section:Isotopy} that $\Phi_R$ is an isomorphism for
each $(J,R)$ such that $J_R$ achieves transversality.

Taking the limit $R\to\infty$, $W$ splits into two spaces.  One is
the space we denoted $W_{\alpha,\beta,\beta'}$ in 
\fullref{Section:Triangles}.  The other is a copy of $W_\beta$.

Choose $J$ so that: 
\begin{itemize}
\item The complex structure $J_{\alpha,\beta,\beta'}$
induced on $W_{\alpha,\beta,\beta'}$ satisfies
(\textbf{J$'$2}) (the other conditions in the triangles section are automatic).
\item The complex structure $J_{\alpha,\beta,\beta'}$ and the
  complex structure $J_\beta$ induced on $W_\beta$ achieve
  transversality for curves of index $0$ (and hence so does $J_R$ for large $R$).
\end{itemize}

For such a $J$ it follows from \fullref{Gluing:Split} that, for $R$ large, 
$$\Phi^\infty_R([\vx,i])=\sum_{\vz}\sum_{\substack{A\in\pi_2(\vx)\\B\in\pi_2(\vx,\vy,\vz)\\\ind(A+B)=0}}\#\left(\cM^A\right)\#\left(\cM^B\right)[\vy,i-n_\Fz(A+B)]
$$
where the first sum is over all intersection points between $\vb$ and $\vb'$.

By additivity of the index, $\#\left(\cM^A\right)\#\left(\cM^B\right)=0$ unless both $A$ and $B$ have index $0$.  So, by \fullref{Lemma:1gons}, it follows that
$$\Phi^\infty_R=(\#\cM(\vec{\theta}_{\beta,\beta'}))\cdot F_{\alpha,\beta,\delta}(\cdot\otimes\theta_{\beta,\beta'})
$$
where $\cM(\vec{\theta}_{\beta,\beta'})$ denotes the index $0$ holomorphic maps in $W_{\beta}$ asymptotic to $\vec{\theta}_{\beta,\beta'}$.  

It only remains to understand $\cM(\vec{\theta}_{\beta,\beta'})$.  We do this by an
indirect argument.  Suppose that $(\Sigma,\va,\vb)$ is the
standard Heegaard diagram for $S^3$.  Then
$\widehat{CF}(\Sigma,\va,\vb)=\widehat{CF}(\Sigma,\va,\vb')=\bZ$ and the
boundary maps are trivial.  So, for
$\Phi_R$ to be an isomorphism it must be multiplication by $\pm
1$.  Further, by explicitly counting triangles (there is only one) we
know that $F_{\alpha,\beta,\beta'}$ is the identity map.  It follows
that $\#\cM(\vec{\theta}_{\beta,\beta'})=\pm1$.  This finishes the proof.
\endproof

\subsection{Proof of Proposition \ref{Prop:HandleslideInvariance}}

\begin{Lem}\label{Lem:HSAssoc}
In the pointed Heegaard quadruple--diagram
  $(\Sigma,\va,\vb,\vb^H,\vb',\Fz)$, both $\delta
  H^1(Y_{\beta,\beta'})$ and $\delta H^1(Y_{\alpha,\beta^H})$ are
  the trivial subgroup of
  $H^2(X_{\alpha,\beta,\beta^H,\beta'},\bdy
  X_{\alpha,\beta,\beta^H,\beta'})$.
\end{Lem}
\proof
(Compare the proof of~\cite[Lemma 9.6]{OS1})

The following diagram commutes:

{\disablesubscriptcorrection

\cl{\scalebox{0.89}{\small
$\phantom{9}
\xymatrix{
H^2(X_{\alpha,\beta,\beta^H,\beta'},\bdy X_{\alpha,\beta,\beta^H,\beta'})\ar@{<->}[r]^(.6)\cong_(.6){P.D.} &
H_2(X_{\alpha,\beta,\beta^H,\beta'})\ar@{<->}[r]^(.31)\cong & 
\ker\left(
\left({\begin{array}{c}
H_1(\ba) \oplus H_1(\bb) 
\\\oplus\\ 
H_1(\bb^H) \oplus H_1(\bb')
\end{array}}\right)
\to H_1(\Sigma)
\right)\\
H^1(\bdy X_{\alpha,\beta,\beta^H,\beta'})\ar[u]\ar@{<->}[r]^(.53)\cong_(.53){P.D.} &
H_2(\bdy X_{\alpha,\beta,\beta^H,\beta'})\ar[u]\ar@{<->}[r]^(.35)\cong &
\left({
\begin{array}{c}
\ker\left(H_1(\ba)\oplus H_1(\bb)\to H_1(\Sigma)\right)\\
\oplus\\
\ker\left(H_1(\bb)\oplus H_1(\bb^H)\to H_1(\Sigma)\right)\\ 
\oplus\\
\ker\left(H_1(\bb^H)\oplus H_1(\bb')\to H_1(\Sigma)\right)\\ 
\oplus\\
\ker\left(H_1(\bb')\oplus H_1(\ba)\to H_1(\Sigma)\right)
\end{array}}\right)
\ar@{>>}[u]
}
$
}}}

It is clear that the right hand vertical map is surjective.  So, from the long exact sequence of the pair $(X_{\alpha,\beta,\beta^H,\beta'},\bdy X_{\alpha,\beta,\beta^H,\beta'})$, the map
$H_2(X_{\alpha,\beta,\beta^H,\beta'})\to H_1(\bdy X_{\alpha,\beta,\beta^H,\beta'})\cong H^2(\bdy
X_{\alpha,\beta,\beta^H,\beta'})$ is injective.  But the image of $\delta
H^1(Y_{\beta,\beta'})$ or $\delta H^1(Y_{\alpha,\beta^H})$ in the
latter group is obviously $0$.
\endproof

\begin{Cor}{\rm(Compare~\cite[Lemma 9.6]{OS1})}
\label{Cor:HSCor}
{\small$$
\hat{F}_{\alpha,\beta^H,\beta'}\left(\hat{F}_{\alpha,\beta,\beta^H}(\cdot\otimes\theta_{\beta,\beta^H})\otimes\theta_{\beta^H,\beta'}\right)
= \hat{F}_{\alpha,\beta,\beta'}\left(\cdot\otimes
  \hat{F}_{\beta,\beta^H,\beta'}(\theta_{\beta,\beta^H}\otimes\theta_{\beta^H,\beta'})\right).
$$}%
Similar statements hold for $F^\infty$, $F^+$ and $F^-$ with
$[\theta,0]$ (appropriately subscripted) in place of $\theta$.
\end{Cor}
\proof
From the previous lemma, the $\delta H^1(Y_{\beta,\beta'})+\delta
H^1(Y_{\alpha,\beta^H})$--orbit of $\Spin^\bC$--structures which
restrict to $\Ss$ has just one element.  Further, strong
admissibility for the quadruple is equivalent to strong admissibility
of the six Heegaard diagrams
$(\Sigma,\va,\vb,\Fz),\cdots,(\Sigma,\vb^H,\vb',\Fz)$, which in turn
follows from admissibility for\break $(\Sigma,\va,\vb,\Fz)$.  So, the
result follows from \fullref{AssocProp}.
\endproof

\begin{Lem}\label{Lem:HSOrientations}
Let $\frako_{\alpha,\beta}$ be a coherent orientation system
  for $(\Sigma,\va,\vb,\Fz)$.  
  Then there is a coherent orientation
  system for the Heegaard quadruple--diagram $(\Sigma,\va,\vb,\vc,\vd,\Fz)$
  extending $\frako_{\alpha,\beta}$ and the orientation systems
  constructed in \fullref{Prop:HSC}.
\end{Lem}
(In fact, it is not hard to see that the coherent orientation system
for the Heegaard quadruple--diagram is essentially unique (compare
\fullref{Prop:OrientationIndependence}).)

\proof
The
orientation systems $\frako_{\alpha,\beta}$ and $\frako_{\beta,\beta^H}$
extend to orientation systems $\frako_{\alpha,\beta,\beta^H}$ and
$\frako_{\alpha,\beta^H}$ by \fullref{Lemma:TriangleOrientations}.
For $\vb'$ close to $\vb$, $\frako_{\alpha,\beta}$ induces
$\frako_{\alpha,\beta'}$, which extends to an orientation system
$\frako_{\alpha,\beta',\beta^H}$.  It is clear that we can perform the
extensions so that the induced orientations
$\frako_{\beta',\beta^H}$ and $\frako_{\beta,\beta^H}$ are consistent
with $\frako_{\beta,\beta^H,\beta'}$.  The result follows.
\endproof

\begin{Lem}\label{Lemma:HSActionsCommute}
The map $F^\infty_{\alpha,\beta,\beta^H}(\cdot\otimes[\theta_{\beta,\beta^H},0])$ commutes with the actions of $U$ and $H_1(Y_{\alpha,\beta})/Tors$.  Similar statements hold for $F^+$, $F^-$ and $\hat{F}$.
\end{Lem}
\proof
The fact that
$F_{\alpha,\beta,\beta^H}(\cdot\otimes[\vec{\theta}_{\beta,\beta^H},0])$
commutes with the $U$--action is obvious.  We check that it commutes with the $H_1/Tors$--action.

Fix a homology class $\zeta\in H_1(Y)$.  Choose a knot $K\mapsinto\Sigma$
representing $\zeta$ and meeting $\ba$ transversely.  Let
$\cM_1$ denote the
space of Riemann surfaces with $2g$ boundary punctures (as usual)
and one additional marked point $p$ on the boundary.  We give
yet another definition of the map $A_\zeta$ for a given
$(\Sigma,\va,\vb,\Fz)$.  Given a
homology class $B\in\pi_2(\vx,\vy)$, let
$\cM_K^B$ denote the moduli space consisting of holomorphic maps
from surfaces $S\in\cM_1$ to $W_{\alpha,\beta}$ in the homology class $B$
mapping $p$ to 
$K\times\{1\}\times\bR$, and $\hcM_K^B=\cM_K^B/\bR$.  (Note that
this definition is slightly different in form from the one we used
to prove \fullref{ActionProp} and \fullref{Isotopy:PreserveStructure}.  This definition is convenient here, but either would work.)

Now, for $A\in\pi_2(\vx,\vy,\vz)$, let $\cM_K^A$ denote the space of
holomorphic maps from surfaces $S\in\cM_1$ to
$W_{\alpha,\beta,\beta^H}$ in the homology class $A$ mapping $p$
to $K\times e_1$.  (Recall that $e_1$ is the edge of the triangle $T$ corresponding to the $\alpha$--circles.)  Define
a map $H\co CF^\infty(Y_{\alpha,\beta})\otimes
CF^\infty(Y_{\beta,\beta^H})\to CF^\infty(Y_{\alpha,\beta^H})$ by
$$
H([\vx,i]\otimes[\vy,j])=\sum_{\vz}\sum_{\substack{A\in\pi_2(\vx,\vy,\vz)\\\ind(A)=0}}
\#\left(\cM_K^A\right)[\vz,i+j-n_\Fz(A)].
$$

For $A\in\pi_2(\vx,\vy,\vz)$, $\ind(A)=1$, consider the space
$\cM_K^A$.  This is a 
one--dimensional manifold the ends of which are height two holomorphic
buildings with an index $0$ holomorphic curve in
$W_{\alpha,\beta,\beta^H}$ and an index $1$ holomorphic curve
(defined up to translation) in one of $W_{\alpha,\beta}$,
$W_{\beta,\beta^H}$ or $W_{\alpha,\beta^H}$.  For each of these
buildings, either
\begin{enumerate}
\item $p$ is mapped to $W_{\alpha,\beta,\beta^H}$ or
\item $p$ is mapped to $W_{\alpha,\beta}$ or
\item $p$ is mapped to $W_{\alpha,\beta^H}$.
\end{enumerate}
Summing the number of ends over the different choices of $A$ with
$n_\Fz(A)=k$ we obtain, in each case respectively,
\begin{enumerate}
\item the coefficient of $[\vz,i+j-k]$ in $(\bdy\circ H +
  H\circ\bdy)([\vx,i],[\vy,j])$
\item the coefficient of $[\vz,i+j-k]$ in
  $f^\infty_{\alpha,\beta,\beta^H}(A_\zeta([\vx,i]),[\vy,j])$ or
\item the coefficient of $[\vz,i+j-k]$ in
  $A_\zeta(f^\infty_{\alpha,\beta,\beta^H}([\vx,i],[\vy,j]))$.
\end{enumerate}
This proves that on the level of homology the
$f^\infty_{\alpha,\beta,\beta^H}$ commute with $A_\zeta$.  As
usual, the result for $HF^+$ and $HF^{\leq0}$ follow; the result
for $\widehat{HF}$ is proved in an analogous way.
\endproof

\proof[Proof of \fullref{Prop:HandleslideInvariance}]
Having proved \fullref{Cor:HSCor}, \fullref{Prop:HSB} and \fullref{Prop:HSC}, we are completely justified in writing
\begin{eqnarray}
\hat{F}_{\alpha,\beta^H\!,\beta'}\left(\hat{F}_{\alpha,\beta,\beta^H}\!(\cdot\otimes\vec{\theta}_{\beta,\beta^H})\otimes\vec{\theta}_{\beta^H,\beta'}\right) & = & \hat{F}_{\alpha,\beta,\beta'}\left(\cdot\otimes \hat{F}_{\beta,\beta^H,\beta'}( \vec{\theta}_{\beta,\beta^H}\otimes\vec{\theta}_{\beta^H,\beta'})\right)\nonumber\\
&=& \hat{F}_{\alpha,\beta,\beta'}(\cdot\otimes\vec{\theta}_{\beta,\beta'})\nonumber\\
&=&\hat{\Phi}_{\beta,\beta'}(\cdot)\nonumber
\end{eqnarray}
where $\hat{\Phi}_{\beta,\beta'}$ is the isomorphism induced by the isotopy from $\vb$ to $\vb'$.

It follows that $\hat{F}_{\alpha,\beta^H,\beta'}$ is surjective and $\hat{F}_{\alpha,\beta,\beta^H}$ is injective.  The same argument with the roles of $\beta$ and $\beta^H$ exchanged shows that $\hat{F}_{\alpha,\beta,\beta^H}$ is surjective.  It then follows from \fullref{Lemma:HSActionsCommute}, that $\hat{F}_{\alpha,\beta,\beta^H}$ is an isomorphism of $\bZ[U]\otimes_\bZ\Lambda^*H^1(Y)/Tors$--modules, proving that $\widehat{HF}$ is invariant under handleslides.  The proofs for $HF^-$, $HF^+$ and $HF^\infty$ are just the same.
\endproof

\section{Stabilization}
\label{Section:Stabilization}
We show that the homology groups defined in \fullref{Section:ChainComplexes} are invariant under
stabilization, or equivalently under taking connected sum of the
Heegaard diagram $(\Sigma,\va,\vb,\Fz)$with the standard Heegaard
diagram for $S^3$.

Because of handleslide invariance, it is enough to prove the result if
we take the connected sum at the point $\Fz$.    Let
$(\Sigma,\va,\vb,\Fz)$ denote the original (pointed) Heegaard diagram,
$(\cT,\alpha_{g+1},\beta_{g+1})$ the standard Heegaard diagram for
$S^3$, and
$(\Sigma',\va',\vb',\Fz')$ a Heegaard diagram obtained by taking a
connect sum near $\Fz$.

The result for
the hat theories is quite easy.  There is an identification between
intersection points in $\Sigma$ and intersection points in
$\Sigma'$, identifying an intersection point
$\vx=\{x_1,\cdots,x_g\}$ in $\Sigma$ with the intersection point
$\vx'=\{x_1,\cdots,x_g,\alpha_{g+1}\cap\beta_{g+1}\}$.  There is then
an obvious identification of $\pi_2^\Sigma(\vx,\vy)$ with
$\pi_2^{\Sigma'}(\vx',\vy')$.  Fix a homology
class $A\in\pi_2^\Sigma(\vx,\vy)$ with $n_z(A)=0$.  There is
an inclusion map $\hcM^A_{\Sigma}\hookrightarrow\hcM^A_{\Sigma'}$
which takes any holomorphic curve $u\in\hcM^A_{\Sigma}$ to the
disjoint union $$u\coprod
\left((\alpha_{g+1}\cap\beta_{g+1})\times[0,1]\times\bR\right).$$  This map is
clearly onto, and hence a homeomorphism.
With these
identifications, the chain complex $\widehat{CF}$ for the stabilized
diagram is isomorphic to the chain complexes for the unstabilized
one.

The other theories require more work.

Our strategy is to insert a longer and longer neck between the
original Heegaard diagram $\Sigma$ and the torus $\cT$ that has
been spliced in.  We will show that in the limit $(\textrm{neck
  length})\to\infty$, the moduli spaces of curves we consider,
correctly defined, are naturally identified with the moduli spaces for
$(\Sigma,\va,\vb,\Fz)$.
A gluing result then shows that the moduli spaces in the
limit can be identified with those of large neck length.
  Since the Floer homologies are independent
of complex structure, the invariance under stabilization will then
follow.

Actually, although everything we say in this section is correct, it is
somewhat dishonest: we leave most of the work for
\fullref{Gluing:Twins} in the appendix.  In particular,
because of transversality issues, to prove the necessary gluing result
we seem to need to work in a more harshly perturbed setting than we
use in this section.  Much of the work of the present section is
redone in the proof of \fullref{Gluing:Twins} in the more
general context.  For example, we prove there a stronger version of the
compactness result \fullref{Prop:TwinCompactness}.

I have chosen to write this section in this slightly dishonest way for
two reasons.  Firstly, perhaps someone else will see a simpler way to
correct the dishonesty (see the discussion at the beginning of the
proof of \fullref{Prop:TwinCompactness}).  Secondly, as
written this section is quite explicit, so one can actually see what
is going on in the stabilization proof.

Returning to mathematics,
observe that from the formula for the index in
\fullref{Section:Index} we know
that for a homology class $A\in\pi_2^\Sigma(\vx,\vy)$ corresponding
to $A'\in\pi_2^{\Sigma'}(\vx,\vy)$, $\ind(A)=\ind(A')$.  Since we only
consider curves with index $1$ when computing Floer homologies, we shall
generally restrict to curves of index $1$.  This allows us to keep
our definitions simpler and our theorems true.

We make precise what we mean by stretching the neck.  Fix complex
structures $j_\Sigma$ on $\Sigma$ and
$j_\cT$ on $\cT$.  Fix a point 
$\Fz_0\in \cT\setminus(\alpha_{g+1}\cup\beta_{g+1})$.  Choose small disks
$D_\Sigma\supset D'_\Sigma\ni \Fz$ and $D_\cT\supset D'_\cT\ni \Fz_0$ so
that $D_\Sigma\setminus D'_\Sigma$ is conformally identified with
$S^1\times[-1,0]$ and $D_\cT\setminus D'_\cT$ is conformally
identified with $S^1\times[0,1]$.  Here, $\bdy D'_\Sigma$ is
identified with $S^1\times\{-1\}$ and $\bdy D'_\cT$ is identified
with $S^1\times\{1\}$.  (The complex structure on a cylinder
  $S^1\times[a,b]$ is given by identifying the cylinder with
  $\{e^{r+i\theta}\in\bC | a\leq r\leq b\}$ under the identification
$(\theta,t)\mapsto e^{t+i\theta}$.)

Let $j_R$ be the
complex structure on 
$$\Sigma'_R=(\Sigma\setminus D'_\Sigma)\bigcup_{D_\Sigma\setminus D'_\Sigma
  \sim S^1\times[-R-1,-R]}S^1\times[-R-1,R+1]\bigcup_{S^1\times[R,R+1]\sim
  \bdy D_\cT}(\cT\setminus D_{\cT}) $$
induced by the complex structures on $\Sigma$, $\cT$ and
$S^1\times[-R-1,R+1]$, and $\omega_R$ the area form.  We refer to
$S^1\times[-R-1,R+1]$ as the neck, and denote it $N_R$.
Notice that in the limit $R\to\infty$, $j_R$ degenerates to the
complex structure $j_\infty=j_\Sigma\vee j_\cT$ on $\Sigma\vee \cT$.

Fix a complex structure $J$ on
$W'_0=\Sigma'_0\times[0,1]\times\bR$ satisfying
(\textbf{J1})--(\textbf{J5}), which is split (ie,
$j_0\times j_\bD$) near $N_0$.
Let $J_R$ be the complex structure on
$W'_R=\Sigma'_R\times[0,1]\times\bR$ which agrees with $J$ outside
$N_R\times[0,1]\times\bR$ and with $j_R\times j_\bD$ on
$N_R\times[0,1]\times\bR$.  Note that $J_R$
converges to a complex structure $J_\infty$ on $W'_{\infty}=(\Sigma\vee
\cT)\times[0,1]\times\bR$.  The space
$W=\Sigma\times[0,1]\times\bR$ lies inside $W'_\infty$ in an
obvious way.  We choose $J$ so that
\begin{enumerate}
\item There is a collection of $R_i\to\infty$ such that $J_{R_i}$
  achieves transversality for holomorphic curves in $W'_{R_i}$
  satisfying (\textbf{M0})--(\textbf{M6}).
\item The restriction of $J_\infty$ to $W$ achieves transversality
  for holomorphic curves in $W$ satisfying (\textbf{M0})--(\textbf{M6}).
\end{enumerate}
Since by choosing $J$ appropriately we can end up with any 
  $J_\infty$ which is split near $\Fz$, we can always find such a
  $J$.

For a generic choice of $\Fz$ and for any holomorphic curve $u$ with
$\ind(u)=1$,\break $\pi_\bD\left((\pi_{\Sigma}\circ u)^{-1}(\Fz)\right)$
consists of $n_\Fz(u)$ distinct points.  We choose $\Fz$ satisfying
this condition.

The holomorphic curves that we consider in $W'_{\infty}$ are
\emph{holomorphic twin towers}.  That is, fix a homology class
$A\in\pi_2^\Sigma(\vx,\vy)\cong\pi_2^{\Sigma'}(\vx',\vy')$.  

By a holomorphic twin tower $u$ in the homology class $A$ we mean a
collection of holomorphic maps $(u_1,\cdots,u_n,v_1,\cdots,v_n)$ where
$u_1,\cdots,u_n\co S_1,\cdots,S_n\to W$ and $v_1,\cdots,v_n\co
S'_1,\cdots,S'_n\to T^2\times[0,1]\times\bR$ such that:
\begin{enumerate}
\item The $u_1,\cdots,u_n$ are holomorphic curves in $W$
  satisfying (\textbf{M0})--(\textbf{M6}).
\item The $S'_1,\cdots,S'_n$ are all closed surfaces.
\item There is a sequence of intersection points
  $\vx_1,\cdots,\vx_{n+1}$ such that $u_i$ connects $\vx_i$ to
  $\vx_{i+1}$.
\item $[u_1]+\cdots+[u_n]=A$.  Here, $[u_i]$ is the class in
  $\pi_2(\vx_i,\vx_{i+1})$ represented by $u_i$.
\item The $v_1,\cdots,v_n$ are holomorphic curves in
  $\cT\times[0,1]\times\bR$ each connecting $\alpha_{g+1}\cap\beta_{g+1}$
  to itself, satisfying (\textbf{M0}), (\textbf{M1}), (\textbf{M3}) and (\textbf{M5})
\item As sets, $\pi_\bD\left((\pi_\Sigma\circ
    u_i)^{-1}(\Fz)\right)=\pi_\bD\left((\pi_\Sigma\circ
    v_i)^{-1}(\Fz_0)\right)$ for each $i$.  (This is the matching condition
    on the horizontal levels.)
\end{enumerate}
\medskip\textbf{Remark}\qua If we had not restricted to holomorphic curves of index
$1$ and chosen $\Fz$ generically, we would need a slightly more
complicated definition of holomorphic twin towers.  In particular, we
would need to allow pieces of the curves to live in the ``horizontal''
cylinders $S^1\times\bR\times[0,1]\times\bR$.

\begin{Lem}\label{ShortTowers}If $\ind(A)=1$ then any holomorphic
  twin tower in the homology class $A$ has height one (ie, in
  the previous definition, $n=1$).
\end{Lem}
\proof
This is trivial:  the index adds between stories.
\endproof

\begin{Lem}If $(u_1,v_1)$ is a
  holomorphic twin tower in a homology class of index 1 then $v_1$
  consists of a trivial cylinder
  $\alpha_{g+1}\cap\beta_{g+1}\times[0,1]\times\bR$ and $n_\Fz(u_1)$
  horizontal tori, $T^2\times
  \pi_\bD\left((\pi_\Sigma\circ u_1)^{-1}(\Fz)\right)$.
\end{Lem}
\proof
The restriction of $\pi_\bD\circ v_1$ to the components on which it
is nonconstant is a $1$--fold covering of the disk.  Hence, it must
be a disk itself.  It follows easily that the restriction of
$\pi_\cT\circ v_1$ to this component must be constant.  Since the
condition on $z$ guarantees that the 
set $(\pi_\bD\circ v_1)\left((\pi_\cT\circ v_1)^{-1}(\Fz_0)\right)$
consists of $n_\Fz([u_1])=n_\Fz([v_1])$ distinct points, each component
on which $\pi_\bD\circ v_1$ is constant must be a copy of $T^2$,
and the restriction of $\pi_\cT\circ v_1$ to each component must be a
diffeomorphism.
\endproof

\begin{Cor}The moduli space of holomorphic twin towers in a given
  homology class $A\in\pi_2^\Sigma(\vx,\vy)$ with $\ind(A)=1$ is
  naturally identified with $\hcM^A_W$.
\end{Cor}
Here, $\hcM^A_W$ denotes the moduli space of holomorphic curves in
the homology class $A$ in $W$.

We now need to identify the space of holomorphic twin towers with
$\hcM^A_{W'_R}$ for $R$ large enough.  As usual, doing so requires
two steps:  compactness and gluing.  To avoid work we torture the
compactness argument slightly.

\begin{Prop}\label{Prop:TwinCompactness}
Fix a sequence $\{u_R\}$ of $J_R$--holomorphic curves in
  $W'$ with index $1$.  Then for $\Fz$ chosen generically, there is a
  subsequence of $\{u_R\}$ which converges to a holomorphic twin
  tower (of height one).  
\end{Prop}
\proof
The proof is in three stages.  First we use the fact that
$\pi_\bD\circ u_R$ is holomorphic to extract the vertical level structure and
conformal structure on the limit surface.  Then we cut $\Sigma$ into
two overlapping regions, and view the $u_R$ as maps into each of the
two regions with Lagrangian boundary conditions.  This allows us to
extract a convergent subsequence of maps to $W'$.

By classical symplectic field theory \cite[Theorem 10.1]{Ya2}, we can replace $u_R$ by a
subsequence so that $\pi_\bD\circ u_R$ converges to a holomorphic
building, which we denote $\pi_\bD\circ u_\infty$ (although
$u_\infty$ does not yet make sense on its own).  We can also assume
that for all $R$ the source of $u_R$
is some fixed topological manifold $S$.  Let $S_\infty$ denote the
source of $\pi_\bD\circ u_\infty$.

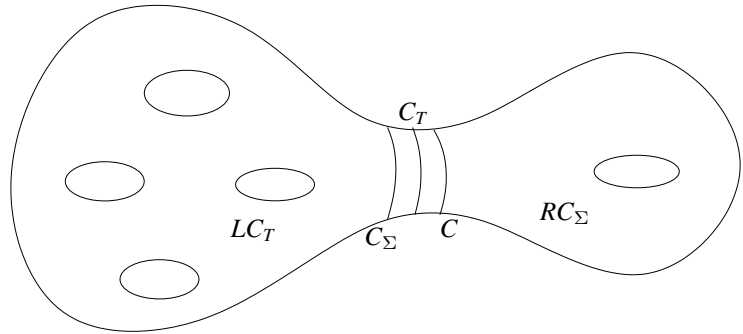
\begin{figure}
\centering
\begin{picture}(0,0)%
\includegraphics[scale=0.7]{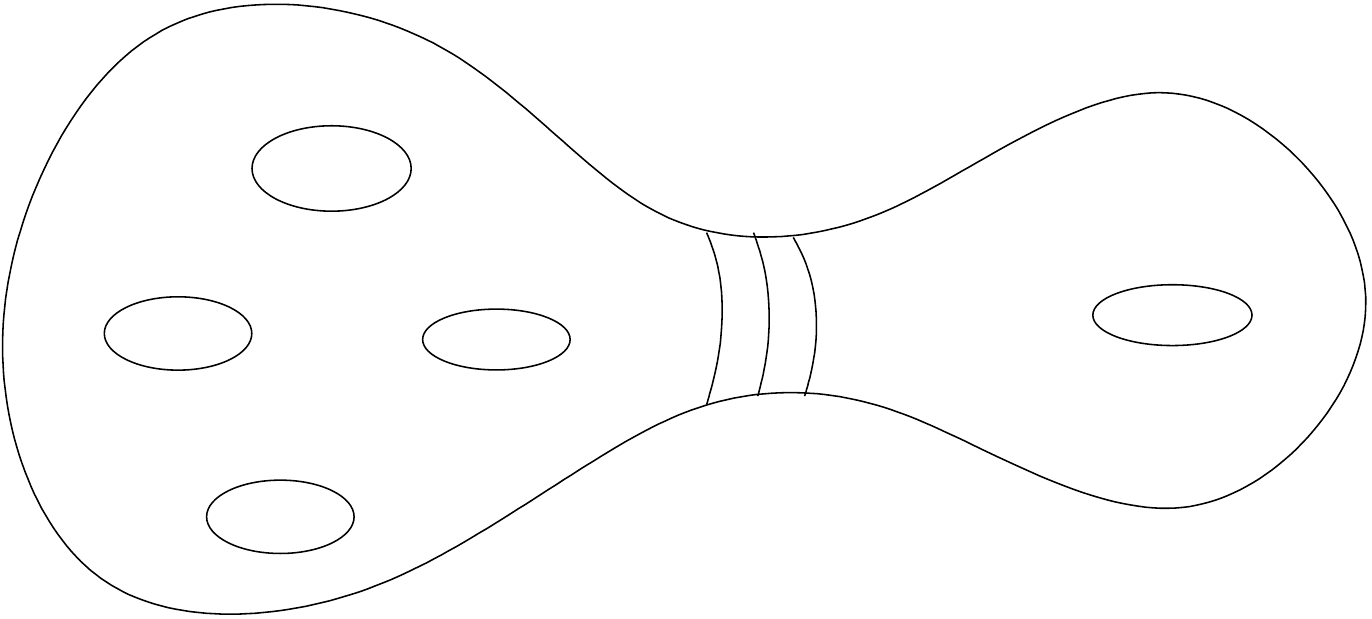}%
\end{picture}%
\setlength{\unitlength}{2763sp}%
\begingroup\makeatletter\ifx\SetFigFont\undefined%
\gdef\SetFigFont#1#2#3#4#5{%
  \reset@font\fontsize{#1}{#2pt}%
  \fontfamily{#3}\fontseries{#4}\fontshape{#5}%
  \selectfont}%
\fi\endgroup%
\begin{picture}(6567,2951)(1185,-3254)
\put(5926,-2236){\makebox(0,0)[lb]{\smash{\SetFigFont{12}{14.4}{\rmdefault}{\mddefault}{\updefault}{\color[rgb]{0,0,0}\(RC_\Sigma\)}%
}}}
\put(3151,-2386){\makebox(0,0)[lb]{\smash{\SetFigFont{12}{14.4}{\rmdefault}{\mddefault}{\updefault}{\color[rgb]{0,0,0}\(LC_T\)}%
}}}
\put(4651,-1336){\makebox(0,0)[lb]{\smash{\SetFigFont{9}{10.8}{\rmdefault}{\mddefault}{\updefault}{\color[rgb]{0,0,0}$C_T$}%
}}}
\put(4351,-2461){\makebox(0,0)[lb]{\smash{\SetFigFont{9}{10.8}{\rmdefault}{\mddefault}{\updefault}{\color[rgb]{0,0,0}$C_\Sigma$}%
}}}
\put(5026,-2386){\makebox(0,0)[lb]{\smash{\SetFigFont{9}{10.8}{\rmdefault}{\mddefault}{\updefault}{\color[rgb]{0,0,0}$C$}%
}}}
\end{picture}
\caption{Interpolating Lagrangian cylinders}
\label{StableFig1}
\end{figure}

Let $C=S^1\times\{0\}\subset\Sigma'_R$ denote the curve in
$\Sigma'$ along which we are splitting.  For
convenience, let us say that $\Sigma$ ``lies to the left of $C$''
while $T$ ``lies to the right of $C$.''  Let
$C_\Sigma=S^1\times\{-R-1\}\subset \Sigma'_R$.  Let
$C_\cT=S^1\times\{-R-1+\rho\}\subset\Sigma'_R$ for some
$\rho<R+1$.  Then,
$C_\Sigma$ and $C_\cT$ lie to the left of $C$, and $C_\cT$
lies to the right of $C_\Sigma$.
Note that the complex structures $J_R$
are split to the right of $C_\Sigma$.
We will choose $\rho$ large enough that $C_\cT$ lies
close to $C$ in a sense we will specify soon.
Let $RC_\Sigma$ denote the region to
the right of $C_\Sigma$ and $LC_\cT$ denote the region to the left
of $C_\cT$.  See \fullref{StableFig1}.

Let $\{C_k,A_\ell\}$ denote the collection of disjoint circles and arcs
in $S$ along which the complex structure degenerates.

Let $S_{RC_\Sigma,R}=(\pi_\Sigma\circ u_R)^{-1}(RC_\Sigma)$.  Since
$J_R$ is split over $RC_\Sigma$, 
$$\pi_\Sigma\circ
u_R\co \left(S_{RC_\Sigma,R},\bdy S_{RC_\Sigma,R}\right)\to\left(RC_\Sigma,
C_\Sigma\cup\ba_{g+1}\cup\bb_{g+1}\right)$$
is a holomorphic map.  So, again by
classical symplectic field theory \cite[Theorem 10.3]{Ya2}, taking a further subsequence we can
assume that $\pi_\Sigma\circ u_R|_{S_{RC_\Sigma,R}}$ converges to a
holomorphic building.  We denote this building by $\pi_\Sigma\circ
u_\infty$, with the understanding that it is, so far, only defined
over $S_{RC_\Sigma,\infty}$.

It now makes sense to talk about the circles in $\{C_k\}$ which
correspond to the degeneration along $C$ -- ie, the level splitting
of $\pi_\Sigma\circ u_R$.  Let $\{p_j\}$ denote the corresponding
points in $S_\infty$.

Since the arcs and circles $\{C_k,A_\ell\}$ are disjoint, we can
choose $C_\cT$ close enough to $C$ (ie, $\rho$ large enough)
that $(\pi_\Sigma\circ u)^{-1}(C_\cT)$ is disjoint from all the
$C_k$ and $A_\ell$.
We can also choose $\rho$ so that $\pi_\Sigma\circ
u_\infty$ is transverse to $C_\cT$.  We do so choose it.
It follows that for $R$ large enough
$\pi_\Sigma\circ u_R$ is transverse to $C_\cT$.

Now, observe that each curve in $(\pi_\bD\circ u_R)(\pi_\Sigma\circ
u_R)^{-1}(C_\cT)$ converges in the $C^\infty$--topology
as $R\to\infty$.
Let $B_R=(\pi_\bD\circ u_R)(\pi_\Sigma\circ u_R)^{-1}(C_\cT)$.  Let
$B_\infty=\lim_{R\to\infty}B_R$.

Let $LC_\cT$ denote the portion of $\Sigma'$ to the left of
$C_\cT$.  View $[0,1]$ as lying inside $S^1=[0,2]/(0\sim 2)$,
say.  We can, thus, consider $\pi_\bD\circ u_R$ as a map to
$S^1\times\bR$.
Consider the symplectic 4--manifold with boundary $LC_\cT\times
S^1\times\bR$, given the obvious (split) symplectic form.  For each
$R$, $C_\cT\times B_R$ is an immersed Lagrangian submanifold, and
these submanifolds $C^\infty$ converge to $C_\cT\times
B_\infty$.

So, applying the compactness theorem \cite[Theorem 10.1]{Ya2} for symplectic field theory to
$u_R|_{S_{LC_T,R}}$, viewed as holomorphic curves with dynamic
  Lagrangian boundary conditions, we can extract a subsequence
  converging to a holomorphic building
  $u_\infty|_{S_{LC_\cT,\infty}}$ defined over $S_{LC_\cT,\infty}$.

Let $q_j=\pi_\bD\circ
u_\infty(p_j)$.  There are $n_\Fz([u_R])$ points $q_j$, all of
them distinct.  It follows that in the (horizontal) cylindrical regions
$S^1\times\bR\times[0,1]\times\bR$ connecting
$\Sigma\times[0,1]\times\bR$ and $T^2\times[0,1]\times\bR$, the
building $u_\infty|_{S_{LC_\cT,\infty}}$ consists of trivial cylinders, so we can ignore
these regions.

The holomorphic building $u_\infty|_{S_{LC_\cT,\infty}}$ must agree
with $(\pi_\Sigma\circ
u_\infty)\times(\pi_\bD\circ u_\infty)$ where both are defined, so we
can patch the two together to obtain a holomorphic twin tower to which
(the subsequence of) the sequence $u_R$ converges.
\endproof

\begin{Prop} \label{Prop:StabilizationInvariance}
The Floer homologies $HF^\infty$, $HF^+$ and $HF^-$
  are invariant under stabilization.
\end{Prop}
\proof
This follows from the previous proposition and \fullref{Gluing:Twins}.
That is, for given intersection points the space of holomorphic twin
towers (of height one, with index 1)
connecting $\vx'$ to $\vy'$ is identified with $\hcM(\vx,\vy)$.
The admissibility criteria ensure that only finitely many homology
classes matter.  Thus, by the previous proposition and
\fullref{Gluing:Twins}, we can use the space of holomorphic
twin towers to compute the boundary maps in $W'$.

Using our original definition of the action of $H_1(Y)/Tors$ it is
immediate that the two actions are the same.  Obviously the
$U$--actions correspond.
\endproof

\textbf{Proof of \fullref{Theorem:Main}.}  
By \fullref{Prop:MaintainAdmis}, there is an admissible Heegaard diagram for $(Y,\Ss)$, and any two admissible Heegaard diagrams can be connected by a sequence of pointed Heegaard moves, ie, Heegaard moves supported in the complement of $\Fz$.  Thus, the theorem is immediate from \fullref{Prop:IsotopyInvariance}, \fullref{Prop:HandleslideInvariance}, and \fullref{Prop:StabilizationInvariance}.
\endproof

\section{Comparison with Heegaard Floer homology}
\label{Section:Comparison}
In this section we prove the equivalence of the theory described in this paper
with Heegaard Floer homology as originally defined by Ozsv\'ath and Szab\'o in~\cite{OS1}.  For notational convenience, we will phrase the argument in terms of $HF^\infty$, but the same proof works for all four theories.  This section assumes familiarity with~\cite{OS1}.

By $CF_{\textrm{ours}}$ we mean the chain complexes defined in this
paper; by $CF_{\textrm{theirs}}$ we mean the chain complex defined
in~\cite{OS1}.  We extend this notation functorially to
$HF_{\textrm{ours/theirs}}$, $\pi_2^{\textrm{ours/theirs}}$, etc.  When we
have successfully identified two corresponding objects we drop the
decorations ``$\textrm{ours}$'' or ``$\textrm{theirs}$''.

Observe that there is an identification between
our intersection points and those of~\cite{OS1}, and so between the generators of
$CF^\infty_{\textrm{ours}}$ and the generators of $CF^\infty_{\textrm{theirs}}$.
Similarly, for any intersection points $\vx$ and $\vy$,
$\pi_2^{\textrm{ours}}(\vx,\vy)$ and $\pi_2^{\textrm{theirs}}(\vx,\vy)$ are
naturally identified, by considering domains for example.

Now we deal with a simple case.  Suppose that the Heegaard diagram
$(\Sigma,\va,\vb,\Fz)$ is such that the split complex structure
$j_{\Sigma}\times j_{\bD}$ achieves transversality for all homology
classes of index $1$ in our theory and such that
$\Sym^g(j_\Sigma)$ achieves transversality for domains of index $1$ in their theory.  Let
$u\co S\to W$ be a holomorphic
curve with respect to $j_\Sigma\times j_\bD$.  Define a map
$u'\co \bD\to \Sym^g(\Sigma)$ as follows.  For $p\in \bD$, let
$(\pi_\bD\circ u)^{-1}(p)$ denote the $g$ preimages of $p$,
listed with multiplicities.  Then $u'(p)=\pi_\Sigma\circ u\left(
(\pi_\bD\circ u)^{-1}(p)\right)$ is a point in $\Sym^g(\Sigma)$.  It
is easy to check that $u$ being $(j_\Sigma\times j_\bD)$--holomorphic implies that $u'$ is
holomorphic with respect to $\Sym^g(j_\Sigma)$.  So, for any
$A\in\pi_2(\vx,\vy)$ we have a map
$\Phi\co \hcM^A_{\textrm{ours}}\to\hcM^A_{\textrm{theirs}}$, which is
clearly injective. 
 
In~\cite[Lemma 3.6]{OS1}, Ozsv\'ath and Szab\'o construct an inverse
for $\Phi$ as follows.  Fix a $\Sym^g(j_\Sigma)$--holomorphic map
$u'\co \bD\to \Sym^g(\Sigma)$.  
Let  $(g-1)!u_\bD\co (g-1)!S\to\bD$ and $(g-1)!u_\Sigma\co (g-1)!S\to
\Sigma^g$ denote the pullback via $u$ of the branched cover
$\Sigma^g\to \Sym^g(\Sigma)$, as in the following diagram:
$$\xymatrix@=36pt{
(g-1)!S\ar[r]^(.54){(g-1)!u_\Sigma} \ar[d]^{(g-1)!u_\bD} & \Sigma^g\ar[d]^{p}\\
\bD\ar[r]^{u'} & \Sym^g(\Sigma)
}$$
Since we are in the holomorphic category, pullbacks of branched
covering maps are well--defined.

The symmetric group $S_g$ acts on $\Sigma^g$ by permuting the
factors, and hence acts on $(g-1)!S$.  The maps $(g-1)!u_\Sigma$ and
$(g-1)!u_\bD$ are $S_g$--equivariant.  Let $S_{g-1}$ denote
the permutations fixing the first factor in
$\Sigma^g=\Sigma\times\Sigma^{g-1}$.  Let $\pi\co \Sigma^g\to\Sigma$
denote projection onto the first factor.  Let $S=(g-1)!S/S_{g-1}$.
Let $u_\Sigma=\left(\pi\circ(g-1)!u\right)/S_{g-1}\co S\to\Sigma$
denote the map induced by $\pi\circ (g-1)!u\co (g-1)!S\to\Sigma$, and
$u_\bD=\left((g-1)!u_\bD\right)/S_{g-1}\co S\to\bD$ the map induced by
$(g-1)!u_{\bD}$.  Then, $u=u_\Sigma\times u_\bD$ is a holomorphic
map $S\to W$.  Define $\Phi^{-1}(u')=u$.

\medskip\textbf{Remark}\qua  In the context of multivalued sections of Lefschetz pencils, the maps $\Phi$ and $\Phi^{-1}$ are called the \emph{tautological correspondence} by M Usher \cite{Usher}, who attributes the term to I Smith.  We will sometimes use this terminology below.

We check that $\Phi\circ \Phi^{-1}$ is the identity map.  It
suffices to show that $\Phi\circ\Phi^{-1}(u')$ agrees with $u'$
away from the diagonal.  Commutativity of the diagram $$
\xymatrix@=36pt{
S\ar[r]^{u_\Sigma}\ar@/_/[ddr]_{u_\bD} & \Sigma\\
&(g-1)!S\ar[r]^(.53){(g-1)!u_\Sigma} \ar[d]^{(g-1)!u_\bD}\ar[ul] &
\Sigma^g\ar[d]^{p}\ar[ul]_\pi\\
&\bD\ar[r]^{u'} & \Sym^g(\Sigma)
}
$$
shows that $\Phi\circ\Phi^{-1}(u')(x)$, viewed as a set of $g$
distinct points in $\Sigma$, agrees with the set $\pi\circ
p^{-1}(u'(x))$.  (Here, it is important that elements of sets do
\emph{not} have multiplicities.)  But this set is exactly $u'(x)$.

Observe that both $\Phi$ and $\Phi^{-1}$ are continuous.  So,
since $\Phi$ was injective, we have proved that $\Phi$ is a
homeomorphism.  This is enough to prove the equivalence of the two theories with $\bZ/2$--coefficients, under our assumption that $j_\Sigma\times j_\bD$ achieves transversality.  We will deal with the issue of orientations presently; first, we discuss the case when $j_\Sigma\times j_\bD$ does not achieve transversality.
The difficulty in the case when $j_\Sigma\times j_\bD$ does not achieve transversality is that the
nearly--symmetric almost complex structures used in~\cite{OS1}
are required to be split near the diagonal,
while compactness is proved in~\cite{Ya2} for almost complex
structures which are cylindrical.  No non--split almost complex
structure satisfies both conditions.

To address this problem, we define a class of complex structures on
$\Sym^g(\Sigma)$ which includes both ($\Sym^g$ of) the complex
structures we consider and the complex structures considered
in~\cite{OS1}.  Specifically, fix an open neighborhood $V_1$ of
$\{\Fz_i\}\times \Sym^{g-1}(\Sigma)\subset \Sym^g(\Sigma)$ and an open
neighborhood $V_2$ of the diagonal $\Delta$ in
$\Sym^{g}(\Sigma)$.  Let $\pi\co \Sigma^g\to \Sym^g(\Sigma)$ denote
projection, and $\omega_0=(dA)^g$ the product symplectic form on
$\Sigma^g$.  As usual, we also fix a complex structure $j_\Sigma$
on $\Sigma$.
\begin{Def} By a \emph{quasi--nearly--symmetric almost complex
    structure} on $\Sym^g(\Sigma)$ we mean an almost complex
    structure $\tilde{J}$ on $\Sym^g(\Sigma)$ such that
\begin{itemize}
\item $\tilde{J}$ is tamed by $\pi_*(\omega_0)$ on
  $\Sym^g(\Sigma)\setminus V_2$.
\item $\tilde{J}$ agrees with $\Sym^g(j_\Sigma)$ on $V_1$.
\item There is some complex structure $j$ on $\Sigma$ such that
  $\tilde{J}$ agrees with $\Sym^g(j)$ on $V_2$.
\end{itemize}
\end{Def}

Observe that the nearly--symmetric almost
complex structures of~\cite[Section 3.1]{OS1} are a special case of
the preceding definition.  So are complex structures of the form
$\Sym^g(j)$ for \emph{any} complex structure $j$ on $\Sigma$
which is tamed by $dA$ and agrees with $j_\Sigma$ near the
$\Fz_i$.  So, a path of nearly--symmetric almost complex structures
is a path of quasi--nearly--symmetric almost complex structures.
Also, an almost complex structure $J$ on $W$ satisfying
(\textbf{J1})--\textbf{J5}) corresponds to a path $J_s$ of complex structures on $\Sigma$, which in turn specifies a path of quasi--nearly--symmetric almost complex structures on
$\Sym^g(\Sigma)$, which we denote $\Sym^g(J_s)$.

For a path $\tilde{J}_t$ of quasi--nearly--symmetric almost complex
structures, let $j_t$ be the complex structure on $\Sigma$ such
that $\tilde{J}_t$ agrees with $\Sym^g(j_t)$ on $V_2$.

Now, we prove the equivalence of our theory with that of~\cite{OS1} in
three steps.  First we show that the compactness proof of~\cite{OS1}
extends to $\tilde{J}_t$--holomorphic curves in $\Sym^g(\Sigma)$,
where $\tilde{J}_t$ is a path of
quasi--nearly--symmetric almost complex structures.  Then we observe
that the class of paths of almost complex structures of the form
$\Sym^g(J_s)$ is sufficient to achieve transversality for holomorphic
disks in $\Sym^g(\Sigma)$.  Finally we show that our moduli spaces of
$J_s$--holomorphic curves can
be identified with the Ozsv\'ath--Szab\'o moduli spaces of
$\Sym^g(J_s)$--holomorphic curves.

After proving the following proposition, the first and third steps are
essentially immediate.
\begin{Prop}\label{Prop:HoloLifts}
 Fix a path $\tilde{J}_t$ of quasi--nearly--symmetric
  almost complex structures on $\Sym^g(\Sigma)$.  Fix a holomorphic
  disk $u'\co (\bD,\bdy\bD)\to(\Sym^g(\Sigma),T_\alpha\cup T_\beta)$.
  Then there is a $g!$--fold branched covering
  $\tilde{u}_\bD\co \tilde{S}\to\bD$ and a map
  $\tilde{u}_\Sigma\co \tilde{S}\to \Sigma^g$ such that the following diagram
  commutes.
$$\xymatrix@=36pt{
\tilde{S}\ar[r]^{\tilde{u}_\Sigma} \ar[d]^{\tilde{u}_\bD} & \Sigma^g\ar[d]^{p}\\
\bD\ar[r]^{u'} & \Sym^g(\Sigma)
}.$$
The map $\tilde{u}_\Sigma$ is holomorphic with respect to the path
of almost complex structures on $\Sigma^g$ induced by $\tilde{J}_t$, in
the obvious sense, and is $S_g$--equivariant.
\end{Prop}
\proof
The idea of the proof is that even though the complex structure is
allowed to vary near the diagonal, since it varies in the class of
split structures, locally near the diagonal we are
still in the integrable case.  The model for this argument is used to
prove the following
\begin{Lem} \label{Lemma:IntersectsDiscretely}
Under the assumptions of \fullref{Prop:HoloLifts},
  $u'$ intersects the diagonal $\Delta$ in a discrete collection
  of points.
\end{Lem}
\proof
Suppose that $u'$ intersects the diagonal $\Delta$ in a collection
of points $p_j$ with limit point $p$.
Throwing out some of the points we may assume all of the $p_j$ lie
in the same stratum of $\Delta$.  For concreteness we will assume
all of the $p_j$ lie in the top dimensional stratum of $\Delta$,
but there is nothing special about this case.  

Write
$p=\{a_1,a_1,a_3,\cdots,a_g\}$.  Choose pairwise disjoint disk
neighborhoods $U_i$ of the $a_i$ such that $U_1\times U_1\times
U_3\times\cdots\times U_g$ is contained in $V_2$.  Choose three
points on the boundary of each $U_i$.  Then the Riemann mapping
theorem gives a well--defined holomorphic identification of $(U_i,j_t)$ and
$(\bD,j_\bD)$ for each $t$.

So in a neighborhood $V$ of ${u'}^{-1}(p)$ we can view $u'$ as a map 
to $(\Sym^g(\bD),\Sym^g(j_{\bD}))$.  The diagonal $\Delta\subset
\Sym^g(\bD)$ is an analytic subvariety, so by elementary complex analysis,
$u'|_V\co V\to\Delta$.

Globally we know that
the image of $u'$ is not entirely contained in $\Delta$.  So, a
standard open--closed argument gives a contradiction.
\endproof

We return to the proof of \fullref{Prop:HoloLifts}.  Away from the diagonal $\Delta$, to obtain $\tilde{S}$, $\tilde{u}_\Sigma$ and $\tilde{u}_\bD$ we simply pullback the
holomorphic covering space $\Sigma^g\to \Sym^g(\Sigma)$.
Near the diagonal,
the argument used in the preceding lemma's proof shows that we can
pull back the branched covering (which locally looks like a piece of
$((\bD^2)^g,j_\bD^g)\to (\Sym^g(\bD),\Sym^g(j_\bD))$) by $u'$.  This
proves the proposition.
\endproof

\begin{Prop}The moduli spaces considered in~\cite{OS1}, computed with
  respect to any quasi--nearly--symmetric almost complex structure,
  are compact.
\end{Prop}
\proof
The only place in their proof that Ozsv\'ath and Szab\'o use the
condition that their complex structures are standard near $\Delta$
is in the proof of their energy estimate \cite[Lemma 3.5]{OS1}.
The only time they use it in that proof is to observe that the
previous proposition holds.  So, compactness is immediate from their work.
\endproof

\begin{Prop}The class of paths of complex structures of the form
  $\Sym^g(J_s)$ is sufficient to achieve transversality for disks
  $u\co (\bD,\bdy \bD)\to (\Sym^g(\Sigma),T_\alpha\cup T_\beta)$.
\end{Prop}
The proof, which we omit, is a simple adaptation of the one in \fullref{Section:Transversality}.

\begin{Prop}Calculated with respect to $J_t$ and $\Sym^g(J_t)$ respectively,
$\hcM^A_{\rm{ours}}$ and $\hcM^A_{\rm{theirs}}$ agree.
\end{Prop}
The proof of this proposition is the same as the proof in the split
case, using \fullref{Prop:HoloLifts} where appropriate.

This is sufficient to prove that the two theories are equivalent with
$\bZ/2$--coefficients.  To prove the equivalence with
$\bZ$--coefficients, we need somewhat more.  Specifically, we need
to check that
\begin{Prop}\label{Prop:OrId}
\begin{enumerate}
\item \label{item:or1}If $u_1$ and $u_2$ are two curves in the moduli space
  $\cM^A$, $\ind(A)=1$, and $\sign_{\rm ours}(u_1)=\sign_{\rm theirs}(u_1)$
  then $\sign_{\rm ours}(u_2)=\sign_{\rm theirs}(u_2)$.  (Here, $\sign$
  denotes the sign of a rigid curve induced by the coherent orientation.)  That is, if the
  orientations agree at one curve then they agree at all curves in
  that moduli space.
\item \label{item:or2}The coherence conditions for the two theories agree.  That is, there is an identification of the determinant line bundles $\cL_{\rm ours}$ and $\cL_{\rm theirs}$ with respect to which any coherent orientation $\frako_{\rm ours}$ of $\cL_{\rm ours}$ specifies a coherent orientation $\frako_{\rm theirs}$ of  $\cL_{\rm theirs}$.
\end{enumerate}
\end{Prop}

To keep the exposition clean, we will assume that the complex structure is split near the
diagonal.  As noted earlier, this assumption is quite
restrictive.  It can be removed by using, at appropriate times, a parametrized version of the Riemann mapping theorem to identify a neighborhood of $u_\Sigma\left(\textit{branch points of $u_\bD$}\right)$ with a standard disk, as in the proof of \fullref{Lemma:IntersectsDiscretely}.  (Alternately, given a curve $u$, we can choose a family $J'_s:T\Sigma\to T\Sigma$ of almost complex structures on $\Sigma$ constant near $\Phi(u)\cap\Delta$ and such that the $(\Sigma,J'_s)$ are conformally isomorphic to the original family $(\Sigma,J_s)$, and then work with the family $(\Sigma,J'_s)$.)  We leave further details to the interested reader.

From now on, when we want to discuss both our and their theories at
once, we will drop the ``ours'' or ``theirs'' from the notation, even
if we have not yet identified the corresponding objects.

Recall from \fullref{Section:Transversality} that the tangent space at $u$ to $\cB_{\rm ours}$ is
$\bR^{2g}\oplus L^{p,d}_1\left(u^*TW,\bdy\right)$, where the $\bR^{2g}$
includes into $\Gamma\left(u^*TW\right)$ as
$\mathrm{Span}\left(\{v_i^\pm\}\right)$ where the $v_i^\pm$ are
some fixed smooth sections with $v_i^\pm$ equal to
$\frac{\partial}{\partial t}$ in a small neighborhood of the $i^{th}$ positive or
negative puncture and zero outside a slightly larger neighborhood of
that puncture and $L^{p,d}_1\left(u^*TW,\bdy\right)$ denotes
sections tangent to $C_\alpha\cup C_\beta$ over $\bdy S$.  The tangent space at $\phi$ to $\cB_{\rm theirs}$ is
just $L^{p,d}_1\left(\phi^*T\Sym^g(\Sigma),\bdy\right)$ of
$L^{p,d}_1$ sections tangent to $T_\alpha\cup T_\beta$ over $\bdy\bD$
and we take as a model for $\bD$ the strip $[0,1]\times\bR$.

There is the minor complication that we are working in
$\bR$--invariant settings, so rigid curves have $\ind=1$.  For
convenience, we introduce the operator $P_{u,ours}\co \bR^{2g}\oplus L^{p,d}_k\left(u^*TW,\bdy\right)\to \bR$ defined by
$P_{u,ours}(v_1,v_2)=\langle v_2,\partial/\partial t\rangle_{L^2}$.  (Here,
$\langle\cdot,\cdot\rangle_{L^2}$ denotes the 
$L^2$ inner product
induced by some Riemannian metric on $W$ and a canonical metric on
the source $S$ of $u$.  Since $L^{p,d}_k$ is finer than $L^2$, the operator $P_{u,ours}$ is continuous.)  Then we replace the linearized
$\dbar$--map $D\dbar_{\rm ours}$ by $\tD\dbar_{\rm ours}=D\dbar_{\rm ours}\oplus P_{\rm ours}$.  Similarly, for $\dbar_{\rm theirs}$, let
$P_{\phi,theirs}(v)=\langle v,d\phi(\partial/\partial t)\rangle_{L^2}$,
where $\partial/\partial t$ generates the one--parameter group of
automorphisms of $(\bD^2,i,-i)$.  Then, replace $D\dbar_{\rm theirs}$
with $\tD\dbar_{\rm theirs}=D\dbar_{\rm theirs}\oplus P_{\rm theirs}$.  We
retain the old meaning of $\ind$, so the index of
$\tD_u\dbar_{\rm ours/theirs}$ is $\ind(u)-1$.

We tackle point~\eqref{item:or1} of the proposition first.  We begin by recalling how the
coherent orientation specifies signs of curves.  Let $u$ be a curve
(not necessarily holomorphic), $\ind(u)=1$, at which
$D\dbar$ (or equivalently $\tD\dbar$) is surjective.  We can view
the coherent orientation $\frako$ as a nonvanishing section of
$\cL$, defined up to multiplication by a positive scalar function.
At $u$, there is also a canonically defined section
{\small$$
1\otimes
1^*\in\bR\otimes\bR^*=\left(\Lambda^0\bR^0\right)\otimes\left(\Lambda^0\bR^0\right)^*=\left(\Lambda^{top}\ker(\tD_u\dbar)\right)\otimes\left(\Lambda^{top}\coker(\tD_u\dbar)\right)^*=\cL_u.
$$}%
The sign $\sign(u)$ is $+1$ (respectively $-1$) if $\frako(u)$ is
a positive (respectively negative) multiple
of $1\otimes 1^*$.

The section $1\otimes 1^*$ does not extend continuously to the whole
configuration space.  Indeed, for a generic path $\{u_a\}$  (ie,
one with only \emph{normal crossings}) between
curves $u_0$ and $u_1$, with $D_{u_0}\dbar$ and $D_{u_1}\dbar$
surjective, the section $1\otimes 1^*$ switches sign at each $a$ for
which $D_{u_a}\dbar$ is not surjective.  So, if $c(\{u_a\})$ denotes
the number of $a$ for which $D_{u_a}\dbar$ is not surjective then
$\sign(u_0)=(-1)^{c(\{u_a\})}\sign(u_1)$.

It follows that to prove part~\eqref{item:or1} of \fullref{Prop:OrId} if suffices to show that
$\ker(\tD\dbar_{\rm ours})$ and $\ker(\tD\dbar_{\rm theirs})$ are nontrivial
at the same curves.  For this to make sense, we first need to identify
the configuration spaces.  Actually, we will define subspaces
$B_{\rm ours}\subset\cB_{\rm ours}$ and $B_{\rm theirs}\subset\cB_{\rm theirs}$
such that
\begin{itemize}
\item $B_{\rm ours/theirs}\supset \cM_{\rm ours/theirs}$ and
\item $B_{\rm ours/theirs}$ contains a path with only normal crossings
  between any two curves at with $D\dbar$ is surjective.  (In
  particular, $B_{\rm ours/theirs}$ is connected.)
\end{itemize}
Then, we will construct an identification $\Phi$ between
$B_{\rm ours}$ and $B_{\rm theirs}$ extending the identification $\Phi$
between $\cM_{\rm ours}$ and $\cM_{\rm theirs}$.  Finally, we will
construct an identification of $\ker(\tD_u\dbar_{\rm ours})$ and
$\ker(\tD_{\Phi(u)}\dbar_{\rm theirs})$.  This suffices to
prove part~\eqref{item:or1} of \fullref{Prop:OrId}.

None of the steps involved are particularly intricate.  The space
$B_{\rm ours}$ consists of embedded curves $u$ such that $u_\bD$ is
holomorphic with only order $2$ branch points, near which
$u_\Sigma$ is holomorphic.  The space $B_{\rm theirs}$ consists of
curves $\phi$ intersecting the diagonal $\Delta$ only in the top--dimensional
stratum, transverse to $\Delta$, and holomorphic in a neighborhood
of $\Delta$.  It is clear that
$B_{\rm ours}$ and $B_{\rm theirs}$ have the requisite properties.  The
identification $\Phi\co B_{\rm ours}\to B_{\rm theirs}$ is given in exactly the
same way as $\Phi\co \cM_{\rm ours}\to\cM_{\rm theirs}$ was defined at the
beginning of the section.

The identification of kernels is somewhat more difficult.  The key
point is that near the branch points, the kernel of
$\tD\dbar_{\rm ours}$ itself consists of holomorphic sections of
$T\Sigma$.  Such sections correspond by the same tautological
correspondence used to define $\Phi$ to holomorphic sections of $T\Sym^g(\Sigma)$, which
comprise the kernel of $\tD\dbar_{\rm theirs}$ near the diagonal.  More
details follow.

Recall \cite[page 28]{MS2}
that the linearized
$\dbar$--operator 
$$(D_u\dbar)\co \bR^{2g}\oplus L^{p,d}_1\left(u^*TW,\bdy)\right)\to L^{p,d}(\Lambda^{0,1}u^*TW) $$
at a curve $u\co (S,\bdy S)\to (W,C_\alpha\cup C_\beta)$, with fixed complex structure on $S$, is given by
$$
D_u\dbar(\xi)(v)=\frac{1}{2}\left(\nabla_v\xi+J\nabla_{j(v)}\xi\right)-\frac{1}{2}J\left(\nabla_\xi J\right)\partial_J(u).
$$
Here, $\nabla$ is the metric connection on $TW$ of the metric
induced by $\omega$ and $J$, and $\xi\in \bR^{2g}\oplus
L^{p,d}_1(u^*TW,\bdy)\subset \Gamma(u^*TW)$.  
The space $L^{p,d}_1\left(u^*TW,\bdy\right)$ consists of
sections of $u^*TM$ tangent to $C_\alpha\cup C_\beta$ over $\bdy S$ with one
derivative in $L^{p,d}$; see \fullref{Section:Transversality}
for more details.

If we allow the complex structure on $S$ to vary, then there is an
additional term.  The tangent space to Teichm\"uller space at $S$ is
a finite--dimensional subspace of $C^\infty\left(\End(TS,j)\right)$, where
$\End(TS,j)$ is the space of endomorphisms of $TS$ anticommuting
with $j$.  The linearized $\dbar$--operator is then a restriction
of the map
$$
D_u\dbar\co \bR^{2g}\oplus L^{p,d}_1\left(u^*TW,\bdy\right)\oplus C^\infty\left(\End(TS,j)\right)\to L^{p,d}(\Lambda^{0,1}u^*TW)
$$
defined by
$$
D_u\dbar(\xi,Y)(v)=\frac{1}{2}\left(\nabla_v\xi+J\nabla_{j(v)}\xi\right)-\frac{1}{2}J\left(\nabla_\xi
  J\right)\partial_{J(v)}(u)+\frac{1}{2}J\circ du\circ Y(v).
$$

For maps $\phi:\bD\to \Sym^g(\Sigma)$ the formulas are the same
except that the $\bR^{2g}$--factors are absent.

In the future, we will suppress the Lagrangians from the notation.

Fix $\epsilon>0$.
Let $\ol{B}_{\rm ours}$ denote the collection of maps $u$ in $\cB_{\rm ours}$
for which
\begin{itemize}
\item $u_\bD$ is holomorphic,
\item $u_\bD$ has only simple branch points $p_1,\cdots,p_\ell$,
\item  $d(p_i,p_j)>2\epsilon$ for $i\neq j$, and
\item $u_\Sigma$ has no branch points inside the
$B_{\epsilon}(p_i)$, $i=1,\cdots,\ell$.
\end{itemize} 

Let $B_{\rm ours}\subset\ol{B}_{\rm ours}$ denote the space of
maps $u\in\ol{B}_{\rm ours}$ for which $u_\Sigma$ is holomorphic over
$\cup_i u_{\bD}^{-1}\left(B_\epsilon(p_i)\right)$.
For $\epsilon$ sufficiently small, $J$ generic,
and $\ind(A)=1$, $\cM^A_{\rm ours}\subset B_{\rm ours}$.  Note that $B^A_{\rm ours}$ is nonempty if and only if the intersection number $\Delta\cdot A$ of the diagonal in $\Sym^g(\Sigma)$ with $A\in H_2(\Sym^g(\Sigma),T_\alpha\cup T_\beta)$ is non--negative.

Let $u\co S\to W$ be a map in $\ol{B}_{\rm ours}$ such that $u_\bD$ has branch points
$p_1,\cdots,p_\ell\in(0,1)\times\bR$.

\begin{Lem} \label{Lemma:TanSpace} The tangent space to
  $\ol{B}_{\rm ours}$ at $u$ is given by
$$
T_u\ol{B}_{\rm ours}=L^{p,d}_k\left((\pi_\Sigma\circ u)^*T\Sigma,\bdy\right)\oplus\left(\oplus_{i=1}^\ell\bC\right).
$$
\end{Lem}
\proof
This is clear.  A point in $\ol{B}_{\rm ours}$ is determined by a
complex structure on $S$ and a map from $S$ to $\Sigma$.  The
complex structure is determined by the branch points
$p_1,\cdots,p_\ell$.
\endproof

An inclusion of $T_u\ol{B}_{\rm ours}$ into $T_u\cB_{\rm ours}$ can be
given as follows.  The inclusion\break $L^{p,d}_1(u_\Sigma^*T\Sigma,\bdy)\mapsinto
L^{p,d}_1(u^*TW,\bdy)$ is obvious.  We include $\bC^\ell$ into the
space of infinitesimal deformations of the almost complex structure on the source.  The
$\bC^\ell$ corresponds to moving the branch points
$p_1,\cdots,p_g$ in $\bD$.  This, in turn, corresponds to
deforming the almost complex structure $j_\bD$ on $\bD$.  But any almost
complex structure $j_\bD$ on $\bD$ specifies an almost complex structure
$j_S$ on $S$ via $j_s=(du_\bD)^{-1}\circ j_\bD\circ du_\bD$.

We choose a
family of infinitesimal deformations parametrizing the $\bC^\ell$ vanishing over the branch points $p_1,\cdots,p_\ell$, ie, such that the
almost complex structure remains fixed near the branch points.

\begin{Lem}
$\ker(D\dbar)|_{\ol{B}_{\rm ours}}=\ker(D(\dbar|_{\ol{B}_{\rm ours}}))$.
\end{Lem}
\proof
Again, this is clear.
\endproof

To discuss the linearized $\dbar$--operator for maps to
$\Sym^g(\Sigma)$, we must first fix a connection on $T\Sym^g(\Sigma)$.  Away from the diagonal we choose the metric connection of the metric induced by the split symplectic form and our almost complex structure.  We extend this connection arbitrarily over the diagonal.  (By~\cite[Corollary 2]{Tim}, we could in fact extend the symplectic form over the diagonal, and work with the induced metric connection.)  Since we will work with curves which are holomorphic near the diagonal, the choice of linearization near the diagonal is unimportant.

Let $\ol{B}_{\rm theirs}$ denote the collection of maps
$\phi\co \bD\to \Sym^g(\Sigma)$ (with boundary on the $\alpha$-- and
$\beta$--tori) intersecting the diagonal $\Delta$ transversely and
only in its top stratum.  Notice that $\ol{B}_{\rm theirs}$ is an open
subset of $\cB_{\rm theirs}$.  Let $B_{\rm theirs}\subset \ol{B}_{\rm theirs}$
denote the subspace of maps which are holomorphic near the diagonal.

\begin{Lem} The tangent space at $\phi$ to $\ol{B}_{\rm theirs}$ is
  $$T_\phi\ol{B}_{\rm theirs}=L^{p,d}_1\left(\phi^*T\Sym^g(\Sigma),\bdy,\Delta\right)\oplus\bC^\ell$$
where
$L^{p,d}_1(\phi^*T\Sym^g(\Sigma),\bdy,\Delta)$ is the space of
$L^{p,d}_1$ vector fields along $\phi$ which are tangent to the diagonal over
$\phi^{-1}(\Delta)$ and to the Lagrangian tori $T_\alpha$ and
$T_\beta$ over $\bdy\bD$, and $\ell$ is the number of intersections of
$\phi$ with $\Delta$.
\end{Lem}
\proof
The tangent space to the space of maps 
$$(\bD^2,\{p_1,\cdots,p_\ell\})\to(\Sym^g(\Sigma),\Delta)$$
is $L^{p,d}_1\left(\phi^*T\Sym^g(\Sigma),\bdy,\Delta\right).$  The
$\bC^\ell$ corresponds to allowing the $p_i$ to move.
\endproof

Again, we can identify the $\bC^\ell$ with the tangent space to an
$\ell$--dimensional family of deformations of $j_\bD$.  We take
this to be the same $\ell$--dimensional family used before.  In
particular, the family of almost complex structures is constant near
the $p_i$.

There is a map
$\Phi\co (\ol{B}_{\rm ours},B_{\rm ours})\to(\ol{B}_{\rm theirs},B_{\rm theirs})$
defined just as $\Phi\co\cM_{\rm ours}\to\cM_{\rm theirs}$ was defined at the
beginning of this section.

Now, suppose that $u\in B_{\rm ours}$ and 
$$(v,w)\in\ker(D_u\dbar|_{B_{\rm ours}})\subset
L^{p,d}_1(u_\Sigma^*T\Sigma,\bdy)\times \bC^\ell$$
We construct from $(v,w)$ an element of
$\ker\left(D_{\Phi(u)}\dbar|_{B_{\rm theirs}}\right)$.  First we construct a vector field $v'\in
L^{p,d}_1(\Phi(u)^*T\Sym^g(\Sigma),\bdy)$.  
For $p\in\bD$, $p\notin \cup_iB_\epsilon(p_i)$, let
$u_\bD^{-1}(p)=\{q_1,\cdots,q_g\}$ and 
$\Phi(u)(p)=\{x_1,\cdots,x_g\}=\{u_\Sigma(q_1),\cdots,u_\Sigma(q_g)\}$.
Then define
$v'(p)=\{v(q_1),\cdots,v(q_g)\}\in T_{\{x_1,\cdots,x_g\}}\Sym^g(\Sigma)$.

Since $\Phi(u)(p_i)$ lies in the top--dimensional stratum of the diagonal, near $\Phi(u)(p_i)$, $\Sym^g(\Sigma)$ decomposes as $\Sym^2(\Sigma)\times\Sigma^{g-2}$.  The projection of $v'$ to $T\Sigma^{g-2}$ is given by the previous construction.  This leaves us to define the projection $v'_1$ of $v'$ to $T\Sym^2(\Sigma)$.  

Inside the $B_{\epsilon}(p_i)$, $w$ is constant so the term
$J\circ du\circ Y$ in $D_u\dbar$ is zero.  Let $U_i$ be the component of $u_\bD^{-1}\left(B_{\epsilon}(p_i)\right)$ on which $du_\bD$ is singular.  Identify
$U_i$ and $V_i=u_\Sigma(U_i)$ holomorphically with
$\bD$.  Then, $TV_i$ is identified with $V_i\times\bC$ and
$u_\Sigma^*TV_i$ with $U_i\times\bC$.  Further, $v$ becomes a
map $U_i\to \bC$, and the statement that $D_u\dbar v=0$ becomes
the statement that $v\co U_i\to \bC$ is holomorphic.

Now, $\Phi(u)(B_\epsilon(p_i))\subset
\Sym^2(V_i)\subset\Sym^2(\bC)\cong\bC^2$.  So,
$\Phi(u)^*T\Sym^2(\Sigma)$ is identified with
$B_{\epsilon}(p_i)\times\bC^2$.  The holomorphic map
$(u_\bD,v)\co U_i\to B_{\epsilon}(p_i)\times\bC$ specifies
a holomorphic map $v'_1\co B_{\epsilon}(p_i)\to\Sym^2(\bC)\cong\bC^2$
by the same tautological correspondence used to define $\Phi$.
We view $v'_1$ as a section of
$\Phi(u)^*T\Sym^2(\Sigma)$.  Observe that $v'_1$ and hence $v_1$ is tangent to the diagonal.

It is easy to check that the two definitions of $v'$ agree over
$\bdy B_{\epsilon}(p_i)$.

Finally, by its definition, $w$ corresponded to an infinitesimal deformation of $j_\bD$.

\begin{Lem}The pair $(v',w')$ lies in the kernel of $D\dbar_{\rm theirs}$.
\end{Lem}
\proof
This is direct from the definitions.  Away from the diagonal, this
follows from the fact that the complex
structure and symplectic form on $\Sym^g(\Sigma)$ have the form
$\Sym^g(j_\Sigma)$ and $\Sym^g(\omega_\Sigma)$ respectively.  Near
the diagonal, on the $\Sym^2(\Sigma)$ factor the kernel of $D\dbar_{\rm theirs}$ corresponds under the
trivializations used above to holomorphic maps $B_{\epsilon}(p_i)\to
\Sym^2(\bC)$.  The map $v'$ is holomorphic near the diagonal by the usual
tautological correspondence.
\endproof

This proves that $\ker(D\dbar_{\rm ours})\subset\ker(D\dbar_{\rm theirs})$.
The reverse inclusion can be proved similarly, using the opposite
direction $\Phi^{-1}$ of the tautological correspondence.

This proves part~\eqref{item:or1} of \fullref{Prop:OrId}.  It remains to check
part~\eqref{item:or2}, ie, that the coherence conditions for the two
theories agree.  Before doing so, we prove another
\begin{Lem}\label{Lemma:ExtendId} Fix $A$ with $\ind(A)>0$ and $A\cdot \Delta\geq 0$.  Let $\frako_{\rm ours}^A$ and $\frako_{\rm theirs}^A$ be sections of $\cL_{\rm ours}$ and $\cL_{\rm theirs}$ respectively.  Let $u_1,u_2\in B_{\rm ours}^A$ with $D_{u_i}\dbar$ surjective for $i=1,2$.  At $u_i$ we constructed an identification of $\cL_{\rm ours}$ and $\cL_{\rm theirs}$.  Suppose that with respect to the identification of determinant liens at $u_1$, $\frako_{\rm ours}^A(u_1)$ is a positive multiple of $\frako_{\rm theirs}^A(u_1)$.  Then with respect to the identification of determinant lines at $u_2$, $\frako_{\rm ours}^A(u_2)$ is a positive multiple of $\frako_{\rm ours}^A(u_2)$.
\end{Lem}
\proof
Choose a path $u_a$ in $B_{\rm ours}$ from $u_1$ to $u_2$.  Choose also a family of subspaces $H$ of $L^{p,d}\left(\Lambda^{0,1}u_a^*T\Sigma\right)$ such that the operator 
$$\ol{D}_{u_a}\dbar\co L^{p,d}_1(u_a^*T\Sigma,\bdy)\oplus \bC^\ell\oplus H\to L^{p,d}\left(\Lambda^{0,1}u_a^*T\Sigma\right)$$
given by
$$
\ol{D}_{u_a}\dbar(v,w,x)=D_{u_a}\dbar(v,w)+x
$$
is surjective for all $a$.  Further, choose $H$ so that all sections in $H$ vanish near the branch points of $u_\bD$.  This is possible for a generic path $u_a$.

Choose an orientation $\frako(H)$ of $H$, ie, a section of $\Lambda^{top}H$.  Then there is a canonical isomorphism between the determinant lines of $D\dbar$ and $\ol{D}\dbar$.

Since the sections in $H$ vanish near the branch points of $u_\bD$, the space $H$ specifies a subspace $H'$ of $L^{p,d}\left(\Lambda^{0,1}\Phi(u_a)^*T\Sym^g(\Sigma)\right)$ as follows:  for $x\in H$ and $y\in T_p\bD$ let $y_1,\cdots,y_g$ be the preimages of $y$ under $u_\bD$.  Then define 
$$x(y)=\{x(y_1),\cdots,x(y_g)\}\in T_{\Phi(u)(p)}\Sym^g(\Sigma).$$

The vanishing of the sections in $H$ near the branch points of $u_\bD$ also means that we have an identification of $\ker(\ol{D}\dbar_{\rm ours})$ and $\ker(\ol{D}\dbar_{\rm theirs})$, defined in the same way as the identification of $\ker(\tD\dbar_{\rm ours})$ and $\ker(\tD\dbar_{\rm theirs})$ above.  Further, the following diagram commutes:
$$
\xymatrix{ Det(D_{u_i}\dbar) \ar[d]^{\otimes\frako(H)}_\cong\ar[r]^(.45)\cong & Det(D_{\Phi(u_i)}\dbar) \ar[d]^{\otimes\frako(H)}_\cong\\
Det(\ol{D}_{u_i}\dbar) \ar[r]^(.45)\cong & Det(\ol{D}_{\Phi(u_i)}\dbar)}
$$
for $i=1,2$.  But this implies the result.
\endproof

Note that as a porism of \fullref{Lemma:ExtendId} we also have an identification of determinant lines for $A$ with $A\cdot\Delta\geq 0$ but $\ind(A)\leq0$.

Next we recall the definition of a coherent
orientation.  A coherent orientation consists of a choice of nonvanishing section
$\frako(A)$ of the determinant line bundle $\cL$ over the
configuration space in the homotopy class $A$ such that the
following coherence condition is satisfied:

Let $u_1\in\cB^{A_1}(\vx,\vy)$ and $u_2\in\cB^{A_2}(\vy,\vz)$ be
curves for which $D_{u_i}\dbar$ is surjective.  In Appendix~\ref{Section:Gluing} we
construct a family of preglued curves $u_1\natural_r
u_2\in\cB^{A_1+A_2}(\vx,\vz)$.  For $r$ large, the kernel of
$D_{u_1\natural_r u_2}\dbar$ is identified with
$\ker(D_{u_1}\dbar)\oplus\ker(D_{u_2}\dbar)$.  So,
$\frako^{A_1}\otimes\frako^{A_2}$ specifies a section of
$\cL^{A_1+A_2}$.  The \emph{coherence condition} is that
$\frako^{A_1+A_2}$ be a positive multiple of
$\frako^{A_1}\otimes\frako^{A_2}$.  One must check that this condition is
independent of $u_1$ and $u_2$; see for instance~\cite[Corollary 7]{BM}.
(In the case that $D_{u_i}\dbar$ is
not surjective, one stabilizes $D_{u_i}\dbar$ by a
finite--dimensional oriented
subspace of $L^{p,d}\left(\Lambda^{0,1}u_i^*TW\right)$, as in the proof of \fullref{Lemma:ExtendId}.)

(From now on, when we write a sum $A+B$ of homotopy classes of curves we implicitly assume that the asymptotics of $A$ at $\infty$ agree with those of $B$ at $-\infty$.)

Returning to our situation, fix a coherent orientation $\frak{o}_{\rm ours}$.  By the previous lemma, this specifies an orientation $\frako_{\rm theirs}^A$ for each $A$ with $\ind(A)>0$ and $A\cdot\Delta\geq 0$.  We will refer to the collection of $A$ with $\ind(A)>0$ and $A\cdot\Delta\geq 0$ as the \emph{positive cone} and denote it $C_+$.  Note that $[\Sigma]\in C_+$ and for any $A$, $A+N[\Sigma]\in C_+$ for $N$ sufficiently large.

\begin{Lem}The orientation $\frako_{\rm theirs}$ is coherent over the positive cone.  That is, for $A,B\in C_+$, with respect to the identification induced by gluing, $\frako_{\rm theirs}^A\otimes\frako_{\rm theirs}^B$ is a positive multiple of $\frako_{\rm theirs}^{A+B}$.
\end{Lem}
\proof
This follows from the commutativity of
$$
\xymatrix{ \cL^{A_1}_{\rm ours}\otimes\cL^{A_2}_{\rm ours}\ar[r]\ar[d] &
  \cL^{A_1}_{\rm theirs}\otimes\cL^{A_2}_{\rm theirs}\ar[d]\\
\cL^{A_1+A_2}_{\rm ours}\ar[r] & \cL^{A_1+A_2}_{\rm theirs}
}
$$
up to positive scaling, which follows easily from the definitions of the various maps involved.
\endproof

\begin{Lem}Any coherent orientation $\frako$ over the positive cone can be extended uniquely to a coherent orientation over all configuration spaces.
\end{Lem}
\proof
Fix $A$, and $N$ large enough that $A+N[\Sigma]\in C_+$.  Then the map $$\xymatrix{\cL^A\ar[r]^(.45){\otimes \frako^{[N\Sigma]}}&\cL^{A+N[\Sigma]}}$$ and the orientation $\frako^{A+N[\Sigma]}$ specifies an orientation $\frako^A$ of $\cL^A$.  By the coherence of $\frako$ over the positive cone, the orientation $\frako^A$ is independent of $N$, and obviously agrees with $\frako$ over the positive cone.  Commutativity of the following diagram, up to positive scaling, implies that $\frako$ is coherent:
$$
\xymatrix{
\cL^A\otimes\cL^B
\ar[rrr]^(.4){\otimes\frako^{[N_1\Sigma]}\otimes\frako^{[N_2\Sigma]}} \ar[d]
 & & &
\cL^{A+N_1[\Sigma]}\otimes\cL^{B+N_2[\Sigma]}\ar[d]\\ 
\cL^{A+B}
\ar[rrr]^(.4){\otimes\frako^{[N_1\Sigma]}\otimes\frako^{[N_2\Sigma]}}
 & & &
\cL^{A+B+(N_1+N_2)[\Sigma]}.
}
$$
(Here, the vertical maps are induced by gluing.  This requires stabilizing the spaces involved, as discussed above.)
\endproof

We have now proved part~\eqref{item:or2} of \fullref{Prop:OrId}.  This concludes the proof of \fullref{Theorem:Equivalent}.

\section{Other Remarks}
\label{Section:OtherRemarks}
\subsection{Elaborations of Heegaard Floer}
The 4--dimensional approach that we have used suggests possible
elaborations of the Floer homology groups considered in this paper.  For the
first one, we will
describe briefly the $\widehat{HF}$--case, the only case that I am confident
works.  Our new chain complex $\widehat{CF}_{big}$ is freely generated
over $\bZ[[t]]$ by the intersection points.  The differential is
$\bdy_{big}=\bdy_0+t\bdy_1+t^2\bdy_2+\cdots$.  Here, $\bdy_0$ is
boundary operator for $\widehat{HF}$ that appears throughout this
paper.  The linear term $\bdy_1$
counts holomorphic curves with exactly $1$ double
point, in homology classes $A$ with $\ind(A)=3$.  (This is the
appropriate condition for the resulting moduli space to be
zero--dimensional.)  In general, $\bdy_i$ counts holomorphic curves
with singularity equivalent to $i$ double points, in homology classes $A$ with
$\ind(A)=2i+1$.  (This idea is inspired by the so--called Taubes
series described by Ionel--Parker in~\cite{IP}.  The analog for disks
in $\Sym^g(\Sigma)$ is to consider only disks with a prescribed
number of tangencies to the diagonal.)

The resulting filtered chain complex, up to filtered chain homotopy, is an invariant of the three--manifold.
The proof is similar to the proof of invariance of $\widehat{HF}$ above.  Undesirable codimension--one degenerations are ruled out just as in \fullref{CompactProp}.  The fact that $\bdy_{big}^2=0$ now expresses the additional fact that when a holomorphic curve with $k$ double points splits, $i$ of them go to one level and $k-i$ to the other.

Isotopy invariance follows from the arguments of \fullref{Section:Isotopy}.  Like $\bdy_{big}$, the chain maps $\Phi$ constructed there now take the form $\Phi_{big}=\Phi_0+t\Phi_1+t^2\Phi_2+\cdots$, where $\Phi_i$ counts holomorphic curves with $i$ double points or equivalent.  Again, the statement that $\Phi_{big}$ is a chain map expresses the observation that when a holomorphic curve with double points splits, some double points go to one level and the rest to the other level.  Otherwise, the argument is unchanged.  Similarly, the triangle maps $F_{\alpha,\beta,\gamma}$ now have the form $\sum_{i\geq 0} F_{\alpha,\beta,\gamma,i}$ where $F_{\alpha,\beta,\gamma,i}$ counts curves with $i$ double points.

\begin{figure}
\begin{picture}(0,0)%
\includegraphics[scale=0.8]{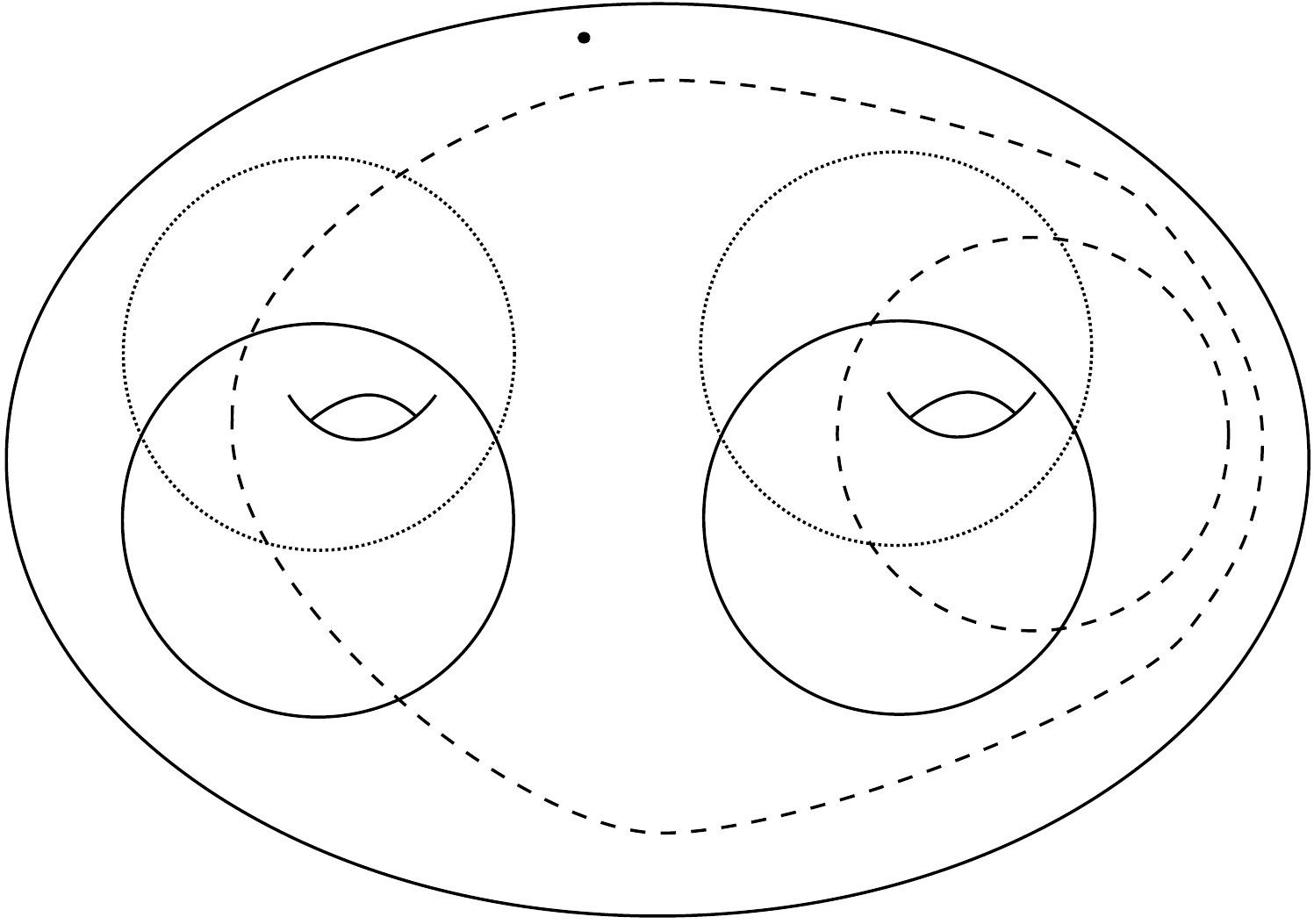}
\end{picture}%
\setlength{\unitlength}{0.08bp}%
\begin{picture}(4310,3030)(0,0)%
\small
\put(941,565){\makebox(0,0)[lb]{\smash{\(\beta_1\)}}}%
\put(286,1880){\makebox(0,0)[lb]{\smash{\({\beta_1}'\)}}}%
\put(2581,251){\makebox(0,0)[lb]{\smash{\(\beta_1^H\)}}}%
\put(3706,949){\makebox(0,0)[lb]{\smash{\(\beta_2^H\)}}}%
\put(2455,707){\makebox(0,0)[lb]{\smash{\({\beta_2}'\)}}}%
\put(2403,2369){\makebox(0,0)[lb]{\smash{\({\beta_2}\)}}}%
\put(1807,2878){\makebox(0,0)[lb]{\smash{\(\mathfrak{z}\)}}}%
\put(358,1563){\makebox(0,0)[lb]{\smash{\(\theta_1\)}}}%
\put(2264,1570){\makebox(0,0)[lb]{\smash{\(\theta_2\)}}}%
\put(1668,1587){\makebox(0,0)[lb]{\smash{\(\eta_1\)}}}%
\put(3567,1608){\makebox(0,0)[lb]{\smash{\(\eta_2\)}}}%
\put(3501,863){\makebox(0,0)[lb]{\smash{\(\eta_2^H\)}}}%
\put(1251,605){\makebox(0,0)[lb]{\smash{\(\eta_1^H\)}}}%
\put(716,1987){\makebox(0,0)[lb]{\smash{\(\theta_1^H\)}}}%
\put(2747,2012){\makebox(0,0)[lb]{\smash{\(\theta_2^H\)}}}%
\put(2766,1137){\makebox(0,0)[lb]{\smash{\({\theta_2}'\)}}}%
\put(788,1165){\makebox(0,0)[lb]{\smash{\({\theta_1}'\)}}}%
\put(1264,2511){\makebox(0,0)[lb]{\smash{\({\eta_1}'\)}}}%
\put(3484,2272){\makebox(0,0)[lb]{\smash{\({\eta_2}'\)}}}%
\put(1331,1103){\makebox(0,0)[lb]{\smash{\(D_1\)}}}%
\put(735,970){\makebox(0,0)[lb]{\smash{\(D_2\)}}}%
\put(596,1586){\makebox(0,0)[lb]{\smash{\(D_3\)}}}%
\put(735,2220){\makebox(0,0)[lb]{\smash{\(D_4\)}}}%
\put(1258,2109){\makebox(0,0)[lb]{\smash{\(D_5\)}}}%
\put(1979,1739){\makebox(0,0)[lb]{\smash{\(D_6\)}}}%
\put(3269,1174){\makebox(0,0)[lb]{\smash{\(D_7\)}}}%
\put(2740,951){\makebox(0,0)[lb]{\smash{\(D_8\)}}}%
\put(2555,1565){\makebox(0,0)[lb]{\smash{\(D_9\)}}}%
\put(2601,2158){\makebox(0,0)[lb]{\smash{\(D_{10}\)}}}%
\put(3289,2024){\makebox(0,0)[lb]{\smash{\(D_{11}\)}}}%
\put(3739,1543){\makebox(0,0)[lb]{\smash{\(D_{12}\)}}}%
\put(1476,273){\makebox(0,0)[lb]{\smash{\(D_{13}\)}}}%
\end{picture}%
\caption{A handleslide in genus two, with notation}
\label{Figure:dblhndl}
\end{figure}
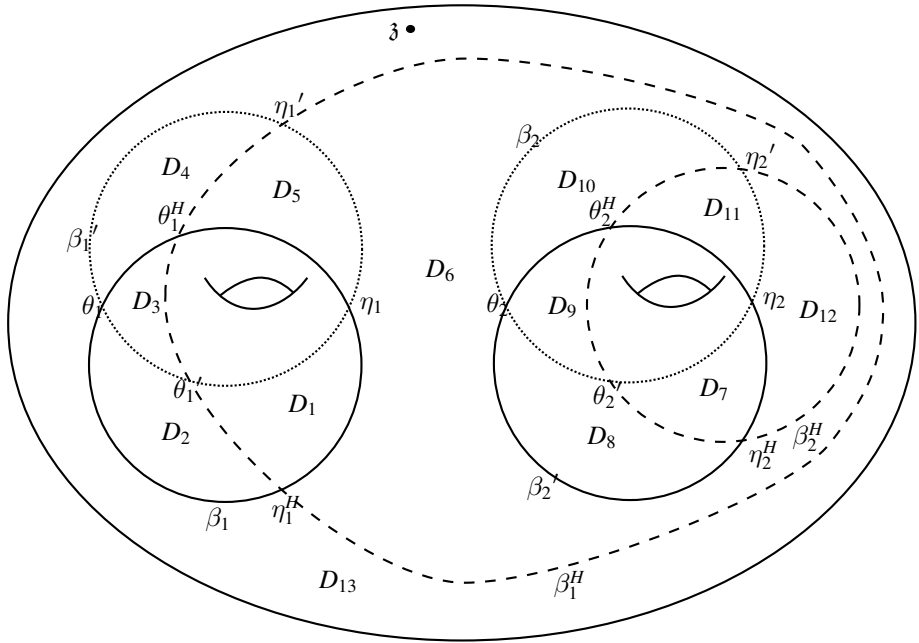

Handleslide invariance is proved as in \fullref{Section:Handleslides}, except that some of the model calculations are more complicated.  We give those computations here.  Since we are considering $\widehat{HF}$, it suffices to study the case $g=2$.  Let $\vb$, $\vb'$ and $\vb^H$ be the curves shown in \fullref{Figure:dblhndl}.  (This notation agrees with the notation of \fullref{Section:Handleslides}.)  The proofs that for an appropriate choice of orientation system the generators $\vec{\theta}_{\beta,\beta'}=\{\theta_1,\theta_2\}$, $\vec{\theta}_{\beta^H,\beta'}=\{\theta'_1,\theta'_2\}$ and $\vec{\theta}_{\beta,\beta^H}=\{\theta_1^H,\theta_2^H\}$ shown in \fullref{Figure:dblhndl} are cycles are the same as in \fullref{Section:Handleslides}, with the additional observation that immersed curves can only contribute to $\bdy_{big}$ if the grading difference is at least three.

However, the proof that $F_{\beta,\beta',\beta^H}\left(\vec{\theta}_{\beta,\beta^H}\otimes\vec{\theta}_{\beta^H,\beta'}\right)=\vec{\theta}_{\beta,\beta'}$  (\fullref{Prop:HSB}) is somewhat more involved.  Since triangle maps count curves in even--dimensional moduli spaces, it is conceivable that the coefficient of $t\{\eta_1,\eta_2\}$ is nonzero in $F_{\beta,\beta',\beta^H}\left(\vec{\theta}_{\beta,\beta^H}\right.$\break $\left.\otimes\vec{\theta}_{\beta^H,\beta'}\right)$.  Again, with notation as in \fullref{Figure:dblhndl}, we need to show that none of the following domains correspond to moduli spaces containing immersed curves:  
$D_1+D_2+D_3+D_7+D_8+D_9$, $D_1+D_2+D_3+D_9+D_{10}+D_{11}$, $D_1+D_2+D_3+D_7+D_{11}+D_{12}$, $D_3+D_4+D_5+D_7+D_8+D_9$, $D_3+D_4+D_5+D_9+D_{10}+D_{11}$, $D_3+D_4+D_5+D_7+D_{11}+D_{12}$, $D_1+D_5+D_6+D_7+D_8+D_9+D_{10}+D_{11}+D_{12}$, $D_1+D_5+D_6+D_7+D_{10}+2D_{11}+2D_{12}$, $D_1+D_5+D_6+D_7+2D_8+2D_9+D_{10}$, $D_1+D_5+D_6+D_9+2D_{10}+2D_{11}+D_{12}$, $D_1+D_5+D_6+D_8+2D_9+2D_{10}+D_{11}$, $D_1+D_5+D_6+2D_7+2D_8+D_9+D_{12}$, and $D_1+D_5+D_6+2D_7+D_8+D_{11}+2D_{12}$.  (These are the only positive domains in $\pi_2(\vec{\theta}_{\beta,\beta^H},\vec{\theta}_{\beta^H,\beta'}, \vec{\eta})$, where $\vec{\eta}=\{\eta_1,\eta_2\}$.)

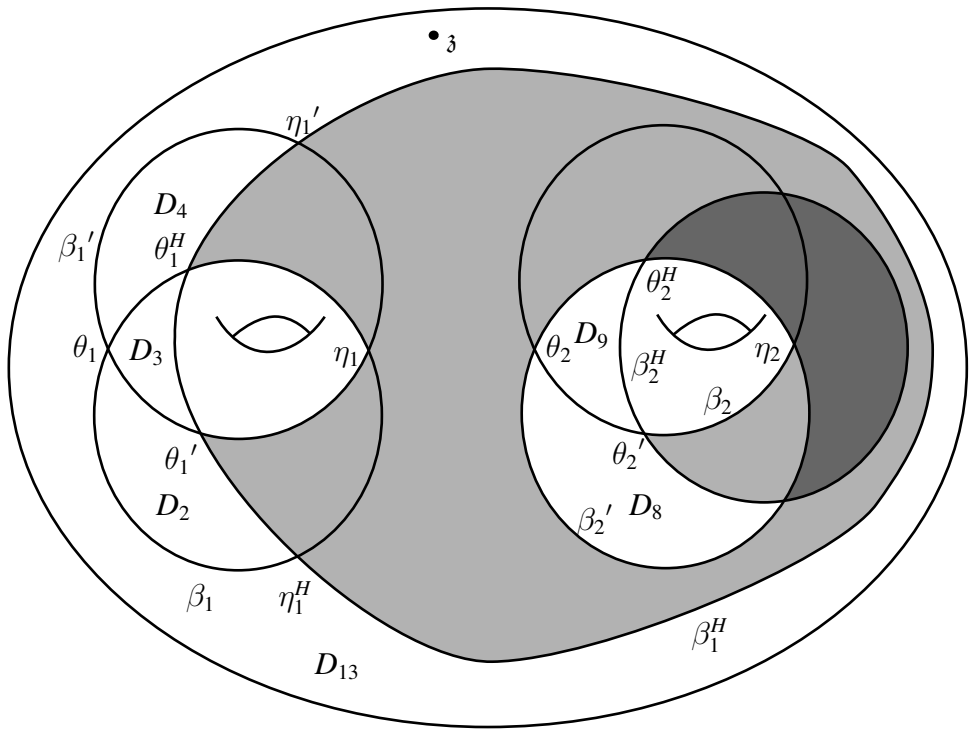
\begin{figure}
\begin{center}
\begin{picture}(0,0)%
\includegraphics{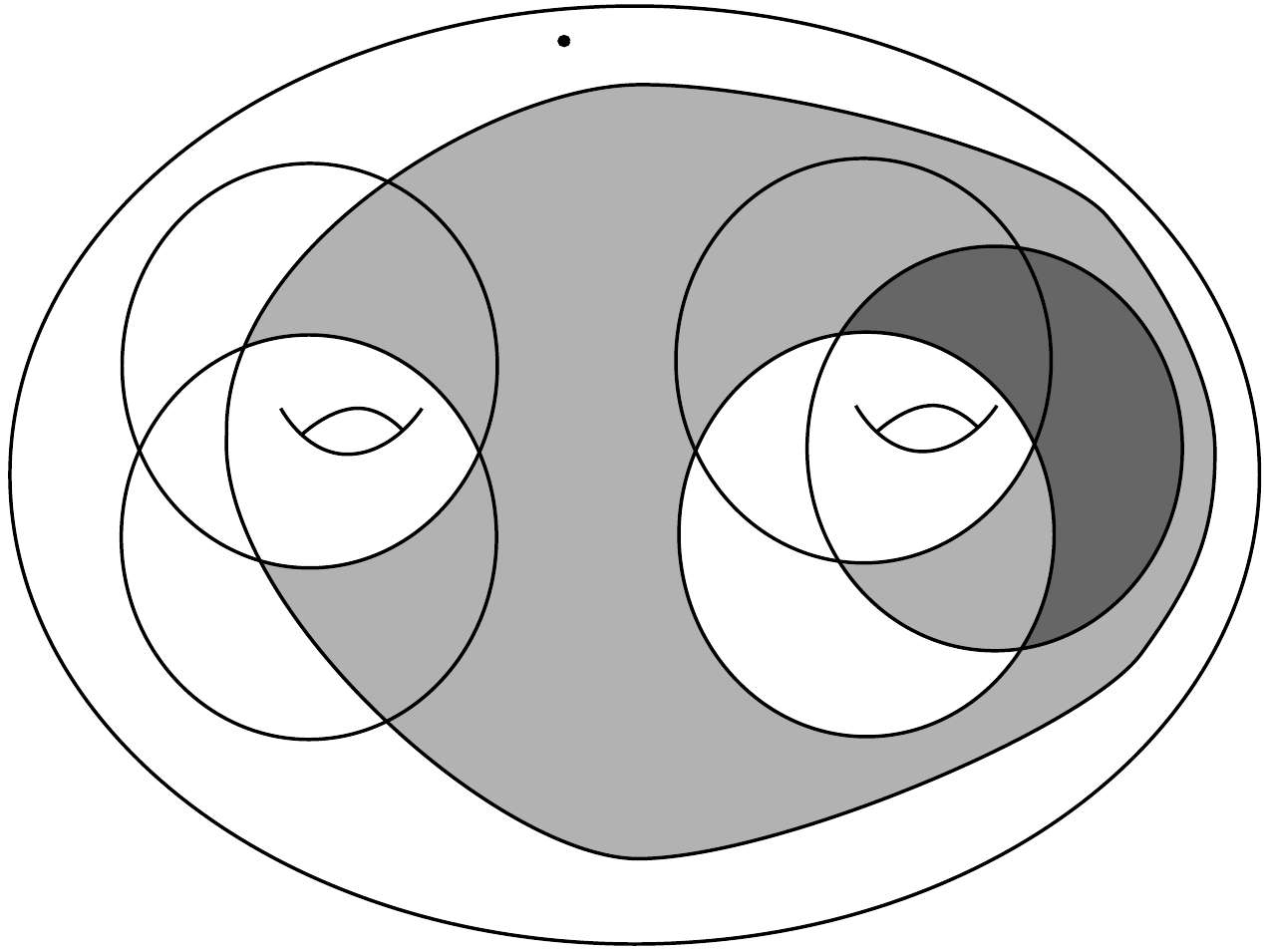}%
\end{picture}%
\setlength{\unitlength}{0.1bp}%
\begin{picture}(3650,2750)(0,0)%
\small
\put(1673,2588){\makebox(0,0)[lb]{\smash{\(\mathfrak{z}\)}}}%
\fontsize{12}{14.4}\selectfont%
\put(689,480){\makebox(0,0)[lb]{\smash{\(\beta_1\)}}}%
\put(207,1798){\makebox(0,0)[lb]{\smash{\({\beta_1}'\)}}}%
\put(2594,336){\makebox(0,0)[lb]{\smash{\(\beta_1^H\)}}}%
\put(2365,1349){\makebox(0,0)[lb]{\smash{\(\beta_2^H\)}}}%
\put(2160,767){\makebox(0,0)[lb]{\smash{\({\beta_2}'\)}}}%
\put(262,1419){\makebox(0,0)[lb]{\smash{\(\theta_1\)}}}%
\put(2045,1418){\makebox(0,0)[lb]{\smash{\(\theta_2\)}}}%
\put(1244,1404){\makebox(0,0)[lb]{\smash{\(\eta_1\)}}}%
\put(2831,1420){\makebox(0,0)[lb]{\smash{\(\eta_2\)}}}%
\put(1037,478){\makebox(0,0)[lb]{\smash{\(\eta_1^H\)}}}%
\put(569,1789){\makebox(0,0)[lb]{\smash{\(\theta_1^H\)}}}%
\put(2423,1662){\makebox(0,0)[lb]{\smash{\(\theta_2^H\)}}}%
\put(2294,1015){\makebox(0,0)[lb]{\smash{\({\theta_2}'\)}}}%
\put(607,995){\makebox(0,0)[lb]{\smash{\({\theta_1}'\)}}}%
\put(1065,2276){\makebox(0,0)[lb]{\smash{\({\eta_1}'\)}}}%
\put(576,821){\makebox(0,0)[lb]{\smash{\(D_2\)}}}%
\put(476,1407){\makebox(0,0)[lb]{\smash{\(D_3\)}}}%
\put(568,1959){\makebox(0,0)[lb]{\smash{\(D_4\)}}}%
\put(2153,1471){\makebox(0,0)[lb]{\smash{\(D_9\)}}}%
\put(1174,221){\makebox(0,0)[lb]{\smash{\(D_{13}\)}}}%
\put(2358,824){\makebox(0,0)[lb]{\smash{\(D_8\)}}}%
\put(2639,1222){\makebox(0,0)[lb]{\smash{\({\beta_2}\)}}}%
\end{picture}%
\end{center}
\caption{The domain $D$}
\label{Figure:HSdblShaded}
\end{figure}

In the first seven of these domains, there are no coefficients larger than one.  It follows that there can be no immersed curve in the corresponding moduli space.  The last six of these domains exhibit enough symmetry that it suffices to consider one of them.  We will focus on  $D=D_1+D_5+D_6+D_7+D_{10}+2D_{11}+2D_{12}$.  It is shown shaded in \fullref{Figure:HSdblShaded}.

We compute the Euler characteristic of an embedded curve representing $D$.  Observe that $D=(D_5 + D_{11}+D_{12}+D_7)+(D_1+D_6+D_{10}+D_{11}+D_{12})=E+F$.  Here, $E=(D_5+D_{11}+D_{12}+D_7)\in\pi_2(\vec{\theta}^H,\{\eta'_1,\theta'_2\},\vec{\eta})$ and $F=(D_1+D_6+D_{10}+D_{11}+D_{12})\in\pi_2(\{\theta'_2,\eta'_1\},\vec{\theta}')$.  The domain $E$ is represented by a pair of embedded disks.  At least for some (nonempty, open set of) almost complex structures, the domain $F$ is represented by an embedded annulus.  Gluing representatives for $E$ and $F$, it follows that the Euler characteristic for an embedded representative of $D$ is $0$.

It follows that the source of a curve with double points representing $D$ must have $\chi\geq 2$.

It is clear from considering the boundary of $D$ that any representative for $D$ satisfying (\textbf{M0})--(\textbf{M5}) must be connected.  But there are no connected surfaces with boundary with $\chi\geq 2$.  This proves the non--embedded moduli spaces are empty, so 
$$F_{\beta,\beta',\beta^H}\left(\vec{\theta}_{\beta,\beta^H}\otimes\vec{\theta}_{\beta^H,\beta'}\right)=\vec{\theta}_{\beta,\beta'}$$
as desired.

Plugging this computation into the proof of handleslide invariance in \fullref{Section:Handleslides} proves handleslide invariance for $\widehat{HF}_{big}$.  Stabilization invariance of $\widehat{HF}_{big}$ follows in the same simple way as stabilization invariance of $\widehat{HF}$.  Note that all of the maps we have used are maps of $\bZ[[t]]$--modules.  It follows that the chain complex $\widehat{CF}_{big}$, up to chain homotopy equivalence over $\bZ[[t]]$, is an invariant of $Y$.

Unfortunately, I have been unable to compute a single case in which $\widehat{CF}_{big}$ is not homotopy equivalent to a complex in which all higher differentials vanish.  

I suspect that one could similarly elaborate $HF^\infty$ and
$HF^\pm$, but have not done the computations necessary to
establish handleslide invariance, and am mildly concerned that there
may be subtleties in the proof of stabilization invariance.  If evidence appears to suggest that these variants of $HF^\infty$, $HF^\pm$ or, for that matter, $\widehat{HF}$ would be new or interesting then they will become the subject of a future paper.

Here are several other elaborations, also inspired by Gromov--Witten
theory.  For convenience, I will formulate only the
$\widehat{HF}$--analogs, but in these cases analogs of the other
theories present few added difficulties.  (One does have to
  deal with annoying curves, but these can be addressed similarly to
  the way we did in \fullref{Section:ChainComplexes}.)  To start, fix a homology class
$[K]\in H_1(Y)$, a knot $K\mapsinto\Sigma$ representing
$[K]$, and a point $s_0\in[0,1]$.  Let $\widehat{CF}_{[K]}$ be
freely generated over $\bZ[t]$
by the intersection points.  Define
$\bdy_{[K]}=\bdy_{0}+\frac{1}{1!}t\bdy_1+\frac{1}{2!}t^2\bdy_2+\cdots$,
where the
coefficient of $\vy$ in $\bdy_i\vx$ counts holomorphic curves
with $i$ marked points in
homology classes $A$ with $ind(A)=1$ so that each marked point is
mapped to $K\times\{s_0\}\times\bR$.  Then the standard proof shows
that $\bdy_{[K]}^2=0$.

The resulting chain complexes, up to chain homotopy equivalence over
$\bZ[t]$, are indeed invariants of $(Y,[K])$.
However, as pointed out to me by M. Hutchings, the resulting chain complex
can be reconstructed from the chain complex with ``totally twisted''
(group ring) coefficients.  (This is not at all surprising; the number
of times a holomorphic curve intersects $K\times\{s_0\}\times\bR$
depends only on the homology class of the curve.  Since this
construction imitates pulling back a 2--dimensional cohomology class
from $W$, this is somewhat analogous to the divisor equation in
Gromov--Witten theory.)

One could try an analogous construction by forcing points to be mapped
to $\{p_0\}\times\{s_0\}\times\bR$ for some choice of
$p_0\in\Sigma$.  The result is again an invariant, and independent
of $p_0$.  However, taking
$p_0\in\Sigma\setminus(\ba\cup\bb)$ and considering $s_0=0$, one sees that all of
the higher differentials vanish.

This exhausts the obvious cohomology classes one could pull back, so
the next idea is to try descendent classes.  For example, let
$\widehat{CF}_{desc}$
be freely generated over $\bZ[[t]]$ by the intersection points.
Define
$\bdy_{desc}=\bdy_0+\frac{1}{1!}t\bdy_1+\frac{1}{2!}t^2\bdy_2+\cdots$
where $\bdy_i$ counts holomorphic curves $u$ from Riemann surfaces $S$
with $i$ marked points $p_1,\cdots,p_i$ such that
\begin{itemize}
\item The homology class of $u$ has index $2i+1$.
\item For each $i$, $u(p_i)\in K\times\{s_0\}\times\bR$.
\item For each $i$, $(\pi_\Sigma\circ u)'(p_i)=0$.
\end{itemize}
Again, the proof that $\bdy_{desc}^2=0$ is standard.  That the chain
complex up to chain homotopy equivalence (over $\bZ[[t]]$) is an
invariant of $(Y,[K])$ is almost identical to the proof of invariance of the $H_1/Tors$--action on Heegaard
Floer homology.  There are no new computations that need to be done.

Of course, this is only one piece of a much bigger chain complex one
could consider, where one keeps track of higher branching of
$\pi_\Sigma\circ u$ and several different elements of $H_1(Y)$.
One could also allow curves with branching of $\pi_\Sigma\circ u$
at prescribed $s_0$ and arbitrary point in $\Sigma$.
Again, invariance of the bigger complex is free.

(One could also consider curves with prescribed
branching of $\pi_\bD\circ u$, but taking $s_0=0$ our proof that
bubbling is impossible forces higher differentials of this kind to be trivial.)

All of these deformed complexes should have roughly the same formal
properties as Heegaard Floer homology.
Unfortunately, even with this apparent wealth of additional
information I have been unable to find nontrivial examples.  That is,
while there are examples where higher differentials are nontrivial, I do
not know examples which are not chain homotopy equivalent to complexes
in which all higher differentials vanish.  Hopefully this is for lack of creativity or perseverance on my part.

One might hope to construct further elaborations by pairing with cohomology
classes of the space of maps $S\to W$, but by
\fullref{Prop:Contractible} this space does not have
interesting topology.
Finally, one might be able to obtain invariants by pulling back
cohomology classes from the space of holomorphic maps to a disk.  Doing so in a useful way, however, would require a better understanding of the cohomology of the moduli space of maps to a disk than I presently posses.

\subsection{Relationship with Taubes' program}
We conclude with a few remarks about a likely relationship between
Heegaard Floer homology as formulated in this paper and Taubes'
program to understand holomorphic curves in 4--manifolds with
singular symplectic forms.  

First, a one paragraph sketch of Taubes' idea.
Any 4--manifold $M$ with $b_2^+>0$ can be endowed with a closed
two--form $\omega$ which is nondegenerate in the complement of a
collection of 
circles, and degenerates in a controlled way near the circles;
see~\cite{Taubes} for further information and references.  The program
is then to fix a complex structure $J$ in the
complement of the singular
circles, adjusted to $\omega$, and obtain smooth invariants by
studying $J$--holomorphic curves of
finite $\omega$--energy in $M$.

The Floer homology associated to Taubes' program would be structured as
follows.  For a three--manifold $Y$, one chooses a closed 2--form
$\omega$ on $Y\times\bR$ which is nondegenerate on the complement
of certain lines $\{p_i\}\times\bR$.  One would then choose a
translation invariant complex structure $J$ on
$Y\times\bR$ adjusted to $\omega$ and study $J$--holomorphic
curves with some specified asymptotics at the singular lines and at $Y\times\{\pm\infty\}$.

Fix a Morse function $f$ on $Y$, and a metric on $Y$.  Then, the
form $df\wedge dt+\star df$
is a closed 2--form on $Y\times\bR$, which is nondegenerate on the
complement of $\{\textrm{critical points of $f$}\}\times\bR$.
(Here, $\star$ means the Hodge star on $Y$, not on
$Y\times\bR$.)  One possible compatible complex structure pairs
$dt$ with $\nabla f$ and is given by rotation by $\pi/2$ on
$\ker(df)$.  We will call this the complex structure induced by $f$.

But this setting bears a close resemblance to $(W,C_\alpha\cup C_\beta)$.  Assume
$f$ was self--indexing.  Then, view $\Sigma\times[0,1]$ as the
slice $f^{-1}(3/2-\epsilon,3/2+\epsilon)$ of $Y$.  Then the
complex structure induced by $f$ satisfies
(\textbf{J1})--(\textbf{J5}).  Further, in the limit
$\epsilon\to1/2$, the boundary conditions we impose become certain
asymptotic conditions on the holomorphic curves at the singular
lines.  So, our setup fits quite nicely in Taubes' picture.

A serious difficulty in studying holomorphic curves in $Y\times\bR$
is understanding the asymptotics at the singular lines.  By studying
only a middle slice of $Y$, Heegaard Floer homology neatly avoids
this issue.  Still, it would be nice to be able to work in
$Y\times\bR$; for instance, this would probably illuminate the
proof of  handleslide invariance.

\appendix
\section{Gluing lemmas}\label{Section:Gluing}
\subsection{Statement of results}
As seems conventional in the subject, we relegate the misery
called ``gluing lemmas'' to the appendix.

In the following, by ``symplectic manifold'' we mean a symplectic
manifold with cylindrical ends, and by ``complex structure'' we mean a
symmetric almost complex structure, as defined in~\cite[Section 2]{Ya2}.  We will
always assume that Reeb orbits and chords corresponding to the almost
complex structure and Lagrangian submanifolds in question are isolated, as the Morse--Bott case requires extra
work, and all Reeb chords in this paper are isolated.  When we refer
to a holomorphic curve we always mean ones with finite energy in the
sense of~\cite[Section 6.1]{Ya2}.
With these conventions, the gluing results used in this paper are
\begin{Prop}\label{Gluing:Cylindrical}
Let $(M_1,M_2)$ be a chain of symplectic manifolds
  (cf~\cite[Section 1.6]{Ya1}) with $M_2$ cylindrical.
  Let $(u_1,u_2)$ be a 
  height two holomorphic building in $(M_1,M_2)$.  Assume that the
  complex structures on $M_1$ and $M_2$ are chosen so that the $\dbar$
  operator is transverse to the zero--section at $u_1$ and
  $u_2$.  Then there is a neighborhood of $(u_1,u_2)$ in the space
  of (height one or two) holomorphic buildings diffeomorphic to
  $\bR^{\ind(u_1)+\ind(u_2)-1}\times(0,1].$  If $M_1$ is also
  cylindrical then the same statement holds with
  $\bR^{ind(u_1)+ind(u_2)-1}$ replaced by $\bR^{ind(u_1)+ind(u_2)-2}$.
  Both statements remain true in the relative case.
\end{Prop}

\begin{Prop}\label{Gluing:Split}
Let $(M_1,M_2)$ be a chain of symplectic manifolds
  (cf~\cite[Section 1.6]{Ya1}) obtained by splitting
  a symplectic manifold $M$ along a co--oriented hypersurface $H$
  as in~\cite[Section 3.4]{Ya2}.
  Let $(u_1,u_2)$ be a 
  height two holomorphic building in $(M_1,M_2)$.  Assume that the
  complex structure on $M$ is chosen so that the $\dbar$
  operator is transverse to the zero--section at $u_1$ and
  $u_2$.  Let $J_R$ denote the complex structure on $M$ obtained
  by inserting a neck of length $R$ along $H$.  Let $\cM_R$
  denote the space of $J_R$--holomorphic curves in $M$  Then there is a
  neighborhood of $(u_1,u_2)$ in the space $\overline{\bigcup_R\cM_R}$
  diffeomorphic to
  $\bR^{\ind(u_1)+\ind(u_2)}\times(0,1].$
  The statement remains true in the relative case.
\end{Prop}

\begin{Prop}\label{Gluing:Twins}\fullref{Gluing:Split} remains
  true in the case of height one holomorphic twin towers studied
  in \fullref{Section:Stabilization}.  That is, with notation as
  in \fullref{Section:Stabilization}, for $\ind(A)=1$,
  we have $\#\hcM^A_{W'_R}=\#\hcM^A_W$.
\end{Prop}

The special case of \fullref{Gluing:Cylindrical} when $M_1$
is the symplectization of a contact manifold follows from work of Bourgeois \cite[Corollary 5.7]{bourgeois}.  His argument, however, extends
without essential changes to prove the more general results stated
here.  To demonstrate this, we imitate his argument to prove \fullref{Gluing:Split}
in the relative case.  \fullref{Gluing:Cylindrical} is similar
but marginally less complicated.  After proving
\fullref{Gluing:Split} we discuss how the proof needs to be
modified to prove \fullref{Gluing:Twins}.

\subsection{The splitting process}\label{Subsection:Splitting}
In this subsection we describe the process of splitting a symplectic
manifold along a hypersurface.  More details can be found
in~\cite[Section 3.4]{Ya2}.  We describe the splitting process (tersely) here
partly because the relative case is hardly discussed in~\cite{Ya2} but
mostly just to fix notation.

Let $(M,\omega)$ be a
symplectic manifold, $J$ a complex structure on $M$ compatible
with $\omega$, and $H\subset M$ a co--oriented compact
hypersurface.  We also assume that we are given a vector field
$\vec{R}\in \ker(\omega|_H)$ on $H$ so that the associated
cylindrical structure on $H\times \bR$ is symmetric
(see~\cite[page 802]{Ya2}).  Let $\alpha$ be the $1$--form
$\omega(J\vec{R},\cdot)$ on $H$.  Then a neighborhood of $H$ is
symplectomorphic to $\left((-\epsilon,\epsilon)\times
  H,\pi_H^*\omega|_H+d(t\alpha)\right)$ so that $J\vec{R}$ is
identified with $\frac{\partial}{\partial t}|_{H=0}$
(cf~\cite[page 806]{Ya2}).  We will assume that $J$ is preserved by
the flow in the $t$--direction.  This can always be arranged by
perturbing $J$ near $H$.

In this case, let $M^\circ$ denote the manifold $M\setminus
(-\epsilon/2,\epsilon/2)\times H$ with $\bdy M^\circ=H'\cup H''$.
Form 
$$M_R=M^\circ\bigcup_{\substack{H'=\{-R\}\times
      H\\ H''=\{R\}\times H}}[-R,R]\times H.
$$
The complex structure on $M_R$ is given by $J$ on $M^\circ$ and
on $[-R,R]\times H$ agrees with $J$ on $H=\{0\}\times H$ and is
translation invariant.  The symplectic form on $M_R$ is given by
$\frac{2R}{\epsilon}\omega$ on $M^\circ$ and by $\pi_H^*\omega|_H
+ d(t\alpha)$ on $[-R,R]\times H$.

Similarly, let
$$M_\infty = (-\infty,0]\times H\bigcup_{\{0\}\times H=H'} M^\circ
\bigcup_{H''=\{0\}\times H} [0,\infty)\times H.
$$
As before, the complex structures on the ends $(-\infty,0]\times H$
and $[0,\infty)\times H$ are induced by translation invariance and
$J|_H$; on $M^\circ$ the complex structure is given by $J$.  On
$M^\circ$, the symplectic form is just $\omega$.  On
$(-\infty,0]\times H$ and $[0,\infty)\times H$ it is given by
$\pi_H^*\omega|_H+d(t\alpha).$

Suppose $L\subset M$ is a Lagrangian submanifold intersecting $H$ transversally.  Perturbing $L$ slightly we can assume that in the tubular neighborhood $(-\epsilon,\epsilon)\times H$ used to perform the stretching, $L$ has the form $(-\epsilon,\epsilon)\times (L \cap H)$.  (In fact, this can be achieved by a Hamiltonian deformation of $L$.)  We will call this property being ``cylindrical'' near $H$.   Then, the splitting procedure just described gives in an obvious way Lagrangian submanifolds $L_R$ in $M_R$ and $L_\infty$ in $M_\infty$.

Now, fix 
\begin{itemize}
\item a symplectic manifold $(M,\omega)$ with cylindrical ends
\item a symmetric, cylindrical almost complex structure $J$ on $M$ adjusted to $\omega$
\item a compact, co--oriented hypersurface $H\subset M$ and
\item a closed Lagrangian submanifold $L\subset M$ which is
  cylindrical near $H$, such that in each end
  $[0,\infty)\times\tilde{E}$ (respectively
  $(-\infty,0]\times\tilde{E}$) of $M$, $L$ has the form
  $[0,\infty)\times\tilde{L}_E$ (respectively
  $(-\infty,0]\times\tilde{L}_E$) with $T\tilde{L}_E\subset
  T\tilde{E}\cap JT\tilde{E}$. 
\end{itemize}

Let $\tilde{E}_1,\cdots,\tilde{E}_N$ denote the ends of $M$, so
near each $\tilde{E}_i$, $M$ is modeled on $I_i\times \tilde{E}_i=:E_i$,
where $I_i$ is a semi--infinite interval.  The isomorphism $E_i=\tilde{E}_i\times
I_i$ is considered fixed as part of the definition of a cylindrical
complex structure.  If $I_i$ is
$[r,\infty)$ we call the end $E_i$ positive; otherwise we call
$E_i$ negative.  Let $\sigma_i$ be $+1$ if $E_i$ is positive
and $-1$ if $E_i$ is negative.  

Let $M_R$ and $M_\infty$ denote the manifolds obtained by
splitting $M$ along $H$.  The ends of $M_R$ are
$\tilde{E}_1,\cdots,\tilde{E}_N$.  The ends of $M_\infty$ are
$\tilde{E}_0^+,\tilde{E}_0^-,\tilde{E}_1,\cdots,\tilde{E}_N$, where
$\tilde{E}_0^\pm$ correspond to $H$.  Note that
$\tilde{E}_0^+=\tilde{E}_0^-$, where this identification respects
the image of $L$ and also the Reeb field.  We sometimes write
$\tilde{E}_0$ for $\tilde{E}_0^\pm$.

For $i=0,\cdots,N$, let $L\cap
E_i=\tilde{L}_i\times I_i$.
  Let $\tau_i$ denote the coordinate on $I_i$,
and $\vec{R}_i=J\bdy_{\tau_i}$.  Let $\Gamma_i^o$ denote the space
of closed orbits of $\vec{R}_i$, $\Gamma_i^c$ the space of
$\vec{R}_i$--chords, ie, flows of $\vec{R}_i$ starting and
ending on $\tilde{L}_i$.  Let
$\Gamma_i=\Gamma_i^o\cup\Gamma_i^c$.

We will assume throughout that the $\Gamma_i$ are discrete.

Fix Reeb chords/orbits
$\gamma_{0,1},\cdots,\gamma_{0,n_0}\in\Gamma_0$,
$\gamma_{i,1}^+,\cdots,\gamma_{i,n_i^+}^+\in\Gamma_i$ ($i=1,\cdots,N$),
and $\gamma_{i,1}^-,\cdots,\gamma_{i,n_i^-}\in\Gamma_i$
($i=1,\cdots,N$).  For notational convenience, let
$\gamma_{0,j}^+=\gamma_{0,j}^-=\gamma_{0,j}$. 
Assume
$\gamma_{i,1}^\pm,\cdots,\gamma_{i,m_i^\pm}^\pm\in\Gamma_i^{o}$ are closed Reeb
orbits and
$\gamma_{i,m_i^\pm+1},\cdots,\gamma_{0,n_i^\pm}\in\Gamma_i^{c}$ are
Reeb chords.
For convenience and because it is the only case relevant to this paper
we will assume that all the $\gamma_{0,k}$ and $\gamma_{i,k}^\pm$
are simple Reeb chords / orbits.

Fix smooth surfaces $\Sigma^\pm$ with boundary and labeled
punctures
$$p_{0,1}^\pm,\cdots,p_{0,n_0}^\pm,\cdots,p_{N,1}^\pm,\cdots,p_{N,n_N^\pm}^\pm.$$
Let $\mathcal{T}^\pm$ be smooth families of conformal structures on $\Sigma$
so that $\mathcal{T}^\pm$ surjects onto an open set
in the moduli space of Riemann surfaces.  (That
  is, points in $\mathcal{T}^\pm$ are honest surfaces with conformal
  structures and marked points, not equivalence classes of such.  The
  map from $\mathcal{T}^\pm$ to moduli space need not be injective.)
  Choose
$\mathcal{T}$ so that there are small neighborhoods of the punctures
in $\Sigma$ in which the conformal structure is constant.  Fix
$p>2$, $k\geq 1$.  For convenience later, we will also fix
$S_0^\pm\in\mathcal{T}^\pm$.

For $S^\pm\in\mathcal{T}^\pm$, let
$\mathcal{B}^\pm_{S^\pm}=\mathcal{B}^{p,d,\pm}_k$ denote the Banach
manifold comprising $W^p_k$--maps\break $(S^\pm,\bdy S^\pm)\to (M_\infty,L)$
converging  to $\gamma_{i,k}^\pm$ at $p_{i,k}^\pm$ ($i=0,\cdots,N$,
$k=1,\cdots,n_i$) in $W^{p,d}_k$.  That is: 

Any $\vec{R}$--chord or orbit $\gamma$ has some period $T_\gamma$ defined by $T_\gamma=\frac{1}{2\pi}\int_\gamma\omega(\bdy_\tau,\cdot)$.  If $\gamma$ is a Reeb orbit choose a tubular neighborhood $N_\gamma$ of $\gamma$ invariant under the $\vec{R}$--flow and a diffeomorphism $\phi$ from a neighborhood of $S^1\times \vec{0}$ in $S^1\times\bR^{2n-2}$ to $N_\gamma$ such that
\begin{enumerate}
\item $\phi(S^1\times \vec{0})=\gamma$
\item the pushforward under $\phi$ of $T_\gamma$ times the unit tangent vector to $S^1\times \vec{0}$ is $\vec{R}$ and
\item $\phi^*\omega|_{S^1\times\vec{0}}=\omega_0|_{S^1\times\vec{0}}$, where $\omega_0$ denotes the standard symplectic form on $\bR^{2n-2}$.
\end{enumerate}
If $\gamma$ is a Reeb chord choose a tubular neighborhood $N_\gamma$ of $\gamma$ invariant under the $\vec{R}$--flow and a diffeomorphism $\phi$ from a neighborhood of $[0,1]\times\vec{0}$ in $[0,1]\times\bR^{2n-2}$ to $N_\gamma$ such that
\begin{enumerate}
\item $\phi([0,1]\times \vec{0})=\gamma$
\item the pushforward under $\phi$ of $T_\gamma$ times the unit tangent vector to $[0,1]\times \vec{0}$ is $\vec{R}$ and
\item $\phi^*\omega|_{[0,1]\times\vec{0}}=\omega_0|_{[0,1]\times\vec{0}}$, where $\omega_0$ denotes the standard symplectic form on $\bR^{2n-2}$.
\end{enumerate}
(Compare~\cite[Lemma A.1]{Ya2}, but note that our $N_\gamma$ is not the same as their $N$.)  Let $\theta_\gamma$ parametrize $S^1$ in the first case, or $[0,1]$ in the second case, by arc length.  Let $z=(x_1,\cdots,x_{n-1},y_1,\cdots,y_{n-1})$ be the standard coordinates on $\bR^{2n-2}$.

Via $\phi$, $\theta_\gamma$ and $z$ become coordinates on $N_\gamma$.

For $S^\pm\in\mathcal{T}^\pm$ we fix a holomorphic identification of
a neighborhood of each $p_{i,j}^\pm$ with $[0,\infty)\times S^1$
if $p_{i,j}^\pm$ is in the interior of $S^\pm$ and with
$[0,\infty)\times[-\pi,\pi]$ if $p_{i,j}^\pm$ is in the boundary
of $S^\pm$.  Let $(t,s)$ denote the coordinates near
$p_{i,j}^\pm$ induced by this decomposition.

The statement that a map $u_\pm\co S^\pm\to M_\infty$ is asymptotic to
$\gamma_{i,j}^\pm$ at $p_{i,j}^\pm$ means that for some $\tau_{i,j,0}^\pm, \theta_{i,j,0}^\pm\in\bR$, $z\circ u_\pm(t,s)\to
0$ as $t\to\infty$, $\tau\circ u_\pm(t,s)-\sigma_iT_{\gamma_{i,j}}t
-\sigma_i\tau_{i,j,0}^\pm\to 0$ as $t\to\infty$ and $\theta\circ
u_\pm(t,s)-s-\theta^\pm_{i,j,0}\to 0$ as $t\to\infty$.  We say that the map $u_\pm$ is
$L_k^{p,d}$ near $p_{i,j}^\pm$ if $z\circ u(t,s)$,  $\tau\circ
u(t,s)-\sigma_iT_{\gamma_{i,j}}t -\sigma_i\tau_{i,j,0}$ and  $\theta\circ u(t,s)-s
-\theta_{i,j,0}^\pm$ are all in $L^{p,d}_k=\{f|f(t,s)e^{d|t|/p}\in
L^p_k\}$.  Note that if $\gamma_{i,j}$ is a Reeb chord then
$\theta_{i,j,0}^\pm=0$ while if $\gamma_{i,j}$ is a closed Reeb orbit
then after replacing $s$ with $s-\theta_{i,j,0}^\pm$ we can assume
that $\theta_{i,j,0}^\pm=0$.  In the future we will assume that $s$
has been chosen so that $\theta_{i,j,0}^\pm=0$.  Replacing some of the
neighborhoods of the punctures with smaller ones if necessary we can
assume that for each $i$ the constants $\tau_{i,j,0}^\pm$ are all the
same, $\tau_{i,j,0}^\pm=\tau_{i,0}$. 

The compactness result~\cite[Theorem 10.2]{Ya2} and asymptotic convergence
result~\cite[Proposition
6.2]{Ya2} imply that there exists a $d>0$ such that any holomorphic
curve of finite energy with the specified asymptotics is in
$\mathcal{B}^{p,d}_k$.  Fix such a $d$, and let
$\mathcal{B}^\pm=\mathcal{B}^{p,d,\pm}_k$.

Similarly, for any Riemannian vector bundle $E$ over $S^\pm$ we can
consider the space of $L^{p,d}_k$--sections of $E$,
$L^{p,d}_k(E)$.  The complex structure $J$ and symplectic form
$\omega$ induce a metric on $M_\infty$, so it makes sense to talk
about the space
$L^{p,d}_{k-1}(\Lambda^{0,1}T^*S^\pm\otimes_Ju^*TM_\infty)$ of
$L^{p,d}_{k-1}$ $(0,1)$--forms on $S^\pm$ with values in $T^*M$.

Now, the spaces $\mathcal{B}_{S^\pm}^\pm$ fit together into a fiber
bundle $\mathcal{B}^\pm$ over $\mathcal{T}^\pm$.  In turn, the spaces
$L^{p,d}_{k-1}(\Lambda^{0,1}T^*S^\pm\otimes_Ju^*TM_\infty)$ fit together
into a vector bundle $\mathcal{E}^\pm$ over $\mathcal{B}^\pm$.
The $\dbar$--operator gives a section $\mathcal{B}^\pm\to
\mathcal{E}^\pm$.  We will assume that this $\dbar$ map is
transversal to the $0$--section.  We let $\mathcal{M}^\pm$ denote
the intersection of $\dbar$ with the $0$--section. 

Recall that we fixed reference Riemann surfaces
$S_0^\pm\in\mathcal{T}^\pm$.  As all Riemann surfaces in $\mathcal{T}^\pm$
have fixed smooth source $\Sigma^\pm$, there is an obvious
projection map $\mathcal{B}^\pm\to\mathcal{B}^\pm_{S_0^\pm}$.  Identify
$\mathcal{B}^\pm$ with
$\mathcal{T}\times\mathcal{B}^\pm_{S_0^\pm}$. 

\medskip\textbf{Remark}\qua  While we have not yet
chosen metrics or norms on the spaces under consideration, when we do,
a figment of the details of our definition is that
the identification of $\mathcal{B}^\pm$ with
$\mathcal{T}\times\mathcal{B}^\pm_{S_0^\pm}$ will not induce
isometries.  However, by choosing $\mathcal{T}$ small we can make the
induced maps arbitrarily close to isometries.

Choosing a connection on $M_\infty$, for instance the Levi--Civita
connection, we obtain also a linearized $\dbar$--operator
$D\dbar\co T_u\mathcal{B}^\pm\to
L^{p,d}_{k-1}(\Lambda^{0,1}T^*S\otimes_J u^*TM)=\mathcal{E}_u$.  (The
space $\mathcal{E}_u$ is identified with the vertical tangent space
to $\mathcal{E}_u$.  The map $D\dbar$ is the composition of
projection onto the vertical tangent space -- defined with the help of
the connection -- with the derivative of the section $\dbar$.
See~\cite[Section 3.3]{MS2}.)   We will also sometimes be interested
in the restriction of the $\dbar$--operator to maps with fixed
source; we denote the $\dbar$--operator on maps with source $S$ by
$\dbar_S$ and its linearization by $D\dbar_S$.

Let $L^{p,d}_k(u_\pm^*TM_\infty,\bdy)$ denote those $L^{p,d}_k$
sections of $u_\pm^*TM_\infty$ which lie in $TL$ over $\bdy S$.

The tangent space at $(S_0^\pm,j)$ to $\mathcal{T}^\pm$ is a
finite--dimensional 
space $V^\pm$.  
The tangent space to $\mathcal{B}^\pm_{S^\pm}$ at
some map $u_\pm$ is 
$$\bC^{\sum_{i=0}^Nm_i^\pm}\oplus\bR^{\sum_{i=0}^N(n_i^\pm-m_i^\pm)}\oplus
L^{p,d}_k(u_\pm^*TM_\infty,\bdy)=F^\pm\oplus L^{p,d}_k(u_\pm^*TM_\infty,\bdy).$$
Fix a vector field $v_{\tau,i,j}^\pm$ (respectively
$v_{\theta,i,j}^\pm$) which agrees with
$\frac{\partial}{\partial\tau}$ (respectively
$\frac{\partial}{\partial\theta}$) near $p_{i,j}^\pm$ and lies in
$L^{p,d}_k\left(u_\pm^*TM_\infty,\bdy\right)$ away from $p_{i,j}^\pm$ for each
interior puncture $p_{i,j}^\pm$ of $S_0^\pm$.
Fix a vector field $v_{\tau,i,j}^\pm$ which agrees with
$\frac{\partial}{\partial\tau}$ near $p_{i,j}^\pm$ and lies in
$L^{p,d}_k\left(u_\pm^*TM_\infty,\bdy\right)$ away from $p_{i,j}^\pm$ for each
boundary puncture $p_{i,j}^\pm$ of $S_0^\pm$.
Then, $F^\pm$ is
$\Span\{v_{\tau,i,j}^\pm,v_{\theta,i,j}^\pm\}\subset\Gamma u^*TM_\infty$.
(Equivalently, $F^\pm$ corresponds to the constants $\tau_{i,j,0}^\pm$
and $\theta_{i,j,0}^\pm$ varying.)  Fix a norm on $F^\pm$
and use this norm to extend the $L^{p,d}_k$--norm
to $F^\pm\oplus L^{p,d}_k(u^*TM_\infty,\bdy)$.

From the decomposition
 $\mathcal{B}^\pm=\mathcal{T}\times\mathcal{B}^\pm_{S_0^\pm}$ we
 obtain an identification of $T_{u_\pm}\mathcal{B}^\pm$ with
 $V^\pm\oplus F^\pm \oplus L^{p,d}_k(u_\pm^*TM_\infty,\bdy)$.  Let
 $W^\pm_{u_\pm}=D\dbar(V^\pm\oplus\{0\}\oplus\{0\})\subset
 \Gamma(\Lambda^{0,1}T^*S^\pm\otimes_J u_\pm^*TM_\infty)$.  (We use $\Gamma$ to denote $C_0^\infty$ sections.)  Later, we
 will assume that the map $D\dbar$ is surjective at $u_\pm$.  This
 is equivalent to the statement that
 $W^\pm_{u_\pm}+D\dbar_S(F^\pm\oplus L^{p,d}_k(u_\pm^*TM_\infty,\bdy))$ is all of
 $$L^{p,d}_{k-1}(\Lambda^{0,1}T^*S^\pm\otimes_J u_\pm^*TM_\infty)=L^{p,d}_{k-1}\left(\Lambda^{0,1}u_\pm^*TM_\infty\right).$$
 (Compare~\cite[pages 55--56]{bourgeois}.)   

As the particulars of our asymptotics are irrelevant to linear
statements,  the following lemma is completely standard; see, for
instance,~\cite[Theorem 4, page 797]{FloerGradFlow}.
\begin{Lem}For appropriate (small) $d>0$, the linearized $\dbar$--operator $D\dbar$ is Fredholm.
\end{Lem}

We will use a gluing construction to define spaces $\mathcal{T}_R$
of glued surfaces; then
 $\mathcal{B}^R$, $\mathcal{E}^R$, and $\mathcal{M}^R$ are defined
 similarly to $\mathcal{B}^\pm$, $\mathcal{E}^\pm$, and $\mathcal{M}^\pm$.  As described later, it will be important that we modify the metric on
 $\mathcal{B}^R$ and norm on $\mathcal{E}^R$ by adding exponential
 weights in the necks, but the topologies on the spaces
 will not be affected. 
We define $\mathcal{B}$ (respectively $\mathcal{E}$,
$\mathcal{M}$) to be the union over $R\in[0,\infty]$ of
$\mathcal{M}^R$ (respectively $\mathcal{E}^R$, $\mathcal{M}^R$).  

\subsection{Notation}
In summary, so far we have a symplectic manifold $M$, hypersurface
$H\rembds M$, and Lagrangian submanifold $L\rembds M$.  Splitting
$M$ along $H$ we obtain manifolds $M_R$ and $M_\infty$, each
with Lagrangian submanifolds obtained from $L$, which we still
denote $L$.  The split manifold $M_\infty$ has two new ends
$E_0^\pm$ corresponding to $H$.

We have families $\mathcal{T}^\pm$ of complex curves and spaces
$\mathcal{B}^\pm$ of $W^{p,d}_k$ maps of curves in
$\mathcal{T}^\pm$ into $M_\infty$, asymptotic to Reeb orbits
$\gamma_{0,1}=\gamma_{0,1}^+=\gamma_{0,1}^-,\cdots,\gamma_{0,m_0}=\gamma_{0,m_0}^+=\gamma_{0,m_0}^-$
and Reeb chords
$\gamma_{0,m_0+1}=\gamma_{0,m_0+1}^+=\gamma_{0,m_0+1}^-,\cdots,\gamma_{0,n_0}=\gamma_{0,n_0}^+=\gamma_{0,n_0}^-$
in  the ends $E_0^\pm$ at punctures $p_{0,j}^\pm$ and various
other fixed Reeb chords in the other ends.  We have defined bundles
$\mathcal{E}^\pm$ over $\mathcal{B}^\pm$ of which $\dbar$ is a
section, and $\mathcal{M}^\pm$ to be the intersection of $\dbar$
with the $0$--section of $\mathcal{E}^\pm$.  The linearized
$\dbar$--operator is denoted 
$$D\dbar:W_{u_\pm}^\pm\oplus F^\pm\oplus
L^{p,d}_k\left(u_\pm^*TM_\infty,\bdy\right)\to L^{p,d}_{k-1}\left(\Lambda^{0,1}u_\pm^*TM_\infty\right)$$
where $F^\pm\cong
\bC^{\sum_{i=0}^Nm_i^\pm}\oplus\bR^{\sum_{i=0}^N(n_i^\pm-m_i^\pm)}$
is generated by sections $v_{\tau,i,j}^\pm$ and
$v_{\theta,i,j}^\pm$ which agree with $\partial/\partial\tau$ and
$\partial/\partial\theta$ respectively near the puncture $p_{i,j}^\pm$.  When restricting to a
fixed source $S^\pm$ those objects are replaced by
$\mathcal{B}_{S^\pm}^\pm$, $\mathcal{E}_{S^\pm}^\pm$,
$\mathcal{M}_{S^\pm}^\pm$, $\dbar_{S^\pm}$ and $D\dbar_{S^\pm}$ respectively.

To perform the gluing we will need cutoff functions.  Fix a smooth function $\beta\co \bR\to[0,1]$ such that
$$
\left\{\begin{array}{ll}
\beta(t)=0 & \qquad\textrm{if $t\leq0$}\\
\beta(t)=1 & \qquad\textrm{if $t\geq 1$}\\
0\leq \beta^{(\ell)}(t)\leq 2^\ell & \qquad\textrm{for all $t\in\bR$ and $\ell\in\bN$}.
\end{array}\right.
$$

For $R/2>T_{\gamma_{0,i}^\pm}$, define the cutoff function $\beta_{i,R}(t)=\beta\left(\frac{T_{\gamma_{0,i}^\pm}(t-1)}{R/2-T_{\gamma_{0,i}^\pm}}\right)$, so that 
$$
\left\{\begin{array}{ll}
\beta_{i,R}(t)=0 &\textrm{\qquad for $t\leq 1$}\\
\beta_{i,R}(t)=1  &\textrm{\qquad for $t\geq \frac{R}{2T}$}\\
0\leq \beta_{i,R}^{(\ell)}(t)\leq\left(\frac{2T}{R/2-T}\right)^\ell & \textrm{\qquad for all $t\in\bR$ and $\ell\in\bN$}.
\end{array}
\right.
$$
\medskip\textbf{Remark}\qua  For our proofs, the explicit bounds on the higher derivatives will not be important.  All that will matter is that all derivatives of $\beta_{i,R}$ are uniformly bounded in $R$, and that $\beta'_{i,R}\to 0$ as $R\to\infty$.

\subsection{Gluing estimates}
Choose elements $u^\pm\co S^\pm\to M_\infty$ in $\mathcal{M}^\pm$.  Let $ind^\pm$ denote the index of $D\dbar$ at $u^\pm$.  Since we are assuming transversality, there is a neighborhood of $u^\pm$ in $\mathcal{M}^\pm$ diffeomorphic to $\bR^{ind^\pm}$.  We want to show that there is a neighborhood of the two--story holomorphic building $(u^-,u^+)$ in $\mathcal{M}$ diffeomorphic to $\bR^{ind^-}\times\bR^{ind^+}\times(0,1]$.  To do this, we will use the
\begin{Lem}(Implicit Function Theorem)\label{Lemma:IFT}  Let $f\co E\to
  F$ be a smooth map of Banach spaces with a Taylor expansion
$$
f(\xi)=f(0)+Df(0)\xi+N(\xi).
$$
Assume $Df(0)$ has a finite--dimensional kernel and a right inverse $Q$ satisfying
$$
\left\|QN(\xi)-QN(\eta)\right\|\leq C\left(\|\xi\|+\|\eta\|\right)\|\xi-\eta\|,
$$
for some constant $C$.  Assume also that
$\|Qf(0)\|\leq\frac{1}{8C}.$  Then for $\delta=1/(4C)$,
$f^{-1}(0)\cap B_\delta(\xi)$
is a smooth manifold of dimension $\dim\ker Df(0)$.  In fact, there
is a smooth function $\pi\co \ker Df(0)\cap B_\delta(\xi)\to Q(F)$ such
that $f(\xi+\phi(\xi))=0$ and all zeroes of $f$ in
$B_\delta(\xi)$ have the form $\xi+\phi(\xi)$.
\end{Lem}
This result is~\cite[Proposition 24, page 25]{FloerMonopoles}.  The proof is
essentially the same as the finite--dimensional case.  The result,
which A. Floer refers to as Newton's method, is called ``Floer's Picard
Lemma'' by some authors.  In the literature, I have found at least five incorrect
references to its location.

We will apply the implicit function theorem to $\dbar\co
\mathcal{B}\to\mathcal{E}$.  Recall that for $u:S\to M$,
$T_u\mathcal{B}=V\oplus F\oplus L^{p,d}_k(u^*TM,\bdy)$ for some
finite--dimensional space $V\oplus F$.  Choose a norm on $V\oplus F$ arbitrarily;
together with the $L^{p,d}_k$--norm this makes
$T\mathcal{B}$ into a Banach space.  

We give two necessary general results before turning to the estimates
needed by the implicit function theorem.
\begin{Lem}\label{Lemma:SomeBounds}At $(u\co S\to M)\in W^{p,d}_k$ the remainder $N$ in the Taylor expansion $$\dbar(u+\Delta u)=\dbar(u)+D\dbar(u)\Delta u + N(\Delta u)\in L^{p,d}_{k-1}(\Lambda^{0,1}T^*S\otimes_Ju^*TM_\infty)$$
for the $\dbar$--operator satisfies 
\begin{equation}\label{Nbound}
\|N(\xi)-N(\eta)\|\leq C'\left(\|\xi\|+\|\eta\|\right)\|\xi-\eta\|
\end{equation}
for $\xi$, $\eta$ in $V\oplus F\oplus L^{p,d}_k(u^*TM,\bdy)$ and some constant $C$ (depending on $\|u\|_{W^{p,d}_{k-1}}$).
(Here, we identify a neighborhood of the map $u$ with a neighborhood
of the $0$--section in $V\oplus F\oplus L^{p,d}_k(u^*TM,\bdy)$.)
\end{Lem}
\proof
The argument is essentially the same as the one used by Floer to prove~\cite[Theorem 3a]{FloerGradFlow}, and we omit it.  Floer finds a relatively explicit formula \cite[Lemma 3.2]{FloerGradFlow} for $N(\xi)$, using charts, and simple bounds on the terms in the formula.
There are several differences between our setup and the setup
of~\cite{FloerGradFlow}.  We claim the bound~\eqref{Nbound} for all
$L^{p,d}_k$, while Floer only states it for $L^{p,d}_1$; however,
Floer in fact proves the result for all $L^{p,d}_k$.  The
holomorphic curves in~\cite{FloerGradFlow} all have source a strip,
but this is again essentially irrelevant for his proof.  The extra
finite--dimensional spaces cause no additional complications.  Finally, the asymptotics considered in~\cite{FloerGradFlow} are somewhat different from the ones we consider, but as the estimate follows from uniform pointwise bounds, this is yet again irrelevant to the proof.

For slightly weaker estimates, still sufficient to apply the implicit function theorem, see~\cite[Section 3.3]{MS2}
\endproof

\begin{Cor}\label{Cor:IFTApplies} If $Q$ is a bounded right inverse
  for $D\dbar(u)$ then there is a constant $C$, linear in $\|Q\|$, so that
$$
\left\|QN(\xi)-QN(\eta)\right\|\leq C\left(\|\xi\|+\|\eta\|\right)\|\xi-\eta\|
$$
for $\xi, \eta\in V\oplus F\oplus L^{p,d}_k(u^*TM,\bdy)$.  Thus, the inverse function theorem applies to find nearby solutions of the $\dbar$--equation.
\end{Cor}
\proof
Take $C=\|Q\|C'$, for $C'$ as in \fullref{Lemma:SomeBounds}.
\endproof

So, the strategy to prove the gluing lemma is to construct a family of pre--glued maps $(u_+\natural_R u_-)\co (S^+\natural_R S^-)\to M_R$ so that
\begin{enumerate}
\item the maps $u_+\natural_R u_-$ converge to the height two building $(u_+,u_-)$ and
\item at each $u_+\natural_R u_-$ there is a right inverse $Q_R$ to $D\dbar$ so that the $Q_R$ are uniformly bounded (in $R$).
\end{enumerate}
The first condition then implies that $D\dbar(u_+\natural_R u_-)\to 0$ as $R\to\infty$.  It follows from this and the second condition that for large enough $R$ the implicit function theorem applies to give families of solutions of the $\dbar$--equation near $u_+\natural_R u_-$, and hence near $(u_+,u_-)$.  Further, any other solution of the $\dbar$--equation near $(u_+,u_-)$ would lie in a small neighborhood of $u_+\natural_R u_-$ for appropriate large $R$, and hence be one of the solutions given by the implicit function theorem.  This then proves the gluing lemma.

Now we define the pre--glued maps $u_+\natural_R u_-$.  Choose
holomorphic coordinates $(t_i^\pm,s_i^\pm)$ ($s_i^\pm\in S^1$ for
$i=1,\cdots,m_0$, $s_i^\pm\in [0,1]$ for $i=m_0+1,\cdots,n_0$,
$t_i^\pm\in\bR$) near the puncture $p_{0,i}^\pm$ so that
$t_i^\pm\to\mp\infty$ as $p\to p_{0,i}^\pm$.  For
$i=1,\cdots,m_0$ we further require that
$\lim_{t_i^+\to-\infty}u_+(t_i^+,s_i^+)=\lim_{t_i^-\to\infty}
u_-(t_i^-,s_i^-).$  For $i=m_0+1,\cdots,n_0$ this is automatic.  We
will call such
coordinates \emph{cylindrical coordinates}.

\medskip\textbf{Remark}\qua  Our
  convention for which coordinate is denoted $s$ and which is
  denoted $t$ is exactly the opposite from~\cite{bourgeois}, but
  agrees with the convention used in the rest of this paper.

Let $\tau^\pm$ denote the $\bR$--coordinate on $E_0^\pm$ chosen
earlier in this section, let $\pi_{\tilde{E}_0}$ denote projection of $E_0^\pm$ onto
$\tilde{E}_0$ and fix some metric on $\tilde{E}_0$ so that
$L\cap\tilde{E}_0$ is totally geodesic (this is possible by the
Lagrangian neighborhood theorem).  Near $p_{0,i}^\pm$ the map
$u_\pm$ has the form 
\begin{eqnarray*}
\tau^\pm\circ
u_\pm(t_i^\pm,s_i^\pm)&=&T_{\gamma_{0,i}}t_i^\pm\pm\tau_{0,0}^\pm+\eta_{i}^\pm(t_i^\pm,s_i^\pm)\\
\pi_{\tilde{E}_0}\circ u &=&
\exp_{\gamma_{0,i}(s_i^\pm)}(U_i^\pm(t_i^\pm,s_i^\pm))
\end{eqnarray*}
where $\eta_i^\pm$ and $U_i^\pm$ decay exponentially in $t_i$.

Choose $R_0$ large enough that for $R>R_0$ and $i=1,\cdots,n_0$,
$\left(\mp\frac{R\pm\tau_{0,0}^\pm}{T_{\gamma_{0,i}}},s_i^\pm\right)$
lies in the neighborhood of $p_{0,i}^\pm$ on which the coordinates in use are defined.  From now on we will assume that $R>R_0$. 
Then $\tau^\pm\circ u_\pm\left(\mp\frac{R\pm\tau_{0,0}^\pm}{T_{\gamma_{0,i}}},s_i^\pm\right)=\pm R$.  So, 
define $S^+\natural_R S^-$ to be the surface obtained by deleting
$\left\{\mp t_i^\pm>\frac{R\pm\tau_{0,0}^\pm}{T_{\gamma_{0,i}}}\right\}$ from
$S^\pm$ and gluing the resulting surfaces along the newly created
boundary in the obvious way.  We will refer to the image of the
neighborhoods of the punctures $p_{0,i}^\pm$ in $S^+\natural_R
S^-$ as the \emph{necks}.  Define coordinates $(t^0_i,s_i^0)$ on the necks by
$$
\left\{\begin{array}{lll}
s^0_i=s^+_i, & t^0_i=t^+_i+\frac{R+\tau_{0,0}}{T_{\gamma_{0,i}}} &\qquad\textrm{for $t^+_i\geq-\frac{R+\tau_{0,0}}{T_{\gamma_{0,i}}}$}\\
s^0_i=s^-_i, & t^0_i=t^-_i-\frac{R-\tau_{0,0}}{T_{\gamma_{0,i}}}&\qquad\textrm{for $t^-_i\leq\frac{R-\tau_{0,0}}{T_{\gamma_{0,i}}}$}
\end{array}\right. .
$$

We define the pre--glued map $u_+\natural_R u_-\co S_+\natural_R
S_-\to M_R$ to agree with $u_+$ or $u_-$ outside the necks.
Recall that $M_R$ has a neck $[-R,R]\times \tilde{E}_0$; let
$\tau^R$ denote the $[-R,R]$--coordinate on the neck of $M_R$,
and $\pi_{\tilde{E}_0}$ projection of the neck of $M_R$ onto
$\tilde{E}_0$.  
Then, on the $i^{\rm th}$ neck of $S_+\natural_R S_-$, $u_+\natural_R u_-$ is defined by
$$
\left\{\begin{array}{ll}
\left.\begin{array}{l}
\tau^R\circ (u_+\natural_R u_-)(s^0_i,t^0_i)=T_{\gamma_{0,i}}t^0_i+\beta(t^0_i-1)\eta_i^+(s^0_i,t^0_i)\\
\pi_{\tE_0}\circ (u_+\natural_R u_-)(s^0_i,t^0_i) = exp_{\gamma_{0,i}(s_i)}\left(\beta(t_i^0-1)U_i^+(s^0_i,t^0_i)\right)
\end{array}\right\} & \hspace{22pt}\textrm{if $t_i^0\geq1$}
\\ & \\
\left.\begin{array}{l}
\tau^R\circ(u_+\natural_R u_-)(s^0_i,t^0_i) = T_{\gamma_{0,i}}t^0_i\\
\pi_{\tE_0}\circ (u_+\natural_R u_-)(s^0_i,t^0_i) = \gamma_{0,i}(s_i)
\end{array}\right\} & \hspace{-4pt}\textrm{if $-1\leq t_i^0\leq 1$}
\\ & \\
\left.\begin{array}{l}
\tau^R\circ (u_+\natural_R u_-)(s^0_i,t^0_i)=T_{\gamma_{0,i}}t^0_i+\beta(-t^0_i-1)\eta_i^-(s^0_i,t^0_i)\\
\pi_{\tE_0}\circ (u_+\natural_R u_-)(s^0_i,t^0_i) = exp_{\gamma_{0,i}(s_i)}\left(\beta(-t_i^0-1)U_i^-(s^0_i,t^0_i)\right)
\end{array}\right\} & \hspace{15pt}\textrm{if $t_i^0\leq -1$}
\end{array}
\right.
$$
(Compare~\cite[page 54]{bourgeois}.)  Note that since $L\cap \tE_0$
is totally geodesic and $U_i^\pm|_{\bdy S^\pm}$ is tangent to $L\cap
\tE_0$, this formula makes perfectly good sense in the relative
setting.

We extend the gluing construction to a neighborhood of $(u_+,u_-)$
in $\mathcal{B}^+\times\mathcal{B}^-$.  The details of this
extension are unimportant, but for completeness we give
them anyway.  Extend the
coordinates $(t_i^\pm,s_i^\pm)$ smoothly to holomorphic coordinates on surfaces in a neighborhood of $(S^+,S^-)$.  Let $\theta_1,\cdots,\theta_{m_0}$ be in some interval around $0\in S^1=\bR/\bZ$ and $r_1,\cdots,r_{n_0}$ in
$\bR$, $R>R_0$.  Let $\theta_{m_0+1}=\cdots=\theta_{n_0}=0$.
Then, for $(u'_+\co {S'}^{+}\to M_\infty, u'_-\co {S'}^{-}\to M_\infty)$, define
${S'}^{+}\natural_{r_1,\cdots,r_{n_0},\theta_1,\cdots,\theta_{m_0}}{S'}^{-}$ by
deleting the
disks $\left\{\mp
  t_i^\pm>\frac{r_i\pm\tau_{0,0}}{T_{\gamma_{0,i}}}\right\}$ from
${S'}^{\pm}$ and identifying the newly created boundaries by
$s_i^+\leftrightarrow s_i^-+\theta_i$.  Define coordinates
$(t_i^0,s_i^0)$ on the necks of the glued surfaces by
$$
\left\{
\begin{array}{lll}
s^0_i=s^+_i-\beta(t_i^0)\theta_i \hspace{.2in}&
t_i^0=\frac{R}{r_i}t_i^++\frac{r_i+\tau_{0,0}}{T_{\gamma_{0,i}}}
&\qquad\textrm{for $t_i^+\geq-\frac{r_i+\tau_{0,0}^+}{T_{\gamma_{0,i}}}$}\\
s^0_i=s^-_i, &
t^0_i=t^-_i-\frac{R}{r_i}\frac{r_i-\tau_{0,0}^-}{T_{\gamma_{0,i}}}&\qquad\textrm{for
  $t^-_i\leq\frac{r_i-\tau_{0,0}^-}{T_{\gamma_{0,i}}}$}
\end{array}
\right.
$$
Then, define
$u'_+\natural_{r_1,\cdots,r_{n_0},\theta_1,\cdots,\theta_{m_0}}u'_-$
by the same formula used to define $u_+\natural_R u_-$ above.  Let
$\mathcal{T}_R$ denote the space of conformal structures given by
the gluing construction just explained.  Note
that $\mathcal{T}_R$ projects onto an open neighborhood of $S^+\natural_R
S^-$ in moduli space.  Let $\mathcal{B}_R$ denote the
space of $W^{p,d}_k$ maps from surfaces in $\mathcal{T}_R$ with
the specified boundary conditions and asymptotics.  Note that the
glued maps fill out a neighborhood of $u_+\natural_R u_-$ in
$\mathcal{B}_R$.

Recall that we defined spaces $W^\pm_{u^\pm}=:W^\pm$ as the image
under $D\dbar$ of $V^\pm=T_{S^\pm}\mathcal{T}$.  Since the
conformal structures in $\mathcal{T}$ agree near the punctures, we
can view $V^+\oplus V^-$ as a subspace of
$\Gamma(\Lambda^{0,1}T^*(S^+\natural_R S^-))$ and $W^+\oplus W^-$
as a subspace of $(u_+\natural_R u_-)^*TM_R$ for any $R$.

Recall also that we defined spaces
$F^\pm=\bC^{\sum_{i=0}^Nm_i^\pm}\oplus\bR^{\sum_i=0^N(n_i^\pm-m_i^\pm)}$
corresponding to the span of a fixed collection of ${\sum_{i=0}^Nm_i^\pm}$ vector
fields $v^\pm_{\tau,i,j}$ and $v^\pm_{\theta,i,j}$ given by $\frac{\partial}{\partial\tau}$ and
$\frac{\partial}{\partial\theta}$ near the ${\sum_{i=0}^Nm_i^\pm}$
interior punctures and a fixed collection of $\sum_{i=0}^N
(n_i^\pm-m_i^\pm)$ vector fields $v^\pm_{\tau,i,j}$ given by $\frac{\partial}{\partial\tau}$
near the $\sum_{i=0}^N(n_i^\pm-m_i^\pm)$ boundary punctures of
$u^\pm$.  Let
$$F^0=\Span\{v_{\tau,0,j},v_{\theta,0,j}\}=\bC^{\sum_{i=1}^N(m_i^+m_i^-)}\oplus\bR^{\sum_{i=1}^N(n_i^++n_i^--m_i^+-m_i^-)}$$
be the sections corresponding to punctures which are not being glued.

The tangent space $T\mathcal{T}_R$ is
$ \bC^{m_0}\oplus\bR^{n_0-m_0}\oplus V^+\oplus V^-$.  (The $\bC$--
and $\bR$--summands correspond to the gluing parameters $\theta_i$
and $r_i$.)  So, the tangent
space to $\mathcal{B}_R$ is
\begin{equation}\label{equation:gluingparamsnotation}
T_{u_+\natural_R u_-}\mathcal{B}_R= \bC^{m_0}\oplus\bR^{n_0-m_0}\oplus
V^+\oplus V^-\oplus F^0\oplus L^{p,d}_{k}\left(
(u_+\natural_R u_-)^*TM_R,\bdy\right).
\end{equation}

The map $V^+\oplus V^-\to W^+\oplus W^-$ induced by $D\dbar$ is
a surjective map of finite--dimensional vector spaces, and independent
of $R$.  It therefore has a uniformly bounded right inverse.  It
therefore suffices to construct a uniformly bounded right inverse to
the map
{\setlength\arraycolsep{0pt}
\begin{eqnarray*}
\bC^{m_0}\oplus\bR^{n_0-m_0}\oplus W^+\oplus W^-&\oplus& F^0\oplus
L^{p,d}_{k} \left((u_+\natural_R u_-)^*TM_R,\bdy\right)\\
&\to&
L^{p,d}_{k-1}\left(\Lambda^{0,1}(u^+\natural_R u^-)^*TM_R\right)
\end{eqnarray*}}
given by 
$$(v_\bC,v_\bR,v_+,v_-,\xi_0,\xi)\mapsto v_+ + v_- + (D\dbar)(v_\bC+v_\bR)+
\left(D\dbar_{S^+\natural_R
    S^-}\right)(\xi_0+\xi).$$  
Here, $v_\bC$ and $v_\bR$ correspond to infinitesimal variations
of the almost complex structure $j_S$, and
$(D\dbar)(v_\bC+v_\bR)$ the image of the sum of these variations
under $D\dbar$.  If one views $v_+$ and $v_-$ as sections of
$\End(TS,j)$ as defined in \fullref{Section:Transversality} then
$D\dbar(v_\bC+v_\bR)$ is given by $J\circ d(u_+\natural_R
u_-)\circ(v_\bC+v_\bR) $.

Constructing the right inverse and proving its boundedness is what we shall do for most of the remainder of this section.

On spaces of sections over $u_+\natural_R u_-$ instead of the Banach
norm specified earlier in this section we use a Banach norm with additional weights of
$e^{d\left(R/T-s_i^0\right)}$ on the necks.  This paragraph and the
next make this precise.  Choose metrics $\langle\cdot,\cdot\rangle_R$ on $S^+\natural_R
S^-$ and diffeomorphisms $\phi_{R_0,R}\co S^+\natural_{R_0} S^-\to
S^+\natural_{R}S^-$ for so that $\langle\cdot,\cdot\rangle_{R_0}\leq \phi^*\langle\cdot,\cdot\rangle_R\leq
\frac{2R}{R_0}\left( \langle\cdot,\cdot\rangle\right)_{R_0}$, pointwise.  Integrals over $S_R$ will be
with respect to the volume forms induced by $\langle\cdot,\cdot\rangle_R$.  Choose metrics
$\langle\cdot,\cdot\rangle^M_R$ on $M_R$ and diffeomorphisms $\phi^M_{R_0,R}\co M_{R_0}\to
M_R$ such that $\langle\cdot,\cdot\rangle^M_{R_0}\leq \phi^*\langle\cdot,\cdot\rangle^M_R\leq \frac{2R}{R_0}\left(
\langle\cdot,\cdot\rangle^M_{R_0}\right)$, pointwise.  The bundles $(u_+\natural_R u_-)^*TM_R$ and
$\Lambda^{0,1}(u_+\natural_R
u_-)^*TM_R$ inherit metrics from $\langle\cdot,\cdot\rangle_R$ and $\langle\cdot,\cdot\rangle^M_R$.  Norms of
elements of these vector bundles will be taken with respect to the
induced metrics. 

Choose cylindrical coordinates $(t_{i,j}^\pm,s_{i,j}^\pm)$ near the
punctures $p_{i,j}^\pm$.  Let $\xi$ be a section of
$(u_+\natural_R u_-)^*TM_R$.  For $1\leq i\leq n_0$ define
$\pi_{\tau,i}(\xi)=\int_{t_i^0=0}\left\langle\frac{\partial}{\partial\tau^R},
\xi\right\rangle$ and for $1\leq i\leq m_0$ define
$\pi_{\theta,i}(\xi)=\int_{t_i^0=0}\left\langle\frac{\partial}{\partial\theta},
\xi\right\rangle$
Then, define
{\setlength\arraycolsep{0pt}}
\begin{eqnarray*}
\ol{\xi}&=&\sum_{i=1}^{n_0}\pi_{\tau,i}(\xi)\left(1-\beta_{i,R}(t_i^0)\right)\left(1-\beta_{i,R}(-t_i^0)\right)\frac{\partial}{\partial\tau^R}
\\
&+&\sum_{i=1}^{m_0}\pi_{\theta,i}(\xi)\left(1-\beta_{i,R}(t_i^0)\right)\left(1-\beta_{i,R}(-t_i^0)\right)\frac{\partial}{\partial\theta}.
\end{eqnarray*}
One can think of $\ol{\xi}$ as an approximate projection of $\xi$ to
$\bC^{m_0}\oplus\bR^{n_0-m_0}$ in $F^+/(F^0\cap F^+)$.

Then, the norm of a section $\xi$ of $(u_+\natural_R u_-)^*TM_R$
is given by the sum of the norm of the vector 
$$
\left(\pi_{\tau,1}(\xi),\cdots,\pi_{\tau,n_0}(\xi),\pi_{\theta,1}(\xi),\cdots,\pi_{\theta,m_0}(\xi)\right)\in\bR^{n_0+m_0}
$$
and 

\cl{\scalebox{0.86}{\small$\displaystyle
\sup_{\substack{ |\alpha|\leq k\\ \left|\frac{\partial}{\partial\alpha}\right|\equiv 1}}\left(\int_{S_+\natural_R S_-}\left|\frac{\partial}{\partial\alpha}\left[\left(1+\sum_{i=0}^{n_0}e^{d(R-|t_i^0|)/p}\beta(R-|t_i^0|)+\sum_\pm\sum_{i=1}^N\sum_{j=1}^{n_i}e^{d|t_{i,j}^\pm|/p}\beta(|t_{i,j}^\pm|)\right)(\xi-\ol{\xi})\right]\right|^pdV\right)
$}}

where the $\sup$ is over all partial derivative of order at most
$k$ and norm $1$.

The norm of a section $\eta$ of
$\Lambda^{0,1}(u_+\natural_R u_-)^*TM_R$ is given by

\cl{\scalebox{0.89}{\small$\displaystyle
\sup_{\substack{ |\alpha|\leq k-1\\ \left|\frac{\partial}{\partial\alpha}\right|\equiv 1}}\left(\int_{S_+\natural_R S_-}\left|\frac{\partial}{\partial\alpha}\left[\left(1+\sum_{i=0}^{n_0}e^{d(R-|t_i^0|)/p}\beta(R-|t_i^0|)+\sum_\pm\sum_{i=1}^N\sum_{j=1}^{n_i}e^{d|t_{i,j}^\pm|/p}\beta(|t_{i,j}^\pm|)\right)\eta\right]\right|^pdV\right).
$}}

It is not necessary to split $\eta$ into two parts.

\begin{Lem}{\rm(Compare~\cite[Lemma 5.5]{bourgeois})}\qua
 
With this Banach structure on $\Gamma\left(\Lambda^{0,1}(u_+\natural_R u_-)\right)$, if $u_\pm$ are holomorphic then 
$$
\lim_{R\to\infty} \|\dbar u_+\natural_R u_-\|_R =0.
$$
\end{Lem}
\proof
This is the same as in~\cite{bourgeois}.  Since $u_\pm$ are holomorphic, the section $\dbar u_+\natural_R u_-$ is identically zero except in regions $-2\leq t_i^0\leq -1$ and $1\leq t_i^0\leq 2$ in the necks.  In these regions we crudely bound $\|\dbar u_+\natural_R u_-\|$ by 
{\small$$\sum_i\left(\left\|\chi_{\left[-\frac{R+\tau_{0,0}^-}{T},-\frac{R+\tau_{0,0}^-}{T}+2\right]}(\nabla U_i^-,\nabla \eta_i^-)\right\|+
\left\|\chi_{\left[\frac{R-\tau_{0,0}^-}{T}-2,-\frac{R-\tau_{0,0}^+}{T}\right]}(\nabla U_i^+,\nabla \eta_i^+)\right\|\right).
$$}%
Recall that these norms are weighted by $e^{d\left(R/T-s_i^0\right)}$.  However, since $u_\pm$ are in $\mathcal{B}^{p,d}_k$, the sections $(U_i^\pm,\eta_i^\pm)$ are in $L^{p,d}_k$ and so their derivatives are in $L^{p,d}_{k-1}$.  It follows that the right hand side goes to zero as $R\to\infty$.
\endproof

\begin{Prop}  \label{Prop:RightInverseBound}Suppose that
  $(D\dbar)_{u_{\pm}}$, the linearized $\dbar$--operator at
  $u_\pm$, is surjective.  Then, for large enough $R$, the
  operator $(D\dbar)_{u_+\natural_R u_-}$ has a right inverse
  $Q_R$ which is uniformly bounded in $R$.
\end{Prop}
\proof
As with all of this section, I learned this proof mainly from
Bourgeois~\cite{bourgeois}.  Bourgeois in turn cites
McDuff--Salamon~\cite{MS2}, who say they adopted the argument from
Donaldson--Kronheimer~\cite{DK}.

Once one has seen how it goes, the proof is not particularly hard.
To prove the proposition we define linearized gluing and splitting maps
$$
g_R\co \Delta\oplus L^{p,d}_k(u_+^*TM_\infty,\bdy)\oplus L^{p,d}_k(u_-^*TM_\infty,\bdy)\to
L^{p,d}_k\left((u_+\natural_R u_-)^*TM_R,\bdy\right)$$
and
$$s_R\co 
L^{p,d}_{k-1}\left(\Lambda^{0,1}(u_+\natural_R u_-)^*TM_R\right)
\to
 L^{p,d}_{k-1}(\Lambda^{0,1}u_+^*TM_\infty)\oplus
 L^{p,d}_{k-1}(\Lambda^{0,1}u_-^*TM_\infty).
$$
where $\Delta$ is the diagonal in
$\bC^{m_0^+}\oplus\bR^{n_0^+-m_0^+}\oplus\bC^{m_0^-}\oplus\bR^{n_0^--m_0^-}\subset
F^+\oplus F^-$ (so $\Delta$ corresponds to the sections at the
punctures being glued which agree on the two sides).

Let $G_R=s_R\left(D\dbar(\bC^{m_0}\oplus\bR^{n_0-m_0})\right)$ where
$\bC^{m_0}\oplus\bR^{n_0-m_0}$ denotes the tangent space to the space of
gluing parameters, as in Equation~\eqref{equation:gluingparamsnotation}.  (The map $s_R\circ D\dbar$ is an
isomorphism when restricted to this $\bC^{m_0}\oplus\bR^{n_0-m_0}$,
and so identifies $G_R$ with $\bC^{m_0}\oplus\bR^{n_0-m_0}$.  So, we
will sometimes abuse notation and use $G_R$ when we mean its
preimage under $s_R\circ D\dbar$.)

We will check that 
\begin{align*}
D\dbar_{(u_+,u_-)}\co& W^+
\oplus 
W^-\oplus \Delta\oplus F^0\oplus G_R\oplus
L^{p,d}_k\left(u_+^*TM_\infty,\bdy\right)\\
&\hspace{-12pt}\oplus L^{p,d}_k\left(u_-^*TM_\infty,\bdy\right)
\to
L^{p,d}_{k-1}\left(\Lambda^{0,1}u_+^*TM_\infty\right)\oplus
L^{p,d}_{k-1}\left(\Lambda^{0,1}u_-^*TM_\infty\right).
\end{align*}
is surjective and has a uniformly bounded right inverse $Q_\infty$ for large
$R$.

We will then define an approximate right inverse $\tilde{Q}_R$ for 
{\setlength\arraycolsep{0pt}
\begin{eqnarray*}D\dbar_{u_+\natural_R u_-}\co \bC^{m_0}\oplus\bR^{n_0-m_0}\oplus W^+\oplus W^-\oplus
  F^0 \oplus&&
L^{p,d}_k\left((u_+\natural_R u_-)^*TM_R,\bdy\right)\\ &&\to
L^{p,d}_{k-1}\left(\Lambda^{0,1}(u_+\natural_R u_-)^*TM_R\right)\end{eqnarray*}
}
by the commutative diagram

{\disablesubscriptcorrection

\cl{\scalebox{0.95}{\small$
\xymatrix{
L^{p,d}_{k-1}\left(\Lambda^{0,1}(u_+\natural_R
  u_-)^*TM_R\right)\ar[d]^(.47){s_R}\ar[r]^(.46){\tilde{Q}_R}
& 
{\begin{array}{c}W^+\oplus W^-\oplus F^0\oplus
    \bC^{m_0}\oplus\bR^{n_0-m_0}\\ \oplus 
L^{p,d}_k\left((u_+\natural_R u_-)^*TM_R,\bdy\right)\end{array}}
\\
  L^{p,d}_{k-1}(\Lambda^{0,1}u_+^*TM_\infty)\oplus
 L^{p,d}_{k-1}(\Lambda^{0,1}u_-^*TM_\infty)\ar[r]^(.52){Q_\infty}
& 
 {\begin{array}{c}W^+\oplus W^-\oplus F^0\oplus G^R\oplus \Delta\\ \oplus L^{p,d}_k(u_+^*TM_\infty,\bdy)\oplus
 L^{p,d}_k(u_-^*TM_\infty,\bdy).\ar[u]_(.47){id\oplus id\oplus id\oplus id\oplus g_R}\end{array}}
}
$}}}

We will check that $s_R$ and $g_R$ are uniformly bounded, so that $\tilde{Q}_R$ is, also.
We will then show that $\|(D\dbar)_{u_+\natural_R
  u_-}\tilde{Q}_R-I\|\leq 1/2$, so that $(D\dbar)_{u_+\natural_R
  u_-}\tilde{Q}_R$ is invertible; the inverse of
$(D\dbar)_{u_+\natural_R u_-}$ is then given by
$\tilde{Q}_R\left((D\dbar)_{u_+\natural_R
    u_-}\tilde{Q}_R\right)^{-1}$, which is bounded by
$2\|\tilde{Q}_R\|$.

For 
\begin{align*}
&(\xi_+,\xi_-)\in \Delta\oplus
L^{p,d}_k\left(u_+^*TM_\infty,\bdy\right)\oplus
L^{p,d}_k\left(u_-^*TM_\infty,\bdy\right)\\
&\ \ \subset
\bC^{m^+_0}\oplus\bR^{n^+_0-m^+_0}\oplus
L^{p,d}_k\left(u_+^*TM_\infty,\bdy\right)\oplus\bC^{m^-_0}\oplus\bR^{n^-_0-m^-_0}
\oplus L^{p,d}_k\left(u_-^*TM_\infty,\bdy\right)\hspace{-3pt}
\end{align*}
define $g_R(\xi_+,\xi_-)$ to agree with $\xi_+$ or $\xi_-$
outside the necks and by
$$
g_R(\xi_+,\xi_-)=\left\{\begin{array}{ll}
\xi_+(t_i^0,s_i^0)+\left(1-\beta_{i,R}(t_i^0)\right)\xi_-(t_i^0,s_i^0)& \textrm{if $t_i^0\geq1$}\\
\xi_+(t_i^0,s_i^0)+\xi_-(t_i^0,s_i^0) & \textrm{if $-1\leq t_i^0\leq 1$}\\
\left(1-\beta_{i,R}(-t_i^0)\right)\xi_+(t_i^0,s_i^0)+\xi_-(t_i^0,s_i^0) & \textrm{if $t_i^0\leq -1$}
\end{array}\right.
$$
in the necks.  (Compare~\cite[pages 56--57]{bourgeois}.)  Note that
this formula makes as much sense in the relative case as in the closed
case.  Exactly the same formulas define a linearized gluing map
\begin{align*}
g_R\co L^{p,d}_{k-1}\left(\Lambda^{0,1}u_+^*TM_\infty\right)&\oplus
L^{p,d}_{k-1}\left(\Lambda^{0,1}u_-^*TM_\infty\right)\\&\to
L^{p,d}_{k-1}\left(\Lambda^{0,1}(u_+\natural_R u_-)^*TM_R\right),
\end{align*}
which shall be useful when we estimate $(D\dbar)\tilde{Q}_R-I$.

Define $s_R(\eta)=(\eta_+,\eta_-)$ where $\eta_\pm$ agrees with $\eta$ away from the punctures and near $p_{0,i}^\pm$ is given by
\begin{eqnarray}
\eta_+(t_i^0,s_i^0) &=& \beta(t_i^0)\eta(t_i^0,s_i^0)\nonumber\\
\eta_-(t_i^0,s_i^0) &=& (1-\beta(t_i^0))\eta(t_i^0,s_i^0)\nonumber.
\end{eqnarray}
(Compare~\cite[page 57]{bourgeois}.)

\begin{Lem} The maps $g_R$ and $s_R$ are uniformly bounded in $R$.
\end{Lem}
\proof
Although for notational reasons it may appear involved, the proof is
in fact straightforward.
Let $(\xi_0,\xi_+,\xi_-)\in\Delta\oplus L^{p,d}_k\left(u_+^*TM_\infty,\bdy\right)\oplus
L^{p,d}_k\left(u_-^*TM_\infty,\bdy\right)$.  Recall that, by definition, 
$$
\|g_R(\xi_0,\xi_+,\xi_-)\|=\|\ol{g_R(\xi_0,\xi_+,\xi_-)}\|+\| g_R(\xi_0,\xi_+,\xi_-)-\ol{g_R(\xi_0,\xi_+,\xi_-)}\|_{L^{p,d}_k}
$$
where $\ol{g_R(\xi_0,\xi_+,\xi_-)}$ is a certain projection of
  $g_R(\xi_0,\xi_+,\xi_-)$ to $\bC^{m_0}\oplus\bR^{n_0}$.  
Observe that as $R\to\infty$,
$\ol{g_R(\xi_0,\xi_+,\xi_-)}\to\xi_0$, and hence
$\|\ol{g_R(\xi_0,\xi_+,\xi_-)}\|$ is bounded.  In fact, by the
Sobolev inequalities, 
$$\|\ol{g_R(\xi_0,\xi_+,\xi_-)}\|\leq
C\left(\|\xi_0\|+\|\xi_+\|_{L^{p,d}_k}+\|\xi_-\|_{L^{p,d}_k}\right)$$
where $C$ depends only on the cutoff function $\beta$, and not on
$R$.

On the other hand,
\begin{eqnarray*}
\|g_R(\xi_0,\xi_+,\xi_-)-\ol{g_R(\xi_0,\xi_+,\xi_-)}\|_{L^{p,d}_k}&\leq& 
\|g_R(0,\xi_+,\xi_-)\|+\|g_R(\xi_0,\xi_+,\xi_-)\\ &&\hspace{.2in}-g_R(0,\xi_+,\xi_-)-\ol{g_R(\xi_0,\xi_+,\xi_-)}\|_{L^{p,d}_k}\\
&=& \|g_R(0,\xi_+,\xi_-)\|_{L^{p,d}_k} + \|g_R(\xi_0,0,0)
\\ && \hspace{1.1in}-\ol{g_R(\xi_0,\xi_+,\xi_-)}\|_{L^{p,d}_k}.
\end{eqnarray*}
Now, on the $i^{\rm th}$ neck,
\begin{align*}
g_R(\xi_0,0,0)-\ol{g_R(\xi_0,\xi_+,\xi_-)} =
\big(1&-\beta_{i,R}(t_i^0)\big)\big(1-\beta_{i,R}(-t_i^0)\big)
\Bigg[\xi_0-\xi_0\\
&\hspace{-20pt}-\int_{t_i^0=0}\left(\left\langle\frac{\partial}{\partial\theta},\xi_++\xi_-\right\rangle
  \frac{\partial}{\partial\theta}
+ \left\langle\frac{\partial}{\partial\tau},\xi_++\xi_-\right\rangle
  \frac{\partial}{\partial\tau}\right)\Bigg].
\end{align*}
So, 
$$
 \|g_R(\xi_0,0,0)-\ol{g_R(\xi_0,\xi_+,\xi_-)}\|_{L^{p,d}_k} \leq
 C\left(\|\xi_+\|_{L^{p,d}_k}+\|\xi_-\|_{L^{p,d}_k}\right)
$$
where $C$ depends only on the cutoff function $\beta$, and not on $R$.

Finally,
{\small
{\setlength\arraycolsep{0pt}
\begin{eqnarray*}
&&\hspace{-10pt}\|g_R(0,\xi_+,\xi_-)\|_{L^{p,d}_k}
\leq\|\xi_+\|_{L^{p,d}_k}+\|\xi_-\|_{L^{p,d}_k}\\
&&+\sum_{i=1}^{n_0}\sum_{|\alpha|\leq
  k}\left( \int_{s_i^0\in
    S^1}\int_{t_i^0=-R}^{-1}\left|\frac{\partial}{\partial\alpha}\left(1-\beta_{i,R}(-t_i^0)\xi_+(t_i^0,s_i^0)+\xi_-(t_i^0,s_i^0)\right)e^{d|t_i^0|/p} \right|^p dt_i^0\right.\nonumber\\
&&+\int_{s_i^0\in S^1}\int_{t_i^0=-1}^1\left|\frac{\partial}{\partial\alpha}\left(\xi_+(t_i^0,s_i^0)+\xi_-(t_i^0,s_i^0)\right)e^{d|t_i^0|/p} \right|^p dt_i^0\nonumber\\
&&\left.+\int_{s_i^0\in S^1}\int_{t_i^0=1}^R\left|\frac{\partial}{\partial\alpha}\left(1-\xi_+(t_i^0,s_i^0)+\left(1-\beta_{i,R}(t_i^0)\right)\xi_-(t_i^0,s_i^0)\right)e^{d|t_i^0|/p} \right|^p dt_i^0\right)^{1/p}\nonumber\\
&\leq&\|\xi_+\|_{L^{p,d}_k}+\|\xi_-\|_{L^{p,d}_k}+\sum_{i=1}^{n_0}\sum_{|\alpha|\leq k}C_{\alpha}\left(\|\xi_+\|_{L^{p,d}_k}+\|\xi_-\|_{L^{p,d}_k}\right)\nonumber\\
&\leq&C\left(\|\xi_+\|_{L^{p,d}_k}+\|\xi_-\|_{L^{p,d}_k}\right).\nonumber
\end{eqnarray*}
}}Here, the terms $\|\xi_+\|_{L^{p,d}_k}$ and $\|\xi_-\|_{L^{p,d}_k}$ on the right hand side of the first inequality take care of the contribution of $S\setminus\{\textrm{the necks}\}$ to $\|g_R(0,\xi_+,\xi_-)\|_{L^{p,d}_k}$.
The second inequality follows from the bound on the derivatives of
$\beta_{i,R}$.  This proves boundedness of $g_R$.  The
$C_\alpha$ are universal constants depending only on $\alpha$ and
the cutoff function $\beta$, not on $R$.

Turning to $s_R$, we have
\begin{equation}\label{sRestimateq}
\|s_R(\eta)\|_{L^{p,d}_{k-1}} = \|\beta(t_i^0)\eta(t_i^0,s_i^0)\|_{L^{p,d}_{k-1}}+\left\|\left(1-\beta(t_i^0)\right)\eta(t_i^0,s_i^0)\right\|_{L^{p,d}_{k-1}}.
\end{equation}
Focussing on the first term on the right hand side, we have
\begin{eqnarray*}
&&\hspace{-10pt} \|\beta(t_i^0)\eta(t_i^0,s_i^0)\|_{L^{p,d}_{k-1}}\\ 
 &&\leq \|\eta\|_{L^{p,d}_{k-1}}
 +\sum_{i=1}^{n_0}\sum_{|\alpha|\leq k-1}\Bigg(\int_{s_i^+\in
 S^1}\int_{t_i^+=-\infty}^{-\frac{R+\tau_{0,0}}{T_{\gamma_0,i}^+}}
 \left|\frac{\partial}{\partial\alpha}
 \beta\left(t_i^++\frac{R+\tau_{0,0}}{T_{\gamma_{0,i}^+}}\right)\right.\\
&&\qquad\qquad\qquad\times\left.\eta\left(t_i^++\frac{R+\tau_{0,0}}{T_{\gamma_{0,i}^+}},s_i^+\right)e^{d|t_i^+|/p}\right|^p dt_i^+\Bigg)^{1/p}.\\
 &&\leq
 \|\eta\|_{L^{p,d}_{k-1}} +
 \sum_{i=1}^{n_0}\sum_{|\alpha|\leq k-1} C'_\alpha\|\eta\|_{L^{p,d}_{k-1}}\\
 &&\leq C'\|\eta\|_{L^{p,d}_{k-1}}
\end{eqnarray*}
where the inequalities follow by the same reasoning as for $g_R(0,\xi_+,\xi_-)$.  A similar argument applies to the second term on the right hand side of Equation~\eqref{sRestimateq}, 
so $s_R$ is uniformly bounded.
\endproof

\begin{Lem}The linearized $\dbar$--operator
{\small\setlength\arraycolsep{0pt}\begin{eqnarray*}
D\dbar_{(u_+,u_-)}\co W^+\oplus W^-\oplus \Delta\oplus F^0\oplus &G_R&\oplus
L^{p,d}_k\left(u_+^*TM_\infty,\bdy\right)\oplus L^{p,d}_k\left(u_-^*TM_\infty,\bdy\right)
\\
&\to&
L^{p,d}_{k-1}\left(\Lambda^{0,1}u_+^*TM_\infty\right)\oplus
L^{p,d}_{k-1}\left(\Lambda^{0,1}u_-^*TM_\infty\right).
\end{eqnarray*}}%
is surjective and has a uniformly bounded right inverse $Q_\infty$
for large $R$.
\end{Lem}
\proof
By assumption,
$$
D\dbar_{u_\pm}\co (W^\pm_{u_\pm}\oplus F^\pm\oplus
L^{p,d}_k(u_\pm^*TM_\infty,\bdy))\to L^{p,d}_{k-1}\left(u_\pm^*TM_\infty\right)
$$
is surjective.  Recall that the inclusion of $F^\pm$ in $\Gamma
\left(u_\pm^*TM_\infty\right)$ required choosing particular sections
$v_{\tau,i,j}^\pm$ and $v_{\theta,i,j}^\pm$,
constant near the punctures, but the image of $D\dbar_{u_\pm}$ was
independent of these choices.  For appropriate choices of these
sections, $D\dbar(F^0\cap F^\pm)+D\dbar(G_R)$ is exactly
$F^\pm$.  The desired surjectivity follows.  Uniform boundedness
only requires that one observe that the ``appropriate choices''
converge, which is clear.  (See also~\cite[Corollary 5.7]{bourgeois}.)
\endproof

Now we estimate $(D\dbar) \tilde{Q}_R-I$.  Given a section
$(\xi_+,\xi_-)$ of $\Delta\oplus L^{p,d}_k(u_+^*TM_\infty,\bdy)\oplus L^{p,d}_k(u_-^*TM_\infty,\bdy)$, observe that
{\small$$
\left((D\dbar)g_R-g_R(D\dbar)\right)(\xi_+,\xi_-)=\left\{\begin{array}{ll}
0 & \textrm{\qquad outside the necks}\\
\left((D\dbar)\beta_{i,R}(t_i^0)\right)\xi_-(t_i^0,s_i^0) & \textrm{\qquad if $t_i^0\geq 1$}\\
0 & \textrm{\qquad if $-1\leq t_i^0\leq 1$}\\
-\left((D\dbar)\beta_{i,R}(-t_i^0)\right)\xi_+(t_i^0,s_i^0)& \textrm{\qquad if $t_i^0\leq -1$}.
\end{array}\right.
$$}%
Now, $D\dbar$ is a pure first order differential operator, so $|(D\dbar)\beta_{i,R}(t_i^0)|\leq C_1/(R/2-T_{\gamma_{0,i}^\pm})\leq C_2/R$ for $R$ large and some constant $C_2$.
So, for large $R$,
$$
\|\left((D\dbar)g_R-g_R(D\dbar)\right)(\xi_+,\xi_-)\|\leq \frac{C_2}{R}\|(\xi_+,\xi_-)\|.
$$

Also, outside the necks $g_R\circ s_R(\eta) =\eta$, while in the necks
\begin{align*}
g_R\circ s_R(\eta) =& g_R\left(\beta(t_i^0)\eta(t_i^0,s_i^0),(1-\beta(t_i^0))\eta(t_i^0,s_i^0)\right)\\
=& \left\{\!\!\begin{array}{ll} 
\beta(t_i^0)\eta(t_i^0,s_i^0)+(1-\beta_{i,R}(t_i^0))(1-\beta(t_i^0))\eta(t_i^0,s_i^0) &\textrm{if $t_i^0\geq 1$}\\
(\beta(t_i^0)+(1-\beta(t_i^0)))\eta(t_i^0,s_i^0)=\eta(t_i^0,s_i^0) &\textrm{if $-1\leq t_i^0\leq 1$}\\
(1-\beta_{i,R}(-t_i^0))\beta(t_i^0)\eta(t_i^0,s_i^0)+(1-\beta(t_i^0))\eta(t_i^0,s_i^0) & \textrm{if $t_i^0\leq -1$}.
\end{array}\right.\\
=& 
\eta(t_i^0,s_i^0).
\end{align*}

It follows that
\begin{align*}
\|D\dbar \tilde{Q}_R\eta - \eta\|=& \|(D\dbar) (id\oplus id\oplus
id\oplus id\oplus
g_R) Q_\infty
s_R\eta-\eta\|\\
=&\|\left((D\dbar (id\oplus id\oplus id\oplus id\oplus g_R)-(g_R)
  D\dbar))\right) Q_\infty s_R\eta\\
&\hspace{1in} +g_R(D\dbar)Q_\infty s_R\eta - \eta\|\\
=&\|\left((D\dbar (id\oplus id\oplus id \oplus id\oplus
  g_R)-(g_R)(D\dbar))\right) Q_\infty s_R\eta+\eta-\eta\|\\
\leq&\frac{C_2}{R}\|Q_\infty s_R\eta\|\\
\leq& \frac{C_3}{R}\|\eta\|\\
\to& 0\textrm{\quad as $R\to\infty$.}
\end{align*}

So, for $R$ large enough we have $\|(D\dbar)\tilde{Q}_R-I\|\leq 1/2$, so $(D\dbar)\tilde{Q}_R$ is invertible for large enough $R$.  Thus, $Q_R=\tilde{Q}_R\left((D\dbar)\tilde{Q}_R\right)^{-1}$ is a uniformly bounded right inverse to $D\dbar$.
\endproof

\proof[Proof of \fullref{Gluing:Split}]
By \fullref{Cor:IFTApplies} and \fullref{Prop:RightInverseBound}, the implicit function theorem (\fullref{Lemma:IFT}) proves the result.
\endproof

\proof[Proof of \fullref{Gluing:Twins}]
The proof of \fullref{Gluing:Split} all takes place in a small neighborhood of the split holomorphic curves in question.  There are, therefore, no additional complications because the spaces $W_\infty$ involved in \fullref{Gluing:Twins} have two kinds of ends: near each end the asymptotics are exactly the kind considered in the proof of \fullref{Gluing:Split}.

Therefore, the one point to check is that the split complex structure
achieves transversality for maps
$\bD^2\coprod\left(\amalg_{k=1}^\ell T^2\right)\to
T^2\times[0,1]\times\bR$.  
Unfortunately, this isn't true.  Presumably these maps do achieve
transversality if one considers them as lying in the space of maps
from a torus--with--boundary to $T^2\times[0,1]\times\bR$.  I do not even know how to
formulate this statement properly, however.  So, instead we use a rather
indirect argument.
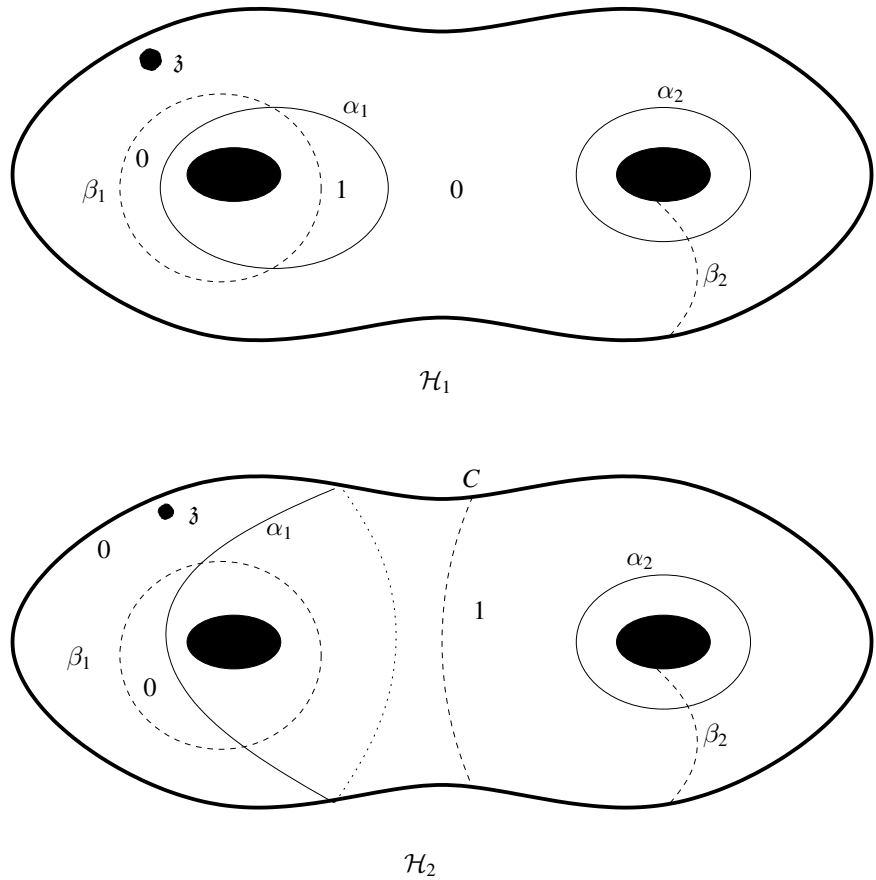
\begin{figure}
\centering
\begin{picture}(0,0)%
\includegraphics[scale=0.8]{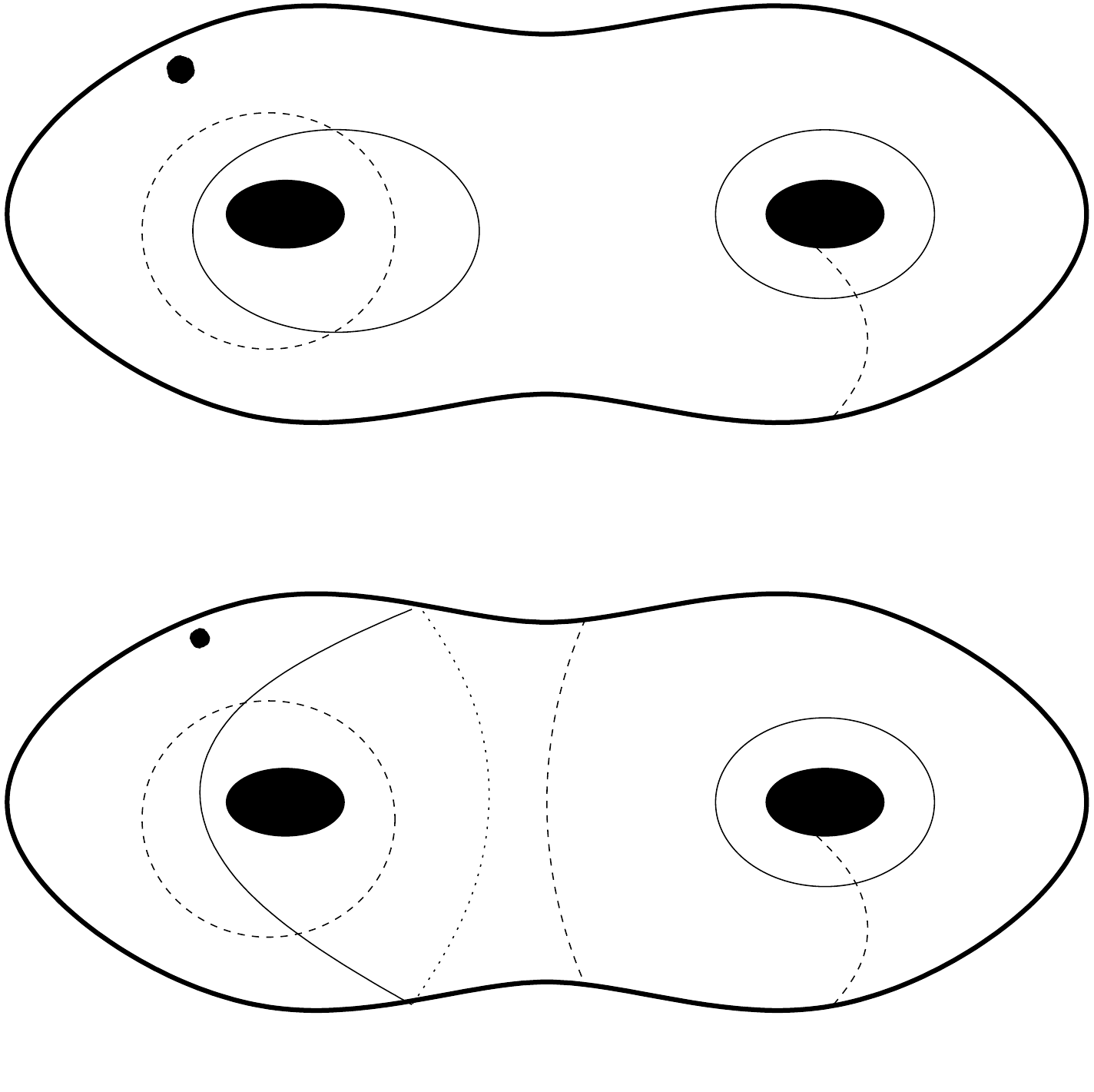}%
\end{picture}%
\setlength{\unitlength}{3157sp}%
\begingroup\makeatletter\ifx\SetFigFont\undefined%
\gdef\SetFigFont#1#2#3#4#5{%
  \reset@font\fontsize{#1}{#2pt}%
  \fontfamily{#3}\fontseries{#4}\fontshape{#5}%
  \selectfont}%
\fi\endgroup%
\begin{picture}(6810,6814)(1996,-7511)
\put(2571,-2213){\makebox(0,0)[lb]{\smash{\SetFigFont{10}{12.0}{\rmdefault}{\mddefault}{\updefault}{\color[rgb]{0,0,0}\(\beta_1\)}%
}}}
\put(4618,-1550){\makebox(0,0)[lb]{\smash{\SetFigFont{10}{12.0}{\rmdefault}{\mddefault}{\updefault}{\color[rgb]{0,0,0}\(\alpha_1\)}%
}}}
\put(7087,-1430){\makebox(0,0)[lb]{\smash{\SetFigFont{10}{12.0}{\rmdefault}{\mddefault}{\updefault}{\color[rgb]{0,0,0}\(\alpha_2\)}%
}}}
\put(7448,-2875){\makebox(0,0)[lb]{\smash{\SetFigFont{10}{12.0}{\rmdefault}{\mddefault}{\updefault}{\color[rgb]{0,0,0}\(\beta_2\)}%
}}}
\put(5220,-3718){\makebox(0,0)[lb]{\smash{\SetFigFont{10}{12.0}{\rmdefault}{\mddefault}{\updefault}{\color[rgb]{0,0,0}\(\mathcal{H}_1\)}%
}}}
\put(5100,-7511){\makebox(0,0)[lb]{\smash{\SetFigFont{10}{12.0}{\rmdefault}{\mddefault}{\updefault}{\color[rgb]{0,0,0}\(\mathcal{H}_2\)}%
}}}
\put(2451,-5885){\makebox(0,0)[lb]{\smash{\SetFigFont{10}{12.0}{\rmdefault}{\mddefault}{\updefault}{\color[rgb]{0,0,0}\(\beta_1\)}%
}}}
\put(4016,-4862){\makebox(0,0)[lb]{\smash{\SetFigFont{10}{12.0}{\rmdefault}{\mddefault}{\updefault}{\color[rgb]{0,0,0}\(\alpha_1\)}%
}}}
\put(6846,-5103){\makebox(0,0)[lb]{\smash{\SetFigFont{10}{12.0}{\rmdefault}{\mddefault}{\updefault}{\color[rgb]{0,0,0}\(\alpha_2\)}%
}}}
\put(7448,-6487){\makebox(0,0)[lb]{\smash{\SetFigFont{10}{12.0}{\rmdefault}{\mddefault}{\updefault}{\color[rgb]{0,0,0}\(\beta_2\)}%
}}}
\put(5642,-5524){\makebox(0,0)[lb]{\smash{\SetFigFont{10}{12.0}{\rmdefault}{\mddefault}{\updefault}{\color[rgb]{0,0,0}\(1\)}%
}}}
\put(2691,-5042){\makebox(0,0)[lb]{\smash{\SetFigFont{10}{12.0}{\rmdefault}{\mddefault}{\updefault}{\color[rgb]{0,0,0}\(0\)}%
}}}
\put(3053,-6126){\makebox(0,0)[lb]{\smash{\SetFigFont{10}{12.0}{\rmdefault}{\mddefault}{\updefault}{\color[rgb]{0,0,0}\(0\)}%
}}}
\put(4558,-2213){\makebox(0,0)[lb]{\smash{\SetFigFont{10}{12.0}{\rmdefault}{\mddefault}{\updefault}{\color[rgb]{0,0,0}\(1\)}%
}}}
\put(5461,-2213){\makebox(0,0)[lb]{\smash{\SetFigFont{10}{12.0}{\rmdefault}{\mddefault}{\updefault}{\color[rgb]{0,0,0}\(0\)}%
}}}
\put(2993,-1972){\makebox(0,0)[lb]{\smash{\SetFigFont{10}{12.0}{\rmdefault}{\mddefault}{\updefault}{\color[rgb]{0,0,0}\(0\)}%
}}}
\put(3294,-1189){\makebox(0,0)[lb]{\smash{\SetFigFont{10}{12.0}{\rmdefault}{\mddefault}{\updefault}{\color[rgb]{0,0,0}\(\Fz\)}%
}}}
\put(3414,-4741){\makebox(0,0)[lb]{\smash{\SetFigFont{10}{12.0}{\rmdefault}{\mddefault}{\updefault}{\color[rgb]{0,0,0}\(\Fz\)}%
}}}
\put(5551,-4486){\makebox(0,0)[lb]{\smash{\SetFigFont{12}{14.4}{\rmdefault}{\mddefault}{\updefault}{\color[rgb]{0,0,0}\(C\)}%
}}}
\end{picture}
\caption{Model stabilization degeneration}
\label{Figure:TwinsGluing}
\end{figure}

Fix $\ell$ distinct points $w_1,\cdots w_\ell$ in
$(0,1)\times\bR$, and a point $z_0\in T^2$.  Fix an almost complex
structure $J$ on $T^2\times[0,1]\times\bR$ satisfying
(\textbf{J1})--(\textbf{J4}) and (\textbf{J5$'$}), which is split near $\{z_0\}\times[0,1]\times\bR$.
Let $\mathcal{N}_J$ denote  the space of maps $u$ from $\bD^2\#_{k=1}^\ell T^2$ with
$\ell$ marked points $p_1,\cdots,p_\ell$ to
$T^2\times[0,1]\times\bR$ in the homology class $\ell[T^2]$, so
that $u(p_i)=(z_0,w_i)$.  For a generic choice of $J$, we have
transversality of the $\dbar$--operator, and the space
$\mathcal{N}_J$ is a smooth, compact, oriented $0$--manifold.  We
want to show that $\#\mathcal{N}_J=1$; then we will use
$\mathcal{N}_J$ to perform the gluing.

To show $\#\mathcal{N}_J=1$, first consider the case $\ell=1$.
We show $\#\mathcal{N}_J=1$ by using a special case of stabilization
invariance.  Specifically, consider the two Heegaard diagrams for
$S^1\times S^2$ shown in \fullref{Figure:TwinsGluing}.  They differ by two
handleslides of $\alpha_1$ over $\alpha_2$.  For an appropriate
choice of coherent orientation system, in the first diagram ($\mathcal{H}_1$),
it is easy to see that $\widehat{HF}=\bZ\oplus\bZ$.  It follows from
handleslide invariance that the same is true for the second diagram ($\mathcal{H}_2$).
It follows that, for $D$ the domain indicated by the numbers in the second diagram,
$\#\mathcal{M}^D=\pm 1$.

Now, stretch the neck in the second diagram along the dark dotted
circle $C$ in \fullref{Figure:TwinsGluing}.
In the limit, $\mathcal{H}_2$ degenerates to
$T^2_{S^1\times S^2}\vee T^2_{S^3}$.  
Choose an almost complex structure on $\mathcal{H}_2$ so that, after
stretching the neck, the corresponding almost complex structure on $T^2_{S^3}$
is $J$.  

\begin{Sublem}Compactness still holds in the current context.  That
  is, let $J_0$ be an almost complex structure on $\mathcal{H}_2$
  satisfying {\rm(\textbf{J1})--(\textbf{J4})} and {\rm(\textbf{J5$'$})}, which is
  split near $C$.  Let $J_R$ denote the almost complex structure
  on $\mathcal{H}_2$ obtained by inserting a neck of modulus $R$
  along $C$.  Let
  $\{u_R\co S_R\to\mathcal{H}_2\times[0,1]\times\bR\}$ be a sequence of
  $J_R$--holomorphic
  curves with $ind(u_R)=1$.  Then there is a subsequence of
  $\{u_R\}$ which converges to a holomorphic twin tower.  (Cf \fullref{Prop:TwinCompactness}.)
\end{Sublem}
\proof
We just sketch the proof.  

Add enough marked points to each component of each $S_R$ to
stabilize it.  Replacing the $u_R$ by a subsequence if necessary, we
may assume that the $S_R$ converge.  Now,
consider the thick--thin decomposition of
$S_R$.  For given $\epsilon>0$ and large $R$,
$(\pi_\Sigma\circ u_R)^{-1}(C)$ lies in $\mathrm{Thin}_\epsilon(S_R)$.  For small
$\epsilon$ and large $R$, the component of
$\mathrm{Thin}_\epsilon(S_R)$ containing $C$ is an
annulus in a neighborhood of which $J_R$ is split.

On the thick part of $S_R$ 
convergence of a subsequence follows from the Gromov--Schwarz lemma
and Arzela--Ascoli theorem, in the standard way (cf~\cite[Section
10.2.2]{Ya2}).  So, from now on assume the $u_R$ converge in the
thick part.

On all of the thin part of $S_R$
except the components intersecting $(\pi_\Sigma\circ u_R)^{-1}(C)$,
convergence of a subsequence also follows in the
standard way (cf~\cite[Section 10.2.3]{Ya2}).  So, from now on assume
the $u_R$ converge in these thin parts.

In components of $\mathrm{Thin}_\epsilon(S_R)$ intersecting
$(\pi_\Sigma\circ u_R)^{-1}(C)$, we extract convergence of
subsequences of 
$\pi_\Sigma\circ u_R$ and then $\pi_\bD\circ u_R$ by viewing them
as sequences of holomorphic curves with converging Lagrangian boundary
conditions, as in the proof of
\fullref{Prop:TwinCompactness}.  That is, fix circles $C_L$
and $C_R$ in $\mathcal{H}_2$ so that the region between $C_L$
and $C_R$ is an annulus $A$ containing the component of
$\mathrm{Thin}_\epsilon(\mathcal{H}_2)$ containing $C$, and so
that $J_R$ is split over $A$ (for all $R$).  Let
$S_R^A=(\pi_\Sigma\circ u_R)^{-1}(A)$.  The maps
$(\pi_\Sigma\circ u_R)|_{S_R^A}\co S_R^A\to A$ are a sequence of maps
with Lagrangian boundary conditions, hence have a convergent
subsequence.  So, from now on assume that the $\pi_\Sigma\circ
u_R|_{S_R^A}$ converge.

Finally, $(\pi_\bD\circ u_R)\left((\pi_\Sigma\circ u_R)^{-1}(C_L\cup C_R)\right)$
is a collection of circles in $[0,1]\times\bR$, which converges as
$R\to\infty$ to a collection of circles in $[0,1]\times\bR$.
Viewing $\pi_\bD\circ u_R|_{S_R^A}$ as a family of maps with these
Lagrangian boundary conditions, we obtain a convergent subsequence.
Replacing the $u_R$ with this subsequence, we finally have a
convergent sequence of $u_R\co S_R\to \mathcal{H}_2\times[0,1]\times\bR$.
\endproof

Let $z_{S^1\times S^2}$ (respectively $z_{S^3}$) denote the wedge point
in $T^2_{S^1\times S^2}$ (respectively $T^2_{S^3}$).  The domain $D$
degenerates to a domain $D_{S^1\times S^2}$ in $T^2_{S^1\times
S^2}$ (respectively $D_{S^3}$ in $T^2_{S^3}$).  Now, there is
clearly a unique holomorphic curve $u_{S^1\times S^2}$ in
$\mathcal{M}^{D_{S^1\times S^2}}$.

By choosing $C$ appropriately, we can force $z_{S^3}$
to be the point $z_0$ and 
$$\pi_\bD\circ u_{S^1\times S^2}\left( (\pi_\Sigma\circ u_{S^1\times
S^2})^{-1}(z_{S^1\times S^2})\right)$$
to be the point $w_1$.
Then, the claim that $\#\mathcal{N}_J=\pm 1$ follows from the fact
that $\#\mathcal{M}^D=\pm 1$, the gluing result
\fullref{Gluing:Split}, and the remark at the beginning of
this proof.  This deals with the case $\ell=1$.

Now, the fact that $\#\mathcal{N}_J=\pm1$ for general $\ell$ follows
from the $\ell=1$ case.  Let $u$ denote a holomorphic curve
constructed for the $\ell=1$ case.  Consider a height $\ell$
holomorphic curve each story of which is $u$.  Then, gluing these
stories together and using the fact that changing the $w_i$ gives
bordant moduli spaces we see that $\#\mathcal{N}_J=\pm1$ for general
$\ell$, $z_0$, and $\{w_1,\cdots,w_\ell\}$.

Since the $\dbar$--operator is transversal for $\mathcal{N}_J$,
stabilization invariance now follows from the gluing
result~\ref{Gluing:Split} and the remark at the beginning of this proof,
just as described in \fullref{Prop:StabilizationInvariance}.
The present proposition is also immediate.
\endproof

\section[Cross--references with {[21]}]{Cross--references
  with~\cite{OS1}}
{\footnotesize\vspace{-.1in}
For the reader's convenience I include a table indicating the
correspondence between results in this paper and those in~\cite{OS1}.
(The correspondence is also indicated in the text.)  The third column
indicates how similar the statements in the two places are, and the
fourth how similar the proofs are.  (Both are on a scale of $0$,
$1$, $2$ or $\infty$, where $\infty$ indicates that I refer
to~\cite{OS1} for the proof.)  The comparison is clearly rather
subjective.

\vspace{-.2in}
\begin{center}
\begin{tabular}{l|l|c|c}
Result of this paper & Result of~\cite{OS1} & \begin{tabular}{c}Similarity\\ of statement\end{tabular}
& \begin{tabular}{c}Similarity\\ of proof\end{tabular}\\ \hline
\hline
\begin{tabular}{@{}l}\fullref{homlem} \\and
  \fullref{Lemma:Epsilon}\end{tabular} & \begin{tabular}{@{}l} Proposition 2.15\\
and Lemma 2.19\end{tabular} & 2 & 1\\ \hline
\fullref{OS:Lemma2.19} & Lemma 2.19 & 2 &$\infty$\\ \hline
\begin{tabular}{@{}l} \fullref{Prop:BdyInjectivity}\\ and
\fullref{Lem:BdyInjectivity} \end{tabular}& Proposition 3.9 & 2 & 2\\ \hline
\fullref{Cor:IndexSigma} & Lemma 2.8 & 1 & 0\\ \hline
\fullref{OSTheorem4.9} & Proposition 7.5 of \cite{OS2} & 2 & $\infty$\\ \hline
\fullref{Index:PeriodicCor} & \begin{tabular}{@{}l}Theorem 4.9\\ and Lemma 2.8\end{tabular} & 2 & 0\\ \hline
\fullref{Lemma:AltAdmiss} & Lemma 4.12 & 1 & 2\\ \hline
\fullref{Lemma:WeakAdmis} & Lemma 4.13 & 2 & 0\\ \hline
\fullref{Lemma:StrongAdmis} & Lemma 4.14 & 2 & 0\\ \hline
\fullref{Prop:MaintainAdmis} & \begin{tabular}{@{}l}Lemma 5.8 and\\ Proposition 7.2\end{tabular} &
2 & $\infty$\\ \hline
\fullref{Prop:Orientable} & Proposition 3.10 & 2 & 1\\ \hline
\fullref{CompactCor} & Theorem 3.18 & 1 & 0\\ \hline
\begin{tabular}{@{}l}\fullref{Lemma:HatDefined} and \\\fullref{Lemma:InfinityDefined} \end{tabular}&
Theorem 4.15 & 2 & 1\\ \hline
\fullref{Lemma:H1ChainMap} & Proposition 4.18 & 2 & 2\\ \hline
\fullref{ActionProp} & Proposition 4.17 & 2 & 1\\ \hline
\fullref{Prop:IsotopyInvariance} & \begin{tabular}{@{}l}Theorem 7.3\\ and Theorem 6.1\end{tabular}
& 2 & 1\\ \hline
Construct~\ref{Construct:TriangleMaps} & Theorem 8.12 & 2 & 1\\ \hline
\begin{tabular}{@{}l}\fullref{Lemma:TriPi}\\ and \fullref{Lemma:TriEpsilon}\end{tabular} & Proposition
8.3 & 2 & 1\\ \hline
\fullref{Lemma:SzWellDefined} & Proposition 8.4 & 2 & 1\\ \hline
\fullref{Lemma:SpinEquiv} & Proposition 8.5 & 2 & 2\\ \hline
\fullref{Lemma:TriangleOrientations} & Lemma 8.7 & 2 & 2\\ \hline
\fullref{Lemma:WeakTriangleAdmisImplies} & Lemma 8.9 & 2 & 0\\ \hline
\fullref{Lemma:StrongTriangleAdmisImplies} & Lemma 8.10 & 2 & 0\\ \hline
\fullref{Prop:TriangleAdmisExists} & Lemma 8.11 & 2 & $\infty$\\ \hline
\fullref{Lemma:TriangleMapsChain} & Theorem 8.12 & 1 & 1\\ \hline
\fullref{Lemma:TriangleCXIndep} & Proposition 8.13 & 2 & 1\\ \hline
\fullref{TrianglesIndependent} & Proposition 8.14 & 2 & 1\\ \hline
\fullref{Lemma:QuadOr} & Proposition 8.15 & 2 & 2 \\ \hline
\fullref{AssocProp} & Theorem 8.16 & 2 & 1\\ \hline
\fullref{Lemma:HandleslideOrientations} & \begin{tabular}{@{}l}Lemma 9.1\\ and Lemma 9.4\end{tabular} &
1 & 2\\ \hline
\fullref{Prop:HandleslideInvariance} & Theorem 9.5 & 2 & 1\\ \hline
\end{tabular}
\end{center}
\eject

\begin{center}
\begin{tabular}{l|l|c|c}
Result of this paper & Result of~\cite{OS1} & \begin{tabular}{c}Similarity\\ of statement\end{tabular}
& \begin{tabular}{c}Similarity\\ of proof\end{tabular}\\ \hline
\hline
\fullref{Prop:HSB} & Proposition 9.8 & 1 & 0\\ \hline
\fullref{Prop:HSC} & Lemma 9.7 & 2 & 2\\ \hline
\fullref{Sublem:Annuli} & Lemma 9.3 & 1 & 1\\ \hline
\begin{tabular}{@{}l}\fullref{Lem:HSAssoc} and\\ \fullref{Cor:HSCor}\end{tabular} & Lemma 9.6 & 2 & 2\\ \hline
\fullref{Prop:StabilizationInvariance} & Theorem 10.1 & 1 & 0\\ \hline
\end{tabular}
\end{center}
}
\bibliographystyle{gtart}
\bibliography{link}

\begin{thebibliography}{}
\providecommand\bibmarginpar{\leavevmode\marginpar}
\def\urlstyle#1{{\tt #1}}

\bibitem{bourgeois}
\textbf{F Bourgeois}, \emph{A {M}orse--{B}ott approach to contact homology},
  from: ``Symplectic and contact topology: interactions and perspectives
  (Toronto, ON/Montreal, QC, 2001)'', Fields Inst. Commun. 35, Amer. Math.
  Soc., Providence, RI (2003)  55--77 \xox{MR}{1969267}

\bibitem{Ya2}
\textbf{F Bourgeois}, \textbf{Y Eliashberg}, \textbf{H Hofer}, \textbf{K
  Wysocki}, \textbf{E Zehnder}, \href{http://dx.doi.org/10.2140/gt.2003.7.799}
  {\emph{Compactness results in symplectic field theory}}, Geom. Topol. 7
  (2003) 799--888 \xox{MR}{2026549}

\bibitem{BM}
\textbf{F Bourgeois}, \textbf{K Mohnke},
  \href{http://dx.doi.org/10.1007/s00209-004-0656--x} {\emph{Coherent
  orientations in symplectic field theory}}, Math. Z. 248 (2004) 123--146
  \xox{MR}{2092725}

\bibitem{DK}
\textbf{S\,K Donaldson}, \textbf{P\,B Kronheimer}, \emph{The geometry of
  four-manifolds}, Oxford Mathematical Monographs, The Clarendon Press Oxford
  University Press, New York (1990) \xox{MR}{1079726}

\bibitem{EES}
\textbf{T Ekholm}, \textbf{J Etnyre}, \textbf{M Sullivan},
  \href{http://dx.doi.org/10.1142/S0129167X05002941} {\emph{Orientations in
  {L}egendrian contact homology and exact {L}agrangian immersions}}, Internat.
  J. Math. 16 (2005) 453--532 \xox{MR}{2141318}

\bibitem{Ya1}
\textbf{Y Eliashberg}, \textbf{A Givental}, \textbf{H Hofer},
  \emph{Introduction to symplectic field theory}, Geom. Funct. Anal.  (2000)
  560--673 \xox{MR}{1826267}

\bibitem{FloerGradFlow}
\textbf{A Floer}, \emph{The unregularized gradient flow of the symplectic
  action}, Comm. Pure Appl. Math. 41 (1988) 775--813 \xox{MR}{948771}

\bibitem{FloerMonopoles}
\textbf{A Floer}, \emph{Monopoles on asymptotically flat manifolds}, from:
  ``The Floer memorial volume'', Progr. Math. 133, Birkh\"auser, Basel (1995)
  3--41 \xox{MR}{1362821}

\bibitem{GS}
\textbf{R\,E Gompf}, \textbf{A\,I Stipsicz}, \emph{4--manifolds and {K}irby
  calculus}, Graduate Studies in Mathematics 20, American Mathematical Society,
  Providence, RI (1999) \xox{MR}{1707327}

\bibitem{HLS}
\textbf{H Hofer}, \textbf{V Lizan}, \textbf{J-C Sikorav}, \emph{On genericity
  for holomorphic curves in four-dimensional almost-complex manifolds}, J.
  Geom. Anal. 7 (1997) 149--159 \xox{MR}{1630789}

\bibitem{IP}
\textbf{E-N Ionel}, \textbf{T\,H Parker}, \emph{The {G}romov invariants of
  {R}uan--{T}ian and {T}aubes}, Math. Res. Lett. 4 (1997) 521--532
  \xox{MR}{1470424}

\bibitem{Ciprian}
\textbf{P Kronheimer}, \textbf{C Manolescu}, \emph{Periodic Floer pro-spectra
  from the Seiberg--Witten equations} \xox{arXiv}{math.GT/0203243}

\bibitem{Liu}
\textbf{C-C\,M Liu}, \emph{Moduli of $J$--Holomorphic Curves with {L}agrangian
  Boundary Conditions and Open {G}romov--{W}itten Invariants for an
  $S^2$--Equivariant Pair} \xox{arXiv}{math.SG/0210257}

\bibitem{MS2}
\textbf{D McDuff}, \textbf{D Salamon}, \emph{$J$--holomorphic curves and
  quantum cohomology}, University Lecture Series 6, American Mathematical
  Society, Providence, RI (1994) \xox{MR}{1286255}

\bibitem{White}
\textbf{M\,J Micallef}, \textbf{B White},
  \href{http://links.jstor.org/sici?sici=0003-486X(199501)2:141:1%3C35:TSOBPI%%
3E2.0.CO%3B2--G} {\emph{The structure of branch points in minimal surfaces and
  in pseudoholomorphic curves}}, Ann. of Math. $(2)$ 141 (1995) 35--85
  \xox{MR}{1314031}

\bibitem{Oh}
\textbf{Y-G Oh}, \emph{Fredholm theory of holomorphic discs under the
  perturbation of boundary conditions}, Math. Z. 222 (1996) 505--520
  \xox{MR}{1400206}

\bibitem{OSTri}
\textbf{P Ozsv{\'a}th}, \textbf{Z Szab{\'o}}, \emph{Holomorphic triangles and
  invariants for smooth four-manifolds} \xox{arXiv}{math.SG/0110169}

\bibitem{OSKnots2}
\textbf{P Ozsv{\'a}th}, \textbf{Z Szab{\'o}},
  \href{http://dx.doi.org/10.2140/gt.2004.8.311} {\emph{Holomorphic disks and
  genus bounds}}, Geom. Topol. 8 (2004) 311--334 \xox{MR}{2023281}

\bibitem{OSKnots}
\textbf{P Ozsv{\'a}th}, \textbf{Z Szab{\'o}},
  \href{http://dx.doi.org/10.1016/j.aim.2003.05.001} {\emph{Holomorphic disks
  and knot invariants}}, Adv. Math. 186 (2004) 58--116 \xox{MR}{2065507}

\bibitem{OS2}
\textbf{P Ozsv{\'a}th}, \textbf{Z Szab{\'o}},
  \href{http://projecteuclid.org/getRecord?id=euclid.annm/1105737569}
  {\emph{Holomorphic disks and three-manifold invariants: properties and
  applications}}, Ann. of Math. $(2)$ 159 (2004) 1159--1245 \xox{MR}{2113020}

\bibitem{OS1}
\textbf{P Ozsv{\'a}th}, \textbf{Z Szab{\'o}},
  \href{http://projecteuclid.org/getRecord?id=euclid.annm/1105737568}
  {\emph{Holomorphic disks and topological invariants for closed
  three-manifolds}}, Ann. of Math. $(2)$ 159 (2004) 1027--1158
  \xox{MR}{2113019}

\bibitem{OSSymplectic}
\textbf{P Ozsv{\'a}th}, \textbf{Z Szab{\'o}},
  \href{http://projecteuclid.org/getRecord?id=euclid.dmj/1072058748}
  {\emph{Holomorphic triangle invariants and the topology of symplectic
  four-manifolds}}, Duke Math. J. 121 (2004) 1--34 \xox{MR}{2031164}

\bibitem{OSContact}
\textbf{P Ozsv{\'a}th}, \textbf{Z Szab{\'o}}, \emph{Heegaard {F}loer homology
  and contact structures}, Duke Math. J. 129 (2005) 39--61 \xox{MR}{2153455}

\bibitem{Tim}
\textbf{T Perutz}, \emph{A remark on K\"ahler forms on symmetric products of
  {R}iemann surfaces} \xox{arXiv}{math.SG/0501547}

\bibitem{Rasmussen}
\textbf{J Rasmussen}, \emph{Floer Homology and Knot Complements}, PhD thesis,
  Harvard University (2003) \xox{arXiv}{math.GT/0306378}

\bibitem{Taubes}
\textbf{C\,H Taubes}, \href{http://dx.doi.org/10.2140/gt.2002.6.657} {\emph{A
  compendium of pseudoholomorphic beasts in {$\Bbb R\times (S^1\times S^2)$}}},
  Geom. Topol. 6 (2002) 657--814 \xox{MR}{1943381}

\bibitem{Turaev}
\textbf{V Turaev}, \emph{Torsion invariants of $\mathrm{Spin}^c$--structures on
  3--manifolds}, Math. Res. Lett. 4 (1997) 679--695 \xox{MR}{1484699}

\bibitem{Usher}
\textbf{M Usher}, \href{http://dx.doi.org/10.2140/gt.2004.8.565} {\emph{The
  {G}romov invariant and the {D}onaldson--{S}mith standard surface count}},
  Geom. Topol. 8 (2004) 565--610 \xox{MR}{2057774}

\end{thebibliography}

\end{document}